\documentclass[12pt, a4paper]{article}

\usepackage[T1]{fontenc}
\usepackage[latin1]{inputenc}
\usepackage[english]{babel}
\usepackage{amsopn}

\usepackage{etoolbox}
\patchcmd{\thebibliography}{\section*}{\section}{}{}

\usepackage{amsmath, amssymb, amsfonts}
\usepackage{amsthm}
\usepackage{latexsym,color,mathrsfs,fancyhdr,enumerate,nicefrac,stmaryrd} 
\usepackage[section]{placeins}
\usepackage[all]{xy}
\usepackage{mathabx}
\usepackage{cancel}
\usepackage{scrextend}
\usepackage{bbm}
\usepackage{enumitem}
\usepackage{tikz}
\usetikzlibrary{arrows}
\usetikzlibrary{arrows.meta, calc}

\usepackage{hyperref}

\theoremstyle{plain}

\newtheorem{sa}{Theorem}[subsection]
\newtheorem{Thm}[sa]{Theorem}
\newtheorem{Lem}[sa]{Lemma}

\newtheorem{Cor}[sa]{Corollary}

\newtheorem{Def}[sa]{Definition}

\newtheorem{Rem}[sa]{Remark}

\newtheorem{Eg}[sa]{Example}

\newcommand{\R}{\mathbb{R}}				%
\newcommand{\F}{\mathbb{F}}				%
\newcommand{\N}{\mathbb{N}}

\newcommand{\x}{\times}						%

\newcommand{\id}{id}

\newcommand{\fa}{\mathfrak{a}}

\newcommand{\Bcal}{\mathcal{B}}
\newcommand{\Ecal}{\mathcal{E}}

\newcommand{\Ccal}{\mathcal{C}}
\newcommand{\Xcal}{\mathcal{X}}
\newcommand{\Ycal}{\mathcal{Y}}

\newcommand{\Scal}{\mathcal{S}}
\newcommand{\Lcal}{\mathcal{L}}

\newcommand{\Mcal}{\mathcal{M}}
\newcommand{\Pcal}{\mathcal{P}}

\newcommand{\sm}{\setminus}						%
\newcommand{\ins}{\subseteq} 					%
\newcommand{\sni}{\supseteq} 					%
\newcommand{\inj}{\hookrightarrow}		%

\newcommand{\ot}{\leftarrow}
\newcommand{\oto}{\leftrightarrow}
\newcommand{\srj}{\twoheadrightarrow}			%

\renewcommand{\Pr}{\mathbb{P}} 		%
\newcommand{\E}{\mathbb{E}}
\newcommand{\I}{\mathbbm{1}}
\newcommand{\Pa}{\mathrm{Pa}} 		%
\newcommand{\Ch}{\mathrm{Ch}} 		%
\newcommand{\Anc}{\mathrm{Anc}} 		
 
\newcommand{\Desc}{\mathrm{Desc}} 	
\newcommand{\Dist}{\mathrm{Dist}} 
\newcommand{\Pred}{\mathrm{Pred}} 
\newcommand{\Sc}{\mathrm{Sc}}
\newcommand{\NonDesc}{\mathrm{NonDesc}}
\DeclareMathOperator*{\Indep}{\perp\!\!\!\perp} 
\DeclareMathOperator*{\nIndep}{\cancel\Indep} 
\DeclareMathOperator*{\given}{|}
\DeclareMathOperator{\doi}{do}
\newcommand{\moral}{\mathrm{mor}}	
\newcommand{\marg}{\mathrm{mar}}
\newcommand{\aug}{\mathrm{aug}}

\newcommand{\acag}{\mathrm{acag}}			
\newcommand{\acy}{\mathrm{acy}}	
	
\newcommand{\HEDG}{HEDG}

\newcommand{\lp}{\left ( }
\newcommand{\rp}{\right ) }

\newcommand{\hyl}[2]{\hyperlink{#1}{#2}}
\newcommand{\hyt}[2]{\hypertarget{#1}{#2}}

\begin{document}

\begin{titlepage}

\title{{\bf Markov Properties for Graphical Models with Cycles and Latent Variables}}

\author{ {\Large Patrick Forré and Joris M. Mooij}\\
Informatics Institute, University of Amsterdam, The Netherlands%
\date{} 
}

\maketitle
\footnotetext{\texttt{P.D.Forre@uva.nl}, \texttt{J.M.Mooij@uva.nl}}
\thispagestyle{empty}
\end{titlepage}

\setcounter{section}{0}

\setcounter{page}{2}

\begin{abstract}

We investigate probabilistic graphical models that allow for both cycles and latent variables.
For this we introduce \emph{directed graphs with hyperedges (\HEDG{}es)}, generalizing and combining both
marginalized directed acyclic graphs (mDAGs) that can model latent (dependent) variables, and %
directed mixed graphs (DMGs) that can model cycles.
We define and analyse several different Markov properties that relate the graphical structure of a \HEDG{} with a probability distribution on
a corresponding product space over the set of nodes, for example factorization properties, structural equations properties, ordered/local/global Markov properties, and marginal versions of these. The various Markov properties for \HEDG{}es are in general not equivalent to each other when cycles or hyperedges are present, in contrast with the simpler case of directed acyclic graphical (DAG) models (also known as Bayesian networks). We show how the Markov properties for \HEDG{}es---and thus the corresponding graphical Markov models---are logically related to each other. 
\end{abstract}

\newpage

\tableofcontents

\newpage

\pagenumbering{arabic}
\pagestyle{headings}
\setcounter{page}{4}

\section{Introduction}

\subsection{Background and related work}

The elegance and simplicity of \emph{Bayesian networks (BN)}, i.e.\ probabilistic graphical models for \emph{directed acyclic graphs (DAGs)}, is rooted in the equivalence of several different versions of \emph{Markov properties} for the corresponding probability distributions (see \cite{Lau90}, \cite{Lau98}, \ref{dag-markov-prop}): 
\begin{enumerate}
\item[i)] the recursive factorization property (\hyl{rFP-DAG}{rFP}),
\item[ii)] the directed local Markov property (\hyl{dLMP-DAG}{dLMP}),
\item[iii)] the ordered local Markov property (\hyl{oLMP-DAG}{oLMP}), 
\item[iv)] the directed global Markov property (\hyl{dGMP-DAG}{dGMP}),
\item[v)] the ancestral undirected global Markov property (\hyl{auGMP-DAG}{auGMP}),
\item[vi)] the structural equations property (\hyl{SEP-DAG}{SEP}).
\end{enumerate}
Unfortunately, such DAG based models have two major shortcomings, namely that these
\begin{enumerate}
\item[I)] are not stable under marginalization/latent projection, and 
\item[II)] do not allow for cycles.
\end{enumerate}
This is in contrast to real world applications where we almost always have incomplete data like unmeasured features
and some kind of cyclic dependences like interactions or feedback loops between the observed variables.
This suggests that the problems (I) and (II) of DAG based models may often lead to misspecified models in practice. 

The problem (I) of marginalization/latent projection has been addressed in depth for \emph{acyclic} graph structures in the literature (see \cite{Lau16}, \cite{Richardson03}, \cite{Richardson09}, \cite{sadeghi2014}, \cite{Eva15}, \cite{Eva15b}, \cite{ER14}, 
\cite{SERR14}, \cite{RERS17}, \cite{RS02}, \cite{Verma93}, a.o.) by introducing 
\emph{acyclic directed mixed graphs (ADMGs)}\footnote{ADMGs are, for instance, used as the underlying graphs for so-called \emph{semi-Markovian models} \cite{Pearl09}.}, i.e.\ DAGs with bidirected edges, or, more generally, \emph{marginalized directed acyclic graphs (mDAGs)}, i.e.\ DAGs with hyperedges.
To preserve the equivalence of versions of the Markov properties from above, \emph{marginal Markov properties} were introduced for mDAGs in \cite{Eva15}. This means that an extended graph including both latent and observed variables was modelled (in an acyclic way) and the probability distribution of the observed variables is the \emph{marginal} distribution of one that satisfies a corresponding (ordinary) Markov property w.r.t.\ the mentioned extended graph.
It was shown that the \emph{marginal} versions of directed global Markov property (\hyl{dGMP}{dGMP}) and structural equations property 
(\hyl{SEP}{SEP}) are equivalent for mDAGs. One of the main results of \cite{Eva15} was then that by summarizing the whole latent space with hyperedges, the corresponding marginal Markov models do not loose or gain any expressiveness or flexibility. This means that in the acyclic case without loss of generality one can always reduce to a model where the observed variables are modelled graphically and explicitly, and where the latent part of the world is represented with hyperedges only.
One can also show (see \cite{Eva15}, \cite{Fritz12}) that i.g.\ it is not possible to reduce to bidirected edges further (e.g.\ consider three  identical copies of \emph{one} fair coin flip). Again note that all these results were proven only for \emph{acyclic} graphical structures, i.e.\ without cycles, but these insights and ideas will be used as a starting point in this paper when we allow for cycles as well.

The problems (II) that arise in directed graphs (DG) with \emph{cycles}, on the other hand, have only been solved partially (and/or under restrictive assumptions and/or with incomplete proofs, see \cite{Spirtes93}, \cite{Spirtes94}, \cite{Spirtes95}, \cite{PearlDechter96}, \cite{Neal00}, \cite{Koster96}, a.o.). In \cite{Spirtes94} an example of a directed graph with cycles and corresponding \emph{non-linear} structural equations was given such that certain expected conditional independence relations were missing. This showed that the structural equations property (\hyl{SEP}{SEP}) i.g.\ does not imply the directed global Markov property (\hyl{dGMP}{dGMP}) in the presence of cycles. Another example will be given in \ref{main-example}. 
So the following two questions arose:
\begin{enumerate}
\item[a)] What kind of conditional independencies, phrased as a general directed global Markov property, %
can we expect in
 general structural equations models (\hyl{SEP}{SEP}) in the presence of cycles?

\item[b)] Under which general assumptions on the structural equations (\hyl{SEP}{SEP}), and/or on the graph structure or variables, etc., do we get back the directed global Markov property (\hyl{dGMP}{dGMP})?
\end{enumerate}
Both questions were analysed for directed graphs in \cite{Spirtes94} under the assumptions of differentiability, reconstructable error terms, the existence of densities w.r.t.\ a Lebesgue measure for both the observed and unobserved variables, invertible Jacobian matrix, etc.,
by analysing the structure of the Jacobian matrix for the density transformation.
For such structural equations attached to a directed graph, question (a) found a solution in \cite{Spirtes94} by introducing a corresponding \emph{collapsed graph}, which then entailed correct conditional independencies. Also question (b) under the additional assumptions of a \emph{constant} Jacobian determinant was solved in \cite{Spirtes94}, which captured the case of \emph{linear} models. So at least \emph{linear} structural equations models (\hyl{SEP}{SEP}) attached to directed graphs (acyclic or not) will lead to the directed global Markov property (\hyl{dGMP}{dGMP}).
Note that all these results under these assumptions only hold if there are \emph{no (dependent) latent variables}.

In \cite{PearlDechter96} an attempt was made to show that at least for \emph{discrete} variables the structural equations property (\hyl{SEP}{SEP}) for a directed graph (acyclic or not) always implies the directed global Markov property (\hyl{dGMP}{dGMP}), but the proof was incomplete and a counterexample was given in \cite{Neal00}.

Other Markov properties for directed graphs (with cycles) were investigated in \cite{Koster96}. There it was shown that under the assumption of \emph{no (dependent) latent variables} and the existence of a strictly \emph{positive density} a certain factorization property (\hyl{aFP}{aFP}) is equivalent to the directed global Markov property (\hyl{dGMP}{dGMP}). It is worth noting that in \cite{Koster96} some undirected edges in the graphical structure were allowed as well, leading to \emph{reciprocal graphs}. %

To the best of our knowledge, the only publications so far addressing both two problems---latent (dependent) variables (I) and cycles (II)---at once restrict to linear Gaussian models for \emph{directed mixed graphs (DMGs)}, i.e.\ directed graphs with bidirected edges, but do not consider the relations to other Markov properties (e.g.\ \cite{foygel2012}).

\begin{figure}[!h]
\centering
\begin{tikzpicture}[transform shape]	
\node [draw,shape=rectangle,rounded corners,align=left] (DAG) at (0,0) {\underline{\bf No latent variables,}\\\underline{\bf no cycles}\\ \emph{BNs, DAGs}\\
	\cite{Pearl86c}, \cite{Lau90},\\ \cite{Lau98}, a.o.
	};
\node [draw,shape=rectangle,rounded corners,align=left] (mDAG) at (6,0) {\underline{\bf With latent variables,}\\\underline{\bf no cycles}\\ 	\emph{ADMGs, mDAGs, etc.}\\ 
		\cite{Verma93}, \cite{Richardson03},\\ \cite{Richardson09}, \cite{Eva15},\\ \cite{Eva15b}, \cite{ER14}, a.o.
		};
\node [draw,shape=rectangle,rounded corners,align=left] (DG) at (0,-4) {\underline{\bf No latent variables,}\\\underline{\bf with cycles}\\
	\emph{DGs, etc.}\\
	\cite{Spirtes93}, \cite{Spirtes94}, \cite{Spirtes95},\\ \cite{PearlDechter96}, \cite{Neal00}, \cite{Koster96},\\ a.o.
	\\ (BUT: incomplete proofs\\
	or restrictive assumptions)};
\node [draw,shape=rectangle,rounded corners,align=left] (HEDG) at (6,-4) {\underline{\bf With latent variables,}\\ \underline{\bf with cycles}\\
	 \emph{DMGs, HEDGes}\\ 
	\cite{foygel2012}, a.o.\\
	 (BUT: only linear Gaussian)\\
	???};
\draw[->] (DAG) to (mDAG);
\draw[->] (DAG) to (DG);
\draw[->] (DG) to (HEDG);
\draw[->] (mDAG) to (HEDG);
	\end{tikzpicture}
  \caption{Bridging the gap between  models with latent (dependent) variables and models with cycles.
	}
	\label{fig:work}
\end{figure}
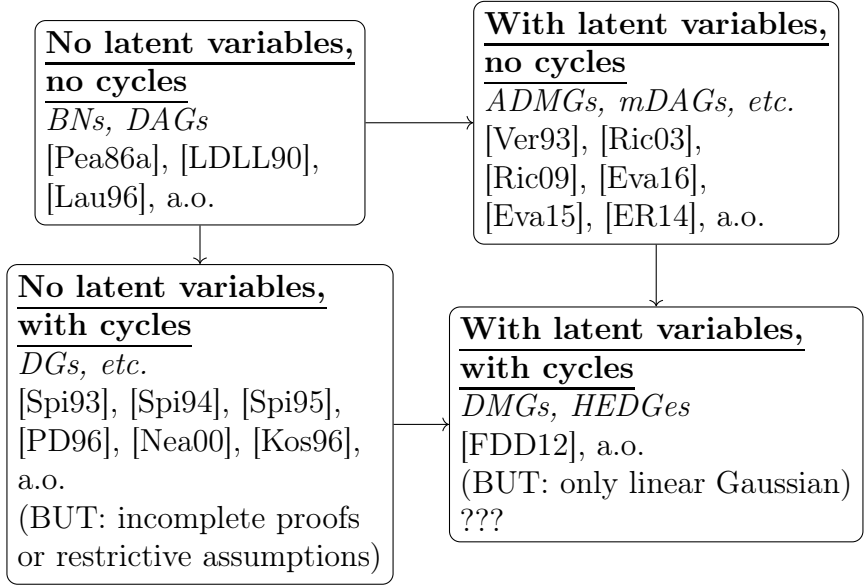

\begin{figure}[!ht]
\centering
\begin{tikzpicture}[scale=.7, transform shape]
\tikzstyle{every node} = [draw,shape=circle]
\node (v1) at (0,0) {$v_1$};
\node (v2) at (0,-2) {$v_2$};
\node (v3) at (-1.5,-4) {$v_3$};
\node (v4) at (1.5,-4) {$v_4$};
\node (v5) at (3.5,-4) {$v_5$};
\node (v6) at (3.5,-2) {$v_6$};
\node (v7) at (2,0) {$v_7$};
\node[fill,circle,red,inner sep=0pt,minimum size=5pt] (v8) at (2,-2) {};
\draw[-{Latex[length=3mm, width=2mm]}, bend right] (v1) to (v2);
\draw[-{Latex[length=3mm, width=2mm]}, bend right] (v2) to (v1);
\foreach \from/\to in {v2/v3, v2/v4, v7/v6, v7/v1, v4/v5, v6/v5}
\draw[-{Latex[length=3mm, width=2mm]}] (\from) -- (\to);
\foreach \from/\to in {v8/v7, v8/v4, v8/v5}
\draw[-{Latex[length=3mm, width=2mm]}, red] (\from) -- (\to);
\end{tikzpicture}
\hspace{1cm}
\begin{tikzpicture}[scale=.7, transform shape]
\tikzstyle{every node} = [draw,shape=circle]
\node (v1) at (-2,-2) {$v_1$};
\node (v2) at (0,-2) {$v_2$};
\node (v3) at (-2,-4) {$v_3$};
\node (v4) at (2,-4) {$v_4$};
\node (v6) at (-2,-6) {$v_6$};
\node (v5) at (-4,-6) {$v_5$};
\node (v7) at (2,-6) {$v_7$};
\node (v8) at (4,-6) {$v_8$};
\draw[{Latex[length=3mm, width=2mm]}-{Latex[length=3mm, width=2mm]}, red, bend left] (v6) to node[fill,circle,,inner sep=0pt,minimum size=5pt] {} (v7);
\draw[-{Latex[length=3mm, width=2mm]},out=135, in=225,looseness=8] (v3) to (v3);
\draw[-{Latex[length=3mm, width=2mm]}, bend right] (v5) to (v6);
\draw[-{Latex[length=3mm, width=2mm]}, bend right] (v6) to (v5);
\draw[-{Latex[length=3mm, width=2mm]}, bend right] (v7) to (v8);
\draw[-{Latex[length=3mm, width=2mm]}, bend right] (v8) to (v7);
\foreach \from/\to in {v1/v2, v2/v3, v3/v4, v4/v2, v3/v6, v4/v8}
\draw[-{Latex[length=3mm, width=2mm]}] (\from) -- (\to);
\end{tikzpicture}
  \caption{Directed graphs with hyperedges (\HEDG{}es) including both directed cycles and hyperedges (in red).}
 \label{fig:hedg}
\end{figure}
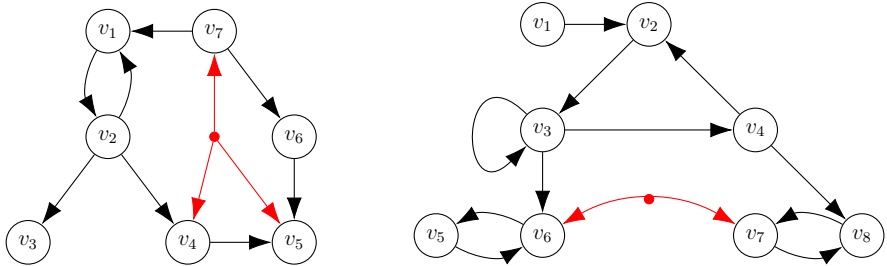

 Our paper aims to fill the gap (see figure \ref{fig:work}) of analysing Markov properties for general probabilistic graphical models for \emph{directed graphs with hyperedges} (\HEDG{}es, see figure \ref{fig:hedg}) in the presence of cycles and latent (dependent) variables while allowing for non-linear functional relations.

We also shortly want to mention other work about Markov properties on other graphical structures that allow for several edge types 
(directed edges, bidirected edges, undirected edges) like 
\emph{chain graphs} (see \cite{Lau89}, \cite{Fry90},  \cite{CW93}, \cite{Lau98}), 
\emph{AMP chain graphs} (see \cite{AMP01}, \cite{Drt09}) and
\emph{marginal AMP chain graphs} (see \cite{Pen14}),
\emph{chain mixed graphs} (see \cite{Sad16}), %
\emph{MC graphs} (see \cite{Kos02}),
\emph{summary graphs} (see \cite{Wer11}), 
\emph{maximal ancestral graphs} (see \cite{RS02}),
 and the general \emph{acyclic graphs} (see \cite{Lau16}), etc..
Since all these papers focus on different aspects, like including undirected edges into the graphical structure, and make some acyclicity assumption on the directed graphical substructure, they are complementary to our paper.

\subsection{Contributions of this paper}
Aiming for solutions to both two problems latent (dependent) variables (I) and cycles (II) at once, the contributions of this paper are as follows:
\begin{enumerate}
\item We will introduce a graphical structure that directly addresses the mentioned problems by allowing for both latent structure (for allowing for marginalizations) and cycles. These will be called \emph{directed graphs with hyperedges (\HEDG{}es)} and each such \HEDG{} will consist of nodes $V$, directed edges $E$ and omni-directed hyperedges $H$. This combines and generalizes mDAGs and directed mixed graphs (DMG) with and without cycles. In this way the nodes of such a \HEDG{} will represent observed variables and appear as a latent projection of a world with many more variables, and, as in the mDAG case (see \cite{Eva15}), the hyperedges will represent a summary of the unobserved latent world.
\item We will analyse several different appoaches of how to merge such a directed graph with hyperedges $G=(V,E,H)$ with a probability distribution $\Pr_V$ over the spaces of observed variables (the nodes). Every such approach will be called a \emph{Markov property} for $(G,\Pr_V)$, e.g.\ directed global Markov property (\hyl{dGMP}{dGMP}), structural equations properties (\hyl{SEP}{SEP}), ancestral factorization property (\hyl{aFP}{aFP}), etc..\\
Note that unlike in the DAG case all these Markov properties for \HEDG{}es might not be equivalent anymore and all induced models might be different. So before working with models as a first layer of theory we need to figure out how these Markov properties are logically related to each other and under which assumptions they might become equivalent. This will be the major contribution of this paper. Since every generalization of a Markov property from the DAG to the \HEDG{} case focuses on a different aspect, the number of Markov properties quickly blows up. To keep track of the found relations an overview is given in figure \ref{fig:overview}.
\item Also, since assumptions like having a density, or even a strictly positive one, having a discrete distribution or linear equations etc.\ might become very restrictive when working with structural equations we make the effort of trying to state the mentioned properties in most generality, be mathematically precise, and only make additional assumptions when necessary. To convince the reader and ourselves of the mentioned generality and for completeness purposes we reproduced, generalized and fixed proofs also of already known and published results. 
\end{enumerate}

\subsection{Results and description of the sections}

The main part of this paper is divided into three sections: \emph{Graph Theory}, \emph{Markov Properties for \HEDG{}es} and \emph{Graphical Models}. \\

In the \emph{Graph Theory} section we will introduce the theory of \emph{directed graphs with hyperedges (\HEDG{}es)} as a generalization of marginalized directed acyclic graphs (mDAGs) from \cite{Eva15} and directed mixed graphs (DMGs) (with cycles). 
A \HEDG{} $G=(V,E,H)$ will consist of a set of nodes $V$, directed edges $E$ and omni-directed hyperedges $H$. As in \cite{Eva15} the hyperedges will be interpreted as a summary of the whole latent world (mathematical justifications for this heuristic will be given in \ref{mMP-lift}, \ref{mdGMP-lift}, \ref{lsSEP-lift}).
To capture the notion of cycles we will see that every \HEDG{} $G$ can be decomposed into its \emph{strongly connected components}, i.e.\ its biggest subsets $S \ins G$ for which every two distinct elements $v,w \in S$ are connected via a ``cycle'', i.e.\ via two directed paths $v \to \dots \to w$ \emph{and} $w \to \dots \to v$. 
Complementary to these and to capture the notion of hyperedges we can also decompose every \HEDG{} $G$ into its \emph{districts}, i.e.\ its biggest subsets $D \ins G$ for which every two distinct elements $v,w \in D$ are connected via a path of hyperedges, e.g.\ $v \oto \cdots \oto w$.

One of the most important concepts for \HEDG{}es are the \emph{marginalizations} (latent projections) w.r.t.\ some subset of nodes.
We will see that the class of directed graphs with hyperedges is closed under such marginalization operations.

Furthermore, we will introduce two different separation concepts for \HEDG{}es that both generalize all the d-/m-/m*-separation notions of DAGs, mDAGs (see \cite{Eva15}) and ADMGs (see \cite{Richardson03}), resp., by also respecting hyperedges and cycles. The ``naive'' generalization will again be called \emph{d-separation}, whereas the ``non-naive'' generalization will be called \emph{$\sigma$-separation}. $\sigma$-separation also generalizes the \emph{collapsed graph} criterion of \cite{Spirtes94} to general \HEDG{}es. We will see that $\sigma$-separation (in contrast to d-separation) will correctly capture conditional independence relations in general (non-linear) structural equations models with cycles and latent (dependent) variables (see \ref{csSEP-smgdGMP}).
We will, furthermore, replace the collapsed graph from \cite{Spirtes94} with a more canonical construction called \emph{acyclification} of a \HEDG{}. We will also provide a version of $\sigma$-separation that is formulated in terms of properties of paths on the \HEDG{} itself, without requiring the construction of another graph.
Further important construction for \HEDG{}es will be the (generalized) \emph{moralization}, a corresponding undirected graph partially encoding d-separation relations, and the \emph{augmentation}, a directed graph that represents hyperedges as nodes with additional directed edges, and the \emph{acyclic augmentation} combining acyclification and augmentation. We will analyse in this section how all these constructions like marginalization, moralization, augmentation, acyclification, etc., interact with each other.

In addition to that, we will generalize and combine the concept of \emph{topological orders} of DAGs and mDAGs and \emph{perfect elimination orders} of undirected graphs in several ways to work for \HEDG{}es (e.g.\ \emph{quasi-topological orders}) and find criteria for their existence, which will be weaker and thus more general than acyclicity assumptions. These will be used to find generalizations of the ordered local Markov property (\hyl{oLMP}{oLMP}) in a later section.\\

In the \emph{Markov Properties for \HEDG{}es} section we will in addition to directed graphs with hyperedges (\HEDG{}es) bring probability distributions into the picture.  
We will introduce and analyse several \emph{Markov properties} between \HEDG{}es and corresponding probability distributions on a product space over the nodes. For comparison and as a guideline we start with the known Markov properties for DAGs and then go over to \HEDG{}es.
The Markov properties for \HEDG{}es will be grouped into the following different categories and an overview between their relations is given in figure \ref{fig:overview}: 
\begin{enumerate}
\setcounter{enumi}{2}
\item \underline{Directed Markov properties}: These will include generalizations of the \emph{directed global Markov property} (\hyl{dGMP}{dGMP}) and \emph{directed local Markov property} (\hyl{dLMP}{dLMP}) based on d-separation in \HEDG{}es and the \emph{general directed global Markov property} (\hyl{gdGMP}{gdGMP}) based on the new $\sigma$-separation in \HEDG{}es. We will see that \hyl{dGMP}{dGMP} will imply the other two (\hyl{dLMP}{dLMP}, \hyl{gdGMP}{gdGMP}, see \ref{dGMP-dLMP}, \ref{dGMP-gdGMP}) and give criteria for the reverse implications (see \ref{gdGMP-dGMP} and \ref{oLMP-dLMP-dGMP} in next subsection).

\item \underline{The ordered local Markov property} (\hyl{oLMP}{oLMP}): We will show (see \ref{oLMP-dGMP}, using semi-graphoid axioms only) that the \hyl{dGMP}{dGMP} and \hyl{oLMP}{oLMP} are equivalent if the underlying \HEDG{} has a \emph{perfect elimination order}. 
This generalizes the DAG cases (see \cite{Lau90} \S6 and \cite{Verma93} Thm. 1.2.1.3), the ADMG cases (see \cite{Richardson03} \S3 Thm. 2), which both only treat topological orders, and in some sense the cases of undirected chordal graphs (see \cite{West01} Def. 5.3.12)  to the much more general cases of \HEDG{}es with perfect elimination orders. We also analyse the relations between \hyl{oLMP}{oLMP} and \hyl{dLMP}{dLMP} (see \ref{dLMP-oLMP}).

\item \underline{Undirected Markov properties}: All these Markov properties will be based on the generalized \emph{moralization} of marginalizations of a \HEDG{}, e.g.\ \emph{ancestral/refined undirected} \emph{pairwise/local/global} \emph{Markov properties} (\hyl{auPMP}{auPMP}, \hyl{auLMP}{auLMP}, \hyl{auGMP}{auGMP}, \hyl{ruPMP}{ruPMP}, \hyl{ruLMP}{ruLMP}, \hyl{ruGMP}{ruGMP}). We will show (see \ref{dGMP-auPMP}, \ref{auPMP-dGMP}) that all these undirected Markov properties are closely related to the directed global Markov property (\hyl{dGMP}{dGMP}) similar to the DAG case.

\item \underline{Factorization properties}: Here we will introduce the most general factorization property for densities, the \emph{ancestral factorization property} (\hyl{aFP}{aFP}). We will see that \hyl{aFP}{aFP} implies the \hyl{dGMP}{dGMP} (see \ref{aFP-dGMP}, \ref{dGMP-auPMP}). Under the assumption of the existence of a density and perfect elimination order or under strict positivity of the density we also have the reverse implication (see \ref{dGMP-aFP}). 
This generalizes the results of \cite{Koster96} for directed graphs to general \HEDG{}es.
We also introduce a generalized \emph{recursive factorization property} (\hyl{rFP}{rFP}) w.r.t.\ some total order, which turns out to be equivalent to the \hyl{oLMP}{oLMP} for \HEDG{}es (see \ref{oLMP-rFP}) similar to the DAG case. If we have a perfect elimination order and density then \hyl{rFP}{rFP} is also equivalent to the \hyl{aFP}{aFP} (see \ref{rFPwd-aFP}). Furthermore, we also introduce the \emph{marginal factorization property} (\hyl{mFP}{mFP}), which also implies \hyl{aFP}{aFP} (see \ref{mFP-aFP}) and has an easier, but more restrictive, factorization structure.

\item \underline{Marginal Markov properties}: Based on the previous ``ordinary'' Markov properties we will introduce several \emph{marginal Markov properties} like the marginal directed global Markov property (\hyl{mMP}{mdGMP}), the marginal general directed global Markov property 
(\hyl{mMP}{mgdGMP}), the marginal ancestral factorization property (\hyl{mMP}{maFP}), etc.. These will basically enforce the ``ordinary'' Markov property on the augmented graph including latent variables. We will prove that for the \hyl{mMP}{mdGMP} this condition will lift to any arbitrary large latent graphical structure showing that w.l.o.g.\ the hyperedges of a \HEDG{} can summarize the latent space without change of expressiveness (see \ref{mdGMP-lift}, \ref{dGMP-stable-marg}).
To indicate how to generalize this to other (marginal) Markov properties we introduce a general construction (see \ref{mMP-lift}, also see \ref{lsSEP-lift}, \ref{lsSEP-marg} in the next section). This result generalizes the corresponding result of \cite{Eva15} for mDAGs to general \HEDG{}es.

\item \underline{Structural equations properties}: Here we will introduce several generalizations of \emph{structural equations properties} (\hyl{SEP}{SEP}) which in addition to the structural equations come with solution functions, e.g.\ \emph{ancestrally (uniquely) solvable structural equations property} (\hyl{asSEP}{asSEP}, \hyl{ausSEP}{ausSEP}), \emph{component-wisely (uniquely) solvable structural equations property} 
(\hyl{csSEP}{csSEP}, \hyl{cusSEP}{cusSEP}) and \emph{loop-wisely (uniquely) solvable structural equations property} 
(\hyl{lsSEP}{lsSEP}, \hyl{lusSEP}{lusSEP}). We will study their relations to each other (see \ref{SEP-equiv}) and show that they are all equivalent in the acyclic, i.e.\ in the mDAG case. 
We will see that \hyl{lsSEP}{lsSEP} for arbitrary \HEDG{}es is stable under marginalizations (see \ref{lsSEP-marg}), has the lifting property (see \ref{lsSEP-lift}, \ref{mMP-lift}) and implies \hyl{mMP}{mdLMP} (see \ref{lsSEP-mdLMP}). The main results of this subsections are as follows:
\begin{enumerate}
\item We will show that the \emph{component-wisely solvable structural equations property} (\hyl{csSEP}{csSEP} or \hyl{lsSEP}{lsSEP}) implies the 
\emph{(strong marginal) general directed global Markov property} \hyl{smgdGMP}{smgdGMP} (see \ref{csSEP-smgdGMP}), roughly stating that random variables satisfying structural equations along a \HEDG{} (under some solvability assumptions) will induce the conditional independence relations coming from the $\sigma$-separation criterion. This generalizes \cite{Spirtes94} Thm.\ 5 about \emph{collapsed graphs}, where a similar statement for directed graphs but under more restrictive assumptions (existence of densities, reconstructible and independent error terms, ancestral solvability, differentiability, invertible Jacobian determinant, etc.) and a more technical proof (e.g.\ the explicit factorization of a Jacobian determinant) was shown, to the general case of \HEDG{}es, which besides cycles also allows for dependent error terms.

\item We show that any distribution of \emph{discrete} variables (not necessarily of finite domain) satisfying the \emph{ancestrally uniquely solvable structural equations property} (\hyl{ausSEP}{ausSEP}) w.r.t.\ some \HEDG{} also satisfies the \emph{marginal ancestral factorization property} (\hyl{mMP}{maFP}) and thus the \emph{marginal directed global Markov property} (\hyl{mMP}{mdGMP}) based on d-separation (see \ref{ausSEP-wmaFP}) and not just \hyl{mMP}{mgdGMP} based on $\sigma$-separation from the point above (see \ref{csSEP-smgdGMP}). This generalizes the statement and fixes the proof of \cite{PearlDechter96}. Note that by imposing  \hyl{ausSEP}{ausSEP} we avoid the counterexample from \cite{Neal00}. Further note that only the observed variables need to be discrete and their domains don't need to be finite.  We explicitely allow for latent (dependent) variables of any type and the result is stronger than in \cite{PearlDechter96} (by having \hyl{mMP}{maFP} instead of only \hyl{dGMP}{dGMP}).

\item We also give a list of conditions formulated in \hyl{SEPwared}{SEPwared} under which general structural equations with (not necessarily discrete) variables attached to a \HEDG{} will lead to the ancestral factorization property (\hyl{aFP}{aFP}) and thus to the directed global Markov property (\hyl{dGMP}{dGMP}) encoding conditional independencies based on d-separation (see \ref{SEPwared-aFP}) and not just \hyl{gdGMP}{gdGMP} based on $\sigma$-separation from above (see \ref{csSEP-smgdGMP}). The conditions in \hyl{SEPwared}{SEPwared} roughly assume that reductions of the error terms are reconstructible from the observed variables and that the Jacobian determinant of the density transformation factorizes according to the \HEDG{}.
Note that any (solvable) linear model with or without cycles and with arbitrary latent (dependent) variables with a density (not necessarily Gaussian) satisfies the conditions of \hyl{SEPwared}{SEPwared} (see \ref{linear-SEP}). The result of \ref{SEPwared-aFP}) generalizes the results  of \cite{Spirtes94} and \cite{Koster96} for linear models to structural equations models with the property \hyl{SEPwared}{SEPwared} in the presence of cycles, latent (dependent) variables and non-linear functional relations.
\end{enumerate}
\end{enumerate}

In the \emph{Graphical Models} section we will translate the introduced \emph{Markov properties} for \HEDG{}es into corresponding \emph{probabilistic graphical Markov models} (or model classes, resp.) by just considering all probability distributions that satisfy the given Markov property w.r.t.\ the \HEDG{}. The overview figure \ref{fig:overview} then can be read as inclusion of models, model classes, resp.. The most important inclusions are highlighted in \ref{prob-model-eg}. 
Since most of the marginal models are different for general \HEDG{}es (in contrast to ADMGs) an analogous analysis of the so called ``\emph{nested Markov models}'' of ADMGs or mDAGs (see \cite{Eva15}, \cite{SERR14}, \cite{RERS17}), which sit in between the ordinary and the marginal Markov models, needs much more elaboration in the \HEDG{} case and was not done in this paper.\\

In the end of this section we will show how to use a more rigid version of the loop-wisely solvable structural equation property 
(\hyl{lsSEP}{lsSEP}) to define a well-behaved class of \emph{causal models} that allow for non-linear functional relations, interactions, feedback loops and latent confounding at once (see \ref{mSCM-def}). 
This definition follows the idea that not only each single variable is the joint effect of its direct causes but that also every subsystem is the joint effect of all its joint direct causes. We will see that with this definition we easily can define arbitrary \emph{marginalizations} (see \ref{mSCM-marg})  and \emph{interventions} (see \ref{mSCM-intv}) and any combinations of these. This was one of the key challenges in causal modelling in the presence of feedback loops and latent (dependent) variables (see also \cite{BPSM16} for a more general treatment of these issues). Furthermore, these models (and any combination of their marginalizations and interventions) will then have all the properties \hyl{lsSEP}{lsSEP} implies, like the \hyl{smgdGMP}{smgdGMP} encoding conditional independencies based on $\sigma$-separation, etc..
Since we will not %
go into the depth of causality, causal modelling and causal inference in this paper we will refer to other literature for these purposes: 
 \cite{Pearl09}, \cite{SGS00}, \cite{PJS17}, \cite{BPSM16}, \cite{Richardson96}, \cite{Richardson96b}, 
\cite{RiS99}, \cite{EHS10}, \cite{HEH10}, \cite{HEH12}, \cite{HHEJ13}, \cite{foygel2012}, a.o.. \\

We also want to mention explicitely that all the theory about Markov properties and Markov model classes for directed graphs with hyperedges (\HEDG{}es) in this paper also applies to all graphical structures that can be interpreted as \HEDG{}es like directed acyclic graphs (DAGs), acyclic directed mixed graphs (ADMGs), marginalized directed acyclic graphs (mDAGs), directed graphs (DGs), directed mixed graphs (DMGs), etc.. Further note that some (but not all) of the Markov properties only depend on the bidirected edges (instead of the whole generality of the hyperedges) in addition to the directed edges. So in these cases replacing a \HEDG{} $G$ with its induced DMG $G_2$ will not change the Markov property and its corresponding model (class). The rule of thumb is that the ``ordinary'' Markov properties only depend on the induced DMG, i.e.\ its bidirected (and directed) edges, whereas the ``marginal'' versions and the structural equations properties in general do depend on the whole hyperedge structure (in addition to the directed edges). In \ref{eg-dGMP-mdGMP}, figure \ref{fig:eg-dGMP-mdGMP}, we give a simple example (taken from \cite{Eva15}) that shows that hyperedges are necessary.

\subsection{Overview over the Markov properties}

\begin{figure}[!h]
\centering
\scalebox{0.85}{
\begin{tikzpicture}[transform shape]
\newcommand\sk{1.25}
  \node [draw,shape=rectangle,rounded corners] (ausSEP) at (2*\sk,7*\sk) {\hyl{ausSEP}{ausSEP}};
  \node [draw,shape=rectangle,rounded corners] (cusSEP) at (5*\sk,7*\sk) {\hyl{cusSEP}{cusSEP}};
  \node [draw,shape=rectangle,rounded corners] (lusSEP) at (8*\sk,7*\sk) {\hyl{lusSEP}{lusSEP}};
\node [draw,shape=rectangle,rounded corners] (SEP) at (3.5*\sk,6*\sk) {\hyl{SEP}{SEP}};
  \node [draw,shape=rectangle,rounded corners] (asSEP) at (2*\sk,5*\sk) {\hyl{asSEP}{asSEP}};
  \node [draw,shape=rectangle,rounded corners] (csSEP) at (5*\sk,5*\sk) {\hyl{csSEP}{csSEP}};
  \node [draw,shape=rectangle,rounded corners] (lsSEP) at (8*\sk,5*\sk) {\hyl{lsSEP}{lsSEP}};
  \node [draw,shape=rectangle,rounded corners] (mdLMP) at (10.5*\sk,6*\sk) {\hyl{mMP}{mdLMP}};
  \node [draw,shape=rectangle,rounded corners] (smgdGMP) at (6*\sk,3.75*\sk) {\hyl{smgdGMP}{smgdGMP}};
  \node [draw,shape=rectangle,rounded corners] (SEPwared) at (0*\sk,3*\sk) {\hyl{SEPwared}{SEPwared}};
  \node [draw,shape=rectangle,rounded corners] (rFP) at (3.5*\sk,3*\sk) {\hyl{rFP}{rFP}};
  \node [draw,shape=rectangle,rounded corners] (mgdGMP) at (7*\sk,2.5*\sk) {\hyl{mMP}{mgdGMP}};
  \node [draw,shape=rectangle,rounded corners] (mdGMP) at (10.5*\sk,1.75*\sk) {\hyl{mMP}{mdGMP}};
  \node [draw,shape=rectangle,rounded corners] (oLMP) at (5*\sk,1.75*\sk) {\hyl{oLMP}{oLMP}};
  \node [draw,shape=rectangle,rounded corners] (rFPwd) at (2*\sk,1.75*\sk) {\hyl{rFPwd}{rFPwd}};
  \node [draw,shape=rectangle,rounded corners] (gdGMP) at (8*\sk,1*\sk) {\hyl{gdGMP}{gdGMP}};
\node [draw,shape=rectangle,rounded corners] (dLMP) at (3.25*\sk,0.5*\sk) {\hyl{dLMP}{dLMP}};
  \node [draw,shape=rectangle,rounded corners] (aFP) at (1*\sk,-1*\sk) {\hyl{aFP}{aFP}};
  \node [draw,shape=rectangle,rounded corners] (dGMP) at (5*\sk,-1*\sk) {\hyl{dGMP}{dGMP}};
  \node [draw,shape=rectangle,rounded corners] (maFP) at (-1*\sk,1*\sk) {\hyl{mMP}{maFP}};
  \node [draw,shape=rectangle,rounded corners] (mFP) at (-1*\sk,-2*\sk) {\hyl{mFP}{mFP}};
  \node [draw,shape=rectangle,rounded corners] (auGMP) at (3.5*\sk,-3*\sk) {\hyl{auGMP}{auGMP}};
  \node [draw,shape=rectangle,rounded corners] (auLMP) at (3.5*\sk,-5*\sk) {\hyl{auLMP}{auLMP}};
  \node [draw,shape=rectangle,rounded corners] (auPMP) at (3.5*\sk,-7*\sk) {\hyl{auPMP}{auPMP}};
  \node [draw,shape=rectangle,rounded corners] (ruGMP) at (6.5*\sk,-3*\sk) {\hyl{ruGMP}{ruGMP}};
  \node [draw,shape=rectangle,rounded corners] (ruLMP) at (6.5*\sk,-5*\sk) {\hyl{ruGMP}{ruLMP}};
  \node [draw,shape=rectangle,rounded corners] (ruPMP) at (6.5*\sk,-7*\sk) {\hyl{ruPMP}{ruPMP}};
\draw[-implies,double equal sign distance] (ausSEP) to node[left] {\ref{SEP-equiv}} (asSEP);
\draw[->, out=-180, in=90] (ausSEP) to node[left,text width=1.3cm] {discrete RVs\\\ref{ausSEP-wmaFP}} (maFP);
\draw[-implies,double equal sign distance] (cusSEP) to node[above] {\ref{SEP-equiv}} (ausSEP);
\draw[-implies,double equal sign distance] (cusSEP) to node[right] {\ref{SEP-equiv}} (csSEP);
\draw[-implies,double equal sign distance] (lusSEP) to node[above] {\ref{SEP-equiv}} (cusSEP);
\draw[-implies,double equal sign distance] (lusSEP) to node[left] {\ref{SEP-equiv}} (lsSEP);
\draw[->,bend left=65] (SEP) to node[above,align=center] {mDAG\\\ref{SEP-equiv}} (lusSEP);
\draw[->] (mdLMP) to node[right] {} node[above,align=center] {\ref{mDAG-mdLMP-SEP}\\mDAG} (lusSEP);
\draw[->, bend right=30] (mdLMP) to node[left,align=right] {a.q.t.ord.\\\ref{oLMP-dLMP-dGMP}} (mdGMP);
\draw[-implies,double equal sign distance] (asSEP) to node[above] {\ref{SEP-equiv}$\quad$ } (SEP);
\draw[->, bend right] (asSEP) to node[left,text width=1.3cm] {linear\\\ref{linear-SEP}} (SEPwared);
\draw[-implies,double equal sign distance] (csSEP) to node[above] {$\qquad\;$\ref{SEP-equiv}} (asSEP);
\draw[-implies,double equal sign distance] (csSEP) to node[right] {\ref{csSEP-smgdGMP}} (smgdGMP);
\draw[-implies,double equal sign distance] (lsSEP) to node[above] {\ref{SEP-equiv}} (csSEP);
\draw[-implies,double equal sign distance] (lsSEP) to node[above]{\ref{lsSEP-mdLMP}$\;$} (mdLMP);
\draw[-implies,double equal sign distance] (smgdGMP) to node[right] {\ref{smgdGMP-mgdGMP}} (mgdGMP);
\draw[-implies,double equal sign distance] (mgdGMP) to node[right] {\ref{mdGMP-dGMP}} (gdGMP);
\draw[->, bend left=40] (mgdGMP) to node[left,align=center] {mDAG\\\ref{mdGMP-mgdGMP}} (mdGMP);
\draw[-implies,double equal sign distance] (SEPwared) to node[right] {\ref{SEPwared-aFP}} (asSEP);
\draw[-implies,double equal sign distance] (SEPwared) to node[right] {\ref{SEPwared-aFP}} (aFP);
\draw[->] (rFP) to node[left,align=right,yshift=-5mm] {DAG,\\top.ord.\\\ref{sep-lemma}} (SEP);
\draw[implies-implies,double equal sign distance] (rFP) to  node[right] {\;\ref{oLMP-rFP}} (oLMP);
\draw[->, bend right=40] (rFP) to node[below] {} node[left,align=right,xshift=-2mm] {\ref{rFPwd-rFP}\\density} (rFPwd);
\draw[-implies,double equal sign distance] (rFPwd) to node[right] {\;\ref{rFPwd-rFP}} (rFP);
\draw[->] (rFPwd) to node[right,align=left,yshift=-2mm] {\ref{rFPwd-aFP}\\p.e.ord.} (aFP);

\draw[->, bend left=40] (oLMP) to node[right,align=left,yshift=5mm] {\ref{oLMP-dGMP}\\p.e.ord.} (dGMP);

\draw[-implies,double equal sign distance] (mdGMP) to node[right] {\ref{dGMP-dLMP}} (mdLMP);
\draw[-implies,double equal sign distance] (mdGMP) to node[below] {\ref{mdGMP-mgdGMP}} (mgdGMP); %
\draw[-implies,double equal sign distance, bend left] (mdGMP) to node[right] {\;\;\ref{mdGMP-dGMP}} (dGMP); %

\draw[->, bend left=30] (gdGMP) to node[right] {*\ref{dGMP-gdGMP}} (dGMP); %

\draw[<->] (dLMP) to node[left,align=right,yshift=3mm] {\ref{dLMP-oLMP}\\a.p.t.ord.}  (oLMP);

\draw[-implies,double equal sign distance] (aFP) to node[left] {\ref{aFP-dGMP}} (auGMP);

\draw[-implies,double equal sign distance] (dGMP) to node[right] {\;\ref{dGMP-gdGMP}} (gdGMP);
\draw[-implies,double equal sign distance] (dGMP) to node[left] {\ref{dGMP-dLMP}} (dLMP);
\draw[-implies,double equal sign distance] (dGMP) to node[left] {\ref{dGMP-oLMP}} (oLMP);
\draw[->, bend right=7] (dGMP) to node[below] {} node[above] {density, p.e.ord.} (aFP);
\draw[->, bend left=10] (dGMP) to node[below] {\;\;\;density $p>0$} node[above] {\ref{dGMP-aFP}} (aFP);
\draw[-implies,double equal sign distance] (maFP) to node[above,yshift=5mm] {\ref{mdGMP-dGMP}} (aFP);
\draw[-implies,double equal sign distance] (mFP) to node[below] {\;\;\;\;\;\ref{mFP-aFP}} (aFP);
\draw[-implies,double equal sign distance] (mFP) to node[right] {\ref{mFP-maFP}} (maFP);
\draw[->,bend right] (maFP) to  node[right] {\#} (mFP);
\draw[implies-implies,double equal sign distance] (dGMP) to node[right] {\ref{dGMP-auPMP}} (ruGMP);
\draw[implies-implies,double equal sign distance] (dGMP) to node[right] {\ref{dGMP-auPMP}} (auGMP);
\draw[implies-implies,double equal sign distance] (auGMP) to node[above] {\ref{dGMP-auPMP}} (ruGMP);
\draw[implies-implies,double equal sign distance] (ruGMP) to node[left] {\ref{dGMP-auPMP}} (ruLMP);
\draw[implies-implies,double equal sign distance] (ruLMP) to node[left] {\ref{dGMP-auPMP}} (ruPMP);
\draw[-implies,double equal sign distance] (auGMP) to node[right] {\ref{dGMP-auPMP}} (auLMP);
\draw[-implies,double equal sign distance] (ruLMP) to node[above] {\ref{dGMP-auPMP}} (auLMP);
\draw[-implies,double equal sign distance] (auLMP) to node[right] {\ref{dGMP-auPMP}} (auPMP);
\draw[-implies,double equal sign distance] (ruPMP) to node[above] {\ref{dGMP-auPMP}} (auPMP);

\draw[->,out=-180,in=-160] (auPMP) to node[above,align=right,yshift=5mm] {intersection\\\ref{auPMP-dGMP}}  (auGMP);
\draw[->,out=-180,in=-95] (auPMP) to  node[left,align=right] {density\\$p>0$\\\ref{auPMP-aFP}} (aFP);

\end{tikzpicture}}
\caption{Overview over the relations between the Markov properties for directed graphs with hyperedges (\HEDG{}es). Implication arrows always hold under the assumption of the property where the arrow starts. Simple arrows only hold under the indicated assumptions written close to it and are usually stronger statements. The numbers correspond to the lemma/theorem where it was proven in this paper.
The * at gdGMP holds if: $\forall v \in V\,\forall w \in \Sc^G(v): \{v,w\} \in H$, e.g.\ for mDAGs.
The \# at maFP holds under the existence of a positive density $p^\aug >0$ on the augmented space. Note that ``a.p.t.ord.'', ``a.q.t.ord.'', ``p.e.ord.'' stand for assembling pseudo/quasi-topological/perfect elimination order, resp.\ (see definition~\ref{peo-hedg}).}
\label{fig:overview}
\end{figure}
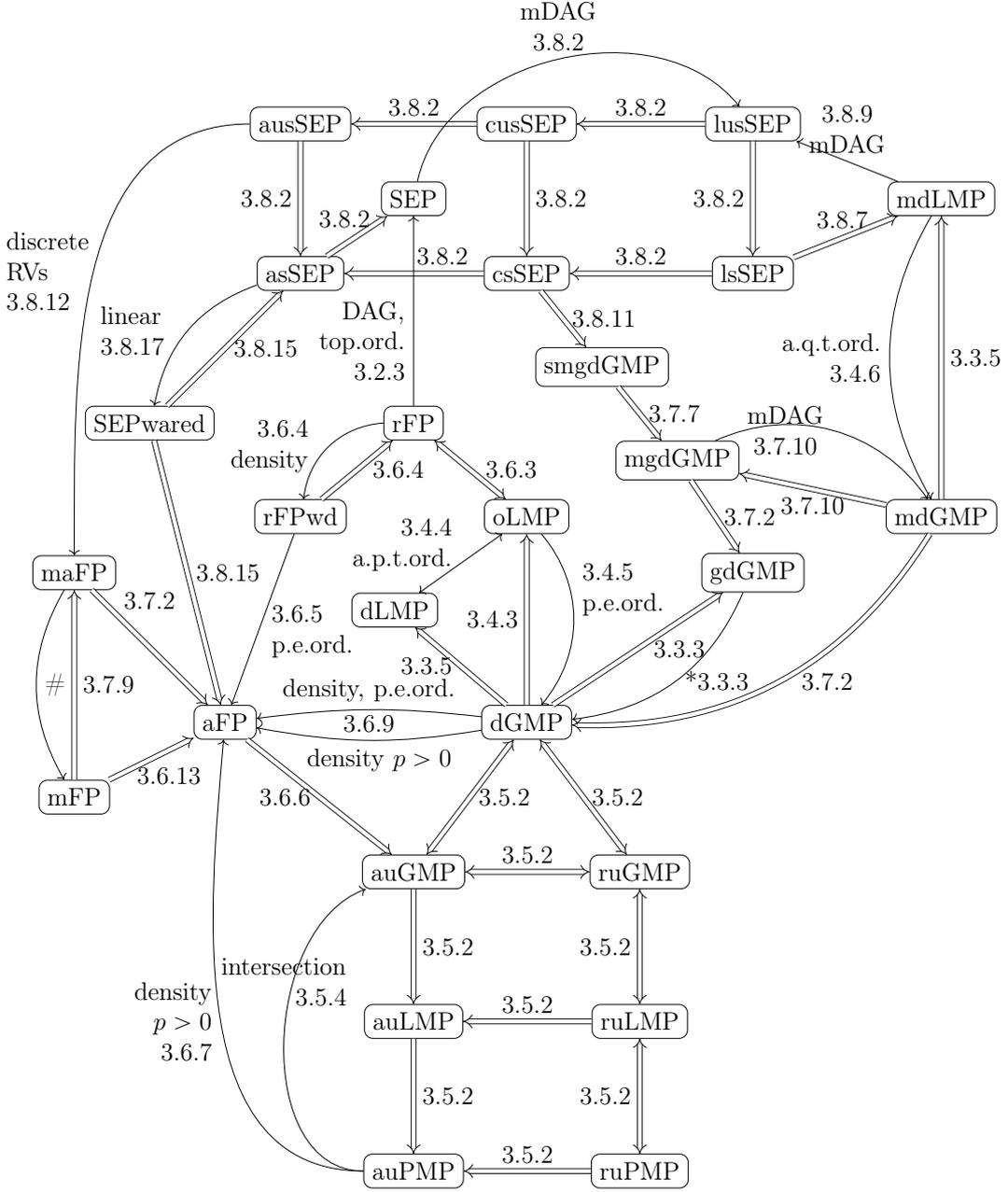

Figure \ref{fig:overview} gives an overview over the Markov properties for directed graphs with hyperedges (\HEDG{}es) and their logical relations between them. These will be analysed in later sections.

\section{Graph Theory}

\subsection{Directed Graphs}

In this subsection we introduce the basic notations needed to understand directed graphs. Of most importance are the notions of \emph{ancestral subgraphs}, \emph{d-separation} and \emph{strongly connected components}
of a directed graph, which are the ``biggest'' directed cycles in there. 

\begin{Def}[Directed Graph]
\begin{enumerate}
\item A \emph{directed graph} is a pair $G=(V,E)$, where $V$ is a finite set of \emph{vertices} or \emph{nodes} and $E$ a subset of $V \x V$ considered as a set of \emph{edges}. We mention that our definition excludes multiple edges between the same vertices, but allows for self-loops of the form $(v,v)$.
\item A directed graph will be represented by a graphical picture. Each node will either be represented by a \emph{dot} or by the name of the node. In a directed graph an edge $(v_1,v_2) \in E$ will be represented as an arrow $v_1 \to v_2$.
\item We will sometimes write $v \in G$ (and $v\to w \in E$ or $v\to w \in G$, resp.) when we mean $v \in V$ (and $(v,w) \in E$, resp.). For a subset of the nodes we will also write $A \ins G$ instead of $A \ins V$.
\item If $v \to w \in E$ then we call $v$ a \emph{parent} of $w$ and $w$ a \emph{child} of $v$ in $G$. 
We define the set of parents of $w$ in $G$ by 
 \[ \Pa^G(w):=\{ v \in V | v \to w \in E\}\]
  and the set of children of $v$ in $G$ by
	\[ \Ch^G(v):=\{w \in V | v \to w \in E\}. \]
\item A \emph{path}\footnote{We explicitly allow for repeated nodes in a path. In some literature the term path refers to ``simple path'',
 i.e.\ without repeated nodes. Paths with repeated nodes are often called ``walks''.} in $G$ is a sequence of $n$ nodes from $V$, $n \ge 1$, and $n-1$ edges from $E$:
    \[ v_1 \ {\to \atop \ot}\ v_2 \ {\to \atop \ot}\  \cdots \ {\to \atop \ot}\ v_{n-1} \ {\to \atop \ot}\  v_n,  \]
where at every position one of the directed edges is chosen. 
		The nodes $v_1$ and $v_n$ are called \emph{endnodes}.
    The definition of path also allows for the case of a \emph{trivial path} consisting only of a single node $v_1$, and also for the case of repeated nodes in the path. 
\item A \emph{directed path} from $v$ to $w$ in $G$ is a path in $G$ of the form:
\[ v \to v_2 \to \cdots \to w,  \] 
where all arrowheads point in the same direction.
\item If there is a directed path from $v$ to $w$ in $G$, $v$ is called an \emph{ancestor} of $w$ in $G$ and $w$ a \emph{descendent} of $w$.
 The set of ancestors of $w$ in $G$ will be denoted by
\[ \Anc^G(w):= \{ v \in V | \exists v \to \cdots \to w \text{ a directed path in } G\}\]
and the set of descendents of $v$ in $G$ by
\[ \Desc^G(v):= \{ w \in V | \exists v \to \cdots \to w \text{ a directed path in } G\}.\]
Note that we have:
\[\{w\} \cup \Pa^G(w) \ins \Anc^G(w)\qquad \text{and}\qquad \{v\} \cup \Ch^G(v) \ins \Desc^G(v).\] 
\item For a subset of the nodes $A \ins V$ of $G$ we define 
\[\begin{array}{lllllll}
\Pa^G(A)&:=&\bigcup_{v \in A} \Pa^G(v), && \Ch^G(A)&:=&\bigcup_{v \in A}\Ch^G(v),\\
\Anc^G(A)&:=&\bigcup_{v \in A}\Anc^G(v),&& \Desc^G(A)&:=&\bigcup_{v \in A}\Desc^G(v). \\
\end{array}\]
\item Furthermore, we define the set of \emph{non-descendents} by 
\[\NonDesc^G(A):= V \sm \Desc^G(A).\]
\item For $v \in V$ we call $\Sc^G(v):=\Anc^G(v) \cap \Desc^G(v)$ the \emph{strongly connected component}\footnote{\label{fn:str-con-comp}In \cite{Spirtes94} the strongly connected components were called \emph{cycle groups}.} %
 of $v$ in $G$ (see \ref{str-conn-comp-def}).
\item A graph $G'=(V',E')$ is called a \emph{subgraph} of $G=(V,E)$ if $V' \ins V$ and $E' \ins E$.
\item If $A \ins V$ is a subset of the nodes of $G$, then $A$ induces a \emph{subgraph} of $G$ given by $G(A):=(V(A),E(A))$,
where $V(A):=A$ and $E(A):= E \cap (A \x A)$, i.e.\ we take all edges with both endpoints in $A$.
We will often just write $A$ instead of $G(A)$ to adress the subgraph.
\item A subgraph $A \ins G$ is called an \emph{ancestral subgraph} in $G$ if $A=\Anc^G(A)$, i.e.\ if $A$ is closed under taking ancestors of $A$ in $G$.
\item If there is no \emph{cycle}, i.e., non-trivial cyclic directed path $v \to \cdots \to v$, in $G$ (i.e.\ with the same endnodes) then $G$ is called a \emph{directed acyclic graph (DAG)}.
\end{enumerate}
\end{Def}

\begin{Def}[Strongly connected components]
\label{str-conn-comp-def}
Let $G=(V,E)$ be a directed graph.
\begin{enumerate}
\item For $v,w \in V$ we define the relation:
 \[ v \sim w :\iff w \in \Anc^G(v) \cap \Desc^G(v) = \Sc^G(v).  \]
This is reflexive, symmetric and transitive, and thus defines an equivalence relation on $V$.
The equivalence classes are exactly of the form $\Sc^G(v)=\Anc^G(v) \cap \Desc^G(v)$ and are called \emph{the strongly connected components of $G$}.
    Let $V/{\sim}{}=\{ \Anc^G(v) \cap \Desc^G(v)| v \in V\}$ be the set of strongly connected components.
\item For edges $v \to w$ and $v' \to w'$ in $E$ we define the relation:
 \[  (v \to w) \sim (v'\to w') :\iff v \sim v' \;\&\; w \sim w'. \]
This defines an equivalence relation on $E$. $E/{\sim}$ can be considered as a subset of $(V/{\sim}) \x (V/{\sim})$.
    Let $$E'/{\sim} := (E/{\sim}) \ \setminus\ \{v' \to v' : v' \in V / {\sim}\}$$
		be $E/{\sim}$  
		without the diagonal	(i.e.\ adjusted for selfloops). %
\item The pair $\Scal(G):=(V/{\sim}, E'/{\sim})$ defines a directed graph and is called the \emph{graph of strongly connected components of $G$}. We will refer to the set of strongly connected components (and its additional directed graph structure) with $\Scal(G)$. %
\end{enumerate}
\end{Def}

\begin{Lem}[DAG of strongly connected components]
\label{sc-dag}
Let $G=(V,E)$ be a directed graph.
Then the directed graph of strongly connected components $\Scal(G)$ of $G$ is a DAG.
\begin{proof}
Every cycle in $\Scal(G)$ would lift to a cycle in $G$, which by definition would lie in one strongly connected component of $G$. 
But since we erased the selfloops in $E'/\!\!\!\sim$ this is a contradiction. Therefore, no cycles in $\Scal(G)$ exist.
\end{proof}
\end{Lem}

\begin{Def}[d-separation in a directed graph]
Let $G=(V,E)$ be a directed graph and $X,Y,Z \ins V$ subsets of the nodes.
\begin{enumerate}
\item Consider a path in $G$: $v_1 \genfrac{}{}{0pt}{}{\to}{\ot}  \cdots \genfrac{}{}{0pt}{}{\to}{\ot} v_n$ with $n \ge 1$.
It will be called \emph{$Z$-blocked} or \emph{blocked by $Z$} if:
 \begin{enumerate}
   \item at least one of the endnodes $v_1$ or $v_n$ is in $Z$, or
	 \item there are two adjacent edges in the path of one of the following forms: 
	 \[ \begin{array}{lllll}
	  \text{ a non-collider: } & \to v_i \to & \text{ with } & v_i \in Z, & \text{ or }\\
		\text{ a non-collider: } & \ot v_i \ot & \text{ with } & v_i \in Z, & \text{ or }\\
		\text{ a non-collider: } & \ot v_i \to & \text{ with } & v_i \in Z, & \text{ or }\\
		\text{ a collider: } &\to v_i \ot & \text{ with } & v_i \notin \Anc^G(Z).
	\end{array}\]
 \end{enumerate}
If none of the above holds then the path is called \emph{$Z$-open} or \emph{$Z$-active}.\footnote{As long as we allow paths to have repeated nodes we can change the collider case to only check if $v_i$ is in $Z$ instead of $\Anc^G(Z)$.}
\item  We say that $X$ is \emph{d-separated}\footnote{This was introduced by J. Pearl, see \cite{Lau98} \S3.2.2, \cite{Pearl86b} Def. 1.2, and was further developed by T. Verma and J. Pearl in \cite{VerPea90a}, \cite{VerPea91}. } 
from $Y$ given $Z$ if every path in $G$ with one endnode in $X$ and the other endnode in $Y$ is blocked by $Z$.
In symbols this will be written as follows:
\[ X \Indep_G^d Y \given Z.\]
\end{enumerate}
\end{Def}

\begin{Rem}
For checking the d-separation $X \Indep_G^d Y \given Z$ it is enough to only check paths where every node occurs at most once. If a node occured twice in a $Z$-open path then there would be a $Z$-open ``short-cut''. 
\end{Rem}

\subsection{Directed Graphs with Hyperedges (\HEDG{}es)}

In the following we will introduce directed graphs with hyperedges (\HEDG{}es). This combines and generalizes the notions of directed graphs and marginalized directed acyclic graphs (mDAG) from \cite{Eva15} (Defs.\ 8 and 9) to the general case that allows for both cycles and hyperedges, which represent unobserved latent structure. In particular, \HEDG{}es also generalize acyclic directed mixed graphs (ADMGs) from e.g.\ \cite{Richardson03}. Of most importance in this subsection is to understand how \emph{hyperedges} are modelled, the notions of \emph{districts}, \emph{marginalizations} of \HEDG{}es and the generalized \emph{d-separation} definition.

\begin{Def}[Directed graphs with hyperedges (\HEDG{}es)]$ $
\begin{enumerate}
\item Let $V$ be a finite set. A \emph{simplicial complex} $H$ over $V$ is a set of %
  subsets of $V$ such that:
 \begin{enumerate}
 \item[(i)] all single element sets $\{v\}$ are in $H$ for $v \in V$, and
 \item[(ii)] if $F \in H$ then also all subsets $F' \ins F$ are elements of $H$.
\end{enumerate}
\item A \emph{directed graph with hyperedges} or %
  \emph{hyperedged directed graph}\footnote{We could also call it \emph{marginalized directed graph}.}
(\emph{\HEDG{}})  is a tuple $G=(V,E,H)$, where $(V,E)$ is a directed graph (with or without cycles) and  $H$ a simplicial complex over the set of vertices $V$ of $G$. 
\item If $G=(V,E,H)$ is a \HEDG{} then the elements $F$ of $H$ will be called the %
    \emph{hyperedges} of $G$. The subset of inclusion-maximal elements of $H$ will be abbreviated with $\tilde{H}$. The elements $F \in \tilde{H}$ are called \emph{maximal hyperedges} of $G$ or again just \emph{hyperedges} if the context is clear.
  \item We can represent a \HEDG{} $G=(V,E,H)$ as an ordinary directed graph $(V,E)$ consisting of vertices $V$ and directed edges $E$, with additional hyperedges corresponding to the maximal hyperedges $F \in \tilde H$ with $\# F \ge 2$ pointing to their target nodes $v \in F$. See figure \ref{fig:hedg}.
\item If for two distinct vertices $v,w \in V$ we have $\{v,w\} \in H$ then we will often represent this with a \emph{bidirected edge} $v \oto w$ and write $v \oto w \in H$ or $v \oto w \in G$.
\end{enumerate}
\end{Def}

\begin{Def}
Let $G=(V,E,H)$ be a \HEDG{}.
\begin{enumerate}
\item A \emph{path} in $G$ is a sequence of $n$ nodes\footnote{We again allow for repeated nodes in the definition of a path.} with $n \ge 1$, and $n-1$ (directed or bidirected) edges:
\[\begin{array}{ccccccc} 
& \ot &&\ot \\[-5pt]
v_1 & \to & \cdots & \to &v_n, \\[-5pt]
& \oto & &\oto 
\end{array}\]
where at every position one of the (directed or bidirected) edges is chosen. 
\item A \emph{directed path} from $v$ to $w$ in $G$ is a path in $G$ of the form:
\[ v \to  \cdots \to w,  \] 
i.e., with only directed edges such that all arrowheads point in the same direction.
\item A \emph{bidirected path} in $G$ is a path of the form $v \oto \cdots \oto w$, only consisting of bidirected edges.
\item The \emph{district}\footnote{A district was also called \emph{C-component} e.g.\ in \cite{Tian2002} or \cite{ShP06}.} of $v$ is the set:
\[\Dist^G(v):=\{ w \in V | w=v \text{ or } \exists w \oto \cdots \oto v \text{ a bidirected path in } G\}.\]
\end{enumerate}
\end{Def}

\begin{Def}[Sub-\HEDG{}es]
Let $G=(V,E,H)$ be a \HEDG{}.
\begin{enumerate}
\item A \HEDG{} $G'=(V',E',H')$ is called a \emph{sub-\HEDG{}} of $G=(V,E,H)$ if $V' \ins V$, $E' \ins E$ and $H' \ins H$.
\item Let $A \ins V$ be a subset. Then $A$ induces a \emph{sub-\HEDG{}} of $G$ as follows: $G(A):=(V(A),E(A),H(A))$ with $V(A):=A$, $E(A):= E \cap(A\x A)$ and $H(A):=\{ F \in H | F \ins A\} = \{ F \cap A | F \in H\}$. We will often write $A$ instead of $G(A)$.
\end{enumerate}
\end{Def}

\begin{Def}[mDAG, DMG, ADMG]
Let $G=(V,E,H)$ be a \HEDG{}.
\begin{enumerate}
\item $G=(V,E,H)$ is called a \emph{marginalized directed acyclic graph (mDAG)}\footnote{See \cite{Eva15} Def. 8, 9.} if the tuple $(V,E)$ is a directed acyclic graph (DAG).
\item A \HEDG{} $G=(V,E,H)$ with $\#F \le 2$ for all $F \in H$ is called \emph{directed mixed graph} or \emph{DMG}. It basically consists only of nodes, directed edges ($\to$) and bidirected edges ($\oto$).
\item If a \HEDG{} $G$ is both an mDAG and a directed mixed graph then it is also called \emph{acyclic directed mixed graph} or \emph{ADMG}.
\item If $G=(V,E,H)$ is a \HEDG{} we can define the \emph{induced directed mixed graph} by $G_2:=(V,E,H_2)$, where $H_2:=\{ F \in H | \# F \le 2\}$. Note that if $G$ is an mDAG then the induced directed mixed graph is an ADMG.
\end{enumerate}
\end{Def}

\begin{Def}[Strongly connected components]
Let $G=(V,E,H)$ be a \HEDG{}.
\begin{enumerate}
\item For $v \in V$ the \emph{strongly connected component} of $v$, abbreviated as $\Sc^G(v)$, is the strongly connected component of $v$ according to the underlying directed graph $(V,E)$ defined in \ref{str-conn-comp-def}. 
\item The \emph{mDAG of strongly connected components} $\Scal(G)$ of $G$ is given by the DAG structure from \ref{sc-dag} plus the additional set of hyperedges $H(\Scal(G))$ given by: $\{S_1,\dots,S_r\} \in H(\Scal(G))$ if and only if there are elements $v_i \in S_i$, $i=1,\dots,r$, with $\{v_1,\dots,v_r\} \in H$. 
\end{enumerate}
\end{Def}

\begin{Def}[d-separation]
\label{d-sep-hedg}
Let $G=(V,E,H)$ be a \HEDG{} and $X,Y,Z \ins V$ subsets of the nodes.
\begin{enumerate}
\item Consider a path in $G$ with $n \ge 1$ nodes:
\[\begin{array}{ccccccc} 
& \ot &&\ot \\[-5pt]
v_1 & \to & \cdots & \to &v_n. \\[-5pt]
& \oto & &\oto 
\end{array}\]
The path will be called \emph{$Z$-blocked} or \emph{blocked by $Z$} if:
 \begin{enumerate}
   \item at least one of the endnodes $v_1$, $v_n$ is in $Z$, or
	 \item there is a node $v_i$ with two adjacent (hyper)edges in the path of one of the following forms:
	 \begin{enumerate}
     \item $v_i \notin 	\Anc^G(Z)$
 and the two (hyper)edges are a \emph{collider}\footnote{Again, as long as we allow a path to have repeated notes we can restrict to $Z$ instead of $\Anc^G(Z)$ here.} at $v_i$, i.e.\ each has an arrowhead towards $v_i$:
 \[\to v_i \ot , \quad  \oto v_i \ot, \quad \to v_i \oto, \quad \oto v_i \oto. \]
 \item $v_i \in Z$ and the two (hyper)edges are a \emph{non-collider} at $v_i$, i.e.\ at least one of the (hyper)edges has no arrowhead towards $v_i$:
\[  \to v_i \to, \quad \oto v_i \to, \quad \ot v_i \ot, \quad \ot v_i \oto, \quad \ot v_i \to.\]
	 \end{enumerate}
 \end{enumerate}
If none of the above holds then the path is called \emph{$Z$-open} or \emph{$Z$-active}.\\
Note that for $Z=\emptyset$ colliders are always $Z$-blocked and non-colliders are always $Z$-open.
\item  We say that $X$ is \emph{d-separated}\footnote{Some authors use the term m-separation for ADMGs instead, e.g.\ in \cite{Richardson03} \S2.1. But since we want to avoid introducing a new letter for every graph type and this definition naturally generalizes Pearl's d-separation for directed graphs by \ref{hedg-d-sep-rem} we will stick with ``d-separation''.}
 from $Y$ given $Z$ if every path in $G$ with one endnode in $X$ and one endnode in $Y$ is blocked by $Z$.
In symbols this will be written as follows:
\[ X \Indep_G^d Y \given Z.\]
\end{enumerate}
\end{Def}

\begin{Rem}
\label{hedg-d-sep-rem}
\begin{enumerate}
  \item For checking the d-separation $X \Indep_G^d Y \given Z$ it is enough to only check paths where every node occurs at most once. If a node occured twice in a $Z$-open path then there would be a $Z$-open ``short-cut''.
\item If we identify a directed graph $(V,E)$ with the \HEDG{} given by $(V,E,H_1)$, where $H_1$ is trivial (i.e.\ only consists of sets $F$ with $\#F \le 1$), then the two definitions of d-separation are equivalent:
\[ X \Indep_{(V,E,H_1)}^d Y \given Z \;\iff\; X \Indep_{(V,E)}^d Y \given Z. \]
\item Since d-separation only uses bidirected edges (rather than the whole hyperedges), we can check it by only using the induced directed mixed graph, i.e.\ we have the equivalence:
\[ X \Indep_{(V,E,H)}^d Y \given Z \;\iff\; X \Indep_{(V,E,H_2)}^d Y \given Z. \]
\end{enumerate}
\end{Rem}

\begin{figure}
\centering
\begin{tikzpicture}[scale=.7, transform shape]
\tikzstyle{every node} = [draw,shape=circle]
\node (v1) at (0,0) {$v_1$};
\node (v2) at (0,-2) {$v_2$};
\node (v3) at (-1.5,-4) {$v_3$};
\node (v4) at (1.5,-4) {$v_4$};
\node (v5) at (3.5,-4) {$v_5$};
\node (v6) at (3.5,-2) {$v_6$};
\node (v7) at (2,0) {$v_7$};
\node[fill,circle,red,inner sep=0pt,minimum size=5pt] (v8) at (2,-2) {};
\draw[-{Latex[length=3mm, width=2mm]}, bend right] (v1) to (v2);
\draw[-{Latex[length=3mm, width=2mm]}, bend right] (v2) to (v1);
\foreach \from/\to in {v2/v3, v2/v4, v7/v6, v7/v1, v4/v5, v6/v5}
\draw[-{Latex[length=3mm, width=2mm]}] (\from) -- (\to);
\foreach \from/\to in {v8/v7, v8/v4, v8/v5}
\draw[-{Latex[length=3mm, width=2mm]}, red] (\from) -- (\to);
\end{tikzpicture}
\hspace{1cm}
\begin{tikzpicture}[scale=.7, transform shape]
\tikzstyle{every node} = [draw,shape=circle]
\node (v1) at (0,0) {$v_1$};
\node (v3) at (-1.5,-4) {$v_3$};
\node (v4) at (1.5,-4) {$v_4$};
\node (v6) at (3.5,-2) {$v_6$};
\node[fill,circle,red,inner sep=0pt,minimum size=5pt] (v8) at (2,-2) {};
\draw[{Latex[length=3mm, width=2mm]}-{Latex[length=3mm, width=2mm]}, red, bend right] (v4) to node[fill,circle,red,inner sep=0pt,minimum size=5pt] {} (v3);
\draw[-{Latex[length=3mm, width=2mm]},out=135, in=225,looseness=8] (v1) to (v1);
\foreach \from/\to in {v1/v3, v1/v4}
\draw[-{Latex[length=3mm, width=2mm]}] (\from) -- (\to);
\foreach \from/\to in {v8/v1, v8/v4, v8/v6}
\draw[-{Latex[length=3mm, width=2mm]}, red] (\from) -- (\to);
\end{tikzpicture}
\caption{A \HEDG{} on the left and its marginalization w.r.t.\ $\{v_2,v_5,v_7\}$ on the right.}
 \label{fig:marginalization}
\end{figure}
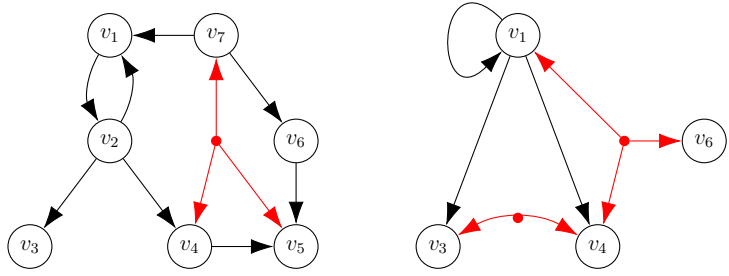

\begin{Def}[Marginalization of a \HEDG{}]
\label{HEDG-marg-def}
Let $G=(V,E,H)$ be a \HEDG{} and $U \ins V$ a subset.
The \emph{marginalization of $G$ with respect to $U$}  or \emph{latent projection}\footnote{This term was used in \cite{Verma93} 3.1.3, \cite{Eva15} Def. 11, 12 and others.} \emph{of $G$ onto $V\sm U$} is the \HEDG{} $G^\marg=(V^\marg,E^\marg, H^\marg)$, where
\begin{enumerate}
\item $V^\marg:=V \sm U$,
\item $v_1\to v_2 \in E^\marg$ iff there exists a directed path in $G$ of the form $v_1 \to u_1 \to \cdots u_r \to v_2$ with $r \ge 0$ and all $u_i \in U$.
\item $F' \ins V \sm U$ is in $H^\marg$ iff there exists an $F \in H$ with $F \ins F' \cup U$ such that
 for every $v \in F'$ we have:
\begin{enumerate}
  \item $v \in F \sm U$, or 
 \item there exists a directed path in $G$ of the form $u_1  \to \cdots \to u_r \to v$ with $r \ge 1$ and all $u_i \in U$
and $u_1 \in F \cap U$. 
\end{enumerate}
So for $F \in H$ we construct $F^\marg \in H^\marg$ by adding all $v$ from (b) and subtracting $U$.
    Every $F' \in \widetilde{H^\marg}$ 
		then is of the form (non-uniquely i.g.): $F'=F^\marg$ with $F \in \tilde{H}$.
\end{enumerate}
We also write $G^{\marg(V \sm U)}$ and $G^{\marg\sm U}$, resp., instead of just $G^\marg$ to emphasize the nodes that are left, or that are marginalized out, resp..
\end{Def}
See Figure~\ref{fig:marginalization} for an example of a marginalization of a \HEDG{}.
\begin{Eg}
Let $G=(V,E,H)$ be a \HEDG{} and $U=\{u\}$ then $G^{\marg\sm U}=(V^\marg,E^\marg,H^\marg)$ is given by:
\begin{enumerate}
\item $V^\marg = V\sm \{u\}$.
\item $v_1 \to v_2 \in E^\marg$ iff $v_1 \to v_2 \in E$ or both edges $v_1 \to u \to v_2$ are in $E$.
\item For $F \in H$ we have: 
\[F^\marg = \left \{ 
\begin{array}{lcl}
 F, &\text { if } &u \notin F,\\
(F \cup \Ch^G(u)) \sm \{u\},& \text{ if }& u \in F.
\end{array}
\right. \]
\item If $F' \in \widetilde{H^\marg}$ then either $F' \in \tilde H$ or there is $F \in \tilde H$ with $(F'\sm \Ch^G(u)) \cup \{u\} \ins F$.
\item $v_1 \oto v_2 \in H^\marg$ iff $v_1 \oto v_2$, $v_1 \oto u \to v_2$, $v_1\ot u \oto v_2$  or $v_1 \ot u \to v_2$ is in $G$.
\end{enumerate}
\end{Eg}

\begin{Rem}
If we are only interested in the marginalization of a directed mixed graph $G=(V,E,H)$ as a directed mixed graph, i.e.\ with only bidirected edges, we can just use
the induced directed mixed graph of the marginalization of $G$ as a \HEDG{}. So we just take the marginalization of $G$ and then restrict to bidirected edges there.
\end{Rem}

\begin{Rem}[Marginalization preserves ancestral relations]
Let $G=(V,E,H)$ be a \HEDG{} and $G^\marg=(V^\marg,E^\marg,H^\marg)$ a marginalization of $G$.
Let $v,w \in V^\marg \ins V$ be two nodes. Then by the definition of marginalization we have:
\[ v \in \Anc^G(w) \; \iff \; v \in \Anc^{G^\marg}(w).  \]
\end{Rem}

\begin{Lem}
\label{anc-sub=mar}
If $G=(V,E,H)$ is a \HEDG{} and $W \ins V$ a subset, then the set $W$ might carry two different structures as a \HEDG{} induced by $G$:
\begin{enumerate}
\item the sub-\HEDG{} structure $G(W)=(W,E(W),H(W))$ of $G$,
\item the marginalized-\HEDG{} structure $G^{\marg(W)}=(W, E^{\marg(W)}, H^{\marg(W)})$.
\end{enumerate} 
These structures do not agree in general, but if
$W$ is an ancestral subset of $G$ then both structures are equal.
\begin{proof}
If $W$ is ancestral, $U:=V \sm W$ and $v_1 \to v_2$ is in $E^{\marg(W)}$, then this edge comes from a directed path
$v_1 \to u_1 \to \dots u_r \to v_2$ in $G$ with $r \ge 0$ and $u_i \in U$. Since $W$ is ancestral and $v_2 \in W$, also $u_i \in W$. But $W \cap U = \emptyset$.
So $r=0$ and $v_1 \to v_2$ already exists in $G$.
If $F \in H^{\marg(W)}$ then there is a $F' \in H$ such that for every $v \in F$ we either have $v \in F' \cap W$ or there is a path in $G$ of the form $u_1 \to u_2 \to \cdots \to v$ with $u_i \in U$ and $u_1 \in F' \cap U$. But again $W$ is ancestral and thus all $u_i \in W \cap U = \emptyset$. So only $F=F'\cap W$ can hold.
So $G^{\marg(W)}$ is a sub-\HEDG{} of $G$ that contains all (hyper)edges of $G$ between nodes in $W$, thus agrees with the sub-\HEDG{}-structure.
\end{proof}
\end{Lem}

\begin{Lem}
\label{HEDG-marg-commute}
Let $G=(V,E,H)$ be a \HEDG{} and $U_1, U_2 \ins V$ be two disjoint subsets of the nodes.
Then the marginalizations w.r.t.\ $U_1$ and $U_2$ commute, i.e.\ we have:
\[ (G^{\marg\sm U_1})^{\marg\sm U_2} = G^{\marg\sm(U_1 \cup U_2)} = (G^{\marg\sm U_2})^{\marg\sm U_1} . \]
\begin{proof}
  The same proof as in \cite{Eva15} Thm. 1 works. 
First all constructions have $V \sm (U_1 \cup U_2)$ as the underlying set of nodes.
By induction by $\#(U_1 \cup U_2)$ we can assume that $U_1=\{u_1\}$, $U_2=\{u_2\}$ and $u_1 \neq u_2$. It is then enough to consider directed paths that contain both $u_1$ and $u_2$. For these the directed edges and hyperedges on all constructions are directly checked to agree.
\end{proof}
\end{Lem}

\begin{Lem}[Ancestral sub-\HEDG{}es in marginalizations]
\label{hedg-anc-marg}
Let $G=(V,E,H)$ be a \HEDG{} and $W \ins V$ a subset. We endow $W$ with the marginalized \HEDG{}-structure $G^{\marg(W)}$ of $G$.
\begin{enumerate}
\item If $A \ins G$ is an ancestral sub-\HEDG{} of $G$ then $G^{\marg(A \cap W)}$ is an ancestral sub-\HEDG{} of $G^{\marg(W)}$.
\item If $B \ins G^{\marg(W)}$ is an ancestral sub-\HEDG{} of $G^{\marg(W)}$ then $A:=\Anc^G(B)$ is an ancestral sub-\HEDG{} of $G$ with $G^{\marg(A \cap W)}=B$.
\end{enumerate}
\begin{proof}
1. If $v \in A \cap W$ and $w \in \Anc^{G^{\marg(W)}}(v)$ then $w \in \Anc^G(v)$ because marginalization preserves ancestral relations by construction. Since $A$ is ancestral in $G$ we have $w \in A$ (and in $W$) and thus $w \in A \cap W$. 
  It follows that $A \cap W$ forms an ancestral sub-\HEDG{} of $G^{\marg(W)}$. By \ref{anc-sub=mar} we have that the sub-\HEDG{} structure of $A \cap W$
equals $(G^{\marg(W)})^{\marg(A\cap W)}$. By \ref{HEDG-marg-commute} this is the same as $G^{\marg{(A\cap W)}}$.\\
2. We clearly have that the nodes of $B$ lie in $A\cap W$. By the first point $G^{\marg(A \cap W)}$ is an ancestral sub-\HEDG{} of $G^{\marg(W)}$.
It follows that $B$ is an ancestral sub-\HEDG{} of $G^{\marg(A \cap W)}$ then.
Now let $v \in A \cap W$ then $v \in A=\Anc^G(B)$ and thus by marginalization $v \in \Anc^{G^{\marg(A\cap W)}}(B)$. But since $B$ is ancestral in 
$G^{\marg(A\cap W)}$ we have $v \in B$. It follows that $A \cap W \ins B$ and thus $B=G^{\marg(A \cap W)}$ (as \HEDG{}es).
\end{proof}
\end{Lem}

\begin{Lem}[Marginalization and d-separation]
\label{marg-d-sep}
Let $G=(V,E,H)$ be a \HEDG{} and $X,Y,Z,U \ins V$ subsets with $U \cap ( X \cup Y \cup Z) =\emptyset$ and $G^{\marg\sm U}$ the marginalization with respect to $U$.
Then we have the equivalence:
\[ X \Indep_G^d Y \given Z \iff X \Indep_{ G^{\marg\sm U} }^d Y \given Z.\]
\begin{proof}
By \ref{HEDG-marg-commute} and induction by $\# U$ we can reduce to the case where $U=\{u\}$ consists of a single node.
First consider a $Z$-open path in $G^{\marg\sm U}$. Every edge $v \to w$ in $G^{\marg\sm U}$ comes from a edge $v\to w$ or $v \to u \to w$ in $G$ and every bidirected edge $v \oto w$ in $G^{\marg\sm U}$ (take $F'=\{v,w\}$) comes from $v \oto w$, $v \ot u \to w$, $v \oto u \to w$ or $v \ot u \oto w$ in $G$, where $u$ is never a collider and not in $Z$.
So the $Z$-open path in $G^{\marg\sm U}$ can be lifted to a $Z$-open path in $G$. \\
Now consider a $Z$-open path between $X$ and $Y$ in $G$.
If this path does not contain a node in $U$ then it exists also in $G^{\marg\sm U}$ and is then also $Z$-open.
If the path contains a node  $u \in U$ it cannot be an endnode by the $Z$-openness.
If $u$ in the path is a collider ($\genfrac{}{}{0pt}{}{\to}{\oto} u \genfrac{}{}{0pt}{}{\ot}{\oto}$) it must be in $\Anc^G(Z)$ (by the $Z$-openness) but not in $Z$ (by the assumption: $U \cap Z =\emptyset$).
So there is a directed path $u \to v_1 \to \cdots \to v_r$ in $G$ with $r \ge 1$ and $v_r \in Z$. Then also $v_1 \in \Anc^G(Z)$ and the $Z$-open collider marginalizes to the $Z$-open collider $\genfrac{}{}{0pt}{}{\to}{\oto} v_1 \genfrac{}{}{0pt}{}{\ot}{\oto}$ in $G^{\marg\sm U}$.\\
If $u$ in the path is a non-collider
it is of one of the forms (or mirrored):
\[v\to u \to w,\qquad v\ot u \to w,\qquad v \oto u \to w.\]
Marginalizing out $u$ gives $v \to w$, $ v \oto w $, $v \oto w$, resp., for these cases, where head and tail of the arrows at $v$ and $w$ stay the same. So in all cases the marginalized path stays $Z$-open.
\end{proof}
\end{Lem}

\begin{Lem}
\label{d-sep-anc-sub-hedg}
Let $G=(V,E,H)$ be a \HEDG{} with subsets $X,Y,Z \ins V$, $A$ an ancestral sub-\HEDG{} of $G$ with $X,Y,Z \ins A$ 
(e.g.\ $A=\Anc^{G}(X \cup Y \cup Z)$).
Then we have:
\[ X \Indep_G^d Y \given Z \iff X \Indep_{A}^d Y \given Z.\]
\begin{proof}
This follows from the fact \ref{anc-sub=mar} that on ancestral sub-\HEDG{}es $A \ins G$ the sub-\HEDG{}-structure $G(A)$ and the marginalized \HEDG{}-structure 
$G^{\marg(A)}$ are the same and that d-separation is stable under marginalization by \ref{marg-d-sep}. 
\end{proof}
\end{Lem}

\subsection{The Augmented Graph of a \HEDG{}}

This short subsection is dedicated to the definition of the \emph{augmented graph} of a \HEDG{}. It is constructed from the \HEDG{} by replacing the (origins of the) hyperedges with explicit nodes, which then give a directed graph (without hyperedges).

\begin{Def}[Augmented directed graph of a \HEDG{}]
\label{augm-def}
Let $G=(V,E,H)$ be \HEDG{}.
The \emph{augmentation} or \emph{augmented (directed) graph}\footnote{In \cite{Richardson03} \S2.2 this terms was used for a different construction for ADMGs, which we later in \ref{moral-hedg-def} will call \emph{(generalized) moralization}.} 
of $G$ is the directed graph $G^\aug=(V^\aug,E^\aug)$, 
where $V^\aug$ is composed of $V$ and the additional nodes $e_F$ for $F \in \tilde{H}$. For $E^\aug$ we take all the edges from $E$ and add the directed edges $e_F \to v$ for $v \in F$.
\end{Def}

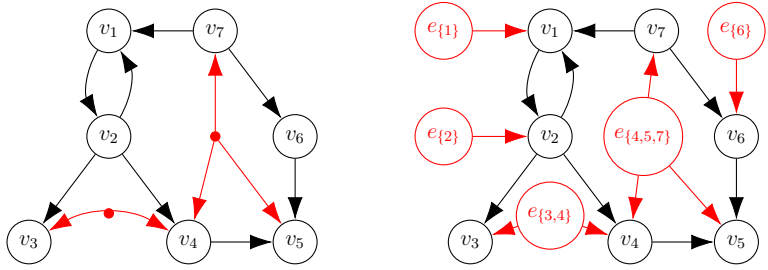
\begin{figure}
\centering
\begin{tikzpicture}[scale=.7, transform shape]
\tikzstyle{every node} = [draw,shape=circle]
\node (v1) at (0,0) {$v_1$};
\node (v2) at (0,-2) {$v_2$};
\node (v3) at (-1.5,-4) {$v_3$};
\node (v4) at (1.5,-4) {$v_4$};
\node (v5) at (3.5,-4) {$v_5$};
\node (v6) at (3.5,-2) {$v_6$};
\node (v7) at (2,0) {$v_7$};
\node[fill,circle,red,inner sep=0pt,minimum size=5pt] (e457) at (2,-2) {};
\draw[-{Latex[length=3mm, width=2mm]}, bend right] (v1) to (v2);
\draw[-{Latex[length=3mm, width=2mm]}, bend right] (v2) to (v1);
\draw[{Latex[length=3mm, width=2mm]}-{Latex[length=3mm, width=2mm]}, red, bend right] (v4) to node[fill,circle,red,inner sep=0pt,minimum size=5pt] {} (v3);
\foreach \from/\to in {v2/v3, v2/v4, v7/v6, v7/v1, v4/v5, v6/v5}
\draw[-{Latex[length=3mm, width=2mm]}] (\from) -- (\to);
\foreach \from/\to in {e457/v7, e457/v4, e457/v5}
\draw[-{Latex[length=3mm, width=2mm]}, red] (\from) -- (\to);
\end{tikzpicture}
\hspace{1cm}
\begin{tikzpicture}[scale=.7, transform shape]
\tikzstyle{every node} = [draw,shape=circle]
\node (v1) at (0,0) {$v_1$};
\node (v2) at (0,-2) {$v_2$};
\node (v3) at (-1.5,-4) {$v_3$};
\node (v4) at (1.5,-4) {$v_4$};
\node (v5) at (3.5,-4) {$v_5$};
\node (v6) at (3.5,-2) {$v_6$};
\node (v7) at (2,0) {$v_7$};
\node[red] (e457) at (1.75,-2) {$e_{\{4,5,7\}}$};
\node[red] (e1) at (-2,0) {$e_{\{1\}}$};
\node[red] (e2) at (-2,-2) {$e_{\{2\}}$};
\node[red] (e34) at (0,-3.5) {$e_{\{3,4\}}$};
\node[red] (e6) at (3.5,0) {$e_{\{6\}}$};
\draw[-{Latex[length=3mm, width=2mm]}, bend right] (v1) to (v2);
\draw[-{Latex[length=3mm, width=2mm]}, bend right] (v2) to (v1);
\foreach \from/\to in {v2/v3, v2/v4, v7/v6, v7/v1, v4/v5, v6/v5}
\draw[-{Latex[length=3mm, width=2mm]}] (\from) -- (\to);
\foreach \from/\to in {e457/v7, e457/v4, e457/v5, e1/v1, e2/v2, e34/v3, e34/v4, e6/v6}
\draw[-{Latex[length=3mm, width=2mm]}, red] (\from) -- (\to);
\end{tikzpicture}
\caption{A \HEDG{} on the left and its augmented graph on the right.}
 \label{fig:augmentation}
\end{figure}

\begin{Lem}
\label{d-sep-aug-hedg}
Let $G=(V,E,H)$ be a \HEDG{} and $X,Y,Z \ins V$ subsets.
Let $G^\aug$ be the augmented graph of $G$.
Then we have:
\[ X \Indep_G^d Y \given Z \quad\iff\quad X \Indep_{G^\aug}^d Y \given Z.  \]
\begin{proof}
Every bidirected edge $v_1 \leftrightarrow v_2$ in a $Z$-open path will be replaced by $v_1 \leftarrow e_F \to v_2$ for an $F \in H$ with $v_1,v_2 \in F$, and vice versa. Since $e_F \notin Z$ the equivalence is clear from the definition of d-separation. Note that the additional edges of the form $e_F \to v$ with $F=\{v\}$ can only extend paths by introducing $v \ot e_F \to v$ and thus do not matter.
\end{proof}
\end{Lem}

\subsection{Undirected Graphs}

The main idea of dealing with the cycles and hyperedges of the \HEDG{}es at the same time is to transform the whole \HEDG{} into an undirected graph and use the well-established theory for undirected graphs (e.g.\ the analysis of (maximal) complete sugraphs etc.). For this purpose we need to introduce some notations and constructions for undirected graphs here.

\begin{Def}[Undirected Graph]
\begin{enumerate}
  \item An \emph{undirected graph} is a pair $G=(V,E)$, where $V$ is a finite set of nodes and $E$ a set of (undirected) edges without self-loops and where $E \ins V \x V$ is symmetric, i.e.\ if $(v_1,v_2) \in E$ then also $(v_2,v_1) \in E$.
  \item An undirected graph will be represented by a graphical picture. Each node will either be represented by a \emph{dot} or by the name of the node. In an undirected graph the pair of edges $(v_1,v_2)$ and $(v_2,v_1) \in E$ will together be represented by $v_1 - v_2$ and we will then write $v_1 - v_2 \in E$.
\item A node $v \in V$ will be called a \emph{neighbour} of another node $w \in V$ in $G$ if there is an edge $v-w \in E$.
\item The set of neighbours of $v$ in $G=(V,E)$ will be denoted by 
 \[\partial_G(v):=\{ w \in V | v-w \in E\}.\]
It will be called the \emph{Markov blanket} of $v$ in $G$.
\item An undirected graph $G'=(V',E')$ is called a \emph{subgraph} of $G=(V,E)$ if $V' \ins V$ and $E' \ins E$.
\item If $A \ins V$ is a subset of the nodes of $G$, then $A$ induces a \emph{subgraph} of $G$ given by $G(A):=(V(A),E(A))$,
where $V(A):=A$ and $E(A):= E \cap (A \x A)$. We again will often just write $A$ instead of $G(A)$.
\item A subgraph $A \ins G$ is called a \emph{complete subgraph}\footnote{Some authors use the term \emph{clique} for complete subgraphs (e.g.\ \cite{Kol09}), whereas some other authors use \emph{clique} for maximal complete subgraphs (e.g.\ \cite{Lau98}). To avoid confusion we will not use the term \emph{clique} in this paper. \label{fn:clique}} of $G$ if for every two distinct $v,w \in A$ there is an edge $v - w \in E$.  
\item A subgraph $A \ins G$ is called a \emph{maximal complete subgraph}\footref{fn:clique} of $G$ if $A$ is complete and there is no complete subgraph $B \ins G$ with 
$A \subsetneq B$.
\item A \emph{path} in $G$ is a sequence of $n$ nodes from $V$ with $n \ge 1$, and $n-1$ edges from $E$:
 \[ v_1 - v_2 - \cdots - v_{n-1} - v_n.  \]
    The nodes $v_1$ and $v_n$ are called \emph{endnodes} of that path.
		The definition of path also allows for \emph{trivial paths} consisting only of a single node $v_1$.
\end{enumerate}
\end{Def}

\begin{Def}[Marginalization of an undirected graph]
Let $G=(V,E)$ be an undirected graph and $W \ins V$ a subset.
The \emph{marginalization of $G$ with respect to $W$} is the undirected graph $G^\marg=(V^\marg,E^\marg)$, where
$V^\marg:=V \sm W$ and for distinct nodes $v_1,v_2 \in V\sm W$ the edge $v_1 - v_2$ is in $E^\marg$ if and only if there exist a finite number of nodes $w_1,\dots,w_n \in W$, $n \ge 0$, such that all the following edges are in $E$:
\[ v_1 - w_1 - w_2 - \cdots - w_n - v_2. \]
In particular, we have an edge $v_1 - v_2$ in $E^\marg$ if $v_1 - v_2$ is already an edge in $E$. 
In other words, we erase the nodes and edges from $W$ and connect the remaining nodes if they had a connecting path through the nodes of $W$ in $E$.
We also write $G^{\marg(V \sm W)}$ and $G^{\marg\sm W}$, resp., instead of just $G^\marg$ to emphasize the nodes that are left, or that are marginalized out, resp..
\end{Def}

\begin{Lem}[Subgraphs and marginalization]
\label{subgr-marg}
Let $G=(V,E)$ be an undirected graph, $G'=(V',E') \ins G$ a subgraph and $W \ins V$ a subset.
Then $(G')^{\marg\sm(W \cap V')}$ is a subgraph of $G^{\marg\sm W}$.
\begin{proof}
We clearly have $V' \sm (W \cap V')=V' \sm W \ins V \sm W$.
Now let $v_1-v_2 \in (E')^{\marg\sm(W \cap V')}$. Then there is an $r \ge 0$ and nodes $u_1,\dots,u_r \in V' \cap W$ such that 
$$ v_1-u_1-\cdots-u_r-v_2$$
is a path in $G'$. Since $G'$ is a subgraph of $G$ the path also exists in $G$.
By the definition of marginalization and $u_i \in W$ get the edge $v_1-v_2$ in $G^{\marg \sm W}$.
\end{proof}
\end{Lem}

\begin{Lem}[Complete subgraphs and marginalization]
\label{compl-subgr-marg}
Let $G=(V,E)$ be an undirected graph, $v \in V$ a node, $W \ins V$ a subset and $G^{\marg\sm W}$ the marginalization of $G$ w.r.t.\ $W$.
We then have:
\begin{enumerate}
\item If $v \in W$ then the nodes from $\partial_G(v) \sm W$ form a complete subgraph in $G^{\marg \sm W}$.
\item If $v \notin W$ and $\partial_G(v)$ is a complete subgraph of $G$ then $\partial_{G^{\marg\sm W}}(v)$ is a complete subgraph of $G^{\marg\sm W}$.
\item If $C \ins G^{\marg\sm W}$ is a complete subgraph and $W \ins G$ is such that for every $v_1,v_2 \in V \sm W$ with a path $v_1-w_1-\cdots-w_r-v_2$, $r \ge 1$, $w_i \in W$, we have $v_1-v_2$ then the nodes of $C$ also build a complete subgraph of $G$. 
\item If $C \ins G^{\marg\sm W}$ is a complete subgraph, $W=\{w\}$ and $\partial_G(w) \ins G$ is complete then also $C \ins G$ is complete.
\end{enumerate}
\begin{proof}
1. Consider $v_1,v_2 \in \partial_G(v) \sm W$ with $v_1 \neq v_2$. Then we have the path $v_1 - v -v_2$ in $G$. Since $v \in W$ we get $v_1-v_2$ in $G^{\marg\sm W}$.\\
2. Consider $v_1,v_2 \in  \partial_{G^{\marg\sm W}}(v)$ with $v_1 \neq v_2$. Then we have paths $v_1-\cdots-w_1-v$ and $v_2-\cdots-w_2-v$ in $G$ with all the intermediate nodes in $W$ (or $v_i=w_i$). Since $\partial_G(v)$ is a complete subgraph of $G$ we have the edge $w_1-w_2 \in E$ leading to a path $v_1-\cdots-w_1-w_2-\cdots-v_2$ in $G$ with again all intermediate nodes in $W$ (or $v_i=w_i$). Marginalizing out $W$ gives the edge $v_1-v_2$ in $G^{\marg\sm W}$. \\
3. If $v_1,v_2 \in C \ins G^{\marg\sm W}$ then $v_1-v_2 \in C$ by the completeness of $C$. So there is a path in $G$ of the form $v_1-w_1-\cdots-w_r-v_2$ with $w_i \in W$. By assumption on $W$ we then have $v_1-v_2 \in G$. So $C$ is also complete in $G$.\\
4. This follows directly from 3.: Since $\partial_G(w)$ is complete the path $v_1-w-v_2$ implies $v_1-v_2 \in G$.
\end{proof}
\end{Lem}

\begin{Def}[Separation in an undirected graph]
Let $G=(V,E)$ be an undirected graph and $X,Y,Z \ins V$ subsets of the nodes.
\begin{enumerate}
\item A path in $G$, $n \ge 1$:
\[ v_1 - v_2 - \cdots - v_{n-1} - v_n  \]
will be called \emph{$Z$-blocked} or \emph{blocked by $Z$} if there exists a node $v_i$ in the path with $v_i \in Z$. \\
If this is not the case then the path is called \emph{$Z$-open} or \emph{$Z$-active}.
\item  We say that $X$ is \emph{separated} from $Y$ given $Z$ if every path in $G$ with one endpoint in $X$ and the other endpoint in $Y$ is blocked by $Z$.
In symbols this will be written as follows:
\[ X \Indep_G Y \given Z.\]
\end{enumerate}
\end{Def}

\begin{Rem}
For checking the separation $X \Indep_G Y \given Z$ it is enough to only check paths where every node occurs at most once. If a node occurred twice there would be a $Z$-open ``short-cut''. 
\end{Rem}

\begin{Lem}[Subgraphs and separation]
\label{sep-subgraph}
Let $G=(V,E)$ be an undirected graph and $G'=(V',E') \ins G$ a subgraph and $X,Y,Z \ins V'$ subsets of the nodes.
Then we have the implication:
\[ X \Indep_G Y \given Z \; \implies \;  X \Indep_{G'} Y \given Z. \]
\begin{proof}
Every $Z$-open path in $G'$ is clearly $Z$-open in $G$.
\end{proof}
\end{Lem}

\begin{Lem}[Separation stable under marginalization]
\label{sep-marg}
Let $G=(V,E)$ be an undirected graph and $X,Y,Z,W \ins V$ subsets with $W \cap (X \cup Y \cup Z) = \emptyset$.
Let $G^\marg$ be the marginalization with respect to $W$.
Then we have:
\[ X \Indep_G Y \given Z \iff X \Indep_{G^\marg} Y \given Z.   \]
\begin{proof}
Every sequence of nodes in $W$ in a $Z$-open path can be replaced by one edge in $G^\marg$ and thus stays $Z$-open.
On the other hand, every edge in a $Z$-open path in $G^\marg$ that was not already in $G$ can be replaced by a sequence of nodes in $W$.
The condition $W \cap (X \cup Y \cup Z) = \emptyset$ ensures that this path stays $Z$-open.
\end{proof}
\end{Lem}

\subsection{The Moralization of a \HEDG{}}

As mentioned before, the main idea of dealing with the cycles and hyperedges of the \HEDG{}es at the same time is to transform the whole \HEDG{} into an undirected graph. Basically, this is done by introducing the \emph{(generalized) moralization} of a \HEDG{}. This construction in general will not capture all information of the \HEDG{} itself, but we will see that the most relevant information of a \HEDG{} (e.g.\ for the purpose of d-separation) is already captured by the family of all moralizations of all ancestral sub-\HEDG{}es. %

\begin{Def}[Moralization of a directed graph]
\label{moral-dg-def}
Let $G=(V,E)$ be a directed graph.
The \emph{moralization} of $G$ is the undirected graph $G^\moral=(V^\moral,E^\moral)$, where $V^\moral:=V$ and $E^\moral$ is constructed from $E$ as follows:
All directed edges from $E$ are replaced by undirected ones and in addition we add the edges $v-w$ for $v \neq w$  if $v$ and $w$ have a common child in $G$.
In other words, for every $v \in V$ we connect every two nodes in $\{v\} \cup \Pa^G(v)$ and make this set a complete subgraph of $G^\moral$.
\end{Def}

\begin{Def}[Moralization of a \HEDG]
\label{moral-hedg-def}
Let $G=(V,E,H)$ be a \HEDG.
The \emph{(generalized) moralization}\footnote{In \cite{Richardson03} 2.2 a similar construction for ADMGs was called \emph{augmentation}. But this term was reserved for a different construction \ref{augm-def} in this paper.}
 of $G$  is given by the undirected graph $G^\moral=(V^\moral,E^\moral)$ which is constructed as follows:
We put $V^\moral:=V$ and for $E^\moral$ we have the edge $v - w$ for $v \neq w$ if and only if 
there exists a sequence of nodes $v_1, \dots, v_n \in V$, $n \ge 1$, with $v \in \{v_1\} \cup \Pa^G(v_1)$ and $w \in \{v_n\} \cup \Pa^G(v_n)$ and the following hyperedges all lie in $H$:
\[ v_1 \leftrightarrow \cdots \leftrightarrow v_n. \]
In other words, for every $v \in V$ we connect every two nodes in $\Dist^G(v) \cup \Pa^G(\Dist^G(v))$ and thus make this set a complete subgraph of $G^\moral$.\\
Note that the particular cases: $v \in \Pa^G(w)$, $w \in \Pa^G(v)$, $\{v,w\} \in H$ and $\{v,w\} \ins \Pa^G(v_1)$ for a node $v_1 \in V$ are all covered by the definition and thus give an edge $v - w$ to $E^\moral$.
\end{Def}

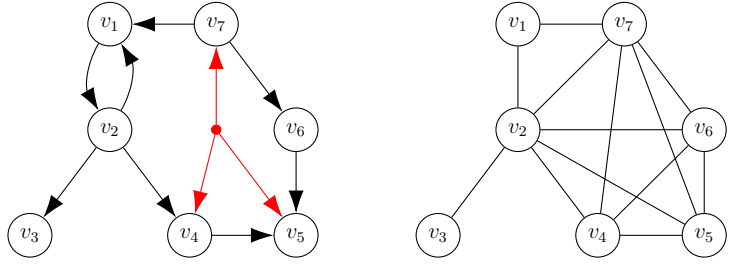
\begin{figure}
\centering
\begin{tikzpicture}[scale=.7, transform shape]
\tikzstyle{every node} = [draw,shape=circle]
\node (v1) at (0,0) {$v_1$};
\node (v2) at (0,-2) {$v_2$};
\node (v3) at (-1.5,-4) {$v_3$};
\node (v4) at (1.5,-4) {$v_4$};
\node (v5) at (3.5,-4) {$v_5$};
\node (v6) at (3.5,-2) {$v_6$};
\node (v7) at (2,0) {$v_7$};
\node[fill,circle,red,inner sep=0pt,minimum size=5pt] (e457) at (2,-2) {};
\draw[-{Latex[length=3mm, width=2mm]}, bend right] (v1) to (v2);
\draw[-{Latex[length=3mm, width=2mm]}, bend right] (v2) to (v1);
\foreach \from/\to in {v2/v3, v2/v4, v7/v6, v7/v1, v4/v5, v6/v5}
\draw[-{Latex[length=3mm, width=2mm]}] (\from) -- (\to);
\foreach \from/\to in {e457/v7, e457/v4, e457/v5}
\draw[-{Latex[length=3mm, width=2mm]}, red] (\from) -- (\to);
\end{tikzpicture}
\hspace{1cm}
\begin{tikzpicture}[scale=.7, transform shape]
\tikzstyle{every node} = [draw,shape=circle]
\node (v1) at (0,0) {$v_1$};
\node (v2) at (0,-2) {$v_2$};
\node (v3) at (-1.5,-4) {$v_3$};
\node (v4) at (1.5,-4) {$v_4$};
\node (v5) at (3.5,-4) {$v_5$};
\node (v6) at (3.5,-2) {$v_6$};
\node (v7) at (2,0) {$v_7$};
\foreach \from/\to in {v1/v2, v2/v3, v2/v4, v7/v6, v7/v1, v4/v5, v6/v5, v2/v7, v4/v7, v5/v7, v2/v6, v4/v6, v2/v5}
\draw (\from) -- (\to);
\end{tikzpicture}
\caption{A \HEDG{} on the left and its moralization on the right.}
 \label{fig:moralization}
\end{figure}

\begin{Rem}
For a \HEDG{} $G = (V,E,H)$ and a node $v \in G$ we have:
  $$\partial_{G^{\moral}}(v) = \Dist^G(\{v\} \cup \Ch^G(v)) \cup \Pa^G(\Dist^G(\{v\} \cup \Ch^G(v))).$$ 
\end{Rem}

\begin{Rem}
\label{moral-hedg-moral-dg}
Note that if we identify a directed graph $(V,E)$ with the \HEDG{} given by $(V,E,H_1)$, where $H_1$ is trivial (i.e.\ only consists of the subsets $F \ins V$ with $\# F\le 1$), then the two definitions of moralization coincide: 
\[(V,E)^\moral = (V,E,H_1)^\moral.\] 
So in this sense the definition of moralization of a \HEDG{} generalizes the one for directed graphs.
\begin{proof}
For trivial $H$ there are only one element sequences $v_1=v_n$. So the edge $v - w$ is in $E^\moral$ if and only if $\{v,w\} \ins \{v_1\} \cup \Pa^G(v_1)$ for some $v_1 \in V$. This is equivalent to the cases $v \to w$, $w \to v$ and $\{v,w\} \ins \Pa^G(v_1)$, which exactly define the edges for the moralization of a directed graph.
\end{proof}
\end{Rem}

\begin{Lem}
\label{moral-hedg}
Let $G=(V,E,H)$ be a \HEDG. Then we have:
\[ G^\moral=((G^\aug)^\moral)^\marg,  \]
where we take the moralization of the augmented graph $G^\aug$ and then marginalize with respect to the set of nodes $V^\aug \sm V=\{e_F| F\in \tilde{H} \}$, which where added in the construction of the augmentation. Alternatively, we could use the canonical graph instead of the augmented graph.
\begin{proof}
First we see that both sides have $V$ as the underlying set of nodes.
$v-w$ is an edge in $((G^\aug)^\moral)^\marg$ if and only if there are $e_{F_1},\dots, e_{F_n}$, $n\ge 0$, with $F_i \in \tilde{H}$ such that the edges 
\[ v - e_{F_1} -\dots - e_{F_n} - w \]
are all in $(G^\aug)^\moral$. But $v - e_{F_1}$ (and $e_{F_n} - w$, resp.) only exists there if and only if $v \in \{v'\} \cup \Pa^G(v')$ for some $v' \in F_1$ (and $w \in \{w'\} \cup \Pa^G(w')$ for some $w' \in F_n$, resp.). Furthermore, $e_{F_i}  - e_{F_{i+1}}$ only exists there if and only if there exists a node $v_i \in V$ with $v_i \in F_i \cap F_{i+1}$.
So the existence of the following edges in $((G^\aug)^\moral)$:
\[ v - e_{F_1} -\dots - e_{F_n} - w \]
 is equivalent to the existence of the following bidirected arrows in $H$:
\[ v' \oto v_1  \oto \dots \oto v_{n-1} \oto w' \]
with $v \in \{v'\} \cup \Pa^G(v')$ and $w \in \{w'\} \cup \Pa^G(w')$ (with $n=0$ meaning that $v'=w'$).
But this is equivalent to the claim that $v - w$ lies in $G^\moral$.
\end{proof}
\end{Lem}

\begin{Lem}
\label{d-sep-moral-dg}
Let $G=(V,E)$ be a directed graph and $X,Y,Z \ins V$ subsets such that $G=\Anc^G(X \cup Y \cup Z)$.
Then we have:
\[ X \Indep_G^d Y \given Z \iff X \Indep_{G^\moral} Y \given Z. \]
\begin{proof}
  This was proven for DAGs in \cite{Lau90} Prop.\ 3 (or \cite{Lau98} Prop.\ 3.25), but as pointed out in \cite{PearlDechter96} Lem. 3 the proof also works for non-acyclic directed graphs. For completeness we add a proof here:\\
($\Longleftarrow$): If $x \genfrac{}{}{0pt}{}{\to}{\ot} \cdots \genfrac{}{}{0pt}{}{\to}{\ot}  y$ is a $Z$-open path from $X$ to $Y$ in $G$, then every non-collider node (e.g.\ $v \to w \to v'$) is not in $Z$ and the corresponding undirected version $v - w - v'$ is then also $Z$-open.
If $v \to w \ot v'$ is a collider on the $Z$-open path then $w \in \Anc^G(Z)$ and $w$ might even be in $Z$, 
 but in the moralized graph we can replace that part with the moralizing edge $v - v'$.
Thus we can construct a $Z$-open path in $G^\moral$.\\
($\Longrightarrow$): %
We now consider a $Z$-open path $x - \cdots - y$ from $X$  to $Y$ in $G^\moral$ with the lowest number of moralized edges (i.e.\ edges that do not correspond to directed edges in $G$). If $v-w$ is such a moralized edge then there is a common child $v'$ in $G$ and we have the directed path $v \to v' \ot w$ in $G$. Since $G=\Anc^G(X \cup Y \cup Z)$ there is a shortest directed path $v' \to \cdots \to w'$ with $w' \in X \cup Y \cup Z$. Moreover, we have $w' \notin X \cup Y$. Otherwise, if $w' \in Y$ the $Z$-open path 
$x - \dots - v' - \dots - w'$ had a lower number of moralized edges. Similar for $w' \in X$ and $w' - \dots - v' - \dots - y$.
It follows that $w' \in Z$ and thus $v' \in \Anc^G(Z)$ and $v \to v' \ot w$ is $Z$-open in $G$.
So replacing all moralized edges in $x - \dots -y$ with such a $Z$-open collider and all other edges with direct edges we get a $Z$-open path in $G$.
\end{proof}
\end{Lem}

We are now ready to prove that also for \HEDG{}es, d-separation can be reduced to a separation criterion in the moralization of ancestral subsets. For DAGs this was shown in \cite{Lau98} Prop.\ 3.25 and for ADMGs in \cite{Richardson03}.

\begin{Thm}
\label{d-sep-moral-hedg}
Let $G=(V,E,H)$ be a \HEDG{}, $X,Y,Z \ins V$ subsets and $\Anc^G(X\cup Y \cup Z)$ the ancestral sub-\HEDG{} generated by $X\cup Y \cup Z$
and $\Anc^G(X\cup Y \cup Z)^\moral$ its generalized moralization.
Then we have:
\[ X \Indep_G^d Y \given Z  \iff X \Indep_{\Anc^G(X\cup Y \cup Z)^\moral} Y \given Z.\]%
\begin{proof}
We have the sequence of equivalences:
\begin{eqnarray*}
X \Indep_G^d Y \given Z 
& \stackrel{\ref{d-sep-anc-sub-hedg}}{\iff} & X \Indep_{\Anc^{G}(X\cup Y \cup Z)}^d Y \given Z. \\
& \stackrel{\ref{d-sep-aug-hedg}}{\iff} & X \Indep_{\Anc^{G}(X\cup Y \cup Z)^\aug}^d Y \given Z. \\
& \stackrel{\ref{d-sep-moral-dg}}{\iff} & X \Indep_{(\Anc^{G}(X\cup Y \cup Z)^\aug)^\moral} Y \given Z. \\
& \stackrel{\ref{sep-marg}}{\iff} & X \Indep_{((\Anc^{G}(X\cup Y \cup Z)^\aug)^\moral)^\marg} Y \given Z \\
& \stackrel{\ref{moral-hedg}}{\iff} & X \Indep_{\Anc^{G}(X\cup Y \cup Z)^\moral} Y \given Z.
\end{eqnarray*}
Note that for $A := \Anc^{G}(X\cup Y \cup Z)^\aug$ we have $A=\Anc^A(X\cup Y \cup Z)$ and \ref{d-sep-moral-dg} applies.
For the application of \ref{sep-marg} we note, that $\{e_F| F \in \widetilde{H(A)} \} \cap (X \cup Y \cup Z) = \emptyset$.\\
Alternatively, we could argue like in \ref{d-sep-moral-dg}: In any $Z$-open path from $X$ to $Y$ replace every non-degenerate unbroken sequence of colliders:
 \[ v  \stackrel{\to}{=} v_1  \oto \dots \oto v_n \stackrel{\ot}{=} w, \] 
$n\ge 0$, with an undirected edge $v-w$. For the converse replacement it can be argued by considering shortest paths with lowest numbers of colliders that $v_1,\dots, v_n \in \Anc^G(Z)$ and $v,w \notin Z$ as it happened in \ref{d-sep-moral-dg}. So every open path in one sense gives another open path in the other sense.
\end{proof}
\end{Thm}

\begin{Cor}
\label{d-sep-moral-hedg-subgraph}
Let $G=(V,E,H)$ be a \HEDG{} and $X,Y,Z \ins V$ be subsets.
Then we have the implication:
\[ X \Indep_{G^\moral} Y \given Z  \;\implies\; X \Indep_{G}^d Y \given Z.\]
\begin{proof}
$\Anc^G(X \cup Y \cup Z)^\moral$ is a subgraph of $G^\moral$. 
So we have: 
\[ X \Indep_{G^\moral} Y \given Z  \;\stackrel{\ref{sep-subgraph}}{\implies}\; X \Indep_{\Anc^G(X \cup Y \cup Z)^\moral} Y \given Z \;
\stackrel{\ref{d-sep-moral-hedg}}{\iff}\; X \Indep_G^d Y \given Z.\]
\end{proof}
\end{Cor}

\begin{Lem}[Marginalization and moralization of a \HEDG{}]
\label{marg-mor1}
Let $G=(V,E,H)$ be a \HEDG{} and $W \ins V$ a subset.
Then $(G^{\marg\sm W})^\moral$ is an (undirected) subgraph of $(G^\moral)^{\marg\sm W}$ with the same set of nodes.
\begin{proof}
First note that the underlying set of nodes in both cases really is $V \sm W$.
So we only need to show that the edges in one graph are also edges in the other.
We will do induction on $\#W$. 
For $\#W=1$ let $W=\{w\}$. 
Let $v_1-v_2$ be an edge in $(G^{\marg\sm W})^\moral$. This will come from a path of the form
\[ v_1  \stackrel{\to}{=} u_1  \oto \dots \oto u_r \stackrel{\ot}{=} v_2\]
in $G^{\marg\sm W}$ where none of the  nodes is $w$. This path comes from a longer one in $G$, where some edges $\to$ might be replaced by $\to w \to$ and some bidirected edges $\oto$ by $\ot w \to$, $\oto w \to$ or $\ot w \oto$.
Moralizing will give a path of the form $v_1-v_2$, $v_1-w-v_2$ or, if more intermediate nodes are $w$, of the form $v_1-w-v_1' \cdots- v_s' - w -v_2$, where we can take the shortcut over $w$ to get $v_1-w-v_2$. 
Marginalizing out $w$ will give us an edge $v_1-v_2$ in  $(G^\moral)^{\marg\sm W}$. This shows the claim for $\#W=1$.
Now let $\#W>1$ and let $W=W' \cup \{w\}$. By induction we already have the inclusion of subgraphs:
\[ \begin{array}{rcl}
(G^{\marg\sm W'})^\moral &\ins&  (G^\moral)^{\marg\sm W'}. %
\end{array}\]
Marginalizing this w.r.t.\ $\{w\}$ we by \ref{subgr-marg} get the inclusion of subgraphs:
\[ \begin{array}{rcl}
((G^{\marg\sm W'})^\moral)^{\marg \sm \{w\}} &\ins&  (G^\moral)^{\marg\sm W}.
\end{array}\]
The case $\#W=1$ applied to $G^{\marg\sm W'}$ gives us the inclusion of subgraphs:
\[ \begin{array}{rcl}
(G^{\marg\sm (W'\cup \{w\})})^\moral &\ins& ((G^{\marg\sm W'})^\moral)^{\marg \sm \{w\}}.
\end{array}\]
Putting both  together we get the inclusion of subgraphs:
\[ \begin{array}{rcl}
(G^{\marg\sm W})^\moral &=& (G^{\marg\sm (W'\cup \{w\})})^\moral \\
&\ins& ((G^{\marg\sm W'})^\moral)^{\marg \sm \{w\}}\\
&\ins&  (G^\moral)^{\marg\sm W}.
\end{array}\]
This shows the claim.
\end{proof}
\end{Lem}

\begin{Lem}[Marginalization and moralization of a \HEDG{}, reverse inclusion]
\label{marg-mor2}
Let $G=(V,E,H)$ be a \HEDG{} and $W=\{w\}$ with $w \in V$.
Then the equality of undirected graphs
\[ (G^{\marg\sm W})^\moral = (G^\moral)^{\marg\sm W}  \]
holds if and only if one of the following cases holds:
\begin{enumerate}
\item $\Ch^G(w) \sm \{w\} \neq \emptyset$.
\item $\Ch^G(w) \sm \{w\} = \emptyset$ and 
for every pair of nodes 
\[ v_1, v_2 \in \lp \Dist^G(w) \cup \Pa^G(\Dist^G(w)) \rp \sm \{w\} \]
there exists an $r \ge 1$, nodes $u_1,\dots,u_r \in V \sm W$ and a path in $G$ of the form:
\[ v_1  \stackrel{\to}{=} u_1  \oto \dots \oto u_r \stackrel{\ot}{=} v_2\]
(including the cases $v_1=v_2$, $v_1 \to v_2$, $v_1 \ot v_2$, $v_1 \oto v_2$, $v_1 \to u_1 \ot v_2$ ($u_1 \neq w$) etc.).
\end{enumerate}
\begin{proof}
By \ref{marg-mor1} we already have the inclusion of subgraphs:
\[ (G^{\marg\sm W})^\moral \ins (G^\moral)^{\marg\sm W}.\]
For the reverse inclusion let $v_1-v_2$ be an edge in $(G^\moral)^{\marg\sm W}$.
The edge then comes from the cases $v_1-v_2$ or $v_1-w-v_2$ in $G^\moral$, which themselves come from paths of the form:
\[ v_1  \stackrel{\to}{=} u_1  \oto \dots \oto u_r \stackrel{\ot}{=} v_2, \qquad \text{or} \]
\[ v_1  \stackrel{\to}{=} u_1  \oto \dots \oto u_r \stackrel{\ot}{=} w \stackrel{\to}{=} u_{r+1}  \oto \dots \oto u_m \stackrel{\ot}{=} v_2. \] 
If $w$ does not occur in the path then marginalizing won't change the path, and hence $v_1 - v_2$ is in $(G^{\marg\sm W})^\moral$.

  In case 1, if $w$ occurs on the path (note: $v_1,v_2 \neq w$), the subpaths $u \stackrel{\to}{\oto} w \stackrel{\ot}{\oto} u'$ can be replaced by $u \stackrel{\to}{\oto} w' \stackrel{\ot}{\oto} u'$ (with the same arrowhead signature) in $G^{\marg\sm W}$, where $w' \in \Ch^G(w)\sm\{w\}$.
  If $w$ still occurs on the path, then it must be in one of the subpaths $u \oto  w \to u'$, $u \ot w \to u'$, $u \ot w \oto u'$. Marginalizing out $w$ will replace these by $u \oto u'$, and we obtain a path of the first form in $G^{\marg\sm W}$. Moralizing will then give an edge $v_1-v_2$ in $(G^{\marg\sm W})^\moral$.

In case 2, if $w$ occurs on the path, we can replace the path by one that contains $w$ only once. The resulting path must be of the first form, because $w$ is childless (up to itself). The second condition means that this path can be replaced by one that does not contain $w$ at all, and this path therefore also exists in $G^{\marg\sm W}$. Moralizing will then give an edge $v_1-v_2$ in 
$(G^{\marg\sm W})^\moral$.

For the reverse implication assume $(G^{\marg\sm W})^\moral = (G^\moral)^{\marg\sm W}$. If $\Ch^G(w) \sm \{w\} \neq \emptyset$ then we have nothing to show. So we may assume $\Ch^G(w) \sm \{w\} = \emptyset$. Then let $v_1, v_2 \in \lp \Dist^G(w) \cup \Pa^G(\Dist^G(w)) \rp \sm \{w\}$ be a pair of distinct nodes. Since $\Dist^G(w) \cup \Pa^G(\Dist^G(w))$ forms a complete subgraph in $G^\moral$ we have the edge $v_1-v_2$ in $G^\moral$. Furthermore, $v_1,v_2 \neq w$ and thus $v_1-v_2$ is also in $(G^\moral)^{\marg\sm W}$. By assumption we have $(G^\moral)^{\marg\sm W} =(G^{\marg\sm W})^\moral$. It follows that $v_1-v_2$ exists in $(G^{\marg\sm W})^\moral$. By the definition of moralization we have an
$r \ge 1$, nodes $u_1,\dots,u_r \in V \sm W$ and a path in $G^{\marg\sm W}$ of the form:
\[ v_1  \stackrel{\to}{=} u_1  \oto \dots \oto u_r \stackrel{\ot}{=} v_2.\]
Again, this path comes from a longer one in $G$, where some edges $\to$ might be replaced by $\to w \to$ and some bidirected edges $\oto$ by $\ot w \to$, $\oto w \to$ or $\ot w \oto$. Since we assumed $\Ch^G(w) \sm \{w\} = \emptyset$ no such replacement would occur and the above path already exists in $G$. This shows the claim from point 2..
\end{proof}
\end{Lem}

\begin{Rem}
Note that for a \HEDG{} and a general set $W \ins V$ we do not have the equality
$(G^{\marg\sm W})^\moral=(G^\moral)^{\marg\sm W}$. This can be seen from figure \ref{fig:subgr-moral}.
\begin{figure}[!htb] 
\centering
\begin{tikzpicture}[scale=.9, transform shape]
\tikzstyle{every node} = [draw,shape=circle]
\node (v2) at (3,2.5) {$v_2$};
\node (v3) at (0,0) {$v_3$};
\node (v4) at (2,0) {$v_4$};
\node (v1) at (-1,2.5) {$v_1$};
\foreach \from/\to in {v2/v4, v1/v3}
\draw[-{Latex[length=3mm, width=2mm]}] (\from) -- (\to);
\draw[-{Latex[length=3mm, width=2mm]}, bend right] (v3) to (v4);
\draw[-{Latex[length=3mm, width=2mm]}, bend right] (v4) to (v3);
\end{tikzpicture} \qquad
\begin{tikzpicture}[scale=.9, transform shape]
\tikzstyle{every node} = [draw,shape=circle]
\node (v2) at (3,2.5) {$v_2$};
\node (v3) at (0,0) {$v_3$};
\node (v4) at (2,0) {$v_4$};
\node (v1) at (-1,2.5) {$v_1$};
\foreach \from/\to in {v2/v4, v1/v3, v4/v3, v1/v4, v2/v3}
\draw[-] (\from) -- (\to);
\end{tikzpicture}\\[20pt]
\begin{tikzpicture}[scale=.9, transform shape]
\tikzstyle{every node} = [draw,shape=circle]
\node (v2) at (3,2.5) {$v_2$};
\node (v1) at (-1,2.5) {$v_1$};
\end{tikzpicture} \qquad
\begin{tikzpicture}[scale=.9, transform shape]
\tikzstyle{every node} = [draw,shape=circle]
\node (v2) at (3,2.5) {$v_2$};
\node (v1) at (-1,2.5) {$v_1$};
\draw[-] (v1) to (v2);
\end{tikzpicture}
\caption{ A \HEDG{} and its moralization (top row). If we marginalize out $W=\{v_3,v_4\}$ on both sides
we get the row on the bottom. So we have $(G^{\marg\sm W})^\moral \subsetneq (G^\moral)^{\marg\sm W}$.
Note that the moralization of the bottom left is not the same as the bottom right even though we have
$\Ch^G(v_i) \sm \{v_i\} \neq \emptyset$ for all nodes. }
\label{fig:subgr-moral}
\end{figure}
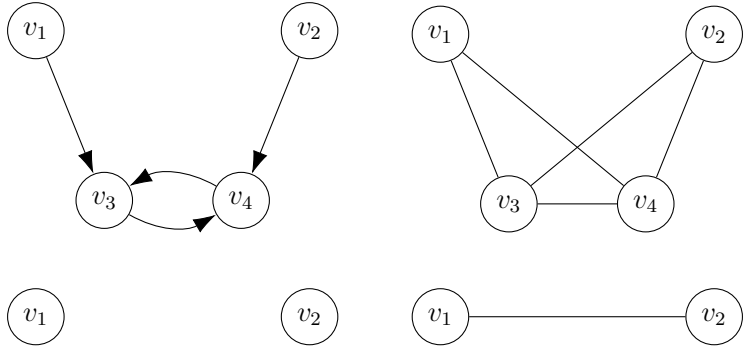
\end{Rem}

\subsection{Total Orders for \HEDG{}es}

One of the most important properties of an acyclic directed graph (DAG) is the existence of a total order that respects parental relations. This is usually called a \emph{topological order} of the DAG. As soon as we allow for cycles, such a total order might not exist anymore. In this subsection we will introduce several generalizations of topological orders on the set of nodes for directed graphs with hyperedges (\HEDG{}es), each highlighting a different property every topological order has. We investigate their mutual relations and provide criteria for their existence.

\begin{Def}[Total Orders for a \HEDG{}]
\label{peo-hedg}
Let $G=(V,E,H)$ be a \HEDG{} and $<$ a total order on the set of nodes $V$. Then:
\begin{enumerate}
\item For a node $v \in V$ we define the set of \emph{predecesssors} of $v$ in $G$ by:
 \[\begin{array}{lll}
\Pred^G_\le(v):=\{ w \in V | w <v \text{ or } w=v\} ,
\end{array}\]
and endow it with the structure of the marginalized \HEDG{} $G^\marg$ induced by $G$, where all other variables are marginalized out.
 \item For every subset $A \ins V$ we have the induced total order on $A$.
\item We call $<$ a \emph{topological order} of $G$, if for every $v \in V$ and every $w \in \Pa^G(v)$ we have $w < v$.
\item We call $<$ a \emph{pseudo-topological order} of $G$, if for every $v \in V$ and 
for every $w \in \Anc^G(v)\sm \Sc^G(v)$ we have $w < v$. %
\item We call $<$ \emph{assembling}, if for every $v,v',w \in V$ with $v' \in \Sc^G(v)$ and $v \le w \le v'$ we also have $w \in \Sc^G(v)$.
\item We call $<$ a \emph{perfect elimination order}\footnote{The definition of a perfect elimination order for \HEDG{}es generalizes the one for undirected graphs, see \cite{West01} Def. 5.3.12.} of $G$, if 
for every $v \in V$ and every ancestral sub-\HEDG{} $A \ins \Pred^G_\le(v)$ with $v \in A$ we have that
\[ \partial_{A^\moral}(v) \cup \{v\}  \]
is a complete subgraph of $A^\moral$.
\item We call $<$ a \emph{quasi-topological order} of $G$, if it is a pseudo-topological order and a perfect elimination order.
\end{enumerate}
\end{Def}

\begin{Rem}
For any \HEDG{} $G = (V,E,H)$, and any node $v \in V$ we have: 
$\partial_{G^{\moral}}(v)$ is a complete subgraph of $G^{\moral}$ if and only if $\partial_{G^{\moral}}(v) \cup \{v\}$ is a complete subgraph of $G^{\moral}$.
\end{Rem}

\begin{figure}[!htb] 
\centering
\begin{tikzpicture}[scale=.7, transform shape]
\tikzstyle{every node} = [draw,shape=circle]
\node (v1) at (0,0) {$v_1$};
\node (v2) at (3,0) {$v_2$};
\node (v3) at (3,-2) {$v_3$};
\node (v4) at (0,-2) {$v_4$};
\node (v5) at (0,-4) {$v_5$};
\draw[{Latex[length=3mm, width=2mm]}-{Latex[length=3mm, width=2mm]}, red, bend right] (v3) to node[fill,circle,red,inner sep=0pt,minimum size=5pt] {} (v4);
\foreach \from/\to in {v2/v1, v3/v2, v4/v3, v1/v4, v4/v5, v3/v5}
\draw[-{Latex[length=3mm, width=2mm]}] (\from) -- (\to);
\end{tikzpicture}
\hspace{1cm}
\begin{tikzpicture}[scale=.7, transform shape]
\tikzstyle{every node} = [draw,shape=circle]
\node (v1) at (0,0) {$v_1$};
\node (v2) at (3,0) {$v_2$};
\node (v3) at (3,-2) {$v_3$};
\node (v4) at (0,-2) {$v_4$};
\node (v5) at (0,-4) {$v_5$};
\foreach \from/\to in {v2/v1, v3/v2, v4/v3, v1/v4, v1/v3, v4/v5, v3/v5}
\draw (\from) -- (\to);
\end{tikzpicture}
\\[0.5cm]
\begin{tikzpicture}[scale=.7, transform shape]
\tikzstyle{every node} = [draw,shape=circle]
\node (v1) at (0,0) {$v_1$};
\node (v2) at (3,0) {$v_2$};
\node (v3) at (3,-2) {$v_3$};
\node (v4) at (0,-2) {$v_4$};
\draw[{Latex[length=3mm, width=2mm]}-{Latex[length=3mm, width=2mm]}, red, bend right] (v3) to node[fill,circle,red,inner sep=0pt,minimum size=5pt] {} (v4);
\foreach \from/\to in {v2/v1, v3/v2, v4/v3, v1/v4}
\draw[-{Latex[length=3mm, width=2mm]}] (\from) -- (\to);
\end{tikzpicture}
\hspace{1cm}
\begin{tikzpicture}[scale=.7, transform shape]
\tikzstyle{every node} = [draw,shape=circle]
\node (v1) at (0,0) {$v_1$};
\node (v2) at (3,0) {$v_2$};
\node (v3) at (3,-2) {$v_3$};
\node (v4) at (0,-2) {$v_4$};
\foreach \from/\to in {v2/v1, v3/v2, v4/v3, v1/v4, v1/v3}
\draw (\from) -- (\to);
\end{tikzpicture}
\\[0.5cm]
\begin{tikzpicture}[scale=.7, transform shape]
\tikzstyle{every node} = [draw,shape=circle]
\node (v1) at (0,0) {$v_1$};
\node (v2) at (3,0) {$v_2$};
\node (v3) at (3,-2) {$v_3$};
\foreach \from/\to in {v2/v1, v3/v2, v1/v3}
\draw[-{Latex[length=3mm, width=2mm]}] (\from) -- (\to);
\end{tikzpicture}
\hspace{1cm}
\begin{tikzpicture}[scale=.7, transform shape]
\tikzstyle{every node} = [draw,shape=circle]
\node (v1) at (0,0) {$v_1$};
\node (v2) at (3,0) {$v_2$};
\node (v3) at (3,-2) {$v_3$};
\foreach \from/\to in {v2/v1, v3/v2, v1/v3}
\draw (\from) -- (\to);
\end{tikzpicture}
\\[0.5cm]
\begin{tikzpicture}[scale=.7, transform shape]
\tikzstyle{every node} = [draw,shape=circle]
\node (v1) at (0,0) {$v_1$};
\node (v2) at (3,0) {$v_2$};
\draw[-{Latex[length=3mm, width=2mm]},bend right] (v1) to (v2);
\draw[-{Latex[length=3mm, width=2mm]}] (v2) to (v1);
\end{tikzpicture} 
\hspace{1cm}
\begin{tikzpicture}[scale=.7, transform shape]
\tikzstyle{every node} = [draw,shape=circle]
\node (v1) at (0,0) {$v_1$};
\node (v2) at (3,0) {$v_2$};
\draw (v1) to (v2);
\end{tikzpicture}
\caption{A \HEDG{} $G$ with an assembling quasi-topological order: On the left the \HEDG{} up to the current node (of highest index) $\Pred^G_\le(v_i)$ and on the right its moralization $\Pred^G_\le(v_i)^\moral$. Note that on the right the neighbours of the node of highest index $\partial_{\Pred^G_\le(v_i)^\moral}(v_i)$ form a complete subgraph of $\Pred^G_\le(v_i)^\moral$ at every step. Also, this holds for every ancestral sub-\HEDG{}.}
 \label{fig:ordering}
\end{figure}
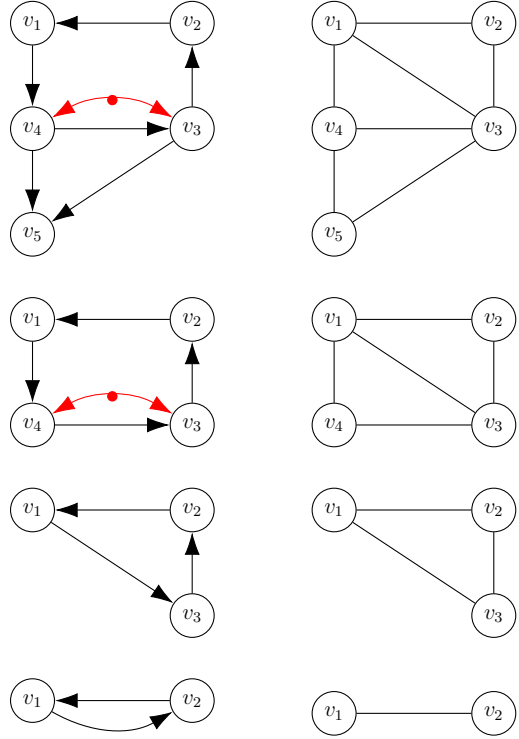

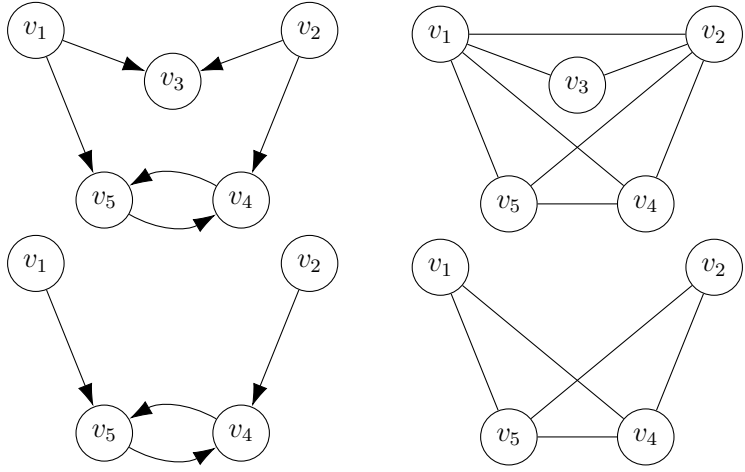
\begin{figure}[!htb] 
\centering
\begin{tikzpicture}[scale=.9, transform shape]
\tikzstyle{every node} = [draw,shape=circle]
\node (v2) at (3,2.5) {$v_2$};
\node (v5) at (0,0) {$v_5$};
\node (v4) at (2,0) {$v_4$};
\node (v1) at (-1,2.5) {$v_1$};
\node (v3) at (1,1.75) {$v_3$};
\foreach \from/\to in {v2/v4, v1/v5, v2/v3, v1/v3}
\draw[-{Latex[length=3mm, width=2mm]}] (\from) -- (\to);
\draw[-{Latex[length=3mm, width=2mm]}, bend right] (v5) to (v4);
\draw[-{Latex[length=3mm, width=2mm]}, bend right] (v4) to (v5);
\end{tikzpicture} \qquad
\begin{tikzpicture}[scale=.9, transform shape]
\tikzstyle{every node} = [draw,shape=circle]
\node (v2) at (3,2.5) {$v_2$};
\node (v5) at (0,0) {$v_5$};
\node (v4) at (2,0) {$v_4$};
\node (v1) at (-1,2.5) {$v_1$};
\node (v3) at (1,1.75) {$v_3$};
\foreach \from/\to in {v2/v4, v1/v5, v2/v3, v1/v3, v4/v5, v1/v4, v2/v5, v1/v2}
\draw[-] (\from) -- (\to);
\end{tikzpicture}\\
\begin{tikzpicture}[scale=.9, transform shape]
\tikzstyle{every node} = [draw,shape=circle]
\node (v2) at (3,2.5) {$v_2$};
\node (v5) at (0,0) {$v_5$};
\node (v4) at (2,0) {$v_4$};
\node (v1) at (-1,2.5) {$v_1$};
\foreach \from/\to in {v2/v4, v1/v5}
\draw[-{Latex[length=3mm, width=2mm]}] (\from) -- (\to);
\draw[-{Latex[length=3mm, width=2mm]}, bend right] (v5) to (v4);
\draw[-{Latex[length=3mm, width=2mm]}, bend right] (v4) to (v5);
\end{tikzpicture} \qquad
\begin{tikzpicture}[scale=.9, transform shape]
\tikzstyle{every node} = [draw,shape=circle]
\node (v2) at (3,2.5) {$v_2$};
\node (v5) at (0,0) {$v_5$};
\node (v4) at (2,0) {$v_4$};
\node (v1) at (-1,2.5) {$v_1$};
\foreach \from/\to in {v2/v4, v1/v5, v4/v5, v1/v4, v2/v5}
\draw[-] (\from) -- (\to);
\end{tikzpicture}
\caption{A directed graph $G$ (top left) and an ancestral subgraph $A$ (bottom left) together with their moralizations on the right. 
We have that $\partial_{\Pred^G_\le(v_i)^\moral}(v_i)$ is a complete subgraph of $\Pred^G_\le(v_i)^\moral$ for all $i=1,\dots,5$.
But $\partial_{\Pred^A_\le(v_5)^\moral}(v_5)$ is not complete in $\Pred^A_\le(v_5)^\moral=A^\moral$.
So the given order is not a perfect elimination order (but an assembling pseudo-topological order).}
\label{fig:wrong-order}
\end{figure}

\begin{Rem}
Let $G=(V,E,H)$ be a \HEDG{} and $<$ a total order on $V$. %
\begin{enumerate}
\item If $v\in V$ is ``childless'' in $\Pred^G_\le(v)$, or more precisely, $\Ch^{\Pred^G_\le(v)}(v) \sm\{v\} = \emptyset$, then for all ancestral $A \ins \Pred^G_\le(v)$ with $v \in A$ we already have that $\partial_{A^\moral}(v) \cup \{v\}$ is a complete subgraph of $A^\moral$. 
So the condition to be a perfect elimination order is a condition on the nodes $v \in V$ with ``proper'' children in $\Pred^G_\le(v)$, e.g.\ nodes in a strongly connected component with two or more nodes. See \ref{anc-complete-neighbourhood} for a more general statement.
\item Not every \HEDG{} has a perfect elimination order (e.g.\ in the directed cycle with 4 nodes every node has two neighbours, which are not connected by an edge in the moralization). 
\item $G$ has a topological order if and only if it is an mDAG, i.e.\ $(V,E)$ is a DAG. In the case of an mDAG every pseudo-topological order is an (assembling) topological order.
\item Every \HEDG{} has an assembling pseudo-topological order. Indeed, from the last to the first index inductively pick an element from one of the smallest strongly connected components $S \in \Scal(G)$ with no outgoing arrows (i.e.\ $\Ch^{\Scal(G)}(S)=\emptyset$ or $\Ch^G(S) \sm S = \emptyset$). %
\item The definition of an (assembling) quasi-topological order of a \HEDG{} merges the definition of a perfect elimination order of an undirected (chordal) graph (see \cite{West01} Def. 5.3.12) with the ancestral structure of a directed graph or \HEDG{}.  
More precisely, in the corner case of a strongly connected \HEDG{} $G$, where we only need to consider one ancestral subset, namely $G$ itself, a perfect elimination order/quasi-topological order of $G$ is the same as a perfect elimination order of its moralization $G^\moral$. 
But as seen above the notion of a quasi-topological order of a \HEDG{} also generalizes the other corner case, the case of a (quasi-)topological order of a DAG.
\end{enumerate}
\begin{proof}
1. If $v \in A$ with $A \ins \Pred^G_\le(v)$ ancestral, then $v$ is also childless in $A$ ($\Ch^A(v) \sm \{v\} = \emptyset)$. We then have $\partial_{A^\moral}(v) \cup \{v\}=\Pa^A(\Dist^A(v)) \cup \Dist^A(v)$, which is a complete subgraph of $A^\moral$.\\
3. It is clear that a topological order does not allow for cycles and every DAG has a topological order (by inductively taking an arbitrary childless node from last to first). A topological order is a perfect elimination order by 1. and clearly a pseudo-topological order.
\end{proof}
\end{Rem}

\begin{Lem}
Let $<$ be a pseudo-topological order on a \HEDG{} $G=(V,E,H)$.
Then for every $v \in V$ we have the inclusion of sets of nodes: 
\[\Pred^G_\le(v) \ins \NonDesc^G(v) \cup \Sc^G(v).\]
\begin{proof}
If $v \ge w$ then $v \notin \Anc^G(w) \sm \Sc^G(w)$.
So $v \in \Sc^G(w)$ or $v \notin \Anc^G(w)$ which implies that $w \in \Sc^G(v)$ or $w \notin \Desc^G(v)$.
It follows $w \in \NonDesc^G(v) \cup \Sc^G(v)$.
\end{proof}
\end{Lem}

\begin{Lem}
\label{anc-complete-neighbourhood}
Let $G=(V,E,H)$ be a \HEDG{} and $v \in V$ a node such that the set $\{v\} \cup \Ch^G(v)$ completely lies in the district $D=\Dist^{\Anc^G(v)}(v)$ of the ancestral sub-\HEDG{} $\Anc^G(v)$ (e.g.\ in the case $\Ch^G(v) \sm \{v\} = \emptyset$). Then for every ancestral sub-\HEDG{} $A \ins G$ with $v \in A$ we have:
\[  \partial_{A^\moral}(v) \cup \{v\} = \Pa^A(\Dist^A(v)) \cup \Dist^A(v),  \]
which is a complete subgraph of $A^\moral$.
\begin{proof}
With $v \in A$ we also have $\{v\} \cup \Ch^G(v) \ins D \ins \Anc^G(v) \ins A$.
It follows that $\Dist^A\lp \{v\} \cup \Ch^G(v) \rp \ins \Dist^A(v)$ and thus:
\[  \partial_{A^\moral}(v) \cup \{v\} = \Pa^A(\Dist^A(v)) \cup \Dist^A(v).\]
\end{proof}
\end{Lem}

\begin{Thm}
\label{peo-general-hedg}
Let $G=(V,E,H)$ be a \HEDG{} and assume that every strongly connected component $S \in \Scal(G)$ of $G$ completely lies in one district of $\Anc^G(S)$
(e.g.\ if $S \in H$ for all $S$). 
Then every pseudo-topological order for $G$ is also a perfect elimination order for $G$ and thus a quasi-topological order.
\begin{proof}
For $v \in V$ we will apply \ref{anc-complete-neighbourhood} to $\Pred^G_\le(v)$: Let $S:=\Sc^{\Pred^G_\le(v)}(v)$. We need to show that $S$ completely lies in one district of $\Anc^{\Pred^G_\le(v)}(S)=\Anc^{\Pred^G_\le(v)}(v)$. 
Let $v_1,v_2 \in S$. Then $v_1,v_2 \in \Sc^G(v)$. By assumption there is a bidirected path in $\Anc^G(v)$ of the form:
\[ v_1 \oto \cdots \oto w \oto \cdots \oto v_2.   \]
Note that $\Anc^{\Pred^G_\le(v)}(v)$ is a marginalization of $\Anc^G(v)$ by \ref{hedg-anc-marg}. %
  For every $w \in \Anc^G(v)$ on the path above consider a directed path $w \to \cdots \to w' \to \cdots \to v$ with $w'$ the first note on that path with $w' \in \Pred^G_\le(v)$. Note that $v \in \Pred^G_\le(v)$ and such a $w'$ always exists.  The corner cases $w'=v$ or $w'=w$ might occur. 
When marginalizing out all nodes outside of $\Pred^G_\le(v)$ the part $\oto w \oto$ becomes $\oto w' \oto$ for every single node $w$ on the bidirected path above. So marginalization gives us a bidirected path:
\[ v_1 \oto \cdots \oto w' \oto \cdots \oto v_2.   \]
It follows that $S$ completely lies in one district $D$ of $\Anc^{\Pred^G_\le(v)}(S)$.
Since $S$ has no outgoing arrows in $\Pred^G_\le(v)$ (by definition of pseudo-topological order) we have that $\{v\} \cup \Ch^{\Pred^G_\le(v)}(v) \ins S \ins D$
 and \ref{anc-complete-neighbourhood} applies. So for every ancestral sub-\HEDG{} $A \ins \Pred^G_\le(v)$ with $v \in A$ we have:
\[  \partial_{A^\moral}(v) \cup \{v\} = \Pa^A(\Dist^A(v)) \cup \Dist^A(v),  \]
which is a complete subgraph of $A^\moral$.
  This shows that the pseudo-topological order $<$ is a perfect elimination order. 
\end{proof}
\end{Thm}

There are perfect elimination orders that do not come from \ref{anc-complete-neighbourhood} or \ref{peo-general-hedg}  (e.g.\ for the directed circle with 3 nodes). Complementary to \ref{peo-general-hedg} we have the following result.

\begin{Thm}
\label{eg-peo}
Let $G=(V,E,H)$ be a \HEDG{} and $<$ a total order on $V$.
For $v \in V$ put $G(v):=\Pred^G_\le(v)$. If for every $v \in V$ we have:
\begin{enumerate}
\item $\Dist^{G(v)}(\{v\} \cup \Ch^{G(v)}(v) ) \ins \Anc^{G(v)}(v)$ as sets of nodes, and
\item $\partial_{\Anc^{G(v)}(v)^\moral}(v)$ is a complete subgraph of $\Anc^{G(v)}(v)^\moral$,
\end{enumerate}
then $<$ is a perfect elimination order for $G$. %
\begin{proof}
Let $A \ins \Pred^G_\le(v)$ be an ancestral sub-\HEDG{} with $v \in A$.
Then we have:
\[v \in \Dist^{G(v)}(\{v\} \cup \Ch^{G(v)}(v) ) \ins \Anc^{G(v)}(v) \ins A \ins \Pred^G_\le(v)=:G(v).\]
  This implies 
\[\Dist^{G(v)}(\{v\} \cup \Ch^{G(v)}(v) ) = \Dist^{A}(\{v\} \cup \Ch^{A}(v) )  \]
und thus 
\[\partial_{G(v)^\moral}(v) = \partial_{A^\moral}(v) =\partial_{\Anc^{G(v)}(v)^\moral}(v)\]
as sets of nodes.
Since $\Anc^{G(v)}(v)$ is an ancestral subgraph of $A$, any complete subgraph of $\Anc^{G(v)}(v)^\moral$ must be a complete subgraph of $A^\moral$.
So the second point implies that $\partial_{A^\moral}(v)$ is a complete subgraph of $A^\moral$.
  This shows that $<$ is a perfect elimination order for $G$.
\end{proof}
\end{Thm}

\begin{Cor}
\label{H-trivial-order}
Let $G=(V,E,H_1)$ be a \HEDG{} with trivial $H_1$ and $<$ a pseudo-topological order on $G$.
If for every $v \in V$ we have that
$\partial_{\Anc^{G(v)}(v)^\moral}(v)$ is a complete subgraph of $\Anc^{G(v)}(v)^\moral$, where $G(v):=\Pred^G_\le(v)$,
then $<$ is also a perfect elimination order and thus a quasi-topological order of $G$.
\begin{proof}
This follows from \ref{eg-peo}. We check the first point:
Put $G(v):=\Pred^G_\le(v)$. 
If $w \in \Ch^{G(v)}(v)$ then $v \in \Anc^G(w)$ by construction of marginalization. Since $<$ is a pseudo-topological order we either have $v < w$ or 
$w \in \Sc^G(v)$. Because $w \le v$ we have $w \in \Sc^G(v)$ and thus $w \in \Sc^{G(v)}(v)$. It follows $\{v\} \cup \Ch^{G(v)}(v) \ins \Anc^{G(v)}(v)$.
Since $H$ is trivial we then get the first point. By \ref{eg-peo} we get the claim.
\end{proof}
\end{Cor}

\begin{Lem}[Total orders under marginalization]
\label{marg-total-ord}
Let $G=(V,E,H)$ be a \HEDG{} and $<$ a total order on $G$. Let $W \ins V$ be a subset and $G^{\marg \sm W}$ the marginalized \HEDG{} w.r.t.\ $W$. 
Then $<$ restricts to a total order on $G^{\marg \sm W}$ with the following property:
If $<$ is a pseudo-topological, quasi-topological, topological, assembling, perfect elimination order, resp., on $G$ so is its restriction $<$ on $G^{\marg \sm W}$.
\begin{proof}
Since marginalization as well as all topological, pseudo-topological, assembling orders respect ancestral relations, the claim is clear for these total orders. So only perfect elimination orders need to be treated.\\
 So lets $<$ be a perfect elimination order of $G$. By induction we can assume that $W=\{w\}$ for one element $w \in V$. Let $v \in V \sm W$ and $A' \ins \Pred_\le^{G^{\marg \sm W}}(v)$ ancestral with $v \in A'$.
By \ref{hedg-anc-marg} $A=\Anc^{\Pred^G_\le(v)}(A')$ is ancestral in $\Pred^G_\le(v)$ and marginalizes to $A'$. By assumption $\partial_{A^\moral}(v)$ is a complete subgraph of $A^\moral$. We now have the two cases: $w \in A$ and $w \notin A$.
If $w \notin A$ then $A=A'$ by \ref{anc-sub=mar} and we have nothing to show. So assume $w \in A$ from now on.
Then $\Ch^A(w) \sm \{w\} \neq \emptyset$. Otherwise $A \sm \{w\}$ would already be ancestral and marginalize to $A'$ and we would be done by the previous argument. So by $\Ch^A(w) \sm \{w\} \neq \emptyset$ and \ref{marg-mor2} we get:
\[ (A')^\moral = (A^\moral)^{\marg \sm W}.  \] 
So we have $\partial_{(A')^\moral}(v)=\partial_{(A^\moral)^{\marg \sm W} }(v)$.
Furthermore, by \ref{compl-subgr-marg} $\partial_{(A^\moral)^{\marg \sm W} }(v)$ is a complete subgraph in $(A^\moral)^{\marg \sm W}$ and thus 
in $(A')^\moral$. So the claim follows. 
\end{proof}
\end{Lem}

\begin{Cor}
Let $G=(V,E,H)$ be a \HEDG{} and $G^\aug$ its augmentation and $<$ a pseudo-topological order on $G^\aug$ such that
for all $v \in V^\aug$ we have that:
$\partial_{\Anc^{G^\aug(v)}(v)^\moral}(v)$ is a complete subgraph of $\Anc^{G^\aug(v)}(v)^\moral$, where $G^\aug(v):=\Pred^{G^\aug}_\le(v)$,
then $<$ is also a perfect elimination order and thus a quasi-topological order for both $G^\aug$ and $G$.
\begin{proof}
This directly follows from \ref{H-trivial-order} and \ref{marg-total-ord}.
\end{proof}
\end{Cor}

\subsection{The Acyclification of a \HEDG{}}

In this subsection we generalize the idea of a \emph{collapsed graph} for directed graphs developed in \cite{Spirtes94} to a version that also works for \HEDG{}es $G$ and that is called \emph{acyclification}. Surprisingly, as soon as we allow for hyperedges the construction becomes much more natural (e.g., it does not depend on a total ordering of the nodes) and later proofs relating structural equations properties (\hyl{lsSEP}{SEP}) to other Markov properties for $G$ become much more general, simple and elegant at the same time (see \ref{csSEP-smgdGMP}).
Furthermore, the d-separation relations in the acyclification $G^\acy$ of $G$ can directly be translated to a new separation criterion in $G$ itself that we introduce in the next subsection and refer to as \emph{$\sigma$-separation}.

\begin{Def}[The acyclification of a \HEDG{}]
Let $G=(V,E,H)$ be a \HEDG{}. The \emph{acyclification} of $G$ is the \HEDG{} $G^\acy=(V^\acy,E^\acy,H^\acy)$ which is defined as follows:
\begin{enumerate}
\item $V^\acy:=V$,
\item Put $v \to w \in E^\acy$ if and only if $v \notin \Sc^G(w)$ and there is a node $w' \in \Sc^G(w)$ with $v \to w' \in E$.
\item $H^\acy:=\{ F' \ins \bigcup_{v \in F} \Sc^G(v)| F \in H\}$,
\end{enumerate}
This means that we draw edges from a node $v$ to all nodes of a strongly connected component $S \in \Scal(G)$ if there was at least one edge from $v$ to a node of $S$, all edges between nodes of a strongly connected component are erased, and hyperedges are extended to cover the union of the strongly connected components of their nodes.
\end{Def}

\begin{figure}
\centering
\begin{tikzpicture}[scale=.7, transform shape]
\tikzstyle{every node} = [draw,shape=circle]
\node (v1) at (0,0) {$v_1$};
\node (v2) at (0,-2) {$v_2$};
\node (v3) at (-2,-4) {$v_3$};
\node (v4) at (2,-4) {$v_4$};
\node (v6) at (-2,-6) {$v_6$};
\node (v5) at (-4,-6) {$v_5$};
\node (v7) at (2,-6) {$v_7$};
\node (v8) at (4,-6) {$v_8$};
\draw[{Latex[length=3mm, width=2mm]}-{Latex[length=3mm, width=2mm]}, red, bend left] (v6) to node[fill,circle,,inner sep=0pt,minimum size=5pt] {} (v7);
\draw[-{Latex[length=3mm, width=2mm]},out=135, in=225,looseness=8] (v3) to (v3);
\draw[-{Latex[length=3mm, width=2mm]}, bend right] (v5) to (v6);
\draw[-{Latex[length=3mm, width=2mm]}, bend right] (v6) to (v5);
\draw[-{Latex[length=3mm, width=2mm]}, bend right] (v7) to (v8);
\draw[-{Latex[length=3mm, width=2mm]}, bend right] (v8) to (v7);
\foreach \from/\to in {v1/v2, v2/v3, v3/v4, v4/v2, v3/v6, v4/v8}
\draw[-{Latex[length=3mm, width=2mm]}] (\from) -- (\to);
\end{tikzpicture}
\hspace{1cm}
\begin{tikzpicture}[scale=.7, transform shape]
\tikzstyle{every node} = [draw,shape=circle]
\node (v1) at (0,0) {$v_1$};
\node (v2) at (0,-2) {$v_2$};
\node (v3) at (-2,-4) {$v_3$};
\node (v4) at (2,-4) {$v_4$};
\node (v6) at (-2,-6) {$v_6$};
\node (v5) at (-4,-6) {$v_5$};
\node (v7) at (2,-6) {$v_7$};
\node (v8) at (4,-6) {$v_8$};
\node[fill,circle,red,inner sep=0pt,minimum size=5pt] (v9) at (0,-4.5) {};
\node[fill,circle,red,inner sep=0pt,minimum size=5pt] (v10) at (0,-3.2) {};
\draw[-{Latex[length=3mm, width=2mm]}, red] (v9) to (v6);
\draw[-{Latex[length=3mm, width=2mm]}, red] (v9) to (v5);
\draw[-{Latex[length=3mm, width=2mm]}, red] (v9) to (v7);
\draw[-{Latex[length=3mm, width=2mm]}, red] (v9) to (v8);
\draw[-{Latex[length=3mm, width=2mm]}, red] (v10) to (v2);
\draw[-{Latex[length=3mm, width=2mm]}, red] (v10) to (v3);
\draw[-{Latex[length=3mm, width=2mm]}, red] (v10) to (v4);
\foreach \from/\to in {v1/v2, v1/v3, v1/v4, v3/v6, v3/v5, v4/v7, v4/v8}
\draw[-{Latex[length=3mm, width=2mm]}] (\from) -- (\to);
\end{tikzpicture}
\caption{A \HEDG{} on the left and its acyclification on the right.}
 \label{fig:acyclification}
\end{figure}
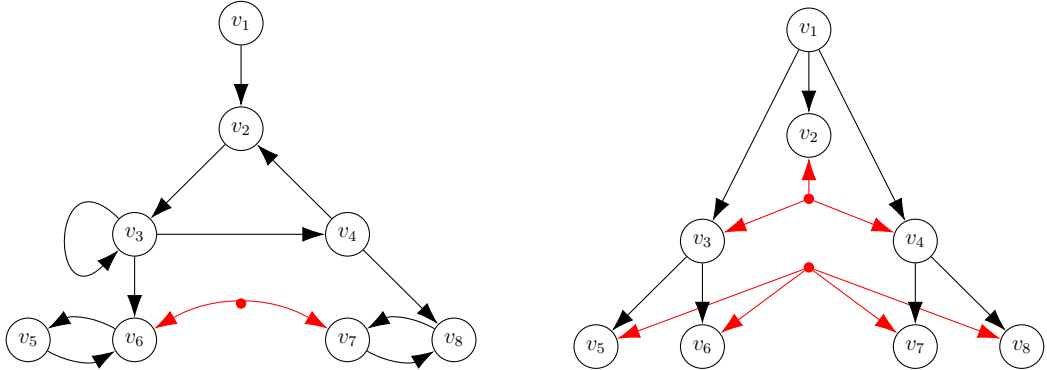

\begin{Lem}
The acyclification $G^\acy=(V,E^\acy,H^\acy)$ of a \HEDG{} $G=(V,E,H)$ is an mDAG, i.e.\ its underlying directed graph $(V,E^\acy)$ is acyclic, i.e.\ a DAG.
\begin{proof}
This follows directly from the fact that the strongly connected components of $G$ form a DAG by \ref{sc-dag} and the edges in $E^\acy$ only go along this DAG.
\end{proof}
\end{Lem}

\begin{Rem}
  Instead of considering the acyclification $G^\acy$ of a \HEDG{} $G=(V,E,H)$ one could consider taking the mDAG of strongly connected components $\Scal(G)$. The reason we focus on the acyclification $G^\acy$ instead is that the acyclification $G^\acy$ carries more information about $G$ than $\Scal(G)$. More precisely, $G^\acy$ still carries the information about which of the nodes of a strongly connected component of $G$ has an edge to another strongly connected component and which not. The mDAG of strongly connected components $\Scal(G)$ of $G$ instead would just identify all nodes of a strongly connected component. This means in perspective of the next chapter that stronger Markov properties will be obtained.
\end{Rem}

\begin{Lem}
Every pseudo-topological order for a \HEDG{} $G$ is a topological order for its acyclification $G^\acy$.
\begin{proof}
This follows directly from the fact that for every $v \in V$ we have that $\Pa^{G^\acy}(v) \ins \Anc^G(v) \sm  \Sc^G(v)$.
\end{proof}
\end{Lem}

\begin{Def}[The acyclic augmentation of a \HEDG{}]
\label{acag-3-def}
Let $G=(V,E,H)$ be a \HEDG{}. Then there are several constructions, where both the augmentation and the acyclification are involved: e.g.\ $(G^\aug)^\acy$ and $(G^\acy)^\aug$. A third one is the \emph{acyclic augmentation} $G^\acag=(V^\aug,(E^\aug)^\acy,H_1)$ of $G$, where we first augment $G$ and then ``break'' the cycles, but do not allow for additional hyperedges, i.e.\ we replace $(H^\aug)^\acy$ with the trivial hyperedges $H_1$. Note that $G^\acag$ is a directed acyclic graph (DAG) and contains all the nodes from the augmentation of $G$.
\end{Def}

\begin{figure}
\centering
\hspace{0.5cm}\\
\vspace{0.5cm}
\begin{tikzpicture}[scale=.7, transform shape]
\tikzstyle{every node} = [draw,shape=circle]
\node (v1) at (0,0) {$v_1$};
\node (v2) at (2,0) {$v_2$};
\node (v3) at (-2,-2) {$v_3$};
\node (v4) at (4,-2) {$v_4$};
\draw[-{Latex[length=3mm, width=2mm]}, bend right] (v1) to (v2);
\draw[-{Latex[length=3mm, width=2mm]}, bend right] (v2) to (v1);
\foreach \from/\to in {v1/v3, v2/v4}
\draw[-{Latex[length=3mm, width=2mm]}] (\from) -- (\to);
\end{tikzpicture}
\hspace{0.5cm}
\begin{tikzpicture}[scale=.7, transform shape]
\tikzstyle{every node} = [draw,shape=circle]
\node (v1) at (0,0) {$v_1$};
\node (v2) at (2,0) {$v_2$};
\node (v3) at (-2,-2) {$v_3$};
\node (v4) at (4,-2) {$v_4$};
\node (e1) at (0,2) {$e_{\{1\}}$};
\node (e2) at (2,2) {$e_{\{2\}}$};
\node (e3) at (-2,0) {$e_{\{3\}}$};
\node (e4) at (4,0) {$e_{\{4\}}$};
\draw[-{Latex[length=3mm, width=2mm]}, bend right] (v1) to (v2);
\draw[-{Latex[length=3mm, width=2mm]}, bend right] (v2) to (v1);
\foreach \from/\to in {e1/v1, e2/v2, e3/v3, e4/v4, v1/v3, v2/v4}
\draw[-{Latex[length=3mm, width=2mm]}] (\from) -- (\to);
\end{tikzpicture}
\hspace{0.5cm}\\
\vspace{0.5cm}
\begin{tikzpicture}[scale=.7, transform shape]
\tikzstyle{every node} = [draw,shape=circle]
\node (v1) at (0,0) {$v_1$};
\node (v2) at (2,0) {$v_2$};
\node (v3) at (-2,-2) {$v_3$};
\node (v4) at (4,-2) {$v_4$};
\foreach \from/\to in {v1/v3, v2/v4}
\draw[-{Latex[length=3mm, width=2mm]}] (\from) -- (\to);
\draw[{Latex[length=3mm, width=2mm]}-{Latex[length=3mm, width=2mm]},bend left, red] (v1) to node[fill,circle,red,inner sep=0pt,minimum size=5pt] {} (v2);
\end{tikzpicture}
\hspace{0.5cm}
\begin{tikzpicture}[scale=.7, transform shape]
\tikzstyle{every node} = [draw,shape=circle]
\node (v1) at (0,0) {$v_1$};
\node (v2) at (2,0) {$v_2$};
\node (v3) at (-2,-2) {$v_3$};
\node (v4) at (4,-2) {$v_4$};
\node (e1) at (0,2) {$e_{\{1\}}$};
\node (e2) at (2,2) {$e_{\{2\}}$};
\node (e3) at (-2,0) {$e_{\{3\}}$};
\node (e4) at (4,0) {$e_{\{4\}}$};
\foreach \from/\to in {e1/v1, e2/v2, e3/v3, e4/v4, v1/v3, v2/v4, e1/v2, e2/v1}
\draw[-{Latex[length=3mm, width=2mm]}] (\from) -- (\to);
\draw[{Latex[length=3mm, width=2mm]}-{Latex[length=3mm, width=2mm]},bend right, red] (v1) to node[fill,circle,red,inner sep=0pt,minimum size=5pt] {} (v2);
\end{tikzpicture}
\hspace{0.5cm}\\
\vspace{0.5cm}
\begin{tikzpicture}[scale=.7, transform shape]
\tikzstyle{every node} = [draw,shape=circle]
\node (v1) at (0,0) {$v_1$};
\node (v2) at (2,0) {$v_2$};
\node (v3) at (-2,-2) {$v_3$};
\node (v4) at (4,-2) {$v_4$};
\node (e12) at (1,2) {$e_{\{1,2\}}$};
\node (e3) at (-2,0) {$e_{\{3\}}$};
\node (e4) at (4,0) {$e_{\{4\}}$};
\foreach \from/\to in {e12/v1, e12/v2, e3/v3, e4/v4, v1/v3, v2/v4}
\draw[-{Latex[length=3mm, width=2mm]}] (\from) -- (\to);
\end{tikzpicture}
\hspace{0.5cm}
\begin{tikzpicture}[scale=.7, transform shape]
\tikzstyle{every node} = [draw,shape=circle]
\node (v1) at (0,0) {$v_1$};
\node (v2) at (2,0) {$v_2$};
\node (v3) at (-2,-2) {$v_3$};
\node (v4) at (4,-2) {$v_4$};
\node (e1) at (0,2) {$e_{\{1\}}$};
\node (e2) at (2,2) {$e_{\{2\}}$};
\node (e3) at (-2,0) {$e_{\{3\}}$};
\node (e4) at (4,0) {$e_{\{4\}}$};
\foreach \from/\to in {e1/v1, e2/v2, e3/v3, e4/v4, v1/v3, v2/v4, e1/v2, e2/v1}
\draw[-{Latex[length=3mm, width=2mm]}] (\from) -- (\to);
\end{tikzpicture}
\caption{ 
A \HEDG{} $G$ (top left), its augmentation $G^\aug$ (top right), its acyclification $G^\acy$ (middle left), also $(G^\aug)^\acy$ (middle right), $(G^\acy)^\aug$ (bottom left) and the acyclic augmentation $G^\acag=(V^\aug, (E^\aug)^\acy)$ (bottom right). 
We have $\{v_3\} \Indep_{G^\acag}^d \{v_4\} \given \{e_{\{1\}},e_{\{2\}}\}$ and $\{v_3\} \Indep_{(G^\acy)^\aug}^d \{v_4\} \given \{e_{\{1,2\}}\}$, 
but $\{v_3\} \nIndep_{(G^\aug)^\acy}^d \{v_4\} \given \{e_{\{1\}},e_{\{2\}}\}$ and $\{v_3\} \nIndep_{G^\aug}^d \{v_4\} \given \{e_{\{1\}},e_{\{2\}}\}$. }
 \label{fig:mdGMP-smgdGMP}
\end{figure}
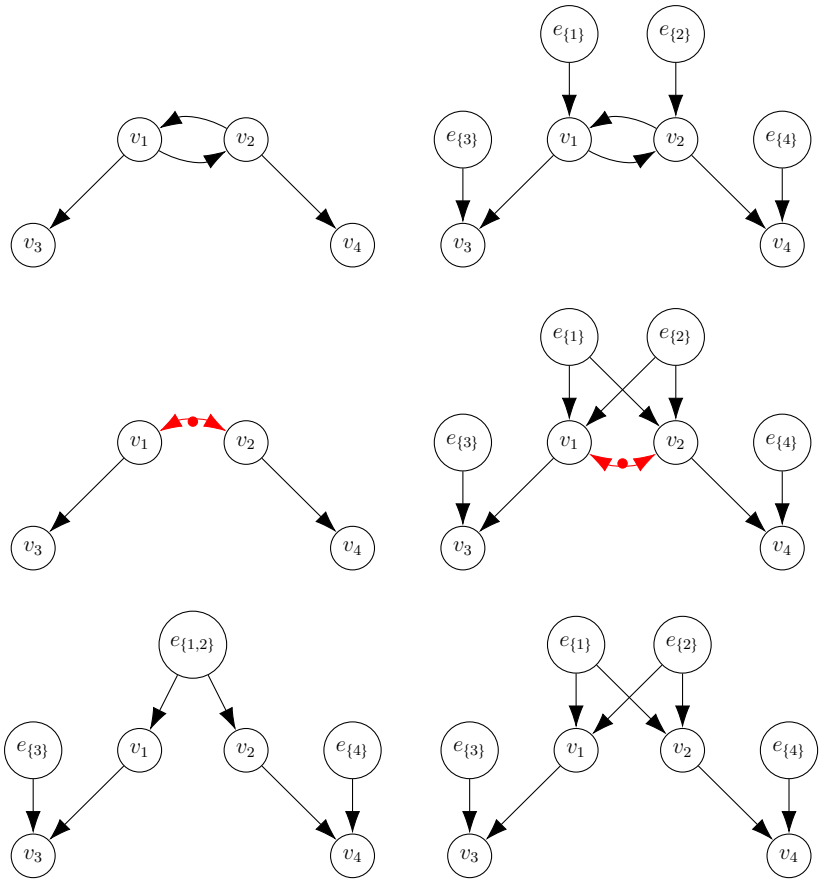

\begin{Rem}
\label{rem-acy-acag}
Note that for a \HEDG{} $G=(V,E,H)$ the three constructions  $(G^\aug)^\acy$, $(G^\acy)^\aug$ and $G^\acag$ are all different in general (see figure \ref{fig:mdGMP-smgdGMP}). But their marginalization to the nodes $V$ is in all cases the acyclification $G^\acy$ of $G$. 
The purpose to analyze three different types of augmentations in this setting is that we are interested in the strongest (marginal) Markov property that holds for structural equations (see \ref{csSEP-smgdGMP}). This is the reason we introduced $G^\acag$ separately.
\end{Rem}

\begin{Lem}
\label{acag-3-d-sep}
Let $G=(V,E,H)$ be a \HEDG{}. Consider the three constructions from \ref{acag-3-def} (see figure \ref{fig:mdGMP-smgdGMP}): $(G^\aug)^\acy$, $(G^\acy)^\aug$ and $G^\acag$.
Then we have:
\begin{enumerate}
\item $G^\acag=(V^\aug, (E^\aug)^\acy,H_1)$ and
\item $(G^\aug)^\acy = (V^\aug,(E^\aug)^\acy,(H^\aug)^\acy)$ with
\item[] $V^\aug = V \dot\cup \tilde H$ and $(E^\aug)^\acy = E^\acy \dot\cup \{ F \to v|v \in F \in \tilde{H} \}$ and 
 $(H^\aug)^\acy = \{ F' \ins \Sc^{G^\aug}(v)| v \in V^\aug\}$.%
\item $(G^\acy)^\aug = ((V^\acy)^\aug, (E^\acy)^\aug, H_1)$ with 
\item[] $(V^\acy)^\aug = V \dot\cup \tilde H^\acy$ and $(E^\acy)^\aug = E^\acy \dot\cup \{ F' \to v|v \in F' \in \tilde{H}^\acy \}$.
\item For $X,Y,Z \ins V^\aug$ we have:
 \[  X \Indep_{(G^\aug)^\acy}^d Y \given Z \;\implies\; X \Indep_{G^\acag}^d Y \given Z.  \]
\item Let $\varphi:\, \tilde H \to \tilde H^\acy$ be a surjective map with $\varphi(F) \supseteq F$ for all $F \in \tilde H$ (see proof for existence) and 
\[\begin{array}{rclcl}
\tilde \varphi =(\id\dot\cup\varphi)&:& V^\aug=V \dot\cup \tilde H &\to& V \dot\cup \tilde H^\acy= (V^\acy)^\aug,\\
&&(E^\aug)^\acy & \to &(E^\acy)^\aug. 
\end{array} \]
 its natural extensions.
For $X',Y',Z' \ins (V^\acy)^\aug$ put $X:=(\tilde \varphi)^{-1}(X')$, $Y:=(\tilde \varphi)^{-1}(Y')$ and $Z:=(\tilde \varphi)^{-1}(Z')$.
Then we have:
\[ %
X' \Indep_{(G^\acy)^\aug}^d Y' \given Z' \;\implies\; X \Indep_{G^\acag}^d Y \given Z,  \]
\end{enumerate}
\begin{proof}
The first three points are clear from the definition of augmentation and acyclification.\\
4.) The only difference between $(G^\aug)^\acy$ and $G^\acag$ lies in the hyperedges. So clearly every $Z$-open path in $G^\acag$ is also a $Z$-open path in $(G^\aug)^\acy$. This implies the claim for all $X,Y,Z \ins V^\aug$:
\[  X \Indep_{(G^\aug)^\acy}^d Y \given Z \;\implies\; X \Indep_{G^\acag}^d Y \given Z.  \]
5.) First fix a map $\varphi: \tilde H \to \tilde H^\acy$ with $\varphi(F) \supseteq F$ for every $F \in \tilde H$. Note that such a map $\varphi$ exists since $H \ins H^\acy$ by construction of the acyclification. Furthermore, $\varphi$ is surjective. Indeed if $F' \in \tilde H^\acy$ then $F'$ is of the form $F'=\bigcup_{v \in \hat F} \Sc^G(v)$ for some $\hat F \in H$. Let $F \in \tilde H$ with $\hat F \ins F$.
Together with $F \ins \varphi(F)$ we get:
\[ F'=\bigcup_{v \in \hat F} \Sc^G(v) \ins \bigcup_{v \in F} \Sc^G(v) \ins \varphi(F). \]
Since $F'$ is inclusion maximal all inclusions are equalities and we get $\varphi(F)=F'$. This shows the surjectivity of $\varphi$.\\
Now let $\pi$ be a $Z$-open path between $x \in X$ and $y \in Y$ in $G^\acag$. %
We now apply $\tilde\varphi$ on every node and edge in $\pi$ to get a path $\tilde \varphi (\pi)$ in $(G^\acy)^\aug$.
Every node in $V$ and every edge between nodes in $V$ (given by $E^\acy$) stay the same.
Every node $F$ in $\pi$ from $V^\aug \sm V=\tilde H$ can only occur as a non-collider with $v_1 \ot F \to v_2$ and $v_1,v_2 \in F$ (or as an endnode, resp.).
So $F \notin Z$ by the $Z$-openness of $\pi$. Since $Z=(\tilde \varphi)^{-1}(Z')$ we have $\varphi(F) \notin Z'$ and $v_1,v_2 \in F \ins \varphi(F)$. So we have the edges $v_1 \ot \varphi(F) \to v_2$ also in $(G^\acy)^\aug$. Furthermore, this will be a $Z'$-open non-collider (and similarly for $F$ an endnode). Also note that $x \in X=(\tilde \varphi)^{-1}(X')$ and $y \in Y=(\tilde \varphi)^{-1}(Y')$ implies $\tilde\varphi(x) \in X$ and $\tilde\varphi(y) \in Y$. So $\tilde\varphi(\pi)$ is a path from $X'$ to $Y'$ in $(G^\acy)^\aug$ and is $Z'$-open. This shows the claim.
\end{proof}
\end{Lem}

\subsection{\texorpdfstring{$\sigma$}{Sigma}-Separation in \HEDG{}es}

In this subsection we introduce a more restrictive version of d-separation for \HEDG{}es, which we call $\sigma$-separation. 
We will later show (see \ref{csSEP-smgdGMP}) that quite general structural equations models (\hyl{csSEP}{csSEP}) for \HEDG{}es always follow a directed global Markov property (\hyl{gdGMP}{gdGMP}) based on $\sigma$-separation.
Here we introduce several approaches to $\sigma$-separation and show that $\sigma$-separation, similarly to d-separation, is stable under marginalization. In addition, on mDAGs (e.g.\ ADMGs and DAGs) we will see that d-separation and $\sigma$-separation are equivalent. From this point of view, $\sigma$-separation can be interpreted as the ``non-naive'' generalization of d-separation from DAGs to \HEDG{}es.
We will also investigate more general relations between $\sigma$-separation in $G$ and d-separation in some graphical transformations of $G$ (e.g.\ its acyclification).

\begin{Def}[$\sigma$-separation (segment version)]
\label{s-sep-def}
Let $G=(V,E,H)$ be a \HEDG{} and $X,Y,Z \ins V$ subsets of the nodes.
\begin{enumerate}
\item Consider a path in $G$ with $n \ge 1$ nodes:
\[\begin{array}{ccccccc} 
& \ot &&\ot \\[-5pt]
v_1 & \to & \cdots & \to &v_n. \\[-5pt]
& \oto & &\oto 
\end{array}\]
Then this path can be uniquely partitioned according to the strongly connected components of $G$:
\[\begin{array}{lllllllllllll} 
&\ot&  & \ot &     & \ot &      & \ot&   & \ot & &\ot \\[-5pt]
v_{i-1}&\to&v_i & \to & v_{i+1} &\to  &\cdots& \to&v_{k-1}& \to & v_{k}&\to & v_{k+1}, \\[-5pt]
&\oto&  & \oto &    &\oto&       &\oto&   &\oto && \oto
\end{array}\]
with  $v_i, \dots, v_{k} \in \Sc^G(v_i)$ and $v_{i-1},v_{k+1} \notin \Sc^G(v_i)$. Note that $v_{i-1}$ or  $v_{k+1}$ might not appear if $v_i$ or $v_k$ is an endnode of the path. We will call the sub-path $\sigma_j$ (given by the nodes $v_i,\dots,v_k$ and its corresponding edges) a \emph{segment} of the path. We abbreviate the left and right endnode of $\sigma_j$ with $\sigma_{j,l}:=v_i$ and $\sigma_{j,r}:=v_k$. The path can then uniquely be written with its segments:
\[\begin{array}{ccccccc} 
& \ot &&\ot \\[-5pt]
\sigma_1 & \to & \cdots & \to &\sigma_m. \\[-5pt]
& \oto & &\oto 
\end{array}\]
We will call $\sigma_1$ and $\sigma_m$ the \emph{end-segments} of the path.
\item Such a path will be called \emph{$Z$-$\sigma$-blocked} or \emph{$\sigma$-blocked by $Z$} if:
\begin{enumerate}
  \item at least one of the endnodes $v_1=\sigma_{1,l}$, $v_n=\sigma_{m,r}$ is in $Z$, or
  \item there is a segment $\sigma_j$ (e.g.\ $\sigma_1$ or $\sigma_m$) with an outgoing directed edge in the path (i.e.\ $\ot \sigma_j \cdots$ or $\cdots \sigma_j \to$, resp.) and its corresponding endnode (i.e.\ $\sigma_{j,l}$ or $\sigma_{j,r}$, resp.) lies in $Z$, or
	 \item there is a segment $\sigma_j$ with two adjacent (hyper)edges that form a \emph{collider} $\genfrac{}{}{0pt}{}{\to}{\oto} \sigma_j \genfrac{}{}{0pt}{}{\ot}{\oto}$ 
     and $\Sc^G(\sigma_j) \cap 	\Anc^G(Z)=\emptyset$.\footnote{As long as we allow for repeated nodes in a path we can relax this to $\Sc^G(\sigma_j) \cap Z = \emptyset$ or $\sigma_j \cap \Anc^G(Z)=\emptyset$ or $\sigma_j \cap Z=\emptyset$. But we then need to make one consistent choice for every path in the definition of $\sigma$-separation below. 
		}
	 \end{enumerate}
If none of the above holds then the path is called \emph{$Z$-$\sigma$-open} or \emph{$Z$-$\sigma$-active}.
\item  We say that $X$ is \emph{$\sigma$-separated} from $Y$ given $Z$ if every path in $G$ with one endnode in $X$ and one endnode in $Y$ is $\sigma$-blocked by $Z$.
In symbols this will be written as follows:
\[ X \Indep_{G}^\sigma Y \given Z.\]
 \end{enumerate}
\end{Def}

\begin{figure}
\centering
\begin{tikzpicture}[scale=.7, transform shape]
\tikzstyle{every node} = [draw,shape=circle]
\node (v1) at (0,0) {$v_1$};
\node (v2) at (3,0) {$v_2$};
\node (v3) at (6,0) {$v_3$};
\node (v4) at (1.5,-1.5) {$v_4$};
\node (v5) at (4.5,-1.5) {$v_5$};
\node (v6) at (3,-3) {$v_6$};
\foreach \from/\to in {v1/v4, v4/v2, v2/v5, v3/v5, v5/v6, v6/v4}
\draw[-{Latex[length=3mm, width=2mm]}] (\from) -- (\to);
\end{tikzpicture}
\caption{We have $\{v_1\} \Indep_G^d \{v_3\} \given \{v_2,v_6\}$ but $\{v_1\} \nIndep_G^\sigma \{v_3\} \given \{v_2,v_6\}$.}
 \label{fig:s-separation}
\end{figure}

\begin{Thm}
\label{s-sep-d-sep}
Let $G=(V,E,H)$ be a \HEDG{}, $G^\acy=(V,E^\acy,H^\acy)$ its acyclification and $X,Y,Z \ins V$ subsets of the nodes.
Let $G'=(V,E',H')$ be a \HEDG{} with the following properties:
\begin{enumerate}
 \item $E' \ins E^\acy \cup E^\mathrm{sc}$ with  $E^\mathrm{sc}:=\{v \to w\,|\, \forall v \in V\, w \in \Sc^G(v) \}$,
 \item $H' \ins H^\acy$.
\end{enumerate}
Then we have the implication:
\[X \Indep_{G}^\sigma Y \given Z \;\implies\; X \Indep_{G'}^d Y \given Z.\]
In particular, we have this for $G'=G$, showing that $\sigma$-separation implies d-separation.
\begin{proof}
Let $P'=v_1 \cdots v_n$ be a $Z$-open path from $X$ to $Y$ in $G'$. Then $v_1,v_n \notin Z$.\\
  Every edge $v \to w$ (or $w \ot v$) in the path $P'$ with $v \notin \Sc^G(w)$ lies in $E^\acy$. So there is an edge $v \to w'$ in $G$ with 
  $w' \in \Sc^G(w)$. So in $G$ we have a directed path $v \to w' \to \cdots \to  w$ (or $w \ot  \cdots \ot w' \ot  v$, resp.).\\
Now consider a (hyper)edge $v \begin{array}{c}\ot\\[-10pt]\to\\[-10pt]\oto\end{array} w$ in the path $P'$ in $G'$ with $v \in \Sc^G(w)$.
  Then we fix an arbitrary directed path in $\Sc^G(w)$ between $v$ and $w$, e.g.\ $v \to\cdots\to w$.\\
In the case we have a bidirected edge $v \oto w$ in the path $P'$ with $v \notin \Sc^G(w)$ then $v \oto w \in H^\acy$. So there is a bidirected edge $v' \oto w'$ in $G$ with $v' \in \Sc^G(v)$ and $w' \in \Sc^G(w')$. In $G$ we then can fix a path
$v  \ot \cdots \ot  v' \oto w' \to \cdots \to w$.\\
If we replace all (hyper)edges in $P'$ with the above mentioned (hyper)edges in $G$ we get a path $P$ in $G$. It is left to show that this path $P$ is $Z$-$\sigma$-open in $G$.\\
First note that in the replacements the left and right nodes stay the same. So the endnodes are not in $Z$ since $P'$ is $Z$-open.\\
  If a segment $\sigma_j$ in the new path $P$ has an outgoing directed edge then the corresponding endnode $v$ of $\sigma_j$ had an outgoing directed edge in $P'$ to a node $w$ with $v \notin \Sc^G(w)$ as well. So we have $v\notin Z$ by the $Z$-openness of $P'$.\\
Now consider a segment $\sigma_j$ in the new path $P$ that is a collider 
$\genfrac{}{}{0pt}{}{\to}{\oto} \sigma_j \genfrac{}{}{0pt}{}{\ot}{\oto}$.
  Then all the corresponding nodes $v_i, \dots, v_k$ in $P'$ that induce $\sigma_j$ by the constructed replacements from above come from a subpath $\genfrac{}{}{0pt}{}{\to}{\oto} v_i \cdots v_k \genfrac{}{}{0pt}{}{\ot}{\oto}$ in $P'$ with the same edges left and right as for $\sigma_j$. 
It is easy to see that a path with two inward pointing arrow heads at its ends must have at least one node $v$ (i.e.\ from $v_i,\dots,v_k$ for $P'$) that is a collider itself on the path. 
  It follows that $v \in \Anc^{G'}(Z)$ by the $Z$-openness of $P'$. Since in the replacements directed paths give directed paths
	we see that $v \in \sigma_j \cap \Anc^G(Z)$. Thus $\sigma_j$ is a $Z$-open collider in the new path $P$.\\ 
Together this shows that $P$ is $Z$-$\sigma$-open and thus the claim.
\end{proof}
\end{Thm}

\begin{Thm}
\label{s-sep-d-sep-general}
Let $G=(V,E,H)$ be a \HEDG{} and $G'=(V,E',H')$ be a \HEDG{} such that for every subset $Z \ins V$ 
and every $Z$-$\sigma$-open path $P$ in $G$ given by its segments
\[\begin{array}{ccccccc} 
& \ot &&\ot \\[-5pt]
\sigma_1 & \to & \cdots & \to &\sigma_m, \\[-5pt]
& \oto & &\oto 
\end{array}\]
at least one of the following replacements (in the table on the right from top to bottom) is always possible in $G'$ for every $\sigma_j$ of the following in $P$ occurring form:
\[\begin{array}{r|llll}
\text{ in } P & \text{replacement in } P' \\
\hline
\begin{array}{c}\to\\[-10pt]\oto\end{array} \sigma_j \begin{array}{c}\to\end{array} & \begin{array}{c}\to\\[-10pt]\oto\end{array} \sigma_{j,r}\begin{array}{c}\to\end{array} &\text{if} &\text{existent},&\text{or} \\
& \begin{array}{c}\to\\[-10pt]\oto\end{array}  \cdot  \oto \cdots \oto  \sigma_{j,r}\begin{array}{c}\to\end{array}&\text{if}&\Sc^G(\sigma_j) \cap \Anc^G(Z) \neq \emptyset, &\text{or} \\
& \begin{array}{c}\to\\[-10pt]\oto\end{array}  \cdot  \to \cdots \to  \sigma_{j,r}\begin{array}{c}\to\end{array}&\text{if}&\Sc^G(\sigma_j) \cap \Anc^G(Z) = \emptyset, \\
\begin{array}{c}\ot\end{array}\sigma_j \begin{array}{c}\ot\\[-10pt]\oto\end{array}  & \begin{array}{c}\ot\end{array}\sigma_{j,l} \begin{array}{c}\ot\\[-10pt]\oto\end{array} &\text{if} &\text{existent},&\text{or} \\
& \begin{array}{c}\ot\end{array}\sigma_{j,l} \oto \cdots \oto  \cdot \begin{array}{c}\ot\\[-10pt]\oto\end{array} &\text{if} & \Sc^G(\sigma_j) \cap \Anc^G(Z) \neq \emptyset,&\text{or}\\
& \begin{array}{c}\ot\end{array}\sigma_{j,l} \ot \cdots \ot  \cdot \begin{array}{c}\ot\\[-10pt]\oto\end{array}&\text{if}&\Sc^G(\sigma_j) \cap \Anc^G(Z) = \emptyset, \\
\begin{array}{c}\ot\end{array} \sigma_j \begin{array}{c}\to\end{array} & \begin{array}{c}\ot\end{array} \sigma_{j,l}  \begin{array}{c}=\\[-10pt]\to\\[-10pt]\ot\\[-10pt]\oto\end{array} \sigma_{j,r} \begin{array}{c}\to\end{array} &\text{if} &\text{existent},&\text{or}  \\
& \begin{array}{c}\ot\end{array}\sigma_{j,l} \oto \cdots \oto  \sigma_{j,r}\begin{array}{c}\to\end{array}&\text{if} & \Sc^G(\sigma_j) \cap \Anc^G(Z) \neq \emptyset,&\text{or}\\
& \begin{array}{c}\ot\end{array}\sigma_{j,l} \to \cdots \to  \sigma_{j,r}\begin{array}{c}\to\end{array}&\text{if} & \Sc^G(\sigma_j) \cap \Anc^G(Z) = \emptyset, &\text{or}\\
& \begin{array}{c}\ot\end{array}\sigma_{j,l} \ot \cdots \ot  \sigma_{j,r}\begin{array}{c}\to\end{array}&\text{if} & \Sc^G(\sigma_j) \cap \Anc^G(Z) = \emptyset, \\
\begin{array}{c}\to\\[-10pt]\oto\end{array} \sigma_j \begin{array}{c}\ot\\[-10pt]\oto\end{array} & \begin{array}{c}\to\\[-10pt]\oto\end{array} \cdot \oto \cdots \oto \cdot \begin{array}{c}\ot\\[-10pt]\oto\end{array} &\text{for} & \Sc^G(\sigma_j) \cap \Anc^G(Z) \neq \emptyset,
\end{array}\]
where all ``$\cdots$''-nodes need to be from $\Sc^G(\sigma_j)$, all non-colliders not in $Z$ and all colliders from $\Sc^G(\sigma_j) \cap \Anc^{G'}(Z)$. Note that to get the correct end-segment replacement from the table one needs to pretend them to be of the form $\ot \sigma_1 \cdots$ and $\cdots \sigma_m \to$.\\
Then every $Z$-$\sigma$-open path $P$ in $G$ can be replaced by these rules with a $Z$-open path $P'$ in $G'$.
In other words, for all subsets $X,Y,Z \ins V$ we have the implication:
\[ X \Indep_{G'}^d Y \given Z \;\implies\; X \Indep_{G}^\sigma Y \given Z.\]
\begin{proof}
It is easy to check that such replacements (if existent) indeed lead to a $Z$-open path in $G'$.
If for example in the $Z$-$\sigma$-open path $P$ we have the segment of the form $\oto \sigma_j \to$ then its endnode $\sigma_{j,r}$ is not in $Z$.
Assume further that $\Sc^G(\sigma_j) \cap \Anc^G(Z) \neq \emptyset$ then the path from the second entry in the table 
$\oto  \cdot  \oto \cdots \oto  \sigma_{j,r}\begin{array}{c}\to\end{array}$ exists in $G'$ by assumption with all intermediate nodes from $\Sc^G(\sigma_j) \cap \Anc^{G'}(Z)$. Since the intermediate nodes are all colliders and $\sigma_{j,r}$ is a non-collider with $\sigma_{j,r} \notin Z$ this subpath is clearly $Z$-open in $G'$ at every point.
\end{proof}
\end{Thm}

\begin{Cor}
\label{s-sep-d-sep-eg}
The statements \ref{s-sep-d-sep-general} and \ref{s-sep-d-sep} both apply to the following cases:
\begin{enumerate}
\item $G'=G^\acy$,
\item $G'=(V,E^\acy \cup E^\mathrm{sc},H^\acy)$,
\item $G'=G^c=(V,E,H^c)$ with $H^c:= H \cup \{ F \ins \Sc^G(v) | \forall v \in V\}$.
\item $G'=G^c_2=(V,E,H^c_2)$ with $H^c_2:= \{ F \in H^c| \# F \le 2 \}$.
\item $G=(V,E,H_1)$ a directed graph and $G'=G^\mathrm{col}=(V,E^\mathrm{col},H_1)$ its \emph{collapsed graph}\footnote{The \emph{collapsed graph} $G^\mathrm{col}$ was introduced in \cite{Spirtes94} for a directed graph $G$.} with $E^\mathrm{col}:=E^\acy \cup E^\mathrm{sc}_<$, where for a given pseudo-topological order $<$ we define $E^\mathrm{sc}_<:=\{v \to w \,|\, v \in V, w \in \Sc^G(v) \text{ with } v<w\}$.
\item $G=(V,E,H_1)$ a directed graph and $G'=(V,E^\acy \cup E^\mathrm{sc},H_1)$.
\end{enumerate}
So in all these cases we have that for all subsets $X,Y,Z \ins V$ the following equivalence holds:
\[ X \Indep_{G}^\sigma Y \given Z \;\iff\; X \Indep_{G'}^d Y \given Z.\]
  Furthermore, if for such $G'$ (satisfying \ref{s-sep-d-sep-general} and \ref{s-sep-d-sep}) we in addition have $G_2=G'_2$ then $\sigma$-separation in $G$ and d-separation in $G$ are equivalent. 
\begin{proof}
  For $G^\acy$ always the first line of the replacements in \ref{s-sep-d-sep-general} already exist (with a bidirected edge or equality in the fork case). In the collider case we have  
one collider element $v_{j,Z}$ in $\Sc^G(\sigma_j)$ that is chosen such that either it is in $Z$ or has a directed path $v_{j,Z} \to w' \to \cdots \to w$ in $G$ 
with $w' \notin \Sc^G(\sigma_j)$ and $w \in Z$. Note that such an element exists since $\sigma_j$ is a collider segment in $P$ in that case and then $\Sc^G(\sigma_j) \ins \Anc^G(Z)$. It then follows that $v_{j,Z} \in \Anc^{G^\acy}(Z)$. \\
  For $G^c$ and $G^c_2$ both the directed paths and collider paths exist. 
	Furthermore, the ancestral relations are the same as for $G$. So colliders are always open if $\Sc^G(\sigma_j) \cap \Anc^G(Z) \neq \emptyset$.\\
For $G^\mathrm{col}$ and the next example the first lines of every non-collider case exists as in $G^\acy$ (with one directed edge or equality in the fork case). The collider case is like for $G^\acy$.\\
Now assume the additional assumption that $G_2=G'_2$. Then we have:
\[\begin{array}{ccccc} X \Indep_{G}^\sigma Y \given Z &\iff& X \Indep_{G'}^d Y \given Z &\iff& X \Indep_{G'_2}^d Y \given Z \\
&\stackrel{G_2'=G_2}{\iff}& X \Indep_{G_2}^d Y \given Z &\iff& X \Indep_{G}^d Y \given Z,
\end{array}\]
showing that d-separation and $\sigma$-separation are equivalent in this case.
\end{proof}
\end{Cor}

\begin{Cor}[$\sigma$-separation stable under marginalization]
\label{s-sep-marg}
Let $G=(V,E,H)$ be a \HEDG{} and $W,X,Y,Z \ins V$ be subsets with $(X \cup Y \cup Z) \cap W =\emptyset$.
Then we have the equivalence:
\[ X \Indep_{G}^\sigma Y \given Z \;\iff\; X \Indep_{G^{\marg\sm W}}^\sigma Y \given Z. \]
\begin{proof}
By induction we can assume $W=\{w\}$. 
We then use $G^c$ from \ref{s-sep-d-sep-eg}. 
Together with \ref{marg-d-sep} we have the equivalences:
\[ X \Indep_{G}^\sigma Y \given Z \;\stackrel{\ref{s-sep-d-sep-eg}}{\iff}\;  X \Indep_{G^c}^d Y \given Z \;\stackrel{\ref{marg-d-sep}}{\iff}\;
 X \Indep_{(G^c)^{\marg\sm W}}^d Y \given Z. \]
and:
\[ X \Indep_{G^{\marg\sm W}}^\sigma Y \given Z \;\stackrel{\ref{s-sep-d-sep-eg}}{\iff}\; X \Indep_{(G^{\marg\sm W})^c}^d Y \given Z. \]
It is then enough to show the equivalence:
\[ X \Indep_{(G^c)^{\marg\sm W}}^d Y \given Z \;\iff\; X \Indep_{(G^{\marg\sm W})^c}^d Y \given Z. \]
First note that $(G^c)^{\marg\sm W}$ and $(G^{\marg\sm W})^c$ have the same underlying set of nodes $V\sm W$ and edges $E^{\marg\sm W}$. So we only need to consider bidirected edges.\\
``$\implies$'': 
Consider a bidirected edge $v_1 \oto v_2$ in $(G^{\marg\sm W})^c$. Then either $v_1 \oto v_2$ lies in $G^{\marg\sm W}$ or $v_1 \in \Sc^{G^{\marg\sm W}}(v_2)$. In the first case we directly get $v_1 \oto v_2$ in $(G^c)^{\marg\sm W}$. In the latter case we get $v_1 \in \Sc^{G}(v_2)$ and then $v_1 \oto v_2$ in $G^c$ and thus in $(G^c)^{\marg\sm W}$.
So we always get $v_1 \oto v_2$ in $(G^c)^{\marg\sm W}$.\\
It follows that every $Z$-open path in $(G^{\marg\sm W})^c$ is a $Z$-open path in $(G^c)^{\marg\sm W}$.\\
``$\Longleftarrow$'':
Now let $v_1 \oto v_2$ be in $(G^c)^{\marg\sm W}$. Then this comes from $v_1 \oto v_2$, $v_1 \ot w \to v_2$, $v_1 \oto w \to v_2$ or $v_1 \ot w\oto v_2$ in $G^c$.  
If the bidirected edges already lie in $G$ (or if $v_1 \ot w \to v_2$) then these cases directly give rise to $v_1 \oto v_2$ in $(G^{\marg\sm W})^c$.
Otherwise $v_1 \in \Sc^G(v_2)$, or $v_1 \in \Sc^G(w)$, or $v_2 \in \Sc^G(w)$, resp.. The first case implies $v_1 \in \Sc^{G^{\marg\sm W}}(v_2)$ and thus $v_1 \oto v_2$ in $(G^{\marg\sm W})^c$.
In the latter cases we can assume 
  $v_1,v_2 \notin \Ch^G(w) \cap \Sc^G(w)$ (since the other cases were covered by $v_1 \ot w \to v_2$ and by $v_1 \in \Sc^G(v_2)$).
Then there is a $w' \in \Ch^G(w) \cap \Sc^G(w) \sm \{w,v_1,v_2\}$. 
For the case $v_1 \oto w \to v_2$ we then get
that $w' \ot w \to v_2$ is in $G$ and thus $w' \oto v_2$ in $G^{\marg\sm W}$.
Since $w' \in \Sc^G(v_1)$ we get both $v_1 \oto w'$ and $v_1 \ot \dots \ot w'$ in $(G^{\marg\sm W})^c$.
If $v_1 \in \Anc^G(Z)$ then also $w' \in \Anc^G(Z)$, and
then $v_1 \oto w' \oto v_2$  is $Z$-open in $(G^{\marg\sm W})^c$.
If $v_1 \notin \Anc^G(Z)$ then all nodes in $v_1 \ot \dots \ot w'$ are not in $\Anc^G(Z)$ and thus the path part
$v_1 \ot \dots \ot w' \oto v_2$ is $Z$-open in $(G^{\marg\sm W})^c$. \\
By symmetry the same arguments apply to the case $v_1 \ot w \oto v_2$.\\
Replacing all bidirected edges with such constructions shows that every $Z$-open path in $(G^c)^{\marg\sm W}$ leads to a $Z$-open path in $(G^{\marg\sm W})^c$.\\
Together this shows the claim.
\end{proof} 
\end{Cor}

\begin{Cor}
\label{s-sep-anc-sub-hedg}
Let $G=(V,E,H)$ be a \HEDG{} with subsets $X,Y,Z \ins V$, $A$ an ancestral sub-\HEDG{} of $G$ with $X,Y,Z \ins A$ 
(e.g.\ $A=\Anc^{G}(X \cup Y \cup Z)$).
Then we have:
\[ X \Indep_G^\sigma Y \given Z \iff X \Indep_{A}^\sigma Y \given Z.\]
\begin{proof}
This follows from the fact \ref{anc-sub=mar} that on ancestral sub-\HEDG{}es $A \ins G$ the sub-\HEDG{}-structure $G(A)$ and the marginalized \HEDG{}-structure 
$G^{\marg(A)}$ are the same and that $\sigma$-separation is stable under marginalization by \ref{s-sep-marg}. 
\end{proof}
\end{Cor}

\begin{Cor}
\label{s-sep-aug-hedg}
Let $G=(V,E,H)$ be a \HEDG{} and $X,Y,Z \ins V$ subsets.
Consider the augmented graph $G^\aug$, the acyclification $G^\acy$ , the acyclic augmentation $G^\acag$  and also $(G^\acy)^\aug$ and $(G^\aug)^\acy$.
Then we have the equivalences:
\[ \begin{array}{clllcl}
 &    X & \Indep_G^\sigma& Y &\given& Z \\ 
\iff& X &\Indep_{G^\aug}^\sigma& Y &\given& Z \\
\iff& X &\Indep_{G^\acy}^d &Y & \given& Z \\
\iff& X &\Indep_{G^\acag}^d &Y &\given& Z \\
\iff& X &\Indep_{(G^\aug)^\acy}^d &Y &\given& Z \\
\iff& X &\Indep_{(G^\acy)^\aug}^d & Y &\given& Z.
\end{array}\]
\begin{proof}
Since all these constructions marginalize either to $G$ or $G^\acy$ the claim follows from \ref{s-sep-d-sep-eg}, \ref{s-sep-marg} and \ref{marg-d-sep}.
\end{proof}
\end{Cor}

\begin{Rem}
Note that even though $\sigma$-separation is stable under marginalization the graphical representations from \ref{s-sep-d-sep-eg} are not necessarily stable under marginalization as can be seen from figure \ref{fig:s-separation-marginalization}.
\begin{figure}[!htb]
\centering
\begin{tikzpicture}[scale=.7, transform shape]
\tikzstyle{every node} = [draw,shape=circle]
\node (v1) at (3,0) {$v_1$};
\node (v2) at (3,-5) {$v_2$};
\node (v3) at (1.5,-1.5) {$v_3$};
\node (v4) at (4.5,-1.5) {$v_4$};
\node (w) at (3,-3) {$w$};
\foreach \from/\to in {v3/v1, v1/v4, w/v2, v4/w, w/v3}
\draw[-{Latex[length=3mm, width=2mm]}] (\from) -- (\to);
\draw[-{Latex[length=3mm, width=2mm]}, dashed, gray] (v1) -- (w);
\end{tikzpicture} 
\begin{tikzpicture}[scale=.7, transform shape]
\tikzstyle{every node} = [draw,shape=circle]
\node (v1) at (3,0) {$v_1$};
\node (v2) at (3,-5) {$v_2$};
\node (v3) at (1.5,-1.5) {$v_3$};
\node (v4) at (4.5,-1.5) {$v_4$};
\foreach \from/\to in {v3/v1, v1/v4, v4/v2, v4/v3}
\draw[-{Latex[length=3mm, width=2mm]}] (\from) -- (\to);
\draw[-{Latex[length=3mm, width=2mm]}, dashed, gray] (v1) -- (v2);
\draw[{Latex[length=3mm, width=2mm]}-{Latex[length=3mm, width=2mm]},bend left, red] (v2) to node[fill,circle,red,inner sep=0pt,minimum size=5pt] {} (v3);
\end{tikzpicture}
\caption{A \HEDG{} $G$ on the left and its marginalization $G^{\marg \sm \{w\}}$ on the right (both without gray dashed arrow).
We have $\{v_1\} \nIndep_G^\sigma \{w\} \given \{v_3,v_4\}$.
 Any graphical representation $G'$ that reduces $\sigma$-separation in $G$ to d-separation in $G'$ needs to achieve
$\{v_1\} \nIndep_{G'}^d \{w\} \given \{v_3,v_4\}$. Assume that for this we introduced the directed edge $v_1 \to w$ (gray dashed arrow on the left) in $G'$. 
Then this will result in a directed edge $v_1 \to v_2$ in $(G')^{\marg \sm \{w\}}$ (gray dashed arrow on the right). But there won't be a directed edge $v_1 \to v_2$ in  $(G^{\marg \sm \{w\}})'$, which is constructed from the right \HEDG{} without the gray dashed arrow by changing the (hyper)edges of $\Sc^{G^{\marg \sm \{w\}}}(v_1)$. So we would have $(G')^{\marg \sm \{w\}} \neq (G^{\marg \sm \{w\}})'$ in this case.
	}
 \label{fig:s-separation-marginalization}
\end{figure}
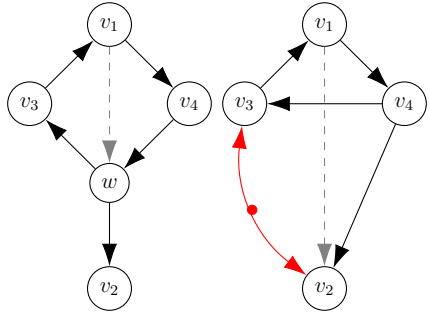
\end{Rem}

In the following we give another two ways how to check $\sigma$-separation in practice. 
One version (see \ref{s-sep-via-marg}) will depend on the context given by a node set $Z$ and mainly uses marginalization operations to reduce the problem of separating two sets in an undirected graph. The second version will use a definition of $\sigma$-separation that only uses triples of nodes instead of triples of segments in a path very similar to the usual d-separation definition. The reason we did not introduce 
$\sigma$-separation via this definition is that the results from \ref{s-sep-d-sep}, \ref{s-sep-d-sep-general}, \ref{s-sep-d-sep-eg}, which reduce $\sigma$-separation to d-separation in another graphical structure, would have become much more lengthy.

\begin{Thm}[Criterion for $\sigma$-separation]
\label{s-sep-via-marg}
Let $G=(V,E,H)$ be a \HEDG{} and $X,Y,Z \ins V$ subsets. Then $\sigma$-separation can be checked by the following construction:
\begin{enumerate}
 \item Let $W=X \cup Y \cup Z$ be the marginalized \HEDG{} structure $G^{\marg(W)}$.
 \item Let $W'$ be the \HEDG{} constructed from $W$, where we deleted all directed edges $z \to v$ for $z \in Z$ and $v \in \Ch^W(z)\sm \Sc^W(z)$.
 \item Replace every directed edge ($ v \to w$) and bidirected edge ($v \oto w$) in $W'$ with an undirected edge $v - w$ to get the \emph{skeleton}  $\tilde W$ of $W'$, which is an undirected graph.
 \item Let $\hat W$ be the marginalization $\tilde W ^{\marg\sm Z}$, which is an undirected graph with nodes $(X \cup Y) \sm Z$.
\end{enumerate}
Then we have:
\[ X \Indep^\sigma_G Y \given Z \;\iff\; X \sm Z \Indep_{\hat W} Y \sm Z. \]
  In other words, if there is an edge $x - y$ in $\hat W$ with $x \in X\sm Z$ and $y \in Y \sm Z$ then $X \nIndep^\sigma_G Y \given Z$ holds, otherwise we have: $X \Indep^\sigma_G Y \given Z$.
\begin{proof}
By \ref{s-sep-marg} $\sigma$-separation is stable under marginalization. So we have:
\[  X \Indep^\sigma_G Y \given Z \;\iff\; X \Indep^\sigma_W Y \given Z. \]
By the $\sigma$-separation rules \ref{s-sep-def} by segments, we clearly have that every $Z$-$\sigma$-open path in $W$ has no directed edge from a node $z \in Z$ to another strongly connected component. So we get:
\[  X \Indep^\sigma_W Y \given Z \;\iff\; X \Indep^\sigma_{W'} Y \given Z. \]
  Since separation in an undirected graph is stable under marginalization by \ref{sep-marg} we also have: 
\[  X\sm Z \Indep_{\tilde W} Y\sm Z \;\iff\; X\sm Z \Indep_{\hat W} Y\sm Z. \]
It is then left to show the equivalence:
\[  X \Indep^\sigma_{W'} Y  \given Z \;\iff\; X \sm Z \Indep_{\tilde W} Y \sm Z. \]
$\Longrightarrow$: Let $x - \cdots - y$ be a shortest path from $X \sm Z$ to $Y \sm Z$ in $\tilde W$ with all intermediate nodes in $Z$ (if any).
  The corresponding path $\pi$ with arrowheads in $W'$ is then $Z$-$\sigma$-open. Indeed, the endnodes are not in $Z$ and 
	all the segments $\sigma_i$ ($1<i<m$) are connected only through bidirected edges since we removed all outgoing directed edges.
	This leaves us only with collider segments whose elements are in $Z$. It follows that the path is $Z$-$\sigma$-open.\\
$\Longleftarrow$: Every $Z$-$\sigma$-open path from $X$ to $Y$ has its endnodes in $X\sm Z$ and $Y \sm Z$, resp.. The undirected version of the path is then open. Note that we do not condition on $Z$ in the undirected version.
\end{proof}
\end{Thm}

In the following we will give an equivalent definition of $\sigma$-separation that only checks for triples of nodes in a path (instead of triple of segments).

\begin{Def}[$\sigma$-separation (node version)]
\label{s-sep-alt-def}
Let $G=(V,E,H)$ be a \HEDG{} and $X,Y,Z \ins V$ subsets of the nodes.
\begin{enumerate}
\item Consider a path in $G$ with $n \ge 1$ nodes:
\[\begin{array}{ccccccc} 
& \ot &&\ot \\[-5pt]
v_1 & \to & \cdots & \to &v_n. \\[-5pt]
& \oto & &\oto 
\end{array}\]
The path will be called \emph{$Z$-$\sigma$-blocked} or \emph{$\sigma$-blocked by $Z$} if:
 \begin{enumerate}
   \item at least one of the endnodes $v_1$, $v_n$ is in $Z$, or
	\item collider case: there is a node $v_i \notin 	\Anc^G(Z)$ with two adjacent (hyper)edges that form a \emph{collider} at $v_i$:
    \[v_{i-1} \to v_i \ot v_{i+1}, \quad  v_{i-1}\oto v_i \ot v_{i+1},\]
   \[ v_{i-1}\to v_i \oto v_{i+1}, \quad v_{i-1}\oto v_i \oto v_{i+1}, \]
 or
  \item non-collider case: there is a node $v_i$ with two adjacent (hyper)edges in the path of the forms:
	 \[v_{i-1}\begin{array}{c}\to\\[-10pt]\ot\\[-10pt]\oto\end{array} v_i \to v_{i+1},\quad v_{i-1} \ot v_i \begin{array}{c}\to\\[-10pt]\ot\\[-10pt]\oto\end{array}v_{i+1},\] 
	 where $v_i \in Z$ has one directed edge pointing from $v_i$ to a node not in $\Sc^G(v_i)$. This summarizes the three cases:
    \begin{enumerate}
		 \item $v_i \in Z$ and at least one $v_{i-1}$ or $v_{i+1} \notin \Sc^G(v_i)$ and the two (hyper)edges are a \emph{fork} at $v_i$:
\[  v_{i-1}\ot v_i \to v_{i+1}.\]
 \item $v_i \in Z$, $v_{i+1} \notin \Sc^G(v_i)$ and the two (hyper)edges are \emph{right-directed} at $v_i$:
\[  v_{i-1}\to v_i \to v_{i+1}, \quad v_{i-1}\oto v_i \to v_{i+1}.\]
\item $v_i \in Z$, $v_{i-1} \notin \Sc^G(v_i)$ and the two (hyper)edges are \emph{left-directed} at $v_i$:
\[  v_{i-1}\ot v_i \ot v_{i+1}, \quad v_{i-1}\ot v_i \oto v_{i+1}.\]
		\end{enumerate}
 \end{enumerate} 
Briefly stated, a path is $Z$-$\sigma$-blocked iff it has an endnode in $Z$, or it contains a collider not in $\Anc^G(Z)$, or it contains a non-collider in $Z$ that points to a node in a different strongly connected component. 
\item If none of the above holds then the path is called \emph{$Z$-$\sigma$-open} or \emph{$Z$-$\sigma$-active}.
\item  We say that $X$ is \emph{$\sigma$-separated} from $Y$ given $Z$ if every path in $G$ with one endnode in $X$ and one endnode in $Y$ is $\sigma$-blocked by $Z$.
\end{enumerate}
\end{Def}

\begin{Cor}
\label{s-sep-alt-s-sep}
Let $G=(V,E,H)$ be a \HEDG{}. 
Then $\sigma$-separation based on triples of segments is equivalent to $\sigma$-separation based on triples of nodes.
I.e. for subsets $X,Y,Z \ins V$ we have:
\[ X \Indep_{G}^{\sigma(\text{segments})} Y \given Z\qquad\iff\qquad X \Indep_{G}^{\sigma(\text{nodes})} Y \given Z. \]
\begin{proof}
For this proof let $\hat \sigma$ indicate $\sigma$-separation based on triple of nodes and $\tilde\sigma$ the one with segments.
Let $\pi$ be a path from $X$ to $Y$ in $G$. By assuming $X \Indep_{G}^{\tilde\sigma} Y \given Z$ the path is then $Z$-$\tilde\sigma$-blocked. If it has an endnode in $Z$ or has an segment of the form $\cdots\sigma_j \to $ with $\sigma_{j,r} \in Z$
or $\ot  \sigma_j \cdots $ with $\sigma_{j,l} \in Z$
 then $\pi$ is also $Z$-$\hat \sigma$-blocked. If there is a collider segment $\genfrac{}{}{0pt}{}{\to}{\oto} \sigma_j \genfrac{}{}{0pt}{}{\ot}{\oto}$ 
with $\Sc^G(\sigma_j) \cap 	\Anc^G(Z)=\emptyset$ then there is at least one collider element $v$ in $\sigma_j$, which then is $v \notin \Anc^G(Z)$.
Together this shows:
\[ X \Indep_{G}^{\tilde\sigma} Y \given Z\qquad\implies\qquad X \Indep_{G}^{\hat \sigma} Y \given Z. \]
For the reverse implication assume $X \Indep_{G}^{\hat \sigma} Y \given Z$ and let $\pi$ be a path from $X$ to $Y$, which then can be assumed to be $Z$-$\hat \sigma$-blocked.
If an endnode or a non-collider with an edge to another component is in $Z$ then clearly $\pi$ is $Z$-$\tilde\sigma$-blocked and we are done.
So we can assume that there is no such endode or non-collider on $\pi$. It follows that the set $T$ of colliders $v_i$ on $\pi$ with $v_i \notin \Anc^G(Z)$ must be non-empty. Otherwise $\pi$ would be $Z$-$\hat \sigma$-open. 
If we had the case that for every $v_i \in T$ the corresponding segment $\sigma_j$ of $v_i$ were a non-collider then we could replace each of these segments $\sigma_j$ with a directed path $\sigma_j' \ins \Sc^G(v_i)$ with the same endnodes as $\sigma_j$ such that the resulting path $\pi'$ would be $Z$-$\hat \sigma$-open from $X$ to $Y$. Indeed $v_i \notin \Anc^G(Z)$ implies $\Sc^G(v_i) \cap \Anc^G(Z) = \emptyset$ and thus $\sigma_j' \cap \Anc^G(Z) = \emptyset$. But by assumption such a $Z$-$\hat \sigma$-open path does not exist.
So there must be a $v_i \in T$ whose segment $\sigma_j$ is also a collider. Since $\sigma_j \cap \Anc^G(Z) = \emptyset$ the path $\pi$ is $Z$-$\tilde\sigma$-blocked showing:
\[ X \Indep_{G}^{\hat \sigma} Y \given Z\qquad\implies\qquad X \Indep_{G}^{\tilde\sigma} Y \given Z. \]
 \end{proof}
\end{Cor}

\section{Markov Properties}

Guided by the theory of Bayesian networks (BN), i.e.\ the DAG case, in this section we will collect several meaningful Markov properties for \HEDG{}es and investigate the logical relations between them. We will see that unlike in the DAG case these properties are not all equivalent and need to be addressed carefully. Furthermore, several different generalizations of the Markov properties for DAGs to \HEDG{}es might exist and the number of Markov properties for \HEDG{}es quickly explodes. We try to stick to the most important Markov properties and group them together in a few subsections.

\subsection{General Notations and Properties}

We will start with the basic notations we will use during the whole section and some general properties about d-separation, $\sigma$-separation and conditional independences, which are part of the \emph{graphoid rules} (see \cite{PeaPaz87}, \cite{GeiVerPea90}) and \emph{separoid axioms} (see \cite{Daw79}).

\begin{Rem}
In the following we fix an index set $V$ of vertices and standard Borel spaces $\Xcal_v$ for $v \in V$.
For a subset $A \ins V$ we will put $\Xcal_A:= \prod_{v \in A} \Xcal_v$ (endowed with the product-$\sigma$-algebra).  
We want to investigate the possible relations between \HEDG{}-structures $G=(V,E,H)$ over $V$ and probability distributions $\Pr_V$ over $\Xcal_V$.
For this reason we will consider the canonical projections $T_v: \Xcal_V \to \Xcal_v$ and $T_A: \Xcal_V \to \Xcal_A$ and view them as random variables as soon as we assign a probability measure $\Pr_V$ to $\Xcal_V$. Since $T_A = (T_v)_{v \in A}$ we will identify subsets $A \ins V$ with the tuple of variables $(T_v)_{v \in A}$. 
For instance, for subsets $X,Y,Z \ins V$ we will write: 
\[ X \Indep_{\Pr_V} Y \given Z,  \]
when we mean that $(T_v)_{v \in X}$ and $(T_v)_{v \in Y}$ are conditionally $\Pr_V$-independent given $(T_v)_{v \in Z}$.
  Furthermore, we will write $\Pr_A$ for the marginal probability distribution of $\Pr_V$ over $\Xcal_A$ 
  (i.e.\ $\Pr_A=T_{A,*}\Pr_V$, the pushforward of $\Pr_V$ under the projection map $T_A$). 
\end{Rem}

\begin{Rem}
Let $\Pr_V$ be a probability distribution over $\Xcal_V$.
In the case that $\Pr_V$ has a density $p$ w.r.t.\ some product measure $\mu_V:=\otimes_{v \in V} \mu_v$ we will define
  the \emph{marginal conditional density} $p(x_B|x_A)$ for sets $A,B \ins V$:
\[ p(x_B|x_A):=\left\{\begin{array}{lll}  \frac{p(x_{B \cup A})}{ p(x_A)  }, & \text{if }\;  p(x_A) > 0 \; \text{ and }\; x_A|_{A \cap B} = x_B|_{A \cap B}, \\
0, & \text{otherwise}.   \end{array}\right.\]
 with the marginal densities for sets $C \ins V$: 
\[p(x_C) :=  \int p(x_V) \,d \mu_{V\sm C}  (x_{V \sm C}),\] 
where we integrate all other variables out w.r.t.\ the product measure $\mu_{V\sm C} :=\otimes_{v \in V\sm C} \mu_v$.
\end{Rem}

\begin{Lem}[(Semi-)graphoid/separoid axioms]
\label{GrAxPr}
Let $V$ be a set and $\Xcal_V$ as before. 
  Then d- and $\sigma$-separation for a \HEDG{} $G=(V,E,H)$, separation in an undirected graph $G=(V,E)$ and $\Pr_V$-independence all satisfy the following \emph{(semi-)graphoid axioms} (or \emph{separoid axioms})\footnote{See \cite{Daw79}, \cite{PeaPaz87} and \cite{GeiVerPea90}.}.\\
For all sets $X,Y,Z,W \ins V$ we have the following rules:
\begin{enumerate}
\item Irrelevance: $X \Indep Y \given Y$ always holds.
\item Symmetry: $X \Indep Y \given Z \implies Y \Indep X \given Z$.
\item Decomposition: $X \Indep Y \cup W \given Z \implies X \Indep Y \given Z$.
\item Weak Union: $X \Indep Y \cup W \given Z \implies X \Indep Y \given W \cup Z$.
\item Contraction: $(X \Indep Y \given W \cup Z) \land (X \Indep W \given Z) \implies X \Indep Y\cup W \given Z$.
\end{enumerate}
Here $\Indep$ either stands for $\Indep_G$, $\Indep_G^d$, $\Indep_{G}^\sigma$ or $\Indep_{\Pr_V}$, resp..\\
Furthermore, in the case of d- and $\sigma$-separation for a \HEDG{} $G=(V,E,H)$, separation in an undirected graph $G$ or $\Pr_V$-independence if $\Pr_V$ has a strictly positive density $p_V >0$ w.r.t.\ a product measure $\mu_V=\otimes_{v \in V} \mu_v$ on $\Xcal_V$ we also have:
\begin{enumerate}[resume]
\item Intersection: $(X \Indep Y \given W \cup Z) \land (X \Indep W \given Y \cup Z) \implies X \Indep Y \cup W \given Z$, whenever
  $X,Y,Z,W$ are pairwise disjoint.
\end{enumerate}
In the case of d- and $\sigma$-separation in a \HEDG{} $G$ or separation in an undirected graph $G$ we further have:
\begin{enumerate}[resume]
\item Composition: $(X \Indep Y \given Z) \land (X \Indep W \given Z) \implies X \Indep Y\cup W \given Z$.
\end{enumerate}
\begin{proof}
For $\Pr_V$ the statements can directly be checked by the definitions and the usual properties of conditional probabilities (e.g.\ see \cite{Kle14} 8.14):
\begin{enumerate}
\item Irrelevance: $\Pr(X \in A |Y,Y)=\Pr(X \in A |Y)$. 
\item Symmetry: $\Pr(X \in A, Y \in B |Z)=\Pr(X \in A|Z) \cdot \Pr(Y \in B |Z)
=\Pr(Y \in B|Z) \cdot \Pr(X \in A |Y)$.
\item Decomposition: $\Pr(Y \in B|X,Z)=\Pr(Y \in B, W \in \Xcal_W|X,Z)= \Pr(Y \in B, W \in \Xcal_W|Z)= \Pr(Y \in B|Z)$.
\item Weak union: $\Pr(X \in A|Y,W,Z)= \Pr(X \in A|Z)= \E[\Pr(X \in A|Z)|W,Z] = \E[\Pr(X \in A|Y,W,Z)|W,Z]=\Pr(X \in A|W,Z)$.
\item Contraction: $\Pr(X \in A|Y,W,Z) =  \Pr(X \in A|W,Z) = \Pr(X \in A|Z)$.
\item Intersection for positive densities $p>0$: %
$X \Indep Y \given W \cup Z$ and $X \Indep W \given Y \cup Z$ can be expressed in the equations:
\[  p(y,w,z) \cdot p(x|w,z) = p(x,y,w,z) = p(x|y,z) \cdot p(y,w,z). \]
Since $p(y,w,z) > 0$ we get the equation: 
\[p(x|w,z) =p(x|y,z)\]
showing that the former is independent of $w$. %
This implies:
\[\begin{array}{lllll}
  && p(x|z) \\
	&=&  \int p(x,w|z) \,d\mu(w) &=&  \int p(x|w,z) \cdot p(w|z) \,d\mu(w)  \\
	&=& \int p(x|y,z) \cdot p(w|z) \,d\mu(w) &=&  p(x|y,z) \cdot \int p(w|z) \,d\mu(w) \\
	&=& p(x|w,z) \cdot 1 &=& p(x|w,z).
\end{array}\]
With this we have: $X \Indep W \given Z$. Together with $X \Indep Y \given W \cup Z$ and contraction we get: $X \Indep Y \cup W \given Z$.
\end{enumerate}
The proofs of the statements for separation in undirected graphs are similar to (but easier than) d-separation in \HEDG{}es.
For the $\sigma$-separation in $G$ one can reduce to d-separation in its acyclification $G^\acy$ by \ref{s-sep-d-sep-eg}.\\
So it is left to show the statements for d-separation in \HEDG{}es. For this it will be easier to work with the equivalent definition of d-separation where colliders need to avoid $Z$ (instead of $\Anc^G(Z)$).\\
Irrelevance follows from the rules about the endnodes in a path. Symmetry and composition are obvious from the definition of d-separation. For decomposition note that a $Z$-open path from $X$ to $Y$ (assuming $X \nIndep Y \given Z$) is also a $Z$-open path from $X$ to $Y \cup W$ (implying $X \nIndep Y \cup W \given Z$).\\
For weak union let $\pi$ be a $W\cup Z$-open path from $X$ to $Y$ (assuming $X \nIndep Y \given W \cup Z$). So every endnode and non-collider of $\pi$ is not in $W \cup Z$ and every collider is in $W \cup Z$. If $\pi$ does not contain nodes from $W \sm Z$ it is already $Z$-open 
(implying $X \nIndep Y \cup W  \given Z$). If a node from $W \sm Z$ occurs in $\pi$ then the shortest subpath of $\pi$ from $X$ to a node $w \in W \sm Z$ is $Z$-open, which shows the claim ($X \nIndep Y \cup W \given Z$).\\
For contraction and intersection let $\pi$ now be a shortest $Z$-open path from $X$ to $Y \cup W$ (assuming $X \nIndep Y \cup W\given Z$). 
Then $\pi$ ends in $Y\sm Z$ or $W\sm Z$. W.l.o.g. we can assume that $\pi$ ends in a node from $Y \sm (W \cup Z)$. Indeed, in the case of contraction ending in $W$ would already give the claim ($X \nIndep W \given Z$) and in the case of intersection $Y$ and $W$ are symmetric in the statement and disjointness of the sets was assumed. 
Then, since $\pi$ was also assumed to be a shortest path, no node from $W$ occurs in $\pi$. So $\pi$ is $W \cup Z$-open and the claim follows ($X \nIndep Y \given W \cup Z$).
\end{proof}
\end{Lem}

\subsection{Markov Properties for DAGs}

The main source of inspiration for the Markov properties for \HEDG{}es come from the Markov properties for DAGs  and the proofs that show that these are all equivalent in the DAG case. %
 So we will give a short overview of the Markov properties for DAGs. Proofs can be found in \cite{Lau90} and \cite{Lau98},  but these results will also follow from the more general versions in this paper later on. An overview over the Markov properties for the \HEDG{} case is given in \ref{fig:overview}, which then also applies for DAGs.

\begin{Thm}
\label{dag-markov-prop}
Let $G=(V,E)$ be a directed acyclic graph (DAG) and $\Pr_V$ a probability distribution on the product space $\Xcal_V=\prod_{v \in V} \Xcal_v$.
Then the following properties, which relate the graphical structure of $G$ to the probability distribution $\Pr_V$, are equivalent:
\begin{enumerate}
 \item The \emph{directed local Markov property (\hyt{dLMP-DAG}{dLMP})}: For every $v \in V$ we have:
 \[ \{v\} \Indep_{\Pr_V} \NonDesc^G(v) \given \Pa^G(v).  \]
\item The \emph{ordered local Markov property (\hyt{oLMP-DAG}{oLMP})}: There exists a topological order $<$ for $G$ such that for every $v \in V$ we have:
\[ \{v\} \Indep_{\Pr_V} \Pred^G_\le(v) \sm \{v\} \given \Pa^G(v).  \]
\item The \emph{directed global Markov property (\hyt{dGMP-DAG}{dGMP})}: For all subsets $X,Y,Z \ins V$ we have the implication:
 \[ X \Indep_G^d Y \given Z \;\implies\; X \Indep_{\Pr_V} Y \given Z.\]
\item The \emph{ancestral undirected global Markov property (\hyt{auGMP-DAG}{auGMP})}: For all subsets $X,Y,Z \ins V$ we have the implication:
 \[ X \Indep_{\Anc^G(X \cup Y \cup Z)^\moral} Y \given Z \;\implies\; X \Indep_{\Pr_V} Y \given Z,\]
where $\Anc^G(X \cup Y \cup Z)^\moral$ is the undirected graph given by the moralization of the ancestral closure of $X \cup Y \cup Z$.
\item The \emph{structural equations property (\hyt{SEP-DAG}{SEP})}: There exist:
 \begin{enumerate}
  \item a probability space $(\Omega, \fa, \Pr)$,
	\item standard Borel spaces $\Ecal_v$ for every $v \in V$,
	\item random variables $E_v: (\Omega, \fa, \Pr) \to \Ecal_v$ for $v \in V$, such that the system
	 $(E_v)_{v \in V}$ is jointly $\Pr$-independent,
  \item measurable functions $f_v: \Xcal_{\Pa^G(v)} \x \Ecal_v \to \Xcal_v$ for every $v \in V$,
	\item random variables $X_v: (\Omega, \fa, \Pr) \to \Xcal_v$ for every $v \in V$ with:
	 \[ X_v=f_v(X_{\Pa^G(v)}, E_v) \qquad \Pr\text{-a.s.},  \]
 \end{enumerate}
such that: $\qquad\Pr^{(X_v)_{v \in V}} = \Pr_V$.		
\item The \emph{recursive factorization property (\hyt{rFP-DAG}{rFP})}:
There exists a topological order $<$ of $G$ given by $v_1,\dots,v_n \in V$ such that 
for all measurable sets $B_1 \ins \Xcal_{v_1},\dots, B_n \ins \Xcal_{v_n}$ we have the equality:
\begin{eqnarray*}
&&\int \I_{B_n}(x_n)\cdots \I_{B_1}(x_1)\, d\Pr_V(x_n,\dots,x_1) \\
&=&\int \cdots \int \I_{B_n}(x_n)\cdots \I_{B_1}(x_1) \, d\Pr_V^{v_n|x_{\Pa^G(v_n)}}(x_n) \cdots \\
&& \qquad \cdots d\Pr_V^{v_{i}|x_{\Pa^G(v_i)}}(x_{i}) \cdots d\Pr_V^{v_1}(x_1), 
\end{eqnarray*}
or in short notation:
\[ d \Pr_V(x) = \prod_{i=n}^1 d\Pr_V^{ v_i|x_{\Pa^G(v_i)} }(x_i),  \] 
where $\Pr_V^{v_i|x_{\Pa^G(v_i)}}$ is the regular conditional probability distribution of $\Pr_V$ on $\Xcal_{v_i}$ given the parent nodes 
$\Pa^G(v_i)$ on the value $x_{\Pa^G(v_i)}$ (the corresponding values of $x_1,\dots,x_n$).
\end{enumerate}
If the probability distribution $\Pr_V$ has a density $p=p_V$ w.r.t.\ a product measure $\mu_V=\otimes_{v \in V}\mu_v$ on $\Xcal_V$
  then the recursive factorization property (\hyl{rFP-DAG}{rFP}) for $(G,\Pr_V)$ is also equivalent to:

\begin{enumerate}[resume]
 \item The \emph{recursive factorization property with density (\hyt{rFPwd-DAG}{rFPwd})} for %
 measures $\mu_v$, $v \in V$: There exists a density $p=p_V$ for $\Pr_V$ w.r.t.\ to the product measure $\mu_V$ such that for $\mu_V$-almost-all $x_V=(x_v)_{v \in V} \in \Xcal_V$ we have the equality:
\[ p(x_V) = \prod_{v \in V} p(x_v|x_{\Pa^G(v)}).  \]
\end{enumerate}
\begin{proof}
This was already proved in \cite{Lau90} and \cite{Lau98}.  
Sketch of proofs by arguments in this paper:\\
  ``\hyl{dGMP-DAG}{dGMP} $\implies$ \hyl{dLMP-DAG}{dLMP}'': Follows from $ \{v\} \Indep_G^d \NonDesc^G(v) \given \Pa^G(v)$.\\
``\hyl{dLMP-DAG}{dLMP} $\implies$ \hyl{oLMP-DAG}{oLMP}'': Follows from $\Pred^G_\le(v) \sm \{v \} \ins \NonDesc^G(v)$ for topological orders $<$.\\
  ``\hyl{oLMP-DAG}{oLMP} $\implies$ \hyl{dGMP-DAG}{dGMP}'': A more general proof is given in \ref{oLMP-dGMP} by using the semi-graphoid rules \ref{GrAxPr}. %
  The whole equivalence ``\hyl{dGMP-DAG}{dGMP} $\iff$ \hyl{dLMP-DAG}{dLMP} $\iff$ \hyl{oLMP-DAG}{oLMP}'' in more generality is summarized in \ref{oLMP-dLMP-dGMP}.\\
	``\hyl{dGMP-DAG}{dGMP} $\iff$ \hyl{auGMP-DAG}{auGMP}'': This follows directly from the following equivalence (\ref{d-sep-moral-hedg}, 
	also see \ref{dGMP-auPMP}):
	 \[ X \Indep_G^d Y \given Z \;\iff\; X \Indep_{\Anc^G(X \cup Y \cup Z)^\moral} Y \given Z.\]
  ``\hyl{oLMP-DAG}{oLMP} $\iff$ \hyl{rFP-DAG}{rFP}'': Induction by $n$ in the same topological order works. Also see \ref{oLMP-rFP}.\\
  ``\hyl{rFP-DAG}{rFP} $\implies$ \hyl{SEP-DAG}{SEP}'': This inductively follows from Lemma \ref{sep-lemma} with $X=T_{\Pa^G(v)}$ and $Y=T_v$ and a fixed topological order on $G$.\\
  ``\hyl{SEP-DAG}{SEP} $\implies$ \hyl{dGMP-DAG}{dGMP}'': This follows from the more general statements later:
  \hyl{SEP-DAG}{SEP} implies the \hyl{dLMP-DAG}{dLMP} for $(G^\aug,\Pr_{V^\aug})$ by \ref{SEP-equiv} and \ref{mDAG-mdLMP-SEP}. The latter implies \hyl{dGMP-DAG}{dGMP} for $(G^\aug,\Pr_{V^\aug})$ by the above equivalences or \ref{oLMP-dLMP-dGMP}. By marginalization (see \ref{mdGMP-dGMP}) we get \hyl{dGMP-DAG}{dGMP} for $(G,\Pr_{V})$.\\
  ``\hyl{rFP-DAG}{rFP} $\iff$ \hyl{rFPwd-DAG}{rFPwd}'': Clear by the use of a density. Also see \ref{rFPwd-rFP}.
\end{proof}
\end{Thm}

\begin{Rem}
\begin{enumerate}
\item By the equivalences in Theorem \ref{dag-markov-prop} the phrase ``There exists a topological order'' can be replaced by ``For every topological order'' in Theorem \ref{dag-markov-prop} everywhere.
\item In the \hyl{SEP-DAG}{SEP} in Theorem \ref{dag-markov-prop} we can without loss of generality assume the spaces $\Ecal_v$ to be $[0,1]$ with the Borel $\sigma$-algebra and $E_v$ to be uniformly distributed by the following Lemma \ref{sep-lemma}. 
\end{enumerate}
\end{Rem}

\begin{Lem}
\label{sep-lemma}
  Let $(\Omega, \fa, \Pr)$ be probability space, $(\Xcal,\Bcal_\Xcal)$, $(\Ycal,\Bcal_\Ycal)$ be measure spaces, where $\Ycal$ is assumed to be a standard Borel space and $\iota: \Ycal \inj [0,1]$ a fixed embedding onto a Borel subset of $[0,1]$. Let 
$X: (\Omega, \fa, \Pr) \to (\Xcal,\Bcal_\Xcal)$, $Y: (\Omega, \fa, \Pr) \to (\Ycal,\Bcal_\Ycal)$ be random variables.
Now construct:
\begin{enumerate}
 \item $(\Omega', \fa', \Pr'):=(\Omega\x [0,1], \fa\otimes\Bcal([0,1]), \Pr \otimes \lambda)$, where $\lambda$ is the Lebesgue measure on $[0,1]$.
 \item $\pi: (\Omega', \fa', \Pr') \srj (\Omega, \fa, \Pr)$ the canonical projection onto the first factor,
 \item $X':= X \circ \pi :\; (\Omega', \fa', \Pr') \to (\Xcal,\Bcal_\Xcal) $,
 \item $E: (\Omega', \fa', \Pr') \srj [0,1]$ the canonical projection onto the second factor,
 \item $\iota Y := \iota \circ Y$,
 \item $F_{Y|x}(y):=\Pr(\iota Y \le y|X=x)$,
 \item $g : \Xcal \x [0,1] \to [0,1]$ with $g(x,e):=\inf\{ y \in [0,1]\,|\, F_{Y|x}(y) \ge e \}$,
 \item $Y':= g(X',E): \; (\Omega', \fa', \Pr') \to ([0,1],\Bcal([0,1])) $.
\end{enumerate}
Then $E$ is uniformly distributed and independent of $X'$ under $\Pr'$, $g$ is measurable and we have:
\[  \Pr^{\iota Y|X} = {\Pr'}^{Y'|X'}, \qquad \text{and}\qquad \Pr^{X} = {\Pr'}^{X'},  \]
and in particular:
\[  \Pr^{(X,\iota Y)} = {\Pr'}^{(X',Y')}. \]
Also note that it follows that $\Pr'(Y' \notin \iota(\Ycal))=0$, and thus $Y'$ can be viewed as a random variable mapping to $\Ycal$. 
\begin{proof}
The independence is clear by construction. 
$g$ is measurable because for all $y' \in [0,1]$ we can write:
\[\begin{array}{rcl} g^{-1}([0,y']) &=&\{ (x,e) \in \Xcal \x [0,1]\,|\, F_{Y|x}(y') \ge e \} \\
&=& (\Pr(B|X=\cdot) \x \id )^{-1}(\Delta), \end{array} \]
with $B:=\{\iota Y \le y'\}$ and $\Delta:=\{(h,u) \in [0,1]^2\,|\, h \ge u \}$,
where all occuring maps and sets are measurable on the right.
After the measurability of $g$ is checked we then get:
\[\begin{array}{lclll}
\Pr'(Y' \le y'|X'=x) &=& \Pr'( g(X',E) \le y'|X'=x) \\
&=&\Pr'( g(x,E) \le y'|X'=x) \\
&\stackrel{X' \Indep_{\Pr'} E}{=}&\Pr'( g(x,E) \le y') \\
&=& \lambda (\{e \in [0,1] \,|\, g(x,e) \le y'\}) \\
&=& \lambda (\{e \in [0,1]\,|\, e \le F_{Y|x}(y')\}) \\
&=& \lambda ([0,F_{Y|x}(y')]) \\
&=& F_{Y|x}(y') \\
&=& \Pr(\iota Y \le y'|X=x).\\
\Pr'(X' \in A) 
&=& \Pr'( X^{-1}(A) \x [0,1] ) \\
&=& \Pr(X^{-1}(A)) \cdot \lambda([0,1]) \\
&=& \Pr(X  \in A).
\end{array} \]
Thus we have: $\Pr^{(X,\iota Y)} = {\Pr'}^{(X',Y')}$.
\end{proof} 
\end{Lem}

\subsection{Directed Markov Properties for \HEDG{}es}

In this subsection we will define the most important directed Markov properties for general \HEDG{}es like the \emph{directed global Markov property (\hyl{dGMP}{dGMP})}. In view of the later structural equations properties (\hyl{csSEP}{csSEP}, see \ref{csSEP-smgdGMP}) we also need a weaker (see \ref{dGMP-gdGMP}) but more generally applicable property: the \emph{general directed global Markov property (\hyl{gdGMP}{gdGMP})} based on $\sigma$-separation instead of d-separation. We show that both of these properties are stable under marginalization (see \ref{dGMP-stable-marg}).
Furthermore, a generalized \emph{directed local Markov property (\hyl{dLMP}{dLMP})} for \HEDG{}es is introduced. The definition will be made in a way such that it is still local in some sense (see \ref{def-dGMP-dLMP-gdGMP}), reduces to the usual one in the DAG case (see \ref{DAG-dLMP}), is strong enough to imply the ordered local Markov property (\hyl{oLMP}{oLMP}) for assembling pseudo-topological total orders (see \ref{dLMP-oLMP}), but is weak enough to be implied by the directed global Markov property (\hyl{dGMP}{dGMP}, see \ref{dGMP-dLMP}) and such that a marginal version (\hyl{mMP}{mdLMP}) is implied by some structural equations property (\hyl{lsSEP}{lsSEP}) later on (see \ref{lsSEP-mdLMP}) like it is used in the proofs of \ref{dag-markov-prop} in the DAG case. One drawback remains, namely, that the directed local Markov property (\hyl{dLMP}{dLMP}) in general is not stable under marginalization (as can be seen from example \ref{main-example}).

\begin{Def}
\label{def-dGMP-dLMP-gdGMP}
Let $G=(V,E,H)$ be a \HEDG{} and $\Pr_V$ a probability distribution $\Pr_V$ on 
$\Xcal_V=\prod_{v \in V} \Xcal_v$. We define the following properties relating $G$ to $\Pr_V$:
\begin{enumerate}
  \item The \emph{directed global Markov property (\hyt{dGMP}{dGMP})}: For all subsets $X,Y,Z \ins V$ we have the implication:
 \[ X \Indep_G^d Y \given Z \;\implies\; X \Indep_{\Pr_V} Y \given Z.\]

\item The \emph{directed local Markov property (\hyt{dLMP}{dLMP})}: For every $v \in V$, every ancestral sub-\HEDG{} $A \ins \NonDesc^G(v) \cup \Sc^G(v)$ with $v \in A$ and every
set $S \ins \Sc^G(v) \sm \{v\}$ we have:
 \[ \{v\} \Indep_{\Pr_V} A_{\sm S} \sm \{v\} \given \partial_{A_{\sm S}^\moral}(v),\]
where $A_{\sm S} := A^{\marg\sm S}$ and
\[\begin{array}{lrl}
\partial_{A_{\sm S}^\moral}(v)&=& \big( \Dist^{A_{\sm S}}(\{v\} \cup \Ch^{A_{\sm S}}(v)) \\&&  \cup \Pa^{A_{\sm S}}(\Dist^{A_{\sm S}}(\{v\} \cup
 \Ch^{A_{\sm S}}(v))) \big) \sm \{v\}. 
\end{array}\]

\item The \emph{general directed global Markov property (\hyt{gdGMP}{gdGMP})}: For all subsets $X,Y,Z \ins V$ we have the implication:
\[ X \Indep_{G}^\sigma Y \given Z \;\implies\; X \Indep_{\Pr_V} Y \given Z.\]
\end{enumerate}
\end{Def}

\begin{Lem}[dGMP/gdGMP stable under marginalization]
\label{dGMP-stable-marg}
Let $G=(V,E,H)$ be a \HEDG{} and $\Pr_V$ a probability distribution $\Pr_V$ on 
$\Xcal_V$. Let $W \ins V$ be a subset and $G_W = G^{\marg(W)}$ the marginalized \HEDG{} with nodes $W$ and $\Pr_W$ the marginal distribution on 
$\Xcal_W$. We then have the implications:
 \[ \xymatrix{ 
      \text{\hyl{dGMP}{dGMP} for } (G,\Pr_V) \ar@{=>}[r]& \text{\hyl{dGMP}{dGMP} for } (G_W,\Pr_W), } \]
			 \[ \xymatrix{ 
      \text{\hyl{gdGMP}{gdGMP} for } (G,\Pr_V) \ar@{=>}[r]& \text{\hyl{gdGMP}{gdGMP} for } (G_W,\Pr_W).
			} \]
  In other words, both the directed global Markov property (\hyl{dGMP}{dGMP}) and the general directed global Markov property (\hyl{gdGMP}{gdGMP}) are stable under marginalization.
\begin{proof}
The statements follow from the fact that d- and $\sigma$-separation are stable under marginalization by \ref{marg-d-sep} and \ref{s-sep-marg}.
  So for $X,Y,Z \ins W$ and \hyl{dGMP}{dGMP} we have:
  \[ X \Indep_{G_W}^d Y \given Z \;\stackrel{\ref{marg-d-sep}}{\iff}\; X \Indep_G^d Y \given Z \;\stackrel{\text{\hyl{dGMP}{dGMP}}}{\implies} \; X \Indep_{\Pr_V} Y \given Z \;\stackrel{X,Y,Z\ins W}{\iff} \; X \Indep_{\Pr_W} Y \given Z.   \]
For \hyl{gdGMP}{gdGMP} we use \ref{s-sep-marg} instead of \ref{marg-d-sep}.
\end{proof}
\end{Lem}

\begin{Thm}
\label{dGMP-gdGMP}
For a \HEDG{} $G=(V,E,H)$ and a probability distribution $\Pr_V$ on $\Xcal_V$ we have the implication:
\[ \xymatrix{ 
      \text{\hyl{dGMP}{dGMP}} \ar@{=>}[r]& \text{\hyl{gdGMP}{gdGMP}}.
} \]
If $\{v,w\} \in H$ for all $v \in V$ and all $w \in \Sc^G(v)$ then the other direction of the implication also holds.%
\begin{proof}
Let $X,Y,Z \ins V$ be subsets. We then have the implications:
\[  X \Indep_{G}^\sigma Y \given Z \;\stackrel{\ref{s-sep-d-sep}}{\implies}\; X \Indep_{G}^d Y \given Z 
  \;\stackrel{\text{\hyl{dGMP}{dGMP}}}{\implies}\;
 X \Indep_{\Pr_V} Y \given Z. \]
  For the other direction note that $G_2=G^c_2$ (under the stated assumptions) with $G^c_2$ from \ref{s-sep-d-sep-eg} making d- and $\sigma$-separation equivalent. We then get the implications: 
\[  X \Indep_{G}^d Y \given Z \;\stackrel{\ref{s-sep-d-sep-eg}}{\iff}\;%
 X \Indep_{G}^\sigma Y \given Z \;\stackrel{\text{\hyl{gdGMP}{gdGMP}}}{\implies}\;  X \Indep_{\Pr_V} Y \given Z. \]
\end{proof}
\end{Thm}

\begin{Rem}
\label{gdGMP-dGMP}
If $G$ is an mDAG, ADMG or DAG then we always have $\Sc^G(v) =\{v\}$ for all $v \in V$.
So the condition $\{v,w\} \in H$ for all $v \in V$ and all $w \in \Sc^G(v)$ in Theorem \ref{dGMP-gdGMP} always holds and we get the equivalence:
\[ \xymatrix{ 
      \text{\hyl{dGMP}{dGMP}} \ar@{<=>}[r]& \text{\hyl{gdGMP}{gdGMP}}.
} \]
\end{Rem}

\begin{Lem}
\label{dGMP-dLMP}
Let $G=(V,E,H)$ be a \HEDG{} and $\Pr_V$ a probability distribution $\Pr_V$ on 
$\Xcal_V$. 
For $(G,\Pr_V)$ we then have:
 \[ \xymatrix{ 
      \text{\hyl{dGMP}{dGMP}} \ar@{=>}[r]& \text{\hyl{dLMP}{dLMP}}.
} \]
\begin{proof}
In the undirected graph $A_{\sm S}^\moral$ we clearly have the separation:
\[\{v\} \Indep_{A_{\sm S}^\moral} A_{\sm S} \sm \{v\} \given \partial_{A_{\sm S}^\moral}(v).\]
Furthermore, by \ref{d-sep-moral-hedg-subgraph} we have the implication: 
\[\{v\} \Indep_{A_{\sm S}^\moral} A_{\sm S} \sm \{v\} \given \partial_{A_{\sm S}^\moral}(v)
\;\implies\; \{v\} \Indep_{A_{\sm S}}^d A_{\sm S} \sm \{v\} \given \partial_{A_{\sm S}^\moral}(v).\]
Note that $A_{\sm S}$ is a marginalization of $G$ by \ref{anc-sub=mar}. With \ref{dGMP-stable-marg} we then get the implication:
\[\{v\} \Indep_{A_{\sm S}}^d A_{\sm S} \sm \{v\} \given \partial_{A_{\sm S}^\moral}(v)
\;\implies\; \{v\} \Indep_{\Pr_V} A_{\sm S} \sm \{v\} \given \partial_{A_{\sm S}^\moral}(v).\]
\end{proof}
\end{Lem}

\begin{Rem}
\label{DAG-dLMP}
For a DAG $G$ the \hyl{dLMP-DAG}{dLMP} for DAGs and the \hyl{dLMP}{dLMP} for \HEDG{}es are equivalent.
This can be seen from the fact that in the DAG case we have $\Sc^G(v) = \{v\}$ and $S=\emptyset$. 
Furthermore, for DAGs we can restrict to $A=\NonDesc^G(v)\cup\{v\}$, where we have $\partial_A(v)=\Pa^G(v)$.
Indeed, every ancestral $A' \ins A$ with $v \in A'$ will then also have $\partial_{A'}(v)=\Pa^G(v)$. This then gives the implication:
 \[ \{v\} \Indep_{\Pr_V} \NonDesc^G(v) \given \Pa^G(v) \; \implies \; \{v\} \Indep_{\Pr_V} A'_{\sm S} \sm \{v\} \given \partial_{A'_{\sm S}}(v).  \]  
\end{Rem}

\begin{Rem}
Note, that the directed local Markov property (\hyl{dLMP}{dLMP}) in general is not stable under marginalizations.
A counterexample is difficult to construct since for the acyclic case \hyl{dLMP}{dLMP} is equivalent to \hyl{dGMP}{dGMP} by \ref{dag-markov-prop}, which is stable under marginalizations by \ref{dGMP-stable-marg}. See \ref{dLMP-not-marg} where example \ref{main-example} is used to indirectly show that \hyl{dLMP}{dLMP}  is not stable under marginalizations.
\end{Rem}

\subsection{The Ordered Local Markov Property for \HEDG{}es}

In this subsection we define the \emph{ordered local Markov property (\hyl{oLMP}{oLMP})} for \HEDG{}es together with a total order. It will be defined in a way such that it is weak enough to be implied by the directed global Markov property (\hyl{dGMP}{dGMP}, see \ref{dGMP-oLMP}) for any total order, and by the directed local Markov property (\hyl{dLMP}{dLMP}) for assembling pseudo-topological orders (see \ref{dLMP-oLMP}), reduces to the usual ordered local Markov property for topological orders in the DAG case (see \ref{DAG-oLMP}), but the definition will be strong enough such that for \HEDG{}es with a perfect elimination order it will give the directed global Markov property (\hyl{dGMP}{dGMP}) back (see \ref{oLMP-dGMP}). This will generalize the DAG case to \HEDG{}es with perfect elimination orders in most generality.

\begin{Def}
\label{def-oLMP}
 Let $G=(V,E,H)$ be a \HEDG{}, $<$ a total order on $V$ and $\Pr_V$ a probability distribution $\Pr_V$ on 
$\Xcal_V=\prod_{v \in V} \Xcal_v$. Then we have the following property for the tuple $(G,\Pr_V,<)$:
\begin{enumerate}
\item[] The \emph{ordered local Markov property (\hyt{oLMP}{oLMP})} for total order $<$ on $G$: 
 For all $v \in V$ and all ancestral sub-\HEDG{}es $A \ins \Pred^G_\le(v)$ with $v \in A$ we have that: 
\[ \{v\} \Indep_{\Pr_V} A \sm \{v\} \given \partial_{A^\moral}(v).  \]
\end{enumerate}
\end{Def}

\begin{Rem}
\label{DAG-oLMP}
For a DAG $G=(V,E,H_1)$ the \hyl{oLMP-DAG}{oLMP} for DAGs (\ref{dag-markov-prop}) and the \hyl{oLMP}{oLMP} for \HEDG{}es (\ref{def-oLMP}) together with a topological order $<$ are equivalent.
Indeed, for every ancestral $A \ins \Pred^G_\le(v)$ with $v \in A$ we have $\partial_{A^\moral}(v)=\Pa^G(v)$.
So we have the implication:
\[ \{v\} \Indep_{\Pr_V} \Pred^G_\le(v) \sm \{v\} \given \Pa^G(v) \;\implies\; \{v\} \Indep_{\Pr_V} A \sm \{v\} \given \partial_{A^\moral}(v).  \]
Note that a more general version will be given in \ref{eg-oLMP}.
\end{Rem}

\begin{Lem}
\label{dGMP-oLMP}
Let $G=(V,E,H)$ be a \HEDG{} and $<$ any arbitrary total order on $V$ and $\Pr_V$ a probability distribution on $\Xcal_V$.
Then we have the implication:
\[ \xymatrix{ 
    \text{\hyl{dGMP}{dGMP} for } (G,\Pr_V)	 \ar@{=>}[r]&   \text{\hyl{oLMP}{oLMP} for } (G,\Pr_V,<).
} \]
\begin{proof}
For $A \ins \Pred^G_\le(v)$ ancestral with $v \in A$ we have 
$$\{v\} \Indep_{A^\moral} A\sm \{v\} \given \partial_{A^\moral}(v)$$
 which is equivalent to
  $$\{v\} \Indep_{A}^d A\sm \{v\} \given \partial_{A^\moral}(v)$$
	by \ref{d-sep-moral-hedg}. 
  Since the \hyl{dGMP}{dGMP} is preserved under marginalization by \ref{dGMP-stable-marg}
	we get $$\{v\} \Indep_{\Pr_V} A\sm \{v\} \given \partial_{A^\moral}(v)$$ and thus the claim. Note that $A$ really is a marginalization of $G$ by \ref{hedg-anc-marg}.
\end{proof}
\end{Lem}

\begin{Lem}
\label{dLMP-oLMP}
Let $G=(V,E,H)$ be a \HEDG{} and $\Pr_V$ a probability distribution on $\Xcal_V$. 
For any assembling pseudo-topological order $<$ on $G$ we have the implication:
\[ \xymatrix{ 
      \text{\hyl{dLMP}{dLMP} for } (G,\Pr_V) \ar@{=>}[r]& \text{\hyl{oLMP}{oLMP} for } (G,\Pr_V,<). 			
} \]
Furthermore, if we consider all assembling pseudo-topological orders $<$ on $G$ we also get the reverse implication:
\[ \xymatrix{ 
			\text{\hyl{oLMP}{oLMP} for } (G,\Pr_V,<)\quad \forall \text{ a.p.t.ord. ``$<$''} \ar@{=>}[r]& \text{\hyl{dLMP}{dLMP} for } (G,\Pr_V). 			
} \]

\begin{proof}
For the first statement assume \hyl{dLMP}{dLMP} and let $<$ be a fixed assembling pseudo-topological order $<$ on $G$.
Let $v \in V$ and $A \ins \Pred^G_\le(v)$ be ancestral with $v \in A$. 
Let $A':=\Anc^G(A)$ and $S:=A' \sm A$. By \ref{hedg-anc-marg} we then have $A = {A'}^{\marg \sm S}=:A'_{\sm S}$. 
Let $s \in S$ then $s>v$. Otherwise $s \in \Pred^G_\le(v) \cap A' = A$, which is not possible by definition of $S$. 
Since $s \in \Anc^G(A)$ we have $s \in \Anc^G(a)$ for some $a \in A$. Since $a \le v$ and $<$ is a pseudo-topological order we have
$s \in \Sc^G(a)$. Because $a \le v < s \in \Sc^G(a)$ and $<$ is assembling we have $v \in \Sc^G(a)$ and thus $s \in \Sc^G(v)$. 
It follows that $S \ins \Sc^G(v) \sm \{v\}$.
  By \hyl{dLMP}{dLMP} we then have:
\[ \{v\} \Indep_{\Pr_V} A'_{\sm S} \sm \{v\} \given \partial_{{A'_{\sm S}}^\moral}(v).  \]
Since $A=A'_{\sm S}$ this shows \hyl{oLMP}{oLMP}.\\
For the second statement assume \hyl{oLMP}{oLMP} for every assembling pseudo-topological order $<$ on $G$.
Let $v \in V$ and $A \ins \NonDesc^G(v) \cup \Sc^G(v)$ ancestral with $v \in A$ and $S \ins \Sc^G(v) \sm \{v\}$ a set.
Since $A$ is ancestral in $G$ as well we can find an assembling pseudo-topological order $<$ such that:
\begin{enumerate}
 \item If $w_1 \in A$ and $w_2 \in V \sm A$ then $w_1 < w_2$.
 \item If $w \in S$ then $v < w$.
 \item If $w \in \Sc^G(v)\sm\lp \{v\} \cup S\rp$ then $w < v$. 
\end{enumerate}
In this order we then have that $A_{\sm S} = \Pred^G_\le(v)$, which is ancestral in $\Pred^G_\le(v)$ and contains $v$.
The oLMP applied to $(G,\Pr_V,<)$ then gives:
\[  \{v\} \Indep_{\Pr_V} A_{\sm S}\sm\{v\} \given \partial_{A_{\sm S}}(v).\]
This shows dLMP.
\end{proof}
\end{Lem}

In the next theorem we will show that for every \HEDG{} that has a perfect elimination order the ordered local Markov property (\hyl{oLMP}{oLMP}) and the directed global Markov property (\hyl{dGMP}{dGMP}) are equivalent. The proof is a generalization of the DAG case (see \cite{Lau90} \S6 and \cite{Verma93} Thm. 1.2.1.3) and ADMG case (see \cite{Richardson03} \S3 Thm. 2), which both only treat topological orders, to the much more general case of \HEDG{}es with perfect elimination orders. Note that the proof only uses the general separoid/graphoid rules (1)-(5) from \ref{GrAxPr}.

\begin{Thm} 
\label{oLMP-dGMP}
Let $G=(V,E,H)$ be a \HEDG{} that has a perfect elimination order $<$ and $\Pr_V$ a probability distribution on $\Xcal_V$.
Then we have the implication:
\[ \xymatrix{ 
      \text{\hyl{oLMP}{oLMP} for } (G,\Pr_V,<) \ar@{=>}[r]& \text{\hyl{dGMP}{dGMP} for } (G,\Pr_V).
} \]
\begin{proof} 
For an ancestral sub-\HEDG{} $A \ins G$ and $v \in V$ let $A(v):=\Pred_\le^A(v):=\{ w \in A | w < v \text{ or } w=v\}$ with the marginalized \HEDG{} structure.
Note that $A(v)$ is then an ancestral sub-\HEDG{} in $G(v)=\Pred^G_\le(v)$ by \ref{hedg-anc-marg}. 
If $(G,\Pr_V,<)$ satisfies \hyl{oLMP}{oLMP} then $(G(v),\Pr_{G(v)},<)$ satisfies \hyl{oLMP}{oLMP} (with the induced ordering).
 With this and $G=G(v_n)$, where the elements of $V=\{v_1,\dots,v_n\}$ are indexed according to $<$, we can do induction. \\
Let $X,Y,Z \ins V$ with $X \Indep_G^d Y \given Z$ be given subsets and $A:=\Anc^G(X \cup Y \cup Z)$ with the sub-\HEDG{}-structure.
  Let $v \in V $ be the smallest element (according to $<$) with $X,Y,Z \ins A(v)$. Then $v \in X \cup Y \cup Z$. If $v=v_1$ we trivially have $X \Indep_{\Pr_V}Y \given Z$. \\
  So consider $v=v_m$ with $m>1$. By induction we can assume that \hyl{dGMP}{dGMP} holds for $(A(w),\Pr_{A(w)})$ for all $w < v$. 
  \[X \Indep_{G}^d Y \given Z \stackrel{\ref{marg-d-sep}}{\iff} X \Indep_{G(v)}^d Y \given Z \stackrel{\ref{d-sep-moral-hedg}}{\iff} 
  X \Indep_{A(v)^\moral} Y \given Z \iff X\sm Z \Indep_{A(v)^\moral} Y \sm Z \given Z\]
where we note that $A(v) = \Anc^{G(v)}(X \cup Y \cup Z)$ by definition.
  Since we also have: $X \Indep_{\Pr_{A(v)}} Y \given Z \iff X\sm Z \Indep_{\Pr_{A(v)}} Y \sm Z \given Z$ we can assume that $X,Y,Z$ are pairwise disjoint. 
We now have three cases: \\
First consider: $v \in X$. Then $X=X' \cup \{v\}$ with $v \notin \partial_{A(v)^\moral}(v) \cup X' \cup Y \cup Z$. 
So $ \partial_{A(v)^\moral}(v) \cup X' \cup Y \cup Z \ins A(w)$ with $w=v_{m-1} <v$, for which induction applies.
Then $A(w)=A(v)^{\marg\sm\{v\}}$. By \ref{marg-mor1} we then have that $A(w)^\moral$ is a subgraph of $A(v)^{\moral}_\marg:=(A(v)^\moral)^{\marg\sm\{v\}}$.
By standard separoid/semi-graphoid rules \ref{GrAxPr} we have the implications (symmetry \ref{GrAxPr} (2) is implicitly  used): 
\begin{align*}
X' \cup \{v\} \Indep_{A(v)^\moral} Y \given Z & \stackrel{\ref{GrAxPr}\; (3)}{\implies} & X' \Indep_{A(v)^\moral} Y \given Z \\
& \stackrel{\ref{sep-marg}}{\implies} & X' \Indep_{A(v)^{\moral}_\marg} Y \given Z \\
& \stackrel{\ref{sep-subgraph}}{\implies} & X' \Indep_{A(w)^\moral} Y \given Z \\
& \stackrel{\ref{d-sep-moral-hedg-subgraph}}{\implies} & X' \Indep_{A(w)}^d Y \given Z \\
 & \stackrel{\text{ind.}}{\implies} & X' \Indep_{\Pr_V} Y \given Z. \label{eqn:bli} \tag{\#1}\\
X' \cup \{v\} \Indep_{A(v)^\moral} Y \given Z & \stackrel{\ref{GrAxPr}\; (4)}{\implies} & \{v\} \Indep_{A(v)^\moral} Y \given X' \cup Z   \\
 & \stackrel{(*)}{\implies} & \partial_{A(v)^\moral}(v)  \Indep_{A(v)^\moral} Y \given X' \cup Z \\ 
& \stackrel{\ref{sep-marg}}{\implies} & \partial_{A(v)^\moral}(v)  \Indep_{A(v)^{\moral}_\marg} Y \given X' \cup Z \\
 &\stackrel{\ref{sep-subgraph}}{\implies} & \partial_{A(v)^\moral}(v)  \Indep_{A(w)^\moral} Y \given X' \cup Z \\
& \stackrel{\ref{d-sep-moral-hedg-subgraph}}{\implies} & \partial_{A(v)^\moral}(v)  \Indep_{A(w)}^d Y \given X' \cup Z \\ 
& \stackrel{\text{ind.}}{\implies} &   \partial_{A(v)^\moral}(v)  \Indep_{\Pr_V} Y \given X' \cup Z. \label{eqn:bla} \tag{\#2} 
\end{align*}
\begin{align*}
\{v\} \Indep_{\Pr_V} A(v)\sm\{v\} \given \partial_{A(v)^\moral}(v) &\stackrel{\ref{GrAxPr}\; (3)}{\implies}& \{v\} \Indep_{\Pr_V}  Y \cup X' \cup Z \given \partial_{A(v)^\moral}(v) \\
& \stackrel{\ref{GrAxPr}\; (4)}{\implies} & \{v\} \Indep_{\Pr_V} Y \given X' \cup Z \cup \partial_{A(v)^\moral}(v).  \label{eqn:blub} \tag{\#3}\\
\eqref{eqn:bla} \;\&\; \eqref{eqn:blub}
& \stackrel{\ref{GrAxPr}\; (5)}{\implies} & \{v\} \cup \partial_{A(v)^\moral}(v)  \Indep_{\Pr_V} Y \given X' \cup Z \\
& \stackrel{\ref{GrAxPr}\; (4)}{\implies} & \{v\} \Indep_{\Pr_V} Y \given X' \cup Z.  \label{eqn:blob} \tag{\#4} \\
\eqref{eqn:bli} \;\&\; \eqref{eqn:blob}
& \stackrel{\ref{GrAxPr}\; (5)}{\implies} & X' \cup \{v\} \Indep_{\Pr_V} Y \given Z\\
& \stackrel{=}{\implies} & X \Indep_{\Pr_V} Y \given Z.
\end{align*}
(*) holds since every $X' \cup Z$-open path to $\partial_{A(v)^\moral}(v)$ extends to an $X' \cup Z$-open path to $v$.\\
If $v \in Y$ by symmetry the same arguments hold.\\
Now assume the case when $v \in Z$. Then $X \Indep_{A(v)^\moral} Y \given Z' \cup \{v\}$ and $\partial_{A(v)^\moral}(v) \cup X \cup Y \cup Z' \ins 
A(w)$ hold for $w < v$.
We now argue that $X  \cup \{v\} \Indep_{A(v)^\moral} Y \given Z'$ or $X \Indep_{A(v)^\moral} Y \cup \{v\} \given Z'$ holds:
Assume the contrary. With $X \Indep_{A(v)^\moral} Y \given Z' \cup \{v\}$, we then had $Z'$-open paths $v - w_1 - \cdots - y$ and $x - \cdots - w_2 - v$.
By considering shortest paths we can assume that all the nodes in between are not $v$. We also have $w_1,w_2 \in \partial_{A(v)^\moral}(v)$, which is a complete subgraph by assumption. So we have the edge $w_1-w_2$ (or even equality) in $A(v)^\moral$. 
So we get a $Z' \cup \{v\}$-open path $x - \cdots - w_2 - w_1 - \cdots - y$ in $A(v)^\moral$ in contradiction to
$X \Indep_{A(v)^\moral} Y \given Z' \cup \{v\}$.\\
So $X  \cup \{v\} \Indep_{A(v)^\moral} Y \given Z'$ or $X \Indep_{A(v)^\moral} Y \cup \{v\} \given Z'$ must hold.\\
Now the previous points imply $X  \cup \{v\} \Indep_{\Pr_V} Y \given Z'$ or $X \Indep_{\Pr_V} Y \cup \{v\} \given Z'$, which both imply 
$X \Indep_{\Pr_V} Y \given Z' \cup \{v\}$ by \ref{GrAxPr} (2+4).
\end{proof}
\end{Thm}

\begin{Cor} 
\label{oLMP-dLMP-dGMP}
Let $G=(V,E,H)$ be a \HEDG{} that has an assembling quasi-topological order $<$ (e.g.\ a topological order) and $\Pr_V$ a probability distribution on $\Xcal_V$.
Then we have the equivalences:
\[ \xymatrix{ 
      \text{\hyl{oLMP}{oLMP} for } (G,\Pr_V,<) \ar@{<=>}[r]& \text{\hyl{dLMP}{dLMP} for } (G,\Pr_V) \ar@{<=>}[r]& \text{\hyl{dGMP}{dGMP} for } (G,\Pr_V).
} \]
\begin{proof} 
  \hyl{dGMP}{dGMP} $\implies$ \hyl{dLMP}{dLMP}: Use \ref{dGMP-dLMP}.\\
  \hyl{dLMP}{dLMP} $\implies$ \hyl{oLMP}{oLMP}: Use \ref{dLMP-oLMP} for assembling pseudo-topological orders.\\
  \hyl{oLMP}{oLMP} $\implies$ \hyl{dGMP}{dGMP}: Use \ref{oLMP-dGMP} for perfect elimination orders.
\end{proof}
\end{Cor}

In the following Corollary we will investigate under which circumstances we only need to check one set for the conditional independence per node (instead of all ancestral subsets it is contained in). 

\begin{Cor}
\label{eg-oLMP}
Let $G=(V,E,H)$ be a \HEDG{}, $\Pr_V$ a probability distribution on $\Xcal_V$ and a total order $<$. If 
for every $v \in V$ we have:
\begin{enumerate}
\item $\Dist^{G(v)}(\{v\} \cup \Ch^{G(v)}(v) ) \ins \Anc^{G(v)}(v)$ as sets of nodes, and
\item $\partial_{\Anc^{G(v)}(v)^\moral}(v)$ is a complete subgraph of $\Anc^{G(v)}(v)^\moral$, and%
\item $\{v\} \Indep_{\Pr_V} G(v) \sm \{v\} \given \partial_{G(v)^\moral}(v)$,
\end{enumerate}
with $G(v):=\Pred^G_\le(v)$,
then $<$ is a perfect elimination order and $(G,\Pr_V,<)$ satisfies the \hyl{oLMP}{oLMP} and thus the \hyl{dGMP}{dGMP}.
\begin{proof}
For $A \ins G(v)$ ancestral with $v \in A$ we have seen in \ref{eg-peo} that $\partial_{A^\moral}(v)=\partial_{G(v)^\moral}(v)$. 
So we immediately get by the decomposition rule \ref{GrAxPr}:
\[ \{v\} \Indep_{\Pr_V} A \sm \{v\} \given \partial_{A^\moral}(v).\]
This shows \hyl{oLMP}{oLMP}. \ref{eg-peo} shows that $<$ must be a perfect elimination order.
The claim then follows directly from \ref{oLMP-dGMP}.
\end{proof}
\end{Cor}

\subsection{Undirected Markov Properties for \HEDG{}es}

In this subsection we investigate Markov properties for \HEDG{}es that are related to Markov properties for undirected graphs (mainly local, global and pairwise Markov properties) by the use of the moralization of a \HEDG{}. To capture all the relevant information of the \HEDG{} we need to consider the moralization of every ancestral sub-\HEDG{} or even every marginalization of it. In the latter case (or under the intersection property) we have the equivalence to the directed gobal Markov property (\hyl{dGMP}{dGMP}). The precise statements are given in the theorems \ref{dGMP-auPMP} and \ref{auPMP-dGMP}. Later in \ref{auPMP-aFP} we will see that the undirected Markov properties will imply the ancestral factorization property (\hyl{aFP}{aFP}) under the existence of a strictly positive density $p_V >0$ w.r.t.\ some product measure.

\begin{Def} Let $G=(V,E,H)$ be a \HEDG{} and $\Pr_V$ a probability distribution $\Pr_V$ on 
$\Xcal_V=\prod_{v \in V} \Xcal_v$. Then we have the following properties relating $G$ to $\Pr_V$:
\begin{enumerate}

  \item The \emph{ancestral undirected global Markov property (\hyt{auGMP}{auGMP})}: For all ancestral sub-\HEDG{}es $A \ins G$ and all subsets $X,Y,Z \ins A$ we have the implication:
 \[ X \Indep_{A^\moral} Y \given Z \;\implies\; X \Indep_{\Pr_V} Y \given Z.\] 

\item The \emph{refined undirected global Markov property (\hyt{ruGMP}{ruGMP})}: For all subsets $W \ins V$ and all subsets $X,Y,Z \ins W$ we have the implication:
 \[ X \Indep_{(G^{\marg(W)})^\moral} Y \given Z \;\implies\; X \Indep_{\Pr_V} Y \given Z.\]

\item The \emph{ancestral undirected local Markov property (\hyt{auLMP}{auLMP})}: For all ancestral sub-\HEDG{}es $A \ins G$ and every $v \in A$ we have:
 \[ \{v\} \Indep_{\Pr_V} A \sm \{v\} \given \partial_{A^\moral}(v).\]

\item The \emph{refined undirected local Markov property (\hyt{ruLMP}{ruLMP})}: For all subsets $W \ins V$ and every $v \in W$ we have:
 \[ \{v\} \Indep_{\Pr_V} W \sm \{v\} \given \partial_{(G^{\marg(W)})^\moral}(v).\]

\item The \emph{ancestral undirected pairwise Markov property (\hyt{auPMP}{auPMP})}: For all ancestral sub-\HEDG{}es $A \ins G$ and every pair of nodes
$v,w \in A^\moral$ we have:
\[ w \notin \partial_{A^\moral}(v) \cup \{v\} \;\implies \; \{v\} \Indep_{\Pr_V} \{w\} \given A\sm\{v,w\}. \]

\item The \emph{refined undirected pairwise Markov property (\hyt{ruPMP}{ruPMP})}: For all subsets $W \ins V$ and every pair of nodes
$v,w \in W$ we have:
\[ w \notin \partial_{(G^{\marg(W)})^\moral}(v) \cup \{v\} \;\implies \; \{v\} \Indep_{\Pr_V} \{w\} \given W\sm\{v,w\}. \] 
\end{enumerate}
\end{Def}

\begin{Thm} 
\label{dGMP-auPMP}
For a \HEDG{} $G=(V,E,H)$ and a probability distribution $\Pr_V$ on $\Xcal_V$ we have the implications:
\[ \xymatrix{ 
      \text{\hyl{dGMP}{dGMP}}\ar@{<=>}[r]\ar@{<=>}[dr]& \text{\hyl{ruGMP}{ruGMP}} \ar@{<=>}[d]\ar@{<=>}[r] & \text{\hyl{ruLMP}{ruLMP}} \ar@{<=>}[r]\ar@{=>}[d] & \text{\hyl{ruPMP}{ruPMP}}\ar@{=>}[d] 		\\
                             & \text{\hyl{auGMP}{auGMP}} \ar@{=>}[r]& \text{\hyl{auLMP}{auLMP}} \ar@{=>}[r] & \text{\hyl{auPMP}{auPMP}} 		
} \]
\begin{proof} 
The implications from top down only need that for ancestral sub-\HEDG{}es $A \ins G$ we have $A=G^{\marg(A)}$ by \ref{anc-sub=mar}. \\
  ``\hyl{dGMP}{dGMP} $\implies$ \hyl{ruGMP}{ruGMP}'': Let $W \ins V$ and $X,Y,Z \ins W$ and $G_W:=G^{\marg(W)}$ and $A_W:=\Anc^{G_W}(X \cup Y \cup Z)$. Then $A_W^\moral \ins G_W^\moral$ is a subgraph and we have:
\[ \begin{array}{lllll}
X \Indep_{G_W^\moral} Y \given Z & \stackrel{\ref{d-sep-moral-hedg-subgraph}}{\implies} &
 X \Indep_{G_W}^d Y \given Z 
&\stackrel{\ref{marg-d-sep}}{\iff}& X \Indep_{G}^d Y \given Z \\
  &\stackrel{\text{\hyl{dGMP}{dGMP}}}{\implies}& X \Indep_{\Pr_V} Y \given Z. 
\end{array}\]
  This shows \hyl{ruGMP}{ruGMP}.\\
All the other implications from left to right are clear.\\
  ``\hyl{auGMP}{auGMP} $\implies$ \hyl{dGMP}{dGMP}'': Let $X,Y,Z \ins V$ and $A:=\Anc^G(X \cup Y \cup Z)$. Then we have:
  \[  X \Indep_{G}^d Y \given Z \qquad \stackrel{\ref{d-sep-moral-hedg}}{\iff} \qquad X \Indep_{A^\moral} Y \given Z  \quad\stackrel{\text{\hyl{auGMP}{auGMP}}}{\implies}\quad X \Indep_{\Pr_V} Y \given Z.  \]
  ``\hyl{ruPMP}{ruPMP} $\implies$ \hyl{dGMP}{dGMP}'': Let $X,Y,Z \ins V$ with $X \Indep_G^d Y \given Z$. By the separoid/graphoid rules \ref{GrAxPr} we have the equivalence for $X' = X \sm \{x\}$:
\[ \{x\}  \cup X'  \Indep_G^d Y \given Z \qquad \iff\qquad \{x\} \Indep_G^d Y \given Z \cup X' \qquad\land\qquad X' \Indep_G^d Y \given Z.  \]
Repeating this and together with symmetry we end up with a number of statements of the form:
\[ \{x\} \Indep_G^d \{y\} \given Z'.\]
Take $W:=\{x,y\} \cup Z'$ with the marginalized \HEDG{}-structure of $G$. Noting that $\Anc^W(W)=W$, we get:
  \[ \{x\} \Indep_G^d \{y\} \given Z' \quad\stackrel{\ref{marg-d-sep}}{\iff}\quad \{x\} \Indep_W^d \{y\} \given Z' \quad\stackrel{\ref{d-sep-moral-hedg}}{\iff}\quad \{x\} \Indep_{W^\moral} \{y\} \given Z'. \]
  The \hyl{ruPMP}{ruPMP} now implies: $ \{x\} \Indep_{\Pr_V} \{y\} \given Z'$. By using the same separoid/graphoid rules \ref{GrAxPr} backwards we get:
\[X \Indep_{\Pr_V} Y \given Z. \]

\end{proof}
\end{Thm}

\begin{Lem}[Undirected pairwise Markov property under marginalization]
\label{uPMP-marg}
Let $G=(V,E)$ be an undirected graph, $\Pr_V$ a probability distribution on $\Xcal_V$ and $W \ins V$ a subset.
Assume that $(G,\Pr_V)$ satisfy the \emph{undirected pairwise Markov property} (\hyt{uPMP}{uPMP}): I.e. for every $v_1,v_2 \in V$ we have:
\[  v_1 \notin \partial_G(v_2) \cup \{v_2\} \; \implies \; \{v_1\} \Indep_{\Pr_V} \{v_2\} \given G \sm \{v_1,v_2\}.  \] 
If $\Pr_V$ satisfies the intersection property (\ref{GrAxPr} (6)), then $(G^{\marg \sm W},\Pr_{V\sm W})$ also satisfies the undirected pairwise Markov property (\hyl{uPMP}{uPMP}).
\begin{proof}
By induction we can assume that $W=\{w\}$. Now let $v_1,v_2 \in V \sm W$. We have the sets:
\[\partial_{G^{\marg \sm W}}(v_i) = \left\{ \begin{array}{lll}
  \partial_G(v_i) & \text{ if } & w \notin \partial_G(v_i), \\
  (\partial_G(v_i) \cup \partial_G(w)) \sm \{v_i, w\} & \text{ if } & w \in \partial_G(v_i).
\end{array} \right.\]
  Put $S:=G \sm \{v_1,v_2,w\}$ and let $v_1 \notin \partial_{G^{\marg \sm W}}(v_2) \cup \{v_2\}$. Then $v_1 \notin \partial_G(v_2) \cup \{v_2\}$ and we have $w \notin \partial_G(v_i) \cup \{v_i\}$ for at least one $i=1,2$, say $w \notin \partial_G(v_2) \cup \{v_2\}$. Then we have by \hyl{uPMP}{uPMP} for $(G,\Pr_V)$:
\[\{v_1\} \Indep_{\Pr_V} \{v_2\} \given S \cup \{w\}, \qquad \text{ and } \qquad
\{w\} \Indep_{\Pr_V} \{v_2\} \given S \cup \{v_1\}.\]
By the intersection property we get:
\[\{v_1,w\} \Indep_{\Pr_V} \{v_2\} \given S, \]
which implies the claim:
\[\{v_1\} \Indep_{\Pr_V} \{v_2\} \given G \sm \{v_1,v_2,w\}. \]
\end{proof}
\end{Lem}

\begin{Thm} 
\label{auPMP-dGMP}
  For a \HEDG{} $G=(V,E,H)$ and a probability distribution $\Pr_V$ on $\Xcal_V$ that satisfies the intersection property \ref{GrAxPr} (6) (e.g.\ in the case of a strictly positive density $p_V >0$ w.r.t.\ a product measure $\mu_V$ of %
	measures $(\mu_v)_{v\in V}$) we have the implications:
\[ \xymatrix{ 
      \text{\hyl{auPMP}{auPMP}}\ar@{=>}[r]& \text{\hyl{auGMP}{auGMP}}, 			
} \]
making all the implications from \ref{dGMP-auPMP} equivalences.
\begin{proof} 
See \cite{Lau90} Thm. 3.7 or \cite{PeaPaz87}.
By \ref{dGMP-auPMP} we only need to show \hyl{ruPMP}{ruPMP}.
For this let $W \ins V$ be a subset. We then need to show the \hyl{uPMP}{uPMP} for $(G^{\marg(W)})^\moral$.
Let $A:=\Anc^G(W)$. Then $A^{\marg(W)} = G^{\marg(W)}$ and every node $v \in \Anc^G(W) \sm W$ has a directed path $v \to \cdots \to w$ in $A$ with a $w \in W$. In particular, we have $\Ch^A(v)\sm \{v\} \neq \emptyset$. By \ref{marg-mor2} we then get that:
\[ (A^{\marg\sm\{v\}})^\moral = (A^\moral)^{\marg\sm\{v\}}. \tag{*} \]
Using the same argument for a node in $ v' \in \Anc^{A^{\marg\sm\{v\}}}(W) \sm W$ we get:
\[ (A^{\marg\sm\{v,v'\}})^\moral \stackrel{\ref{marg-mor2}}{=} ((A^{\marg\sm\{v\}})^\moral)^{\marg\sm\{v'\}} \stackrel{(*)}{=} (A^\moral)^{\marg\sm\{v,v'\}}. \]
Inductively we arrive at:
\[ (G^{\marg(W)})^\moral= (A^{\marg(W)})^\moral = (A^\moral)^{\marg(W)}. \]
By the assumption \hyl{auPMP}{auPMP} $A^\moral$ satisfies the \hyl{uPMP}{uPMP}. By \ref{uPMP-marg} and the intersection property also $(A^\moral)^{\marg(W)}=(G^{\marg(W)})^\moral$ satisfies the \hyl{uPMP}{uPMP}. This implies \hyl{ruPMP}{ruPMP} and thus the claim.
\end{proof}
\end{Thm}

\subsection{Factorization Properties for \HEDG{}es}

In this subsection we define the most genereral factorization property for \HEDG{}es with an associated probability density $p$: the \emph{ancestral factorization property (\hyl{aFP}{aFP})}, and variants of it. It basically holds if $p$ factors according to the maximal complete subgraphs of the moralization of every ancestral sub-\HEDG{} of $G$.  We will see that this implies the directed global Markov property (\hyl{dGMP}{dGMP}, see \ref{aFP-dGMP}).
Furthermore, there are basically two ways back (see \ref{dGMP-aFP}):\\
One way is in the case of a strictly positive density $p>0$. Then we will show that the ancestral factorization property (\hyl{aFP}{aFP}) is implied by the ancestral undirected pairwise Markov property (\hyl{auPMP}{auPMP}, see \ref{auPMP-aFP}), which itself is implied by the directed global Markov property (\hyl{dGMP}{dGMP}, see \ref{dGMP-auPMP}). \\
The second way is in the case where we have a perfect elimination order for $G$ (e.g.\ in the case of an mDAGs). Then we will prove that a general version of the recursive factorization property (\hyl{rFPwd}{rFPwd}) for \HEDG{}es will imply the ancestral factorization property (\hyl{aFP}{aFP}, see \ref{rFPwd-aFP}). We will also see that the recursive factorization property (\hyl{rFP}{rFP}) is equivalent to the ordered local Markov property (\hyl{oLMP}{oLMP}, see \ref{oLMP-rFP}), which is again implied by the directed global Markov property (\hyl{dGMP}{dGMP}, see \ref{dGMP-oLMP}), closing the implication loop.\\
The question when a structural equations property (\hyl{ausSEP}{SEP}) will imply the ancestral factorization property (\hyl{aFP}{aFP}) will be handled later in \ref{ausSEP-wmaFP} for discrete variables under ``\hyl{ausSEP}{ausSEP}'' and in \ref{SEPwared-aFP} under the different assumption ``\hyl{SEPwared}{SEPwared}'', which is a generalization of the linear to the non-linear case. Note, that by the combination of ``\hyl{csSEP}{csSEP}'' and the results \ref{csSEP-smgdGMP} and \ref{gdGMP-dGMP} we will get the recursive factorization property (\hyl{rFP}{rFP}) for the acyclification, e.g.\ in the mDAG case itself.  
Also see \ref{dag-markov-prop} for the DAG case and the reverse implication.

\begin{Def} Let $G=(V,E,H)$ be a \HEDG{} and $\Pr_V$ a probability distribution $\Pr_V$ on 
$\Xcal_V=\prod_{v \in V} \Xcal_v$. Then we have the following properties relating $G$ to $\Pr_V$:
\begin{enumerate}

  \item The \emph{ancestral factorization property (\hyt{aFP}{aFP})} for %
measures $\mu_v$ on $\Xcal_v$, $v \in V$: 
There exists a density $p_V$ for $\Pr_V$ w.r.t.\ the product measure $\mu_V:=\otimes_{ v \in V} \mu_v$ such that
    for all ancestral sub-\HEDG{}es $A \ins G$ there exists a finite tuple $\Ccal_A$ of complete subgraphs\footnote{Again, note that complete subgraphs are often called \emph{cliques} in some literature.} $C \ins A^\moral$ and for every such complete subgraph $C \in \Ccal_A$ there exists an integrable function 
\[k_{A,C} : \Xcal_C \to \R_{\ge 0} \]
such that for $\mu_A$-almost-all $x \in \Xcal_A$ the marginal density w.r.t.\ to $A$, 
$ p_A:  \Xcal_A \to \R$, is given by the product:
\[ p_A(x_A) = \prod_{C \in \Ccal_A} k_{A,C}(x_C), \]
where $x_A$ is the tuple $(x_v)_{v \in A}$ and $x_C = (x_v)_{ v \in C}$ denotes the tuple of values $x_v$ corresponding to the nodes in $C$.

\item The \emph{marginal factorization property (\hyt{mFP}{mFP})} for measures $\mu_v$ on $\Xcal_v$, $v \in V$:
 For every $F \in \tilde{H}$ there exists a standard Borel space $\Ecal_F$ with a %
 measure $\nu_F$ on it and for every strongly connected component $S \in \Scal(G^\aug)$ of $G^\aug$ there is a finite tuple $\Ccal_S$ of complete subgraphs 
 $C \ins \Anc^{G^\aug}(S)^\moral$ in the moralization of the ancestral closure of $S$ in the augmented graph $G^\aug$ with $C \ins S \cup \Pa^{G^\aug}(S)$, and for every such complete subgraph $C \in \Ccal_S$ there exists an integrable function 
\[k_{S,C} :  \Xcal_{C \cap V } \x  \Ecal_{C \sm V } \to \R_{\ge 0} \]
such that there is a joint density 
    $p_{V^\aug}:  \,\Xcal_{V^\aug}:=\Xcal_V \x \Ecal_V \to \R_{\ge 0} $ 
		(with $\Ecal_V:=\prod_{F \in \tilde{H}}\Ecal_F$)  
		w.r.t.\ the product measure $\mu_{V^\aug}:=(\otimes_{v\in V} \mu_v) \otimes (\otimes_{F \in \tilde{H}} \nu_F) $ which is given by:
\[ p_{V^\aug}(x,e) =  \prod_{S \in \Scal(G^\aug)} \prod_{C \in \Ccal_S} k_{S,C}(x_C,e_C), \]
and with the property that for all $S \in \Scal(G^\aug)$ we have:
\[\begin{array}{rclll}
\int \prod_{C \in \Ccal_S} k_{S,C}(x_C,e_C) \, d\mu_S(x_S) &=& 1, & \text{for} & S \in \Scal(G),\\
\int \prod_{C \in \Ccal_S} k_{S,C}(e_F)\, d\nu_F(e_F) &=& 1,& \text{for} & S=\{F\}, F \in \tilde{H}.\\
\end{array}\]
Here $x_C,e_C$ are the corresponding components of $x\in \Xcal_V$ and $e \in \Ecal_V$. %
Note that if $S$ is a strongly connected component of $G^\aug$ that is not in $G$ then $S$ consists of a single variable $S=\{F\}$ corresponding to a hyperedge $F \in \tilde{H}$.

\item The \emph{recursive factorization property (\hyt{rFP}{rFP})} for a total order $<$:
For every ancestral sub-\HEDG{} $A \ins G$ and every measurable product set $B:=\prod_{v \in A,\ge} B_v \ins \Xcal_A$, where the product from left to right is taken in decreasing order according to $<$, 
we have the equality
\[ \Pr_A(B) = \int \cdots \int \prod_{v \in A, \ge} \I_{B_v}(x_v)\, d\Pr_V^{v|x_{\partial_\le^A(v)} }(x_v), \]
    where the product and integration of the regular conditional probabilities $\Pr_V^{v|x_{\partial_\le^A(v)}}$ of the $v$-component given the value $x_{\partial_\le^A(v)}$ is taken from the highest $v$ to the lowest $v$ according to $<$, 
		and where we made the abbreviation:
\[ \partial_\le^A(v):= \partial_{\Pred_\le^A(v)^\moral }(v), \]
the neighbours of $v$ in the moralization of the marginalization of all elements in $A$ that are smaller than or equal to $v$.

\item The \emph{recursive factorization property with density (\hyt{rFPwd}{rFPwd})} for %
 measures $\mu_v$, $v \in V$, and total order $<$: There exists a density $p=p_V$ for $\Pr_V$ w.r.t.\ the product measure $\mu_V$ such that
for every ancestral sub-\HEDG{} $A \ins G$ and for $\mu_A$-almost-all $x_A=(x_v)_{v \in A} \in \Xcal_A$ we have the equality:
\[ p_A(x_A) = \prod_{v \in A} p_A(x_v|x_{\partial^A_\le(v)}).  \]

\end{enumerate}
\end{Def}

\begin{Rem}
\label{DAG-rFP}
  Note that the recursive factorization property (\hyl{rFP-DAG}{rFP}) for a DAG with topological order $<$ is equivalent to the recursive factorization property (\hyl{rFP}{rFP}) for \HEDG{}es in that topological order. Indeed for every ancestral $A \ins G$ and every $v \in A$ we have $\partial^A_\le(v)=\Pa^G(v)$ and the factorization for $A$ follows from marginalizing out the other variables $G \sm A$ from the factorization of $G$.
\end{Rem}

\begin{Lem}
\label{oLMP-rFP}
Let $G=(V,E,H)$ be a \HEDG{} and  $\Pr_V$ a probability distribution on $\Xcal_V$ and $<$ a total order on $V$.
Then for $(G,\Pr_V,<)$ we have the equivalence:
\[ \xymatrix{ 
    \text{\hyl{oLMP}{oLMP}}\ar@{<=>}[r] &	\text{\hyl{rFP}{rFP}}. 	
} \]
\begin{proof}
This follows from \ref{hedg-anc-marg}.
  \hyl{oLMP}{oLMP} $\implies$ \hyl{rFP}{rFP}: Let $A \ins G$ be ancestral and $v \in A$.
Then $A(v): = \Pred^A_\le(v)$ is ancestral in $\Pred^G_\le(v)$ by \ref{hedg-anc-marg} with node set $A \cap \Pred^G_\le(v)$.
  By \hyl{oLMP}{oLMP} we then have:
\[ \{v\} \Indep_{\Pr_V} A(v) \sm \{v\} \given \partial_{A(v)^\moral}(v),  \]
which translates to 
\[ \Pr_A( T_v \in B_v |T_{A(v)\sm \{v\}}=x_{A(v)\sm \{v\}}) = \Pr_A( T_v \in B_v |T_{\partial_{A(v)^\moral}(v)}=x_{\partial_{A(v)^\moral}(v)}) \]
for all measurable sets $B_v \ins \Xcal_v$ and almost all values $x_{A(v)\sm \{v\}} \in \Xcal_{A(v)\sm \{v\}}$, or in short:
\[ d\Pr_{A(v)}^{v|x_{A(v)\sm \{v\}}}(x_v) =  d\Pr_{A(v)}^{v|x_{\partial^A_\le(v)}}(x_v).  \]
Since we always have the product rule:
\[ d\Pr_{A(v)}(x) = d\Pr_{A(v)}^{v|x_{A(v)\sm \{v\}}}(x_v) \; d\Pr_{A(v)\sm \{v\}}(x_{A(v)\sm \{v\}})  \]
  we can just do induction from the highest $v \in A$ to the lowest $v \in A$ to get \hyl{rFP}{rFP}.\\
  \hyl{rFP}{rFP} $\implies$ \hyl{oLMP}{oLMP}: Now let $w \in V$ and $A' \ins \Pred^G_\le(w)$ ancestral with $w \in A'$. By \ref{hedg-anc-marg} we have that $A' = G^{\marg(A \cap \Pred^G_\le(w))}$ 
	for an ancestral $A \ins G$. We then get $A'=\Pred^A_\le(w)$. Note that for $v \in A'$,  $\Pred^A_\le(v) = \Pred^{A'}_\le(v)$ and $\partial^A_\le(v) = \partial^{A'}_\le(v)$. Therefore, by the mentioned product rule, by \hyl{rFP}{rFP} and after integrating out the factors corresponding to $v \in A \sm A'$ we get the factorization:
\[  \prod_{v \in A', \ge} d\Pr_{A'}^{v|x_{\Pred^A_\le(v) \sm \{v\}}}(x_v) =  d \Pr_{A'}(x) = \prod_{v \in A', \ge} d\Pr_{A'}^{v|x_{\partial^{A}_\le(v)}}(x_v).  \]
So by the uniqueness of the conditional expectations (see \cite{Kle14} Thm. 8.12) we directly get for the highest node $v=w$:
\[ d\Pr_{A'}^{w|x_{\Pred^A_\le(w) \sm \{w\}}}(x_w) =  d\Pr_{A'}^{w|x_{\partial^A_\le(w)}}(x_w).  \]
Since $A'=\Pred^A_\le(w)$ and $\partial_{{A'}^\moral}(w) = \partial^A_\le(w)$ the above equality encodes the conditional independence:
\[ \{w \} \Indep_{\Pr_V} A'\sm \{w\} \given \partial_{{A'}^\moral}(w). \]
  This shows \hyl{oLMP}{oLMP}.
\end{proof}
\end{Lem}

\begin{Rem}
\label{rFPwd-rFP}
Let $G=(V,E,H)$ be a \HEDG{} and  $\Pr_V$ a probability distribution on $\Xcal_V$ and $<$ a total order on $V$.
Then for $(G,\Pr_V,<)$ we clearly have the implication:
\[ \xymatrix{ 
    \text{\hyl{rFPwd}{rFPwd}}\ar@{=>}[r] &	\text{\hyl{rFP}{rFP}}, 	
} \]
 which is an equivalence if $\Pr_V$ has a density $p_V$ w.r.t.\ some product measure $\mu_V$.
\end{Rem}

\begin{Thm}
\label{rFPwd-aFP}
Let $G=(V,E,H)$ be a \HEDG{} that has a perfect elimination order $<$ and $\Pr_V$ a probability distribution on $\Xcal_V$ that has a density $p_V$ w.r.t.\ some product measure $\mu_V$.
Then we have the implication:
\[ \xymatrix{ 
    \text{\hyl{rFPwd}{rFPwd} for } (G,\Pr_V,\mu_V, <) \ar@{=>}[r] &	\text{\hyl{aFP}{aFP} for } (G,\Pr_V,\mu_V). 	
} \]
\begin{proof}
Let $A \ins G$ be an ancestral sub-\HEDG{}.
  By \hyl{rFPwd}{rFPwd} we then have the factorization:
\[ p_A(x_A) = \prod_{v \in A} p_A(x_v|x_{\partial^A_\le(v)}).  \]
Since $<$ is a perfect elimination order we have that for every $v \in A$ the set 
\[C:=\{v\} \cup \partial^A_\le(v)=\{v\} \cup \partial_{\Pred_\le^A(v)^\moral }(v)\]
 is a complete subgraph of $\Pred_\le^A(v)^\moral$. 
Indeed by \ref{hedg-anc-marg} $\Pred_\le^A(v)=A \cap \Pred^G_\le(v)$ is an ancestral sub-\HEDG{} of $\Pred^G_\le(v)$ which contains $v$.\\
It is left to show that $C$ is also a complete subgraph of $A^\moral$ because then we can just put $k_{A,C}(x_C):=p_A(x_v|x_{\partial^A_\le(v)})$ to get a desired factorization.\\
We inductively show that $C$ is a complete subgraph of $\Pred^A_\le(w)$ for every $w \in A$ with $v <w$. So let $v=w_0<w_1 < \dots <w_r$ be all elements in  $A$ with $v \le w_i$.
By \ref{marg-mor1} we have the inclusion of undirected subgraphs:
\[ \Pred_\le^A(w_{i-1})^\moral \ins (\Pred_\le^A(w_i)^\moral)^{\marg \sm \{w_i\}}. \]
Since $<$ is a perfect elimination order we have that $\partial_{\Pred_\le^A(w_i)^\moral}(w_i)$ is a complete subgraph of $\Pred_\le^A(w_i)^\moral$. By \ref{compl-subgr-marg} then $C$ is also a complete subgraph of $\Pred_\le^A(w_i)^\moral$ (if it was already in $\Pred_\le^A(w_{i-1})^\moral$). So inductively $C$ becomes a complete subgraph of $A^\moral$.
\end{proof}
\end{Thm}

\begin{Thm}
\label{aFP-dGMP}
Let $G=(V,E,H)$ be a \HEDG{} and  $\Pr_V$ a probability distribution on $\Xcal_V$. Assume that there are %
measures $\mu_v$ on $\Xcal_v$ for every $v \in V$ such that $\Pr_V$ has a density $p_V$ w.r.t.\ the product measure $\mu_V=\otimes_{v\in V} \mu_v$.
We then have the following implication for $(G,\Pr_V)$:
\[ \xymatrix{ 
    \text{\hyl{aFP}{aFP}}\ar@{=>}[r] &	\text{\hyl{auGMP}{auGMP}}. 	
} \]
\begin{proof}
Let $A \subseteq G$ be an ancestral sub-\HEDG{} of $G$.
Let $X,Y,Z$ be pairwise disjoint subsets of $A$ with $X \Indep_{A^\moral} Y \given Z$. We have the factorization:
\[p_A(a) = \prod_{C \in \Ccal_A} k_{A,C}(a_C). \]
Let $X'$ be the set of all $j \in A \sm Z$ such that there is a $Z$-unblocked path in $A^\moral$ to a node in $X$. 
Let $Y'$ be the set of nodes $j \in A \sm Z$ for which every chain in $A^\moral$ from $j$ to $X$ is blocked by $Z$.
  Since $X \Indep_{A^\moral} Y \given Z$ we have $Y \ins Y'$. It now follows that $A^\moral=X' \cup Y' \cup Z$ and $X' \Indep_{A^\moral} Y' \given Z$.
So for every complete subgraph $C \ins A^\moral$ with $C \cap X' \neq \emptyset$ we have $C \cap Y'=\emptyset$ (and the other way around). 
Using this for the complete subgraphs $C \in \Ccal_A$ in $A^\moral$ we get the decomposition:
\[p_A(a) = \prod_{\substack{C \in \Ccal_A,\\ X' \cap C \neq \emptyset}} k_{A,C}(a_C) \cdot \prod_{\substack{C \in \Ccal_A,\\ X' \cap C =\emptyset}}  k_{A,C}(a_C), \]
where the first product is a function in variables from $X' \cup Z$ only and the second a function in variables from $Y' \cup Z$ only.
  Integrating both sides over all variables $A \sm( X \cup Y \cup Z)$ we get that the marginal density w.r.t.\ $X \cup Y \cup Z$ has the following form
\[p(x,y,z)=p_{X\cup Y \cup Z}(x,y,z) = g(x,z) h(y,z), \]
with integrable functions $g,h$. 
So we get for all values $x,y,z$:
\begin{eqnarray*}
&& p(x,z) \cdot p(y,z) \\
&=& g(x,z) \int h(y',z) d\mu_Y(y') \cdot  h(y,z) \int g(x',z) d\mu_X(x') \\
&=& g(x,z)h(y,z) \cdot \int h(y',z) d\mu_Y(y')  \int g(x',z) d\mu_X(x')  \\
&=& p(x,y,z) \cdot p(z).
\end{eqnarray*}
But this implies: $X \Indep_{\Pr_V} Y \given Z$.\\
\end{proof}
\end{Thm}

\begin{Thm}
\label{auPMP-aFP}
Let $G=(V,E,H)$ be a \HEDG{} and  $\Pr_V$ a probability distribution on $\Xcal_V$ that has a strictly positive density $p_V > 0$ w.r.t.\ some product measure $\mu_V=\otimes_{v\in V} \mu_v$ on $\Xcal_V$.
Then the following implication for $(G,\Pr_V)$ holds:
\[ \xymatrix{ 
    \text{\hyl{auPMP}{auPMP}}\ar@{=>}[r] &	\text{\hyl{aFP}{aFP}},	
} \]
making all the implications from \ref{dGMP-auPMP} and \ref{aFP-dGMP} equivalences.
\begin{proof}
The claim directly follows from the famous Hammersley-Clifford theorem (see e.g.\ \cite{Lau98} Thm. 3.9) applied to each $A^\moral$ separately, where $A$ runs through all ancestral sub-\HEDG{}es of $G$. 
Note that with $\Pr_V$ also all marginalizations $\Pr_A$ have strictly positive densities w.r.t.\ the corresponding product measures $\mu_A$
by the following Lemma \ref{positive-marginalisation-density}. 
\end{proof}
\end{Thm}

\begin{Lem}
\label{positive-marginalisation-density}
Let $\Xcal = \Xcal_1 \x \Xcal_2$ be the product space of two measurable spaces. Let $\Pr$ be a probability distribution on $\Xcal$ that has a strictly positive density $p >0$ w.r.t.\ some product measure $\mu = \mu_1 \otimes \mu_2$ on $\Xcal$.
Then also the marginalizations $\Pr_1$ (and $\Pr_2$) of $\Pr$ have strictly positive densities w.r.t.\ the corresponding measure $\mu_1$ (and $\mu_2$, resp.).
\begin{proof}
To see this let $p_1$ be the marginalization of $p$ to $\Xcal_1$ (given by $p_1(x_1)=\int_{\Xcal_2} p(x_1,x_2)\, d\mu_2(x_2)$) and $N_1:=\{x_1\in \Xcal_1| p_1(x_1)=0\}$. For the claim it will be enough to show that $\mu_1(N_1)=0$, because then $p_1$ can be changed to an arbitrary positive value on $N_1$.
First note that:
 \[ \begin{array}{lclcl}
0 &=& \int \I_{N_1}(x_1) \cdot p_1(x_1) \,d\mu_1(x_1) \\
 &=& \int \I_{N_1}(x_1) \cdot \int p(x_1,x_2)\,d\mu_2(x_2) \,d\mu_1(x_1)\\
&=& \int \I_{N_1 \x \Xcal_2}(x) \cdot p(x) \,d\mu(x).
\end{array}\]
So we have $\int g(x) \,d\mu(x) =0 $ with $g(x):=\I_{N_1 \x \Xcal_2}(x) \cdot p(x)$.
Put $S:=\{x \in \Xcal|g(x) >0\}$ and $S_n:=\{x \in \Xcal|g(x) > \frac{1}{n}\}$ for $n \in \N$. We then have $S=\bigcup_n S_n$ and:
\[ \begin{array}{lclclcl}
0 &=& \int g(x) \,d\mu(x) \\
 &=& \int \I_S(x) \cdot g(x) \,d\mu(x) \\
 &=& \sup_{n} \int \I_{S_n}(x) \cdot g(x) \,d\mu(x) \\
 &\stackrel{\forall n}{\ge}& \int \I_{S_n}(x) \cdot g(x) \,d\mu(x) \\
 &\ge& \int \I_{S_n}(x) \cdot \frac{1}{n} \,d\mu(x) \\
 &=& \frac{1}{n} \cdot \mu(S_n) \\& \ge &0.
\end{array}\]
This implies $\mu(S_n) = 0$ for all $n \in \N$ and thus $\mu(S)= \sup_n \mu(S_n) =0$. \\
On the other hand we have:
\[ \begin{array}{lcl}
S &=& \{x \in \Xcal | \I_{N_1 \x \Xcal_2}(x) \cdot p(x) > 0\} \\
&\stackrel{p(x) >0}{=} & \{x \in \Xcal | \I_{N_1 \x \Xcal_2}(x) > 0\} \\
&=& N_1 \x \Xcal_2.
\end{array}\]
Together we get:
\[ 0 = \mu(S) =\mu(N_1 \x \Xcal_2) = \mu_1(N_1) \cdot \mu_2(\Xcal_2).\]
Since $\mu_2(\Xcal_2) \neq 0$ (otherwise $0=\Pr_2(\Xcal_2)=1$) we get $\mu_1(N_1) =0$.
\end{proof}
\end{Lem}

\begin{Cor}
\label{dGMP-aFP}
Let $G=(V,E,H)$ be a \HEDG{} and  $\Pr_V$ a probability distribution on $\Xcal_V$ that has a density $p_V$ w.r.t.\ some product measure $\mu_V=\otimes_{v\in V} \mu_v$ on $\Xcal_V$. Furthermore, assume one (or both) of the following:
\begin{enumerate}
\item $p_V$ is strictly positive ($p_V > 0$), or
\item $G$ has a perfect elimination order.
\end{enumerate}
  Then for $(G,\Pr_V)$ the directed global Markov property (\hyl{dGMP}{dGMP}) is equivalent to the ancestral factorization property (\hyl{aFP}{aFP}):
\[ \xymatrix{ 
    \text{\hyl{dGMP}{dGMP}}\ar@{<=>}[r] &	\text{\hyl{aFP}{aFP}}.	
} \]
\begin{proof}
  ``\hyl{aFP}{aFP} $\implies$ \hyl{dGMP}{dGMP}'' was proven in \ref{aFP-dGMP} and \ref{dGMP-auPMP}. For ``\hyl{dGMP}{dGMP} $\implies$ \hyl{aFP}{aFP}'' note:
If $p_V > 0$ the claim follows from the following implication chain (\ref{dGMP-auPMP}, \ref{auPMP-aFP}):
  \[ \text{\hyl{dGMP}{dGMP}} \;\stackrel{\ref{dGMP-auPMP}}{\Longleftrightarrow}\; \text{\hyl{auGMP}{auGMP}} \;\stackrel{\ref{dGMP-auPMP}}{\implies}\; \text{\hyl{auLMP}{auLMP}} 
  \;\stackrel{\ref{dGMP-auPMP}}{\implies}\; \text{\hyl{auPMP}{auPMP}} \;\stackrel{\ref{auPMP-aFP},\, p_V>0}{\implies}\; \text{\hyl{aFP}{aFP}}. \]
If $G$ has a perfect elimination order the claim follows from the following implication chain (\ref{dGMP-oLMP}, \ref{oLMP-rFP}, \ref{rFPwd-rFP}, \ref{rFPwd-aFP}):
  \[ \text{\hyl{dGMP}{dGMP}} \;\stackrel{\ref{dGMP-oLMP}}{\implies}\; \text{\hyl{oLMP}{oLMP}} \;\stackrel{\ref{oLMP-rFP}}{\Longleftrightarrow}\; \text{rFP} 
  \;\stackrel{\ref{rFPwd-rFP},\text{ density}}{\Longleftrightarrow}\; \text{\hyl{rFPwd}{rFPwd}} \;\stackrel{\ref{rFPwd-aFP}, \text{ p.e.ord.}}{\implies}\; \text{\hyl{aFP}{aFP}}. \]
\end{proof}
\end{Cor}

In the next example we will see that without the positivity of the density or perfect elimination order we might not get the ancestral factorization property (\hyl{aFP}{aFP}) back. The example was taken from \cite{Lau98} Ex. 3.10, where it was stated for undirected graphs.

\begin{Eg}[\hyl{dGMP}{dGMP} does not imply \hyl{aFP}{aFP} without strictly positive density or perfect elimination order]
\label{no-p-dGMP-aFP}
Let $X_1,\dots,X_4$ be four discrete $\{0,1\}$-valued random variables that have a joint probability of $p=\frac{1}{8}$ on each of the following points:
\[\begin{array}{cccc}
(0,0,0,0), & (1,0,0,0), & (1,1,0,0), & (1,1,1,0), \\
(0,0,0,1), & (0,0,1,1), & (0,1,1,1), & (1,1,1,1). 
\end{array}\]
For all other possible points we have $p=0$.
Straightforward calculation shows that we have the following independencies:
\[ X_1 \Indep_\Pr X_3 \given \;(X_2,X_4), \qquad \text{and} \qquad X_2 \Indep_\Pr X_4 \given \;(X_1,X_3).\]
So it satisfies the directed global Markov property (\hyl{dGMP}{dGMP}) according to the directed circle with four nodes, figure \ref{fig:directed-cycler}.
\begin{figure}[!htb] 
\centering
\begin{tikzpicture}[scale=.9, transform shape]
\tikzstyle{every node} = [draw,shape=circle]
\node (v1) at (0,0) {$X_1$};
\node (v2) at (2,0) {$X_2$};
\node (v3) at (2,2) {$X_3$};
\node (v4) at (0,2) {$X_4$};
\foreach \from/\to in {v2/v1, v3/v2, v4/v3, v1/v4}
\draw[-{Latex[length=3mm, width=2mm]}] (\from) -- (\to);
\end{tikzpicture}
  \caption{Directed cycle with four nodes corresponding to Example~\ref{no-p-dGMP-aFP}.}
\label{fig:directed-cycler}
\end{figure}
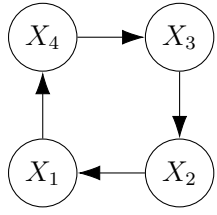
But there is no ancestral factorization of $\Pr$ according to $G$. Otherwise we had the equations:
\[\begin{array}{ccccccc}
0 &\neq  & p(0,0,0,0) & = & \phi_{12}(0,0)\cdot\phi_{23}(0,0)\cdot\phi_{34}(0,0)\cdot\phi_{41}(0,0),\\
0 &\neq  & p(0,0,1,1) & = & \phi_{12}(0,0)\cdot\phi_{23}(0,1)\cdot\phi_{34}(1,1)\cdot\phi_{41}(1,0),\\
0 &\neq  & p(1,1,1,0) & = & \phi_{12}(1,1)\cdot\phi_{23}(1,1)\cdot\phi_{34}(1,0)\cdot\phi_{41}(0,1),\\
0 &=     & p(0,0,1,0) & = & \phi_{12}(0,0)\cdot\phi_{23}(0,1)\cdot\phi_{34}(1,0)\cdot\phi_{41}(0,0),\\
\end{array}\]
  showing that $\phi_{12}(0,0),\phi_{23}(0,1),\phi_{34}(1,0),\phi_{41}(0,0) \neq0$, but their product equals $0$, which is a contradiction. Also, $G$ has no perfect elimination order. Indeed, whatever order one chooses, for $v$ the last element in that order, we have that the only ancestral subgraph of $\Pred^G_\le(v)$ is just $G$ itself, and $\partial_{A^\moral}(v) \cup \{v\}$ is not a complete subgraph of $A^\moral$.
\end{Eg}

\begin{Thm}[\hyl{aFP}{aFP} stable under marginalization]
\label{aFP-marg}
Let $G=(V,E,H)$ be a \HEDG{} and $\Pr_V$ a probability distribution on $\Xcal_V$ and $W \ins V$ a subset.
  If $(G,\Pr_V)$ has the  ancestral factorization property (\hyl{aFP}{aFP}) then also the marginalization $(G^{\marg\sm W}, \Pr_{V\sm W})$ has the corresponding ancestral factorization property (\hyl{aFP}{aFP}):
\[ \xymatrix{ 
    \text{\hyl{aFP}{aFP} for } (G,\Pr_V, \otimes_{v \in V}\mu_v) \ar@{=>}[r] &	\text{\hyl{aFP}{aFP} for } (G^{\marg\sm W}, \Pr_{V\sm W}, \otimes_{v \in V\sm W}\mu_v ).
} \]
\begin{proof}
By induction we can assume $W=\{w\}$ has only one element. 
Let $A' \ins  G^{\marg \sm W}$ be ancestral. The ancestral sub-\HEDG{} $A := \Anc^G(A')$ of $G$ will then by \ref{hedg-anc-marg} marginalize to $A'$ in $G^{\marg\sm W}$.
  By \hyl{aFP}{aFP} we have the factorization:
\[ p_A(x_A) = \prod_{C \in \Ccal_A} k_{A,C}(x_C), \]
for a set of complete subgraphs $\Ccal_A$ of $A^\moral$. 
We distinguish two cases: $w \in A$ and $w \notin A$.
If $w \notin A$ then $A'=A$ is already an ancestral sub-\HEDG{} of $G$. Note that ancestral sub-\HEDG{}es carry the marginalized \HEDG{} structure by \ref{anc-sub=mar}. So $p_{A'}=p_A$ has the same factorization as $p_A$.\\
If $w \in A$ then the marginalization $p_{A'}$ of $p_A$ is:
\[ p_{A'}(x_{A'}) = \prod_{\substack{C \in \Ccal_A\\ w \notin C}} k_{A,C}(x_C) \cdot \int_{\Xcal_w} \prod_{\substack{C \in \Ccal_A\\ w \in C}} k_{A,C}(x_C)\,d\mu_w(x_w).   \]
Since $w \in A$ we can assume that $\Ch^A(w) \sm \{w\} \neq \emptyset$. Otherwise $A \sm \{w\}$ would already be an ancestral sub-\HEDG{} of $G$ marginalizing to $A'$. By \ref{marg-mor2} we then get:
\[  (A')^\moral = (A^\moral)^{\marg\sm\{w\}}. \]
So if $w \notin C$ then $C$ stays a complete subgraph of $(A^\moral)^{\marg\sm\{w\}} = (A')^\moral$. So for these $C$ we can put
$k_{A',C}(x_C):= k_{A,C}(x_C)$.\\
We are left to argue about the integral on the right of the factorization. Note that it is a function in the variables of $D:=\bigcup_{\substack{C \in \Ccal_A\\ w \in C}} C \sm \{w\}$.
Since these $C \ins A^\moral$ are complete and contain $w$ we have $D \ins \partial_{A^\moral}(w)$. By \ref{compl-subgr-marg} we get that 
$\partial_{A^\moral}(w)$ and thus $D$ becomes a complete subgraph of $(A^\moral)^{\marg\sm\{w\}} = (A')^\moral$. So we can put
$$k_{A',D}(x_D):= \int_{\Xcal_w} \prod_{\substack{C \in \Ccal_A\\ w \in C}} k_{A,C}(x_C)\,d\mu_w(x_w).$$
  This shows the required factorization and thus \hyl{aFP}{aFP} for the marginalization.
\end{proof}
\end{Thm}

\begin{Lem}[\hyl{mFP}{mFP} stable under ancestral subsets]
\label{mFP-anc}
Let $G=(V,E,H)$ be a \HEDG{} and $\Pr_V$ a probability distribution on $\Xcal_V$ and $A \ins G$ an ancestral sub-\HEDG{}.
Then we have the implication:
\[ \xymatrix{ 
  \text{\hyl{mFP}{mFP} for } (G,\Pr_V, \otimes_{v \in V}\mu_v) \ar@{=>}[r] &	\text{\hyl{mFP}{mFP} for } (A, \Pr_A, \otimes_{v \in A}\mu_v ).		
} \]
\begin{proof} %
We first show the little stronger claim that the form of the joint density $p_{V^\aug}$ as mentioned in \hyl{mFP}{mFP} is stable under ancestral subgraphs $A'$ of $G^\aug$ (we can then just take $A':=\Anc^{G^\aug}(A)$ later on). 
So let $A' \ins G^\aug$ be ancestral and $S_1,\dots,S_r \in \Scal(G^\aug)$ all strongly connected components of $G^\aug$ ordered by a topological order of $\Scal(G^\aug)$ such that $S_1,\dots, S_k$ are the strongly connected components of $A'$ and $S_{k+1},\dots,S_r$ are the rest.
Since $A'$ is ancestral such a topological order always exists.
  By \hyl{mFP}{mFP} we have the joint density:
\[ p_{V^\aug}(x,e) =   \prod_{S \in \Scal(G^\aug)} \prod_{C \in \Ccal_S} k_{S,C}(x_C,e_C), \]
and for all $m=1,\dots,r$:
\[\begin{array}{lclll}
\int \prod_{C \in \Ccal_{S_m}} k_{{S_m},C}(x_C,e_C) \, d\mu_{S_m}(x_{S_m}) &=& 1, & \text{for} & {S_m} \in \Scal(G),\\
\int \prod_{C \in \Ccal_{S_m}} k_{{S_m},C}(e_F)\, d\nu_F(e_F) &=& 1,& \text{for} & {S_m}=\{F\}, F \in \tilde{H}.\\
\end{array}\]
  So integrating over $x_{S_m}$ and $e_F$ for $F \in V^\aug \sm (A'\cup V)$ for decreasing $m=r, r-1, \dots k+1$ then gives the claim for ancestral $A' \ins G^\aug$:
\begin{eqnarray*}
&&p_{A'}(x_{A'},e_{A'}) \\
&=&   \prod_{m=1}^k \prod_{C \in \Ccal_{S_m}} k_{S_m,C}(x_C,e_C) \cdots \\
&&\cdots \int \int\prod_{m=r}^{k+1} \prod_{C \in \Ccal_{S_m}} k_{S_m,C}(x_C,e_C)\, 
d\mu(x_{V\sm A'})\, d\nu(e_{V^\aug\sm (A'\cup V)}),\\
 &=&\prod_{m=1}^k \prod_{C \in \Ccal_{S_m}} k_{S_m,C}(x_C,e_C)\,\cdot\,1.
\end{eqnarray*}
For the actual claim for $A \ins G$ ancestral we need to construct a factorization for the marginal density over $A^\aug$ out of the above one for $A':=\Anc^{G^\aug}(A)$. Fortunately, these two directed graphs might only differ by the latent nodes outside of $A$ corresponding to hyperedges $F \in \tilde H$ with $F \sm A \neq \emptyset$ and $F \cap A \neq \emptyset$. Since $F \cap A \ins F'$ for an $F' \in  \widetilde{H(A)}$ we can subsume such $F$ in an arbitrary way to some $F'$ and gather the corresponding factors $k_{\{F\},C}(e_C)$ into one product $k_{\{F'\},C}(e_C)$. This then shows the claim.
\end{proof}
\end{Lem}

\begin{Lem}
\label{mFP-aFP}
  For a \HEDG{} $G=(V,E,H)$ and a probability distribution $\Pr_V$ on $\Xcal_V$ with a density $p_V$ w.r.t.\ a product measure $\mu_V$ of %
	measures $\mu_v$, $v \in V$, we have the implication:
\[ \xymatrix{ 
      \text{\hyl{mFP}{mFP}}\ar@{=>}[r]& \text{\hyl{aFP}{aFP}}.
} \]
\begin{proof}
Let $A \ins G$ be ancestral. 
  By \ref{mFP-anc} \hyl{mFP}{mFP} is stable under ancestral subsets. So $(A,\Pr_A)$ satisfies \hyl{mFP}{mFP} as well.
The claim then follows by gathering the terms unter the integral exactly like in the proof of \ref{aFP-marg}
by marginalizing out $W=A^\aug \sm A$. 
\end{proof}
\end{Lem}

The next result is \cite{Koster96} Thm. 3.4 applied to our setting and gives another connection between the form of \hyl{mFP}{mFP} and \hyl{aFP}{aFP}. 

\begin{Lem}
\label{aFP-mFP}
Let  $G=(V,E,H_1)$ be \HEDG{} with trivial $H_1$ and $\Pr_V$ a probability distribution on $\Xcal_V$.
  Assume that $(G,\Pr_V,\mu_V)$ satisfies the ancestral factorization property (\hyl{aFP}{aFP}) with a strictly positive density $p_V>0$ w.r.t.\ some product measure $\mu_V = \otimes_{v\in V} \mu_v$. 
Then for every strongly connected component $S \in \Scal(G)$ there exists a set $\Ccal_S$ of complete subgraphs  $C \ins \Anc^G(S)^\moral$ with $C \ins S \cup \Pa^G(S)$ together with integrable functions:
   \[ k_{S,C}:\; \Xcal_C \to \R_{\ge 0},  \]
	 such that the joint density $p_V$ w.r.t.\ $\mu_V$ is given by:
\[ p_V(x) = \prod_{S \in \Scal(G)} \prod_{C \in \Ccal_S} k_{S,C}(x_C),   \]
and such that for every $S \in \Scal(G)$ and almost all $x_{\Pa^G(S)\sm S}$ we have:
\[ \int_{\Xcal_S} \prod_{C \in \Ccal_S} k_{S,C}(x_C)\,d\mu_S(x_S) = 1.   \]
\begin{proof}
  By \ref{aFP-dGMP} the tuple $(G,\Pr_V)$ satisfies the directed global Markov property (\hyl{dGMP}{dGMP}). From this we get for every $S \in \Scal(G)$:
\[  S \Indep_{\Pr_V}  \NonDesc^G(S) \given \Pa^G(S) \sm S,  \]
and, in particular:
\[  S \Indep_{\Pr_V}  \Anc^G(S)\sm S \given \Pa^G(S) \sm S.  \]
So inductively---using a topological order for the mDAG of strongly connected components $\Scal(G)$---we get that the joint density $p=p_V$ factorizes as:
\[ p(x) = \prod_{S\in \Scal(G)} p(x_S | x_{\Pa^G(S)\sm S}).  \]
Furthermore, we have:
\[ p(x_S | x_{\Pa^G(S)\sm S}) = p(x_S | x_{\Anc^G(S)\sm S}) 
 \stackrel{\ref{positive-marginalisation-density},\;p>0}{=} \frac{p(x_{\Anc^G(S)})}{p(x_{\Anc^G(S)\sm S}) }.\]  
  By the ancestral factorization property (\hyl{aFP}{aFP}) applied to the two ancestral subsets $A:=\Anc^G(S)$ and $A':=\Anc^G(S)\sm S$ we get:
\[ p_A(x_A) = \prod_{\substack{C \ins A^\moral\\\text{complete}}} k_{A,C}(x_C),\qquad \text{ and } \qquad 
   p_{A'}(x_{A'}) = \prod_{\substack{D \ins (A')^\moral\\\text{complete}}} k_{A',D}(x_D).  \]
Since $(A')^\moral$ is a subgraph of $A^\moral$ every $D$ is also complete in $A^\moral$ and we can gather some $D$'s with $D \ins C$ and assign them to the complete $C \ins A^\moral$.
So we can write:
\[\begin{array}{rcl}
 p(x_S | x_{\Pa^G(S)\sm S}) &=&\frac{p_A(x_A)}{p_{A'}(x_{A'})} \\%
&=& \prod_{\substack{C \ins A^\moral\\\text{complete}}} \frac{ k_{A,C}(x_C)}{ k_{A',C}(x_C) }\\%
&=& \prod_{\substack{C \ins A^\moral\\\text{complete}}} k_{S,C}(x_C),
\end{array}\]
for the functions $k_{S,C}(x_C):= \frac{ k_{A,C}(x_C)}{ k_{A',C}(x_C) }$. %
Since the left hand side only depends on the variables $S \cup \Pa^G(S)$ we can (by plugging in fixed values for all other variables) restrict to $C \cap (S \cup \Pa^G(S) )$ on the right hand side as well. These are still complete subgraphs of $A^\moral$ and lie in $S \cup \Pa^G(S)$.
Since $\int p(x_S | x_{\Pa^G(S)\sm S})\, d\mu_S(x_S) = 1$ we also get the remaining claim.
\end{proof}
\end{Lem}

\begin{Rem}
The procedure from \ref{aFP-mFP} also throws light on the difference between \hyl{dGMP}{dGMP} and \hyl{gdGMP}{gdGMP} based on  d-separation and $\sigma$-separation, resp..
For simplicity again assume that $G=(V,E,H_1)$ has a trivial $H_1$ and $\Pr_V$ has a strictly positive density $p>0$ w.r.t.\ some product measure $\mu_V$. Then by \ref{auPMP-aFP}, \ref{dGMP-auPMP}, \ref{aFP-dGMP} and \ref{dGMP-gdGMP}, \ref{s-sep-d-sep-eg} we have the implications:
\[ \xymatrix{ 
      \text{\hyl{aFP}{aFP} for } (G,\Pr_V) \ar@{<=>}[r]& \text{\hyl{dGMP}{dGMP} for } (G,\Pr_V) \ar@{=>}[d]\\
      \text{\hyl{aFP}{aFP} for } (G',\Pr_V) \ar@{<=>}[r]& \text{\hyl{dGMP}{dGMP} for } (G',\Pr_V) \ar@{<=>}[r]& \text{\hyl{gdGMP}{gdGMP} for } (G,\Pr_V), \\
} \]
  where we used the graph with fully connected strongly connected components $G'=(V,E^\acy \cup E^\mathrm{sc},H_1)$ (see \ref{s-sep-d-sep-eg})  instead of $G^\acy$ to not introduce hyperedges for simplicity.
Like in \ref{aFP-mFP} now both d-separation and $\sigma$-separation give us for every strongly connected component $S \in \Scal(G)=\Scal(G')$:
\[  S \Indep  \NonDesc^G(S) \given \Pa^G(S) \sm S, \]
  from which the following factorization follows (assuming at least \hyl{gdGMP}{gdGMP}):
\[ p(x) = \prod_{S\in \Scal(G)} p(x_S | x_{\Pa^G(S)\sm S}).  \]
  For \hyl{gdGMP}{gdGMP} the factorization now stops here because for every ancestral $A \ins G'$ containing $S$ and thus $D:=S \cup \Pa^G(S)=S \cup \Pa^{G'}(S)$ the nodes of $D$ form a complete subgraph of $A^\moral$ and $p(x_S | x_{\Pa^G(S)\sm S})$ cannot be factorized any further.
  If we in contrast assume the stronger \hyl{dGMP}{dGMP} we get the further factorization from \ref{aFP-mFP} :
\[ p(x_S | x_{\Pa^G(S)\sm S}) = \prod_{\substack{ C \ins S \cup \Pa^G(S)\\ C \ins \Anc^G(S)^\moral\\\text{complete} }} k_{S,C}(x_C). \]
  This shows the main difference between \hyl{gdGMP}{gdGMP} and \hyl{dGMP}{dGMP} on a factorization level.
\end{Rem}

\subsection{Marginal Markov Properties for \HEDG{}es}

In this subsection we define ``marginal Markov properties'' for \HEDG{}es $G$ and probability distributions $\Pr_V$ on $\Xcal_V$.
We will see that the marginal Markov properties will be more closely related to the structural equations properties (\hyl{SEP}{SEP}) later on (see \ref{csSEP-smgdGMP}) than the ordinary Markov properties.
Roughly speaking, ``marginal Markov property'' will mean that $\Pr_V$ appears to be the marginal of a distribution $\Pr_{V^\aug}$ on $\Xcal_{V^\aug}$  that satisfies the corresponding ``ordinary Markov property'' on the augmented graph $G^\aug$ of $G$.
This will reflect the idea that the latent space  of all non-observed variables  together with the observed space should follow the same rules as if everything were observed, and the observed variables are just a projection of the joint (observed and unobserved) world. 
To justify the use of hyperedges (represented as nodes in the augmented graph) as a summary of the whole latent space we will
show (see \ref{mMP-lift}, \ref{mdGMP-lift}) that adding more structure to the latent space will not change the expressiveness of the model class for the marginal directed global Markov property (\hyl{mMP}{mdGMP}). So imposing the mentioned conditions only on the augmented graph will be enough.
This result generalizes the corresponding result of \cite{Eva15}, where it was shown for mDAGs, to the case of general \HEDG{}es.
We will also give a general construction (see \ref{mMP-lift}) to extend this result to other (marginal) Markov properties, but we restrict ourselves to give another proof only for a certain structural equations property (\hyl{lsSEP}{lsSEP}) later on (see \ref{lsSEP-lift}).

In contrast to this we will see that as long as the ordinary Markov property is stable under marginalization (e.g.\ \hyl{dGMP}{dGMP}, \hyl{gdGMP}{gdGMP}, \hyl{aFP}{aFP} by \ref{dGMP-stable-marg}, \ref{aFP-marg}) the marginal Markov property implies the ordinary one and imposes restrictions on the set of meaningful probability distributions $\Pr_V$ related to $G$. In other words, the model class of all distributions satisfying a marginal Markov property i.g.\ is smaller than the class of the corresponding ordinary Markov property (as it was already shown in \cite{Eva15}). This makes the introduction of marginal Markov properties meaningful.

\begin{Def}[Marginal Markov properties]
\label{def-marg}
Let $G=(V,E,H)$ be a \HEDG{} and $\Pr_V$ a probability distribution on $\Xcal_V$.
Let \emph{MP} be one of the Markov properties from before (e.g.\ \hyl{dGMP}{dGMP}, \hyl{dLMP}{dLMP}, etc.).
We can define a corresponding \emph{\hyt{mMP}{marginal Markov property} mMP} (e.g.\ mdGMP, mdLMP, etc.) for $(G,\Pr_V)$ as follows:
\begin{enumerate}
	   \item For every $F \in \tilde{H}$ there exists a standard Borel space $\Ecal_F$.
		 \item There exists a probability distribution $\Pr_{V^\aug}$ on 
		       $\Xcal_{V^\aug}:=\prod_{v \in V}\Xcal_v \x \prod_{F \in \tilde{H}} \Ecal_F$, such that:
			   \item The marginal distribution of $\Pr_{V^\aug}$ on $\Xcal_V$ equals $\Pr_V$, i.e.\ we have:
				  \[ \pi_{V,*}\Pr_{V^\aug} = \Pr_V,\]
				 where $\pi_V: \Xcal_{V^\aug} \to \Xcal_V$ is the canonical projection, and:
       \item $(G^\aug,\Pr_{V^\aug})$ satisfies \emph{MP} (e.g.\ \hyl{dGMP}{dGMP}, \hyl{dLMP}{dLMP}, etc.), where 
				  $G^\aug$ is the augmented graph of $G$ viewed as the \HEDG{} $(V^\aug,E^\aug,H^\aug_1)$, where $H^\aug_1=\{\{w\}|w \in V^\aug\}$ is trivial.
    \end{enumerate}
\end{Def}

\begin{Rem}
\label{mdGMP-dGMP}
Let $G=(V,E,H)$ be a \HEDG{} and $\Pr_V$ a probability distribution on $\Xcal_V$.
\begin{enumerate}
\item 
  If a Markov property is \emph{stable under marginalization} (e.g.\ \hyl{dGMP}{dGMP}, \hyl{gdGMP}{gdGMP} by \ref{dGMP-stable-marg} and \hyl{aFP}{aFP} by \ref{aFP-marg}) then the marginalized Markov property implies the ordinary Markov property (by marginalizing out $V^\aug \sm V$), for instance we have:
    \[ \text{\hyl{mMP}{mdGMP} for } (G,\Pr_V) \qquad \implies \quad \text{\hyl{dGMP}{dGMP} for } (G,\Pr_V), \]
    \[ \text{\hyl{mMP}{mgdGMP} for } (G,\Pr_V) \qquad \implies \quad \text{\hyl{gdGMP}{gdGMP} for } (G,\Pr_V), \]
    \[ \text{\hyl{mMP}{maFP} for } (G,\Pr_V) \qquad \implies \quad \text{\hyl{aFP}{aFP} for } (G,\Pr_V). \]
  \item For \hyl{mMP}{mgdGMP} also see \ref{smgdGMP-mgdGMP}. 
\item Most of the Markov properties are not (or only under additional assumptions) stable under marginalization (e.g.\ \hyl{dLMP}{dLMP}, \hyl{auPMP}{auPMP}, etc.).
\item If for a class of \HEDG{}es $G=(V,E,H)$ with at least trivial $H$ we have the implication of ordinary Markov properties:
\[ \text{ordinary Markov property 1} \qquad \implies \qquad  \text{ordinary Markov property 2},  \]
then we also have the implication:
\[ \text{marginal Markov property 1} \qquad \implies \qquad  \text{marginal Markov property 2}  \]
for the corresponding marginalized class of \HEDG{}es with non-trivial $H$.
For example, for DAGs we have the implication:
    \[ \text{\hyl{dLMP}{dLMP}} \qquad \implies \qquad  \text{\hyl{dGMP}{dGMP}}.  \]
So for mDAGs we get the implication:
    \[ \text{\hyl{mMP}{mdLMP}} \qquad \implies \qquad  \text{\hyl{mMP}{mdGMP}}.  \]
\end{enumerate}
\end{Rem}

\begin{Eg}
\label{eg-dGMP-mdGMP}
Figure \ref{fig:eg-dGMP-mdGMP} shows an example from \cite{Eva15} Example 9 and \cite{Fritz12}, resp.. It shows that the ordinary directed global Markov property (\hyl{dGMP}{dGMP}) does not imply the marginal directed global Markov property (\hyl{mMP}{mdGMP}).
\begin{figure}[!htb]
\centering
\begin{tikzpicture}[scale=.7, transform shape]
\tikzstyle{every node} = [draw,shape=circle]
\node (v1) at (0,-1) {$X_1$};
\node (v2) at (-2,-4) {$X_2$};
\node (v3) at (2,-4) {$X_3$};
\node[fill,circle,red,inner sep=0pt,minimum size=5pt] (v0) at (0,-3) {};
\draw[-{Latex[length=3mm, width=2mm]}, red] (v0) to (v1);
\draw[-{Latex[length=3mm, width=2mm]}, red] (v0) to (v2);
\draw[-{Latex[length=3mm, width=2mm]}, red] (v0) to (v3);
\end{tikzpicture}
\hspace{1cm}
\begin{tikzpicture}[scale=.7, transform shape]
\tikzstyle{every node} = [draw,shape=circle]
\node (v1) at (0,-1) {$X_1$};
\node (v2) at (-2,-4) {$X_2$};
\node (v3) at (2,-4) {$X_3$};
\draw[{Latex[length=3mm, width=2mm]}-{Latex[length=3mm, width=2mm]}, red] (v1) to node[fill,circle,,inner sep=0pt,minimum size=5pt] {} (v2);
\draw[{Latex[length=3mm, width=2mm]}-{Latex[length=3mm, width=2mm]}, red] (v1) to node[fill,circle,,inner sep=0pt,minimum size=5pt] {} (v3);
\draw[{Latex[length=3mm, width=2mm]}-{Latex[length=3mm, width=2mm]}, red] (v3) to node[fill,circle,,inner sep=0pt,minimum size=5pt] {} (v2);
\end{tikzpicture}
\caption{A \HEDG{} $G$ on the left and its induced DMG $G_2$ on the right. Three identical copies $X_1=X_2=X_3$ of one fair coin flip $E_{\{1,2,3\}}$ (i.e.\ with $\Pr(E_{\{1,2,3\}}=0)=\Pr(E_{\{1,2,3\}}=1)= \frac{1}{2}$) will satisfy the marginal directed global Markov property w.r.t.\ $G$.
The marginal directed global Markov property w.r.t.\ $G_2$ on the other hand will not hold because we can not find three \emph{independent} variables $E_{\{1,2\}},E_{\{2,3\}},E_{\{1,3\}}$ which induce the distribution of the highly dependent variables $X_1,X_2,X_3$. In contrast to that, the ordinary directed global Markov property holds for both $G$ and $G_2$.}
 \label{fig:eg-dGMP-mdGMP}
\end{figure}
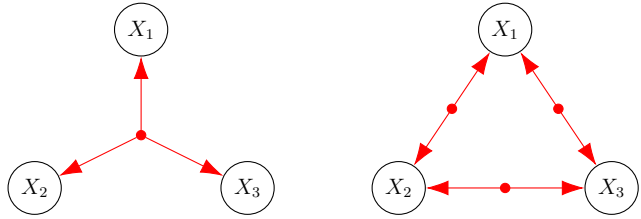

\end{Eg}

As a generalization of the results of \cite{Eva15} we will introduce the following notion to justify the use of hyperedges as summaries of latent variables.

\begin{Rem}[Lifting of marginal Markov properties]
\label{mMP-lift}
To justify the use of hyperedges in a \HEDG{} as a summary of all latent variables w.r.t.\ a marginal Markov property (\hyl{mMP}{mMP}) one does not only need that the ordinary MP is \emph{stable under marginalizations} but also that \hyl{mMP}{mMP} satisfies the following \emph{lifting property}:
\begin{enumerate}
\item[] 
For any \HEDG{} $G=(V,E,H)$, any subset $W \ins V$, $U:=V\sm W$, and any probability distribution $\Pr_U'$ on a product space $\Xcal_U:=\prod_{u \in U} \Xcal_u$ of standard Borel spaces $\Xcal_u$, $u \in U$, such that $(G^{\marg(U)},\Pr_U')$ satisfies mMP there exist standard Borel spaces $\Xcal_w$, $w \in W$, and a probability distribution $\Pr_V$ on $\Xcal_V:=\prod_{ v \in V} \Xcal_v$ such that $(G,\Pr_V)$ satisfies mMP as well and the marginal distribution is given by: $\Pr_{V\sm W} = \Pr_U'$.
\end{enumerate}
For proving such lifting property one can usually do induction on $\#W$ and assume $W=\{w\}$. Identify $U$ with $G^{\marg(U)}=(V\sm \{w\},E^\marg,H^\marg)$. One then first needs to define $\Ecal_F$ for $F \in \tilde{H}$ and $\Xcal_w$ out of $\Ecal'_{F'}$ for $F' \in \tilde{H}^\marg$.
We will sketch a standard strategy for this. 
Since every $F' \in \tilde H^\marg$ is of the form $F'=\hat F^\marg$ for an $\hat F \in \tilde H$ we can fix a map:
\[ \psi: \tilde H^\marg \to \tilde H\qquad \text{with}\qquad \psi(F')^\marg = F'.\]
Note that for $F \in H$ we have defined (see \ref{HEDG-marg-def}):
\[ F^\marg = \left\{ 
\begin{array}{lcl}
  F & \text{ if } & w \notin F,\\
	(F \cup \Ch^G(w))\sm \{w\} &\text{ if } & w \in F.
\end{array}\right. \]
It is clear that $\psi$ is injective. Let $\Psi \ins \tilde H$ be the image of $\psi$ in $\tilde H$, i.e.:
\[\Psi := \{ F \in \tilde H \,|\, \exists F' \in \tilde H^\marg \text{ with } F=\psi(F')\}.\]
We then define $\Ecal_F$ for $F \in \tilde H$ as follows:
\[ \Ecal_F := \left\{
\begin{array}{lll}
\{*\} & \text{ if }& F \notin \Psi,\\
\Ecal'_{\psi^{-1}(F)} & \text{ if } & F \in \Psi,
\end{array}
\right.  \]
where $\{*\}$ denotes a space with one point only. 
For the corresponding random variables $E_F$, $F \in \tilde H$, we similarly define:
\[ E_F := \left\{
\begin{array}{lll}
* & \text{ if }& F \notin \Psi,\\
E'_{\psi^{-1}(F)} & \text{ if } & F \in \Psi,
\end{array}
\right.  \]
where $*$ means the constant map.
We can then also define:
\[ \begin{array}{rcl} 
\Xcal_w &:=& \prod_{v \in \Pa^G(w)\sm \{w\}} \Xcal_v \x \prod_{\substack{F \in \tilde H \\ w \in F}}\Ecal_F,\\
X_w &:=& ((X_v)_{v \in \Pa^G(w)\sm \{w\}}, (E_F)_{\substack{F \in \tilde H \\ w \in F}}).
\end{array}  \]
Note that for the marginalization from $G$ to $G^\marg$ we need to assign each variable $E_F$, $F \in \tilde H$, to a maximal hyperedge 
$F' \in \tilde H^\marg$ with $F' \supseteq F^\marg$, i.e.\ we need to fix a map:
\[ \varphi: \tilde H \to \tilde H^\marg\qquad \text{with}\qquad \varphi(F) \supseteq F^\marg.\]
For such $\varphi$ it is then clear that we have: $\varphi(\psi(F'))=F'$ for every $F' \in \tilde H^\marg$. 
Marginalizing along $\varphi$ gives for every $F' \in \tilde H^\marg$ the isomorphism:
\[  \Ecal'_{F'} = \Ecal_{\psi(F')} \stackrel{\x \{*\}}{\cong} \prod_{\substack{F \in \tilde H \\ \varphi(F) = F'}}\Ecal_F.   \]
So we can identify $\Ecal'_{F'}$ with the right hand side and similarly for $E_{F'}'$.
With these definitions and identifications one can often prove the mMP for $(G,\Pr_V)$ (see \ref{mdGMP-lift}, \ref{lsSEP-lift}).
Note that with these definitions of random variables we clearly have $\Pr_{V^\aug} = \Pr_{U^\aug}'$ and thus $\Pr_{V \sm W} = \Pr_U'$.
\end{Rem}

The following result will generalize results of \cite{Eva15} for mDAGs to the case of all \HEDG{}es.
There it was shown that an unknown and unobserved latent space can be summarized by hyperedges without changing the expressiveness of the underlying marginal Markov model (class). Here we will give a version for the marginal directed global Markov property (\hyl{mMP}{mdGMP}) for \HEDG{}es. A similar result using a structural equations property (\hyl{lsSEP}{lsSEP}) for \HEDG{}es is shown in \ref{lsSEP-lift}.

\begin{Thm}
\label{mdGMP-lift}
The marginal directed global Markov property (\hyl{mMP}{mdGMP}) has the lifting property (\ref{mMP-lift}).
\begin{proof}
Let the definitions and identifications be like in \ref{mMP-lift} with $W=\{w\}$ and $U=G^{\marg\sm W}$.
It is enough to show that $(G^\aug,\Pr_{V^\aug})$ (as defined in \ref{mMP-lift}) satisfies the \hyl{dGMP}{dGMP}.
Let $X,Y,Z \ins V^\aug$ subsets and assume:
\[ X \Indep_{G^\aug}^d Y \given Z.\]
Since every node from $\tilde H \sm \Psi$ is a non-collider on any path the above implies:
\[ X \Indep_{G^\aug}^d Y \given Z \cup (\tilde H \sm \Psi).\]
Furthermore, since every node from $\tilde H \sm \Psi$ corresponds to a constant random variable we have the equivalence:
\[ X \Indep_{\Pr_{V^\aug}} Y \given Z \;\iff\;X \Indep_{\Pr_{V^\aug}} Y \given Z \cup (\tilde H \sm \Psi).\]
So w.l.o.g.\ we can assume $\tilde H \sm \Psi \ins Z$.
By the semi-graphoid axioms \ref{GrAxPr} (1)-(5) we can 
  further reduce to the following three cases: $w \notin X\cup Y \cup Z$ and $X=\{w\}$ and $w \in Z$.\\
``$w \notin X \cup Y \cup Z$'': By replacing $X$ with $X \sm Z$ and $Y$ with $Y \sm Z$ we by \ref{GrAxPr} can be assume that both $X$ and $Y$ are disjoint from $Z$ and thus from $(\tilde H \sm \Psi) \cup W$. So only nodes from $U \cup \Psi$ occur in $X$ and $Y$, i.e.\ we have $X,Y \ins U^\aug$.
We now claim that we have:
\[ X \Indep_{U^\aug}^d Y \given Z' \qquad \text{ with } \qquad Z':=(Z\sm \tilde H) \cup \varphi (Z \cap \Psi).\]
Indeed, let $\pi'$ be a shortest $Z'$-open path from $x \in X$ to $y \in Y$ in $U^\aug$. 
We want to construct a $Z$-open path $\pi$ in $G^\aug$. 
 Every node $v$ in $\pi'$ is either in $U \ins V$ or an $F' \in \tilde H^\marg$. The former nodes will stay for the construction of $\pi$ and the latter ones can be replaced by $\psi(F')$ in $G^\aug$. Furthermore, every edge $v_1 \to v_2$ (or mirrored) in $\pi'$ can be replaced by $v_1\to v_2$ or $v_1 \to w \to v_2$ in $G^\aug$. These replacements lead to a path $\pi$ in $G^\aug$ with nodes only from  $V$ and $\Psi$.
The non-colliders of $\pi$ either correspond to nodes in $\{w\} \cup \tilde H$ or non-colliders of $\pi'$, which are all $Z'$-open.
Since elements from $\tilde H$ can only occur as non-colliders, 
  and $w \notin Z$, the non-colliders of $\pi$ are all $Z$-open in $G^\aug$. Every collider of $\pi$ is also a collider of $\pi'$, which then is an ancestor of $Z\sm \tilde H$. This shows that $\pi$ is $Z$-open in $G^\aug$.\\
By \hyl{dGMP}{dGMP} for $(U^\aug,\Pr_{U^\aug})$ we then get:
\[ X \Indep_{\Pr_{U^\aug}} Y \given Z'.\]
Again by adding the constant maps $\tilde H\sm \Psi$ from $Z$ to $Z'$ we get the claim:
\[ X \Indep_{\Pr_{V^\aug}} Y \given Z.\]
``$X=\{w\}$'': By assumption we have: $\{w\} \Indep_{G^\aug}^d Y \given Z$. 
We first claim that this implies:
\[ \Pa^{G^\aug}(w)\sm\{w\} \Indep_{G^\aug}^d Y \given Z.  \]  
Indeed, consider a $Z$-open path $\pi$ in $G^\aug$ from a node $y \in Y$ to a node $x' \in \Pa^{G^\aug}(w)\sm\{w\}$. We extend $\pi$ to $\pi': \,w \ot x' \cdots y$. Since $w,x' \notin Z$ and $x'$ is not a collider on $\pi'$ the path $\pi'$ is then $Z$-open. 
This reduces to the case: $w \notin X' \cup Y \cup Z$ with $X':=\Pa^{G^\aug}(w)\sm\{w\}$.  By the previous point we get:
\[ \Pa^{G^\aug}(w)\sm\{w\} \Indep_{\Pr_{V^\aug}} Y \given Z.  \]
Since we defined the variable $X_w$ corresponding to the node $w$ exactly as the tuple of variables from $\Pa^{G^\aug}(w)\sm\{w\}$ we get the claim:
\[ \{w\} \Indep_{\Pr_{V^\aug}} Y \given Z.  \]
``$ w \in Z$'': By assumption we have $X \Indep_{G^\aug}^d Y \given Z$.
We claim that this implies:
\[ X \Indep_{U^\aug}^d Y \given Z'  \]
with $Z':= \lp Z \cup \Pa^{G^\aug}(w)\rp\sm(\tilde H \cup \{w\}) \cup \varphi\lp  (Z \cup \Pa^{G^\aug}(w)) \cap \Psi\rp$. 
Let $\pi'$ be a $Z'$-open path in $U^\aug$. As before we can lift this path to a path $\pi$ in $G^\aug$ by replacing nodes from $F' \in \tilde H^\marg$ by $\psi(F')$ and some edges $v_1 \to v_2$ by $v_1 \to w \to v_2$. The latter can not occur because then 
$v_1 \in \Pa^{G^\aug}(w) \cap ( V \cup \Psi) \sm \{w\}$ would be a non-collider on $\pi'$ in $Z'$. So every non-collider on $\pi$ is $Z$-open.
Since every collider on $\pi'$ is an ancestor of either $Z \sm \{w\}$ or $\Pa^{G^\aug}(w) \sm \{w\}$ it follows that every collider on $\pi$ is an ancestor of $Z$ in $G^\aug$. This shows that $\pi$ is $Z$-open.\\
By \hyl{dGMP}{dGMP} for $(U^\aug,\Pr_{U^\aug})$ we then get:
\[ X \Indep_{\Pr_{U^\aug}} Y \given Z'.\]
By adding some constant maps to $Z'$ we get:
\[ X \Indep_{\Pr_{V^\aug}} Y \given \lp Z \cup \Pa^{G^\aug}(w)\rp \sm \{w\}.\]
Again, since $X_w$ is the tuple of variables from $\Pa^{G^\aug}(w)\sm\{w\}$ we get:
\[ X \Indep_{\Pr_{V^\aug}} Y \given Z.  \]
Altogether, these show \hyl{mMP}{mdGMP} for $(G,\Pr_V)$.
\end{proof}

\end{Thm}

Since the general directed global Markov property (\hyl{gdGMP}{gdGMP}) uses the acyclification of a \HEDG{}, which does not commute with the augmentation (see figure \ref{fig:mdGMP-smgdGMP}), we can define an alternative version of the marginal general directed global Markov property as follows.

\begin{Def}[Strong marginal general directed global Markov property]
Let $G=(V,E,H)$ be a \HEDG{} and $\Pr_V$ a probability distribution on $\Xcal_V$.
We define the \emph{strong marginal general directed global Markov property (\hyt{smgdGMP}{smgdGMP})} for $(G,\Pr_V)$ as follows:
\begin{enumerate}
	   \item For every $F \in \tilde{H}$ there exists a standard Borel space $\Ecal_F$.
		 \item There exists a probability distribution $\Pr_{V^\aug}$ on 
		       $\Xcal_{V^\aug}:=\prod_{v \in V}\Xcal_v \x \prod_{F \in \tilde{H}} \Ecal_F$, such that:
     \begin{enumerate}
			   \item for the marginal distribution of $\Pr_{V^\aug}$ on $\Xcal_V$ we have:
				  \[ \pi_{V,*}\Pr_{V^\aug} = \Pr_V,\]
				 where $\pi_V: \Xcal_{V^\aug} \to \Xcal_V$ is the canonical projection, and:
				 \item for all subsets $X,Y,Z \ins V^\aug$ we have the implication:
				\[ X \Indep_{G^\acag}^d Y \given Z \; \implies \; X \Indep_{\Pr_{V^\aug}} Y \given Z,  \]
				where $G^\acag$ is the acyclic augmentation given by the DAG $(V^\aug, (E^\aug)^\acy)$.
    \end{enumerate}
\end{enumerate}
\end{Def}

\begin{Lem}
\label{smgdGMP-mgdGMP}
Let $G=(V,E,H)$ be a \HEDG{} and $\Pr_V$ a probability distribution on $\Xcal_V$.
For $(G,\Pr_V)$ we then have the implications:
\[ \xymatrix{ 
      \text{\hyl{smgdGMP}{smgdGMP}}\ar@{=>}[r]& \text{\hyl{mMP}{mgdGMP}}.%
} \]
\begin{proof}
	\hyl{smgdGMP}{smgdGMP} checks the d-separation in the acyclic augmentation (DAG) $G^\acag$ %
	and \hyl{mMP}{mgdGMP} in the mDAG $(G^\aug)^\acy$. 
By \ref{acag-3-d-sep} we get for all subsets $X,Y,Z \ins V^\aug$:
  \[ X \Indep_{(G^\aug)^\acy}^d Y \given Z \;\stackrel{\ref{acag-3-d-sep}}{\implies}\; X \Indep_{G^\acag}^d Y \given Z \; \stackrel{\text{\hyl{smgdGMP}{smgdGMP}}}{\implies}\; X \Indep_{\Pr_{V^\aug}} Y \given Z.\]
This shows the claim. %
\end{proof}
\end{Lem}

\begin{Lem}
\label{smgdGMP-mdGMP-acy}
Let $G=(V,E,H)$ be a \HEDG{} and $\Pr_V$ a probability distribution on $\Xcal_V$.
We then have the implications:
\[ \xymatrix{ 
      \text{\hyl{smgdGMP}{smgdGMP} for } (G,\Pr_V) \ar@{=>}[r]& \text{\hyl{mMP}{mdGMP} for } (G^\acy,\Pr_V).%
} \]
\begin{proof} 
First fix a surjective map $\varphi: \tilde H \to \tilde H^\acy$ with $\varphi(F) \supseteq F$ for every $F \in \tilde H$
(see \ref{acag-3-d-sep} for existence).
Let \[\begin{array}{rclcl}
\tilde \varphi =(\id\dot\cup\varphi)&:& V^\aug=V \dot\cup \tilde H &\to& V \dot\cup \tilde H^\acy= (V^\acy)^\aug%
\end{array} \]
 its natural extensions. 
By the assumption of \hyl{smgdGMP}{smgdGMP} for $(G,\Pr_V)$ we have random variables $(X_v)_{v \in V}$ and $(E_F)_{F \in \tilde H}$ inducing $\Pr_{V^\aug}$ such that $(G^\acag,\Pr_{V^\aug})$ satisfies \hyl{dGMP}{dGMP}.  
For $F' \in \tilde H^\acy$ we now put $E_{F'}':=(E_F)_{\substack{F \in \tilde{H}\\\varphi(F)=F'}}$, which is well-defined since $\varphi$ is surjective. $(X_v)_{v \in V}$ together with $(E_{F'}')_{F' \in \tilde H^\acy}$ then induce the same distribution $\Pr_{V^\aug}$.
For $X,Y,Z \ins (V^\acy)^\aug$ we then have:
\[ \begin{array}{rcl} &&X \Indep_{(G^\acy)^\aug}^d Y \given Z \\
& \stackrel{\ref{acag-3-d-sep}}{\implies}&
 (\tilde \varphi)^{-1}(X) \Indep_{G^\acag}^d (\tilde \varphi)^{-1}(Y) \given (\tilde \varphi)^{-1}(Z)\\
& \stackrel{\hyl{smgdGMP}{smgdGMP}}{\implies}& (\tilde \varphi)^{-1}(X) \Indep_{\Pr_{V^\aug}} (\tilde \varphi)^{-1}(Y) \given (\tilde \varphi)^{-1}(Z) \\
& \stackrel{E_{F'}'=(E_F)_{F \in \varphi^{-1}(F')}}{\implies}& X \Indep_{\Pr_{V^\aug}} Y \given Z. \\
\end{array}  \]
This shows \hyl{dGMP}{dGMP} for $((G^\acy)^\aug,\Pr_{V^\aug})$ and thus \hyl{mMP}{mdGMP} for $(G^\acy,\Pr_V)$.
\end{proof}
\end{Lem}

\begin{Lem}
\label{mFP-maFP}
Let $G=(V,E,H)$ be a \HEDG{} and $\Pr_V$ a probability distribution on $\Xcal_V$.
Then for $(G,\Pr_V)$ we have the implication:
\[ \xymatrix{ 
    \text{\hyl{mFP}{mFP}}\ar@{=>}[r]&	\text{\hyl{mMP}{maFP}}.
} \]
  If there is a strictly positive joint density $p_{V^\aug}$ on the augmented space $\Xcal_{V^\aug}$ occuring in \hyl{mMP}{maFP} then also the reverse implication holds.
\begin{proof}
The claim follows from the proof of \ref{mFP-anc} and the reverse implication from \ref{aFP-mFP}.
\end{proof}
\end{Lem}

\begin{Rem}
\label{mdGMP-mgdGMP}
\begin{enumerate}
\item By \ref{dGMP-gdGMP} and \ref{aFP-dGMP} we get the implications:
\[ \xymatrix{ 
    \text{\hyl{mMP}{maFP}}\ar@{=>}[r]&	\text{\hyl{mMP}{mdGMP}}\ar@{=>}[r]& \text{\hyl{mMP}{mgdGMP}}.
} \]
The reverse implication on the right holds for mDAGs. The reverse implication on the left holds e.g.\ if there exists a density and perfect elimination order on $G^\aug$ by \ref{dGMP-aFP}, e.g.\ for mDAGs, or if there exists a strictly positive density on $G^\aug$.
\item Note that we do not have the implication ``\hyl{mMP}{mdGMP} $\implies$ \hyl{smgdGMP}{smgdGMP}'' by the example in figure \ref{fig:mdGMP-smgdGMP}. 
\end{enumerate}
\end{Rem}

\subsection{Structural Equations Properties for \HEDG{}es}

In this subsection we will define another kind of ``Markov properties'' that relate a \HEDG{} structure $G=(V,E,H)$ to a probability distribution $\Pr_V$ on $\Xcal_V$: \emph{structural equations properties (SEP)}. Since there are several ways to generalize the structural equations property for DAGs to the case of \HEDG{}es, where we might have cycles and dependent error terms, we need to investigate several meaningful structural equations properties. Some of these are stronger and some weaker, but they will all (except one) be equivalent (see \ref{SEP-equiv}) in the case of mDAGs.\\
One of the main results will be that the \emph{loop-wisely solvable structural equations property (\hyl{lsSEP}{lsSEP})} is stable under marginalizations (see \ref{lsSEP-marg}), has the lifting property (see \ref{mMP-lift}, \ref{lsSEP-lift}), will imply the marginal directed local Markov property (\hyl{mMP}{mdLMP}, see \ref{lsSEP-mdLMP}) for general \HEDG{}es like in the DAG case and will also imply the (strong marginal) general directed global Markov property (\hyl{smgdGMP}{smgdGMP}, see \ref{csSEP-smgdGMP}). The latter result is a generalization of \cite{Spirtes94} Thm. 5, where a similar statement for directed graphs but under more restrictive assumptions (existence of a density, reconstructible and independent error terms, ancestral solvability, differentiability, invertible Jacobian, etc.) and a more technical proof (e.g.\ the explicit factorization of a Jacobian) was shown, to the general case of \HEDG{}es, which besides cycles also allows for dependent error terms.\\
Furthermore, we will investigate the question under which conditions on the structural equations the stronger (see \ref{dGMP-gdGMP}) directed global Markov property (\hyl{dGMP}{dGMP}) based on d-separation and the ancestral factorization property (\hyl{aFP}{aFP}) hold, rather than just the general directed global Markov property (\hyl{gdGMP}{gdGMP}) based on $\sigma$-separation.
An erroneous proof for directed graphs and discrete variables with finite domains was given in \cite{PearlDechter96}. A counterexample was given in \cite{Neal00}. Here, we will correct the proof under the assumptions of ancestrally uniquely solvable structural equations (\hyl{ausSEP}{ausSEP}), which will exclude such counterexamples, and generalize the statement to arbitrary \HEDG{}es with discrete variables (not necessarily of finite domain) in \ref{ausSEP-wmaFP}.\\
For not necessarily discrete variables we will show a result for general \HEDG{}es (see \ref{SEPwared-aFP}) based on structural equations with reconstructible errors and a factorization of a Jacobian determinant (\hyl{SEPwared}{SEPwared}) that implies the ancestral factorization property (\hyl{aFP}{aFP}) and thus the directed global Markov property (\hyl{dGMP}{dGMP}, see \ref{aFP-dGMP}, \ref{dGMP-auPMP}). This result recovers and generalizes the linear case  for directed graphs (see \cite{Spirtes94} Thm 1, \cite{Koster96}) to the non-linear case for \HEDG{}es.

\begin{Def} 
\label{def-sep}
Let $G=(V,E,H)$ be a \HEDG{} and $\Pr_V$ a probability distribution $\Pr_V$ on 
$\Xcal_V=\prod_{v \in V} \Xcal_v$. Then we have the following properties relating $G$ to $\Pr_V$:
\begin{enumerate}

\item The \emph{structural equations property (\hyt{SEP}{SEP})}: There exist:
 \begin{enumerate}
  \item a probability space $(\Omega, \fa, \Pr)$,
	\item random variables $E_F: (\Omega, \fa, \Pr) \to \Ecal_F$ for $F \in \tilde{H}$
	which are jointly $\Pr$-independent, and where $\Ecal_F$ are standard Borel spaces,
  \item measurable functions $f_v: \Xcal_{\Pa^G(v)} \x \Ecal_v \to \Xcal_v$ for $v \in V$,
	   where we put $\Ecal_v:=\prod_{v \in F \in \tilde{H}} \Ecal_F$,
	\item random variables $X_v: (\Omega, \fa, \Pr) \to \Xcal_v$ for $v \in V$ that satisfy:
	 \begin{enumerate}
	   \item $X_v=f_v(X_{\Pa^G(v)}, E_v)$ $\Pr$-a.s., 
			 with $E_v:=(E_F)_{v \in F \in \tilde{H}}$, %
	   \item the joint distribution of $(X_v)_{v \in V}$ under $\Pr$ is given by:
			\[ \Pr^{(X_v)_{v \in V}} = \Pr_V.  \]
	 \end{enumerate}
	In other words, we want the distribution to be explained by variables satisfying structural equations attached to the graph.
 \end{enumerate}

\item The \emph{ancestrally solvable structural equations property (\hyt{asSEP}{asSEP})}: We have the points from \hyl{SEP}{SEP} and:
 \begin{enumerate}[resume]
 \item measurable functions $g_v:  \Ecal_{\Anc^G(v)} \to \Xcal_v$ for $v \in V$ with
			\[  X_v = g_v( E_{\Anc^G(v)}) \qquad \Pr\text{-a.s.},   \]
 with $\Ecal_A:=\prod_{\substack{ F \in \tilde{H} \\ F \cap A \neq \emptyset}} \Ecal_F$ 
 and $E_A:=(E_F)_{\substack{F \in \tilde{H} \\ F \cap A \neq \emptyset}}$ for $A=\Anc^G(v)$.\\
In other words, we want $X_v$ only to be determined by the error terms of it ancestors.
\end{enumerate}

\item The \emph{component-wisely solvable structural equations property (\hyt{csSEP}{csSEP})}: We have the points from \hyl{SEP}{SEP} and:
    \begin{enumerate}[resume]
    \item measurable functions $\tilde{g}_v: \Xcal_{\Pa^G(\Sc^G(v))\sm \Sc^G(v)} \x \Ecal_{\Sc^G(v)} \to \Xcal_v$ for $v \in V$ with
			\[  X_v = \tilde{g}_v(X_{\Pa^G(\Sc^G(v))\sm \Sc^G(v)}, E_{\Sc^G(v)}) \qquad \Pr\text{-a.s.,}   \]
        with $\Ecal_{\Sc^G(v)} := \prod_{ \substack{F \in \tilde{H} \\ F \cap \Sc^G(v) \neq \emptyset}} \Ecal_F$ and
        $E_{\Sc^G(v)} := (E_F)_{\substack{F \in \tilde{H} \\ F \cap \Sc^G(v) \neq \emptyset}}$.
			In other words, we want the variables of each strongly connected component to be determined by the variables of all ingoing arrows.
    \end{enumerate}

  \item The \emph{loop-wisely solvable structural equations property (\hyt{lsSEP}{lsSEP})}: We have the points from \hyl{SEP}{SEP} and:
    \begin{enumerate}[resume]
		\item for every strongly connected induced sub-\HEDG{} $S$ in $G$
		and every $v \in S$ we have a measurable function $\hat{g}_{S,v}: \Xcal_{\Pa^G(S)\sm S} \x \Ecal_{S} \to \Xcal_v$ such that
			\[  X_v = \hat{g}_{S,v}(X_{\Pa^G(S)\sm S}, E_{S}) \qquad \Pr\text{-a.s.},   \]
        where $\Ecal_S:= \prod_{\substack{ F \in \tilde{H} \\ F \cap S \neq \emptyset}} \Ecal_F$ and $E_S:= (E_F)_{\substack{F \in \tilde{H} \\ F \cap S \neq \emptyset}}$.
    \end{enumerate}
		Roughly speaking, we want the variables of every directed cycle/loop of $G$ to be determined by the variables of all ingoing arrows.
		Also note that the existence of functions $f_v$ for $v \in V$ is not needed anymore, as we can put $f_v:=\hat g_{\{v\},v}$ since $\{v\}$ is a strongly connected induced sub-\HEDG{} of $G$ for every $v \in V$.

\end{enumerate}

\noindent For the last property we put $\tilde{g}_v:=\hat g_{\Sc^G(v),v}$. With this we can for the last two properties recursively define the functions $g_v$ only dependent on their ancestral error terms (see \ref{sc-dag}):
\[ g_v(e_{\Anc^G(v)}):=\tilde{g}_v\left( g_{\Pa^G(\Sc^G(v))\sm \Sc^G(v)}(e_{\Anc^G(v)\sm \Sc^G(v)})  , e_{\Sc^G(v)} \right).  \]
We can then also define the following properties:

\begin{enumerate}[resume]

  \item The \emph{ancestrally uniquely solvable structural equations property (\hyt{ausSEP}{ausSEP})}: We have the points from \hyl{asSEP}{asSEP} and:
 \begin{enumerate}[resume]
 \item for all ancestral sub\HEDG{}es $A \ins G$ there exists a 
$\Pr^{E_A}$-zero set 
$N_A \ins \Ecal_A=\prod_{\substack{F \in \tilde{H}\\F \cap A \neq \emptyset}}\Ecal_F$
 such that the induced map $g_A :\Ecal_A \to \Xcal_A$ with components $(g_v)_{v \in A}$ has the property that for all pairs 
$(e_A,x_A) \in N^c_A \x \Xcal_A:$
   \[  g_A(e_A)=x_A \;\iff\; \forall v \in A: x_v = f_v(x_{\Pa^G(v)},e_v).  \]
     In other words, the last equivalence should hold for almost all $e_A$ and all $x_A$. %
\end{enumerate}
		
\item The \emph{component-wisely uniquely solvable structural equations property (\hyt{cusSEP}{cusSEP})}: We have the points from \hyl{csSEP}{csSEP} and:
 \begin{enumerate}[resume]
 \item for all strongly connected components $S \ins G$ there exists a 
   $\Pr^{E_{\Anc^G(S)}}$-zero set $N_S \ins \Xcal_{\Pa^G(S) \sm S} \x \Ecal_S$ w.r.t.\ the map $g_{\Pa^G(S)\sm S}\x \textrm{pr}_S$ (where $\textrm{pr}_S$ is the projection onto the $S$-components)
such that the map $\tilde{g}_S :\Xcal_{\Pa^G(S) \sm S} \x \Ecal_S \to \Xcal_S$ with components $(\tilde{g}_{S,v})_{v \in S}$ has the property that for all tuples 
$(x_{S}, x_{\Pa^G(S)\sm S},e_S) \in \Xcal_S \x N^c_S:$ %
   \[  \tilde{g}_S(x_{\Pa^G(S) \sm S},e_S)=x_S \;\iff\; \forall v \in S: x_v = f_v(x_{\Pa^G(v)},e_v).  \]
	In other words, the last equivalence should hold for almost all $(e_S,x_{\Pa^G(S) \sm S})$ and all $x_S$. %
\end{enumerate}	
		
\item The \emph{loop-wisely uniquely solvable structural equations property (\hyt{lusSEP}{lusSEP})}: We have the points from \hyl{lsSEP}{lsSEP} and:
 \begin{enumerate}[resume]
 \item for every strongly connected induced sub-\HEDG{} $S \ins G$ there exists a 
$\Pr^{E_{\Anc^G(S)}}$-zero set $N_S \ins \Xcal_{\Pa^G(S) \sm S} \x \Ecal_S$ w.r.t.\ the map $g_{\Pa^G(S)\sm S}\x \textrm{pr}_S$
(where $\textrm{pr}_S$ is the projection onto the $S$-components)
such that the map $\hat{g}_S :\Xcal_{\Pa^G(S) \sm S} \x \Ecal_S \to \Xcal_S$ with components $(\hat{g}_{S,v})_{v \in S}$ has the property that for all tuples 
$(x_{S}, x_{\Pa^G(S)\sm S},e_S) \in \Xcal_S \x N^c_S:$
   \[  \hat{g}_S(x_{\Pa^G(S) \sm S},e_S)=x_S \;\iff\; \forall v \in S: x_v = f_v(x_{\Pa^G(v)},e_v).  \]
		In other words, the last equivalence should hold for almost all $(e_S,x_{\Pa^G(S) \sm S})$ and all $x_S$. %
\end{enumerate}			
\end{enumerate}
Another property (\hyl{SEPwared}{SEPwared}) will be defined later in \ref{def-SEPwared}. 
\end{Def}

\begin{Thm} 
\label{SEP-equiv}
For a \HEDG{} $G=(V,E,H)$ and a probability distribution $\Pr_V$ on $\Xcal_V$ we have the implications:
\[ \xymatrix{ 
  \text{\hyl{lusSEP}{lusSEP}}\ar@{=>}[r]\ar@{=>}[d]&		\text{\hyl{cusSEP}{cusSEP}}\ar@{=>}[r]\ar@{=>}[d]& \text{\hyl{ausSEP}{ausSEP}} \ar@{=>}[d]\\%
  \text{\hyl{lsSEP}{lsSEP}}\ar@{=>}[r]&		\text{\hyl{csSEP}{csSEP}}\ar@{=>}[r]& \text{\hyl{asSEP}{asSEP}} \ar@{=>}[r]& \text{\hyl{SEP}{SEP}.} 
} \]
Furthermore, if $G$ is an mDAG (e.g.\ an ADMG or a DAG) then all these properties are equivalent.%
\begin{proof} 
The lower right and the vertical implications are clear.
Since every strongly connected component of $G$ forms a strongly connected induced sub-\HEDG{} of $G$ we have the two very left implications
(with $\tilde{g}_v:=\hat g_{\Sc^G(v),v}$).\\
  Now assume \hyl{csSEP}{csSEP}. Remember that by \ref{sc-dag} the strongly connected components form a DAG and we can recursively insert the functions $\tilde{g}_v$ into each other to get a function $g_v$ only dependent on its ancestral error terms:
\[ g_v(e_{\Anc^G(v)}):=\tilde{g}_v\left( g_{\Pa^G(\Sc^G(v))\sm \Sc^G(v)}(e_{\Anc^G(v)\sm \Sc^G(v)})  , e_{\Sc^G(v)} \right).  \]
  From this \hyl{asSEP}{asSEP} follows recursively.\\
  Now assume \hyl{cusSEP}{cusSEP} and let $g_v$ be defined as above. Let $A \ins G$ be ancestral and $S_1,\dots,S_r$ its strongly connected components ordered according to a topological order of the DAG of strongly connected components of $G$. If $A$ consists only of one strongly connected component, i.e.\ $r=1$ then $g_v=\tilde{g}_v$ and there is nothing to show. We now argue by induction on $r$. Let $S:=S_r$. We then have the following zero sets from before:
$N_{A\sm S} \ins \Ecal_{A\sm S}$, $N_S \ins  \Xcal_{\Pa^G(S)\sm S} \x \Ecal_S $.
The union $N_A$ of their preimages in $\Ecal_A$ under the canonical projection for the first and under the map $g_{\Pa^G(S)\sm S} \x \textrm{pr}_S$ for the second is then a zero set as well. 
We get that for every $e_A \in N_A^c$  the image $(g_{\Pa^G(S)\sm S} \x \textrm{pr}_S)(e_A) \in N_S^c$ and the components 
$e_{A\sm S} \in N_{A\sm S}^c$. %
It follows that for every $(e_A,x_A) \in N_A^c \x \Xcal_A$ we have:
\begin{eqnarray*}
x_A=g_A(e_A) & \iff & 
\left\{ 
\begin{array}{lll}
x_{A\sm S}&=& g_{A\sm S}(e_{A\sm S}),  \\ 
x_S &=& g_{S}(e_{A}) \\ 
    &=& \tilde{g}_S(g_{\Pa^G(S)\sm S}(e_{A\sm S}), e_S)\\
		&=& \tilde{g}_S(x_{\Pa^G(S)\sm S}, e_S),
\end{array} \right \} \\
  & \overset{\text{induction}}{ \underset{\text{\hyl{cusSEP}{cusSEP}}}{\iff} } & 
\left\{ 
\begin{array}{ccccc}
\forall v \in A \sm S : & x_v &=& f_v(x_{\Pa^G(v)}, e_v),  \\ 
\forall v \in S : & x_v &=& f_v(x_{\Pa^G(v)}, e_v), 
\end{array} \right \} \\
&\iff & \forall v \in A : x_v = f_v(x_{\Pa^G(v)}, e_v).
 \end{eqnarray*}
  This shows \hyl{ausSEP}{ausSEP}. \\
Now let $G$ be an mDAG and assume \hyl{SEP}{SEP}.
Since $G$ is acyclic every strongly connected induced sub-\HEDG{} $S$ only consists of a point, i.e.\ $S=\{v\}$ and $\Pa^G(S)\sm S=\Pa^G(v)$. With 
  $\hat{g}_{S,v}:=f_v$, $S=\{v\}$ and $N:=\emptyset$ the property \hyl{lusSEP}{lusSEP} follows trivially.
\end{proof}
\end{Thm}

\begin{Eg}[Marginalization of structural equations]
\label{SEP-no-marg}

\begin{figure}[!h]
\centering
\begin{tikzpicture}[scale=.9, transform shape]
\tikzstyle{every node} = [draw,shape=circle]
\node (X) at (0,2) {$x$};
\node (Y) at (2,0) {$y$};
\node (Z) at (0,-2) {$z$};
\node (W) at (-2,0) {$w$};
\foreach \from/\to in {X/Y, Y/Z, Z/W, W/X}
\draw[-{Latex[length=3mm, width=2mm]}] (\from) -- (\to);
\draw[-{Latex[length=3mm, width=2mm]},out=135, in=225,looseness=8] (W) to (W);
\end{tikzpicture}\qquad
\begin{tikzpicture}[scale=.9, transform shape]
\tikzstyle{every node} = [draw,shape=circle]
\node (X) at (0,2) {$x$};
\node (Y) at (2,0) {$y$};
\node (Z) at (0,-2) {$z$};
\foreach \from/\to in {X/Y, Y/Z, Z/X}
\draw[-{Latex[length=3mm, width=2mm]}] (\from) -- (\to);
\end{tikzpicture}
\caption{A directed graph (left) and its marginalization (right) w.r.t.\ $\{w\}$. A corresponding marginalization operation for structural equations might in general not exist because of the self-loop at $w$ (see \ref{SEP-no-marg}). }
\label{fig:no-marg}
\end{figure}
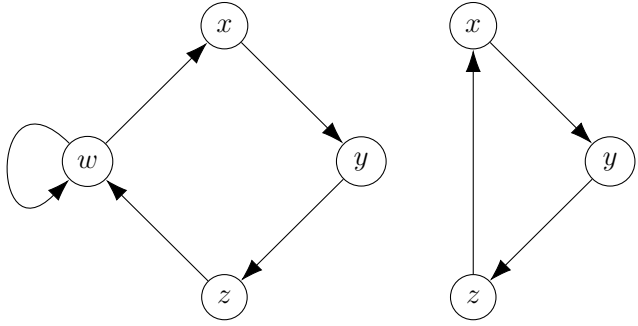
Let $G$ be the directed cycle with set of nodes $V=\{x,y,z,w\}$ and additional selfloop $w \to w$ (see figure \ref{fig:no-marg}).
Assume we have corresponding structural equations of the form:
\[\begin{array}{lllllll}
W & = & f_W(Z,W,E_W) &=& \frac{1}{2}(W+Z+E_W),\\%
X & = & f_X(W,E_X) & =& W \cdot E_X,\\ %
Y & = & f_Y(X,E_Y) &=& X + E_Y, \\
Z & = & f_Z(Y,E_Z) &=& Y \cdot E_Z. %
\end{array}\]
The corresponding variables then induce a distribution $\Pr_V$, which then satisfies the structural equations property (\hyl{SEP}{SEP}) w.r.t.\ $G$.
We now want to marginalize out the node $\{w\}$. This marginalization operation is easily done for the graph $G$ (as in figure \ref{fig:no-marg}).
For the equations, on the other hand, this involves non-trivial assumptions and operations. We need to construct a function of the form:
\[ X=f_X^\marg(Z,E_X^\marg) \qquad \text{ with } \qquad E_X^\marg=(E_W,E_X).\]
The self-dependence of $f_W$ on $W$ forbids us to just ``plug'' $f_W$ into $f_X$ for $W$.
Even the solution functions for the whole loop $S=\{x,y,z,w\}$ as required for \hyl{csSEP}{csSEP} won't directly help:
\[\begin{array}{lllllll}
W & = & \tilde{g}_{S,W}(E_S) &=&\frac{E_W + E_Y \cdot E_Z}{1- E_X \cdot E_Z},\\
X & = & \tilde{g}_{S,X}(E_S) &=& \frac{E_W + E_Y \cdot E_Z }{1- E_X \cdot E_Z} \cdot E_X,\\
Y & = & \tilde{g}_{S,Y}(E_S) &=&\frac{E_X \cdot E_W + E_Y}{1- E_X \cdot E_Z},\\
Z & = & \tilde{g}_{S,Z}(E_S) &=&\frac{E_X \cdot E_W + E_Y}{1- E_X \cdot E_Z} \cdot E_Z,
\end{array}\]
as they are not written in terms of $Z$ and $E_X^\marg=(E_W,E_X)$ as required for $f_X^\marg$, but in terms of all error nodes $E_S=(E_W,E_X,E_Y,E_Z)$ at the same time.
Instead we need to solve for $W$ locally w.r.t.\ the loop $L=\{w\}$, which lies inside the bigger loop $S$. This is a condition of the loop-wisely solvable structural equations property (\hyl{lsSEP}{lsSEP}). In this example we (luckily) get:
\[ W = \hat{g}_{L,W}(Z,E_W) = Z + E_W.   \]
Then we can define the marginalization by plugging $\hat{g}_{L,W}$ into $f_X$:
\[\begin{array}{lllllll}
X &=& f_X^\marg(Z,E_X^\marg) &:=& f_X(\hat{g}_{L,W}(Z,E_W),E_X)\\
&& &=& (Z+E_W)\cdot E_X.
\end{array}\]
The rest of the equations stays unchanged:
\[\begin{array}{lllllll}
Y & = & f_Y^\marg(X,E_Y) &=& X + E_Y, \\
Z & = & f_Z^\marg(Y,E_Z) &=& Y \cdot E_Z.
\end{array}\]
We will see in \ref{lsSEP-marg} that the conditions of \hyl{lsSEP}{lsSEP} ensures the existence of marginalizations w.r.t.\ any set of nodes.
The general idea behind the conditions of \hyl{lsSEP}{lsSEP} is that marginalizations might be done step by step and that marginalizations inside cycles might lead to self-loops, which then also need to be solved to marginalize further. The conditions of \hyl{lsSEP}{lsSEP} anticipate all these possibilities. 
\end{Eg}

\begin{Thm}[\hyl{lsSEP}{lsSEP} stable under marginalization]
\label{lsSEP-marg}
  Let $G=(V,E,H)$ be a \HEDG{} and $\Pr_V$ a probability distribution on $\Xcal_V$ such that $(G,\Pr_V)$ satisfies the loop-wisely solvable structural equations property (\hyl{lsSEP}{lsSEP}).
  Then for every subset $W \ins V$ the marginalizations $(G^{\marg\sm W},\Pr_{V \sm W})$ also satisfy the loop-wisely solvable structural equations property (\hyl{lsSEP}{lsSEP}).
\begin{proof}
By induction we can assume that $W=\{w\}$. Put $G':=G^{\marg\sm W}$. 
  For the \hyl{lsSEP}{lsSEP} for $(G',\Pr_{V\sm W})$ we take the same probability space $(\Omega,\fa,\Pr)$ and the same random variables $X_v$ with $v \in V \sm W$ given by the \hyl{lsSEP}{lsSEP} for $(G,\Pr_{V})$. 
Then clearly we have $\Pr^{(X_v)_{v \in V \sm W}}=\Pr_{V\sm W}$ by marginalization.
For $F \in \tilde{H}$ assign $E_F$ to any fixed $F' \in \tilde{H}^{\marg\sm W}$ with $F^{\marg \sm W} \ins F'$ (see \ref{HEDG-marg-def}).
Every $E_{F'}$ with $F' \in \tilde{H}^{\marg\sm W}$ is then defined to be the tuple of some $E_F$ with $F^{\marg \sm W} \ins F'$.
For $v \in V \sm W$ we put:
\[
f'_v(x_{\Pa^{G'}(v)},e'_v):=  \]
\[\left \{ 
\begin{array}{llll}
f_v(x_{\Pa^{G}(v)},e_v) & \text{if} & w \notin \Pa^G(v), \; \\
f_v(x_{\Pa^{G}(v)\sm \{w\}},f_w(x_{\Pa^G(w)},e_w),e_v) & \text{if} & w \in \Pa^G(v),\\&& w \to w \notin E,  \\
f_v(x_{\Pa^{G}(v)\sm \{w\}},\hat g_{w}(x_{\Pa^G(w)\sm \{w\}},e_w),e_v) & \text{if} & w \in \Pa^G(v),\\&& w \to w \in E, 
\end{array}
\right.
\] 
where $e'_v=(e_{F'})_{v \in F' \in \tilde{H}'}$, $e_v =(e_F)_{v \in F \in \tilde{H}}$.
Note that in the cases where $w \in \Pa^G(v)$ we have $\Pa^{G'}(v) \sni \Pa^G(w)\sm \{w\}$ and that every $F \ni w$ marginalizes to $v$.
  In the last case $\hat g_{w}$ is the function from \hyl{lsSEP}{lsSEP} for $(G,\Pr_V)$ corresponding to the selfloop $w \to w$. \\
We then clearly get for all $v \in V \sm W$ the equations:
\[  X_v = f_v'(X_{\Pa^{G'}(v)},E_v')  \qquad \Pr\text{-a.s.}.   \]
  This shows \hyl{SEP}{SEP} for $(G',\Pr_{V\sm W})$. For \hyl{lsSEP}{lsSEP} now consider a strongly connected induced sub-\HEDG{} $C'$ of $G'$.
Then consider the induced sub-\HEDG{} $\Gamma$ of $C' \cup \{w\}$ in $G$. 
Then $C:=\Sc^\Gamma(C')$ can be viewed %
as a strongly connected induced sub-\HEDG{} 
of $G$ including all nodes from $C'$ and that might include $w$ or not.
For $v \in C'$ we then put:  $$g'_{C',v}(x_{\Pa^{G'}(C')\sm C'}, e_{C'}'):= $$
\[
\left \{ 
\begin{array}{llll}
\hat g _{C,v}(x_{\Pa^{G}(C)\sm C}, e_{C}) & \text{if} & w \notin \Pa^{G}(C)\sm C, \\
\hat g _{C,v}(x_{\Pa^{G}(C)\sm (C \cup \{w\})},f_w(x_{\Pa^G(w)},e_w), e_{C}) & \text{if} & w \in \Pa^{G}(C)\sm C,\\&& w \to w \notin E,  \\
\hat g _{C,v}(x_{\Pa^{G}(C)\sm (C \cup \{w\})}, \hat g_{w}(x_{\Pa^G(w)\sm \{w\}},e_w), e_{C}) & \text{if} & w \in \Pa^{G}(C)\sm C,\\&& w \to w \in E.
\end{array}
\right.
\] 
The right hand side of the equation only makes sense if all occuring arguments on the right also appear on the left.
First note that $\Pa^{G}(C)\sm (C \cup \{w\}) \ins \Pa^{G'}(C')\sm C'$ (and $ w \notin \Pa^{G}(C)\sm C$ in the first case).
So the first argument on the right for all three cases exists also on the left.
Further note that the first case of the equations works for both $w \in C$ and $w \notin C$.
For the last two cases to be well defined we further need to show $\Pa^G(w)\sm \{w\} \ins \Pa^{G'}(C')\sm C'$ (for their second arguments).
We have $\Pa^G(w) \cap C' = \emptyset$ because otherwise $w$ had both a parent and a child in $C$, i.e.\ we had edges $v_i \to w \to v_j$ in $G$ with $v_i,v_j \in C'$. But by the choice of $C$ (as $\Sc^\Gamma(C')$) this would imply $w \in C$, which is a contradiction to $w \in \Pa^G(C)\sm C$. So in these cases we have $\Pa^G(w)\sm \{w\} \ins \Pa^{G'}(C')\sm C'$.
Together this shows that all arguments on the right hand side also exist on the left hand side. Thus the right hand side is well defined.
As before with these definitions we get for all $v \in V \sm W$ the equations:
 \[  X_v = g'_{C',v}(X_{\Pa^{G'}(C')\sm C'}, E_{C'}')  \qquad \Pr\text{-a.s.}.   \]
  This shows the \hyl{lsSEP}{lsSEP} for $(G',\Pr_{V \sm W})$.
\end{proof}
\end{Thm}

The following result will generalize results of \cite{Eva15} from mDAGs to the case of all \HEDG{}es. There it was shown that unknown and unobserved latent spaces can be summarized by hyperedges without changing the expressiveness of the underlying marginal Markov model.
A similar version for the marginal directed global Markov property (\hyl{mMP}{mdGMP}) for \HEDG{}es (instead of \hyl{lsSEP}{lsSEP}) is given in \ref{mdGMP-lift}.

\begin{Thm}[Lifting property of \hyl{lsSEP}{lsSEP}]
\label{lsSEP-lift}
Let $G=(V,E,H)$ be a \HEDG{} and $W \ins V$ a subset, $U:=V\sm W$.
Let $G^{\marg(U)}$ be the marginalization of $G$ w.r.t.\ $W$ and $\Pr_U'$ a probability distribution on $\Xcal_U:=\prod_{u \in U} \Xcal_u$ with standard Borel spaces $\Xcal_u$, $u \in U$, such that $(G^{\marg(U)}, \Pr_U')$ satisfies the loop-wisely solvable structural equations property (\hyl{lsSEP}{lsSEP}). Then there exist standard Borel spaces $\Xcal_w$ for $w \in W$ and a probability distribution $\Pr_V$ on $\Xcal_V=\prod_{v \in V} \Xcal_v$ with marginal $\Pr_U=\Pr_U'$ such that $(G,\Pr_V)$ also satisfies the loop-wisely solvable structural equations property (\hyl{lsSEP}{lsSEP}).
\begin{proof}
By induction we can assume $W=\{w\}$. Identify $U$ with $G^\marg = (V\sm W, E^\marg, H^\marg)$ the marginalization w.r.t.\ $W$.
We use the notations and maps from \ref{mMP-lift}:
\[ \begin{array}{rllllll} 
\varphi &:& \tilde H &\to& \tilde H^\marg & \text{with} & \varphi(F) \supseteq F^\marg,\\
\psi &:& \tilde H^\marg &\to& \tilde H& \text{with}& \psi(F')^\marg = F',
\end{array}\]
with the image $\Psi \ins \tilde H$ of $\psi$ and (see \ref{HEDG-marg-def}):
\[ F^\marg := \left\{ 
\begin{array}{lcl}
  F & \text{ if } & w \notin F,\\
	(F \cup \Ch^G(w))\sm \{w\} &\text{ if } & w \in F.
\end{array}\right. \]
$\varphi$ is surjective, $\psi$ injective and we have $\varphi \circ \psi = \id_{\tilde H^\marg}$.\\
By the assumption of \hyl{lsSEP}{lsSEP} for $(G^\marg, \Pr_{V \sm W}')$ we have random variables $X_v$ for $v \in V \sm \{w\}$ that induce $\Pr_{V \sm W}'$ and $E'_{F'}$ for $F' \in \tilde H^\marg$. 
As in \ref{mMP-lift} we define $\Ecal_F$ for $F \in \tilde H$ as follows:
\[ \Ecal_F := \left\{
\begin{array}{lll}
\{*\} & \text{ if }& F \notin \Psi,\\
\Ecal'_{\psi^{-1}(F)} & \text{ if } & F \in \Psi,
\end{array}
\right.  \]
where $\{*\}$ denotes the space of one point only. %
Similarly we define the variables $E_F$ for $F \in \tilde H$ out of $E_{F'}'$.
We also define (see \ref{mMP-lift}):
\[\begin{array}{rcl}
 \Xcal_w &:=& \Xcal_{\Pa^G(w)\sm\{w\}} \x \prod_{\substack{F \in \tilde H \\ w \in F}}\Ecal_F,\\
 X_w &:=& (X_{\Pa^G(w)\sm\{w\}},(E_F)_{\substack{F \in \tilde H \\ w \in F}}).
\end{array}
\]
Furthermore, we define:
\[\hat g_{\{w\},w}=\id:\quad \Xcal_{\Pa^G(w)\sm \{w\}} \x \prod_{\substack{F \in \tilde H\\w \in F}}\Ecal_F \stackrel{=}{\to} \Xcal_w \]
as the identity map.
For $v \in \Ch^G(w)\sm\{w\}$ we then define $\hat g_{\{v\},v}$ as the composite of the following natural projections ($\srj$) and the given map $\hat g_{\{v\},v}'$ as follows:
\[\begin{array}{rclcccl}
& \hat g_{\{v\},v}: & \Xcal_{\Pa^G(v)\sm \{v\}} &&&\x&\prod_{\substack{F \in \tilde H\\v \in F}}\Ecal_F\\
&  = & \Xcal_{\Pa^G(v)\sm \{v,w\}} &\x&  \Xcal_w &\x&\prod_{\substack{F \in \tilde H\\v \in F}}\Ecal_F  \\
&\srj& \Xcal_{\Pa^G(v)\sm \{v,w\}} &\x& \Xcal_{\Pa^G(w)\sm\{w\}} \x \prod_{\substack{F \in \tilde H\\w \in F\\v \notin F}}\Ecal_F &\x&\prod_{\substack{F \in \tilde H\\v \in F}}\Ecal_F \\
&\srj& \Xcal_{\Pa^{G^\marg}(v)\sm \{v\}} & \x& \prod_{\substack{F \in \Psi\\w \in F\\v \notin F}}\Ecal_F &\x&\prod_{\substack{F \in \Psi\\v \in F}}\Ecal_F \\
&=& \Xcal_{\Pa^{G^\marg}(v)\sm \{v\}} & \x& \prod_{\substack{F' \in \tilde H^\marg\\w \in \psi(F')\\v \notin \psi(F')}}\Ecal_{F'}' &\x&\prod_{\substack{F' \in \tilde H^\marg\\v \in \psi(F')}}\Ecal_{F'}' \\
&\stackrel{(\#)}{=}& \Xcal_{\Pa^{G^\marg}(v)\sm \{v\}} & && \x& \prod_{\substack{F' \in \tilde H^\marg\\v \in F'}}\Ecal_{F'}' \\
&\stackrel{\hat g_{\{v\},v}'}{\to} & \Xcal_v.
\end{array}  \]
Note that in $(\#)$ we used that for $v \in \Ch^G(w) \sm \{w\}$ and $F' \in \tilde H^\marg$ we have the equivalence:
\[ v \in F' \;\stackrel{F'=\psi(F')^\marg}{\iff}\;v \in \psi(F') \;\lor\; w \in \psi(F'). \]
For $v \in V \sm( \Ch^G(w) \cup \{w\})$ we define $\hat g_{\{v\},v}$ similarly (but without the middle column) and the fact that for 
$v \in V \sm( \Ch^G(w) \cup \{w\})$ and $F' \in \tilde H^\marg$ we have the equivalence:
 \[ v \in F' \;\stackrel{F'=\psi(F')^\marg}{\iff}\; v \in \psi(F').\]
This shows \hyl{SEP}{SEP} for $(G, \Pr_V)$ with $\Pr_V:=\Pr^{(X_v)_{v \in V}}$.
Now let $S \ins G$ be a ``loop'' (i.e.\ a strongly connected induced sub-\HEDG{}) with $\# S \ge 2$.\\
First consider the case $w \in S$. Then $S':=S \sm \{w\}$ is a loop of $G^\marg$ and 
$\Pa^{G^\marg}(S')\sm S' = \Pa^G(S) \sm S$.
We then define the map
\[\hat g_{S,w}: \Xcal_{\Pa^G(S)\sm S} \x \prod_{\substack{F \in \tilde H\\F \cap S \neq \emptyset}}\Ecal_F \to \Xcal_w \]
via the equation:
\[\begin{array}{rcl}&& \hat g_{S,w}(x_{\Pa^G(S)\sm S},e_S)\\
  &:=& \hat g_{\{w\},w}(\hat g'_{S',(\Pa^G(w)\cap S)\sm\{w\}}(x_{\Pa^{G^\marg}(S')\sm S'},e_{S'}'),x_{\Pa^G(w)\sm S},e_w). 
\end{array}\]
Note that the components of $e_w$ come from the ones in $\Xcal_w$. Further note for $e'_{S'}$ that with the same arguments as before we have the natural projections:
\[\begin{array}{rclcccl}
\Ecal_{S}&=&\prod_{\substack{F \in \tilde H\\F \cap S \neq \emptyset}}\Ecal_F\\
&\srj& \prod_{\substack{F \in \tilde H\\F \cap (S\sm \{w\}) \neq \emptyset}}\Ecal_F&\x&\prod_{\substack{F \in \tilde H\\w \in F}}\Ecal_F\\
&\srj&\prod_{\substack{F \in \Psi\\F \cap (S\sm \{w\}) \neq \emptyset}}\Ecal_F&\x&\prod_{\substack{F \in \Psi\\w \in F}}\Ecal_F\\
&\stackrel{(\#)}{\srj}& \prod_{\substack{F' \in \tilde H^\marg\\F' \cap (S\sm \{w\}) \neq \emptyset}}\Ecal_{F'}'\\
&=&\Ecal'_{S'}
\end{array}  \]
So the map $\hat g_{S,w}$ is well defined. Similar arguments show
that $\hat g_{S,v}$ for $v \in S$ can be defined via $\hat g_{S',v}'$ in the case $w \in S$.\\
If $w \in \Pa^G(S) \sm S$ and $v \in S$ we again define $\hat g_{S,v}$ via $\hat g_{S',v}'$. For this note that
$\Pa^{G^\marg}(S)\sm S = (\Pa^G(S) \cup \Pa^G(w)) \sm (S \cup \{w\})$. Again $\Xcal_{\Pa^G(w)\sm \{w\}}$ is a component of $\Xcal_w$ and 
$\Ecal'_{S'}$ a projection of $\Ecal_S$. So everything is well defined. \\
With the same arguments we can show the easier cases when $w \notin \Pa^G(S) \cup S$ and $v \in S$.
With all these definitions we clearly get the (almost-surely) equations:
\[  X_S = \hat g_S(X_{\Pa^G(S)\sm S}, E_S).  \]
This shows \hyl{lsSEP}{lsSEP} for $(G,\Pr_V)$.
\end{proof}
\end{Thm}

We will now turn to the question what kind of conditional independencies are encoded in a distribution $\Pr_V$ that satisfies a structural equations property (\hyl{SEP}{SEP}) w.r.t.\ some \HEDG{} $G$. The next example will show that even under the strong assumptions of \hyl{lusSEP}{lusSEP} and strictly positive density ($p>0$) we do not get the directed global Markov property (\hyl{dGMP}{dGMP}). A similar example but with a different graph was given in \cite{Spirtes94}.

\begin{Eg}[\hyl{lusSEP}{lusSEP} does not imply \hyl{dGMP}{dGMP}]
\label{main-example}
Let the notation be like in example \ref{SEP-no-marg} with figure \ref{fig:no-marg}, where we considered a directed cycle with four nodes and the following structural equations:
\[\begin{array}{lllllll}
W & = & f_W(Z,W,E_W) &=& \frac{1}{2}(W+Z+E_W),\\%
X & = & f_X(W,E_X) & =& W \cdot E_X,\\ %
Y & = & f_Y(X,E_Y) &=& X + E_Y, \\
Z & = & f_Z(Y,E_Z) &=& Y \cdot E_Z. %
\end{array}\]
We can then also reconstruct the error terms with the following functions:
\[\begin{array}{llllcll}
E_W & = & h_W(W,Z) &=& W-Z,\\
E_X & = & h_X(X,W) &=& X/W,\\
E_Y & = & h_Y(Y,X) &=& Y-X,\\
E_Z & = & h_Z(Z,Y) &=& Z/Y.
\end{array}\]
This transformation has the Jacobian:
\[ |h'|(w,x,y,z) = \left|\begin{array}{cccc}
1 & 0 & 0& -1\\
- \frac{x}{w^2} & \frac{1}{w} & 0 & 0 \\
0 & -1 & 1 & 0 \\
0 & 0 & -\frac{z}{y^2} & \frac{1}{y} 
\end{array} \right| = \frac{wy - xz}{w^2y^2}.\]
So the probability density for $(W,X,Y,Z)$ is then: $p(w,x,y,z)=$ %
\[ p_{E_W}(h_W(w,z)) \cdot p_{E_X}(h_X(x,w)) \cdot p_{E_Y}(h_Y(y,x)) \cdot p_{E_Z}(h_Z(z,y)) \cdot |h'|(w,x,y,z). \]
With this one can check the equation:
\begin{eqnarray*} 
&&p(w,x,y,z) \cdot p(y,w) - p(w,x,y) \cdot p(w,y,z)\\
&=& \frac{p(w,x,y,z)}{ wy-xz} \cdot \int p_{E_W}(h_W(w,z'))\cdot p_{E_Z}(h_Z(z',y)) \cdot (z'-z)\, dz' \\
&& \quad\cdot \int p_{E_X}(h_X(x',w)) \cdot p_{E_Y}(h_Y(y,x')) \cdot (x-x') \,dx'.
\end{eqnarray*}
Because every factor in the last term equals zero only on a Lebesgue zero-set we do not get the conditional independence: 
\[X \Indep_\Pr Z \given (Y,W),\] which would be implied by the directed global Markov property (\hyl{dGMP}{dGMP}) using the d-separation criterion on the directed graph. \\
Note that this is despite the fact that the density $p$ is strictly positive almost everywhere, the observed variables are unique functions 
($\tilde{g}_S(E_S)$) in the error terms, which are reconstructible from the observed variables, and which can be assumed to be independent, standard Gaussian.\\
This happens because the error terms are entangled in the representation of the variables via the solution functions $\tilde g_S$. Every variable depends on every other variable's error term $E_S=(E_W,E_X,E_Y,E_Z)$ (even on their childrens') even though the error terms were assumed to be independent. The Jacobian $|h'|$ then measures the degree of entanglement of the error terms in a certain sense. Since $|h'|$ is non-constant by the non-linearities in the structural equations and does not factor according to the moralization of the graph of observed variables it obstructs against a needed factorization of the probability density and thus against the directed global Markov property (\hyl{dGMP}{dGMP}).
The definition \hyl{SEPwared}{SEPwared} provides an additional assumption (see \ref{SEPwared-aFP}) that does imply the \hyl{dGMP}{dGMP}.
\end{Eg}

A first positive result about \hyl{lsSEP}{lsSEP} and conditional independencies encoded in \hyl{mMP}{mdLMP} is the following:

\begin{Lem}
\label{lsSEP-mdLMP}
Let $G=(V,E,H)$ be a \HEDG{} and $\Pr_V$ a probability distribution on $\Xcal_V$.
  Assume that $(G,\Pr_V)$ satisfies the loop-wisely solvable structural equations property (\hyl{lsSEP}{lsSEP}).
Then for every subset $W \ins V$ with the marginalized \HEDG{}-structure $G^{\marg(W)}$ we have:
\[ \begin{array}{lllll}
 \{v\} \Indep_{\Pr_{W^\aug}} W^\aug \sm \{v\} \given \Pa^{W^\aug}(v) \sm \{v\} & \text { for every } & v \in W, \\
 \{F\} \Indep_{\Pr_{W^\aug}} \NonDesc^{W^\aug}(F) \given \emptyset & \text { for every } & F \in \tilde{H}^{\marg(W)}.
\end{array}\]
Furthermore, these imply that for every $v \in W^\aug$ we have:
\[ \{v\} \Indep_{\Pr_{W^\aug}} A \sm \{v\} \given \partial_{A^\moral}(v) , \]
  for every 
	$A \ins \NonDesc^{W^\aug}(v) \cup \Sc^{W^\aug}(v)$ with $v \in A$ and in the marginalized structure. 
  In particular, $(G,\Pr_V)$ and all its marginalizations satisfy the marginal directed local Markov property (\hyl{mMP}{mdLMP}).
\begin{proof}
  By \ref{lsSEP-marg} also $(W,\Pr_{V\sm W})$ satisfies \hyl{lsSEP}{lsSEP}. So w.l.o.g.\ let $W=G$.
  By \hyl{lsSEP}{lsSEP} we then have random variables $X_v$, $v \in V$, and $E_F$, $F \in \tilde{H}$ and functional relations $X_v=f_v(X_{\Pa^G(v)\sm \{v\}},E_v)$ with $E_v=(E_F)_{v \in F}$. 
  Note that by \hyl{lsSEP}{lsSEP} $X_v$ can be assumed to be a function in the variables from $\Pa^{G^\aug}(v) \sm \{v\}$ only (i.e.\ no self-dependence). 
So for all sets $A,Z \ins V^\aug$ with $\Pa^{G^\aug}(v) \sm \{v\} \ins Z$ 
  (e.g.\ $Z =\partial_{A^\moral}(v)$) 
	we get:
\[ \{v\} \Indep_{\Pr_{V^\aug}} A \sm \{v\} \given Z.   \] 
  For $F \in \tilde{H}$ and \hyl{asSEP}{asSEP} we can write the variables from $\NonDesc^{G^\aug}(F)$ as functions of $(E_{F'})_{F' \neq F}$. By the joint independence assumption of the error terms we get:
\[ \{F\} \Indep_{\Pr_{V^\aug}} \NonDesc^{G^\aug}(F) \given \emptyset. \]
  Note that $\partial_{A^\moral}(F)=\Pa^{A}(F)=\emptyset$ for the ancestral subgraph $A:=\NonDesc^{G^\aug}(F) \cup \Sc^G(F)$. 
So for every $v \in V^\aug$, every ancestral $A \ins \NonDesc^{G^\aug}(v) \cup \Sc^{G^\aug}(v)$ with $v \in A$ and every $S \ins \Sc^{G^\aug}(v) \sm \{v\}$ we get:
\[ \{v\} \Indep_{\Pr_V} A_{\sm S} \sm \{v\} \given \partial_{A_{\sm S}^\moral}(v),\]
where  $A_{\sm S} := A^{\marg\sm S}$. This shows \hyl{mMP}{mdLMP}.
\end{proof}
\end{Lem}

\begin{Rem}
\label{dLMP-not-marg}
Even though in \ref{lsSEP-mdLMP} we showed for every \HEDG{} $G=(V,E,H)$ and every of its marginalizations the implication:
\[ \xymatrix{ 
      \text{\hyl{lsSEP}{lsSEP} for } (G,\Pr_V)\ar@{=>}[r] & \text{\hyl{mMP}{mdLMP} for every marginalization } (G^{\marg(W)},\Pr_W),
} \]
  i.g.\ we do not get the (ordinary) directed local Markov property (\hyl{dLMP}{dLMP}) for $(G,\Pr_V)$ or any of its marginalizations.
	Otherwise the example \ref{main-example}, which satisfies \hyl{lsSEP}{lsSEP} and thus \hyl{mMP}{mdLMP}, would satisfy \hyl{dLMP}{dLMP} and thus
	$X \Indep_\Pr Z \given (Y,W)$ but this contradicts the fact $X \nIndep_\Pr Z \given (Y,W)$, which was shown there.\\
	If, on the other hand, \hyl{dLMP}{dLMP} is equivalent to \hyl{dGMP}{dGMP} (see e.g.\ \ref{oLMP-dLMP-dGMP}) then we can take marginalizations.
\end{Rem}

\begin{Lem}
\label{mDAG-mdLMP-SEP}
 For an mDAG $G=(V,E,H)$ and a probability distribution $\Pr_V$ on $\Xcal_V$ we have the implication:
\[ \xymatrix{ 
      \text{\hyl{mMP}{mdLMP}}\ar@{=>}[r] & \text{\hyl{lusSEP}{lusSEP}}. 
} \]
\begin{proof}
We have the chain of implications:
\[\begin{array}{rcl}
  \text{\hyl{mMP}{mdLMP} for } (G,\Pr_V) & \stackrel{\ref{def-marg}}{\implies} & \text{\hyl{dLMP}{dLMP} for } (G^\aug,\Pr_{V^\aug}) \\
  & \stackrel{\ref{dLMP-oLMP}, \text{ top.ord.}}{\implies} & \text{\hyl{oLMP}{oLMP} for } (G^\aug,\Pr_{V^\aug}) \\
  & \stackrel{\ref{oLMP-rFP}}{\iff} & \text{\hyl{rFP}{rFP} for } (G^\aug,\Pr_{V^\aug}) \\
  & \stackrel{\ref{DAG-rFP}}{\implies} & \text{\hyl{rFP-DAG}{rFP} for DAG } (G^\aug,\Pr_{V^\aug}) \\
& \stackrel{\ref{sep-lemma}}{\implies} & \text{\hyl{SEP}{SEP} for } (G^\aug,\Pr_{V^\aug}) \\
  & \stackrel{\ref{SEP-equiv}, \text{ mDAG}}{\implies} & \text{\hyl{lsSEP}{lsSEP} for } (G^\aug,\Pr_{V^\aug})\\
  & \stackrel{\ref{lsSEP-marg}}{\implies} & \text{\hyl{lsSEP}{lsSEP} for } (G,\Pr_{V}) \\
  & \stackrel{\ref{SEP-equiv}, \text{ mDAG}}{\implies} & \text{\hyl{lusSEP}{lusSEP} for } (G,\Pr_{V}).
\end{array}\]
\end{proof}
\end{Lem}

\begin{Cor}
\label{mDAG-all}
Let $G=(V,E,H)$ be an mDAG and $\Pr_V$ a probability distribution on $\Xcal_V$. 
Then the following properties for $(G,\Pr_V)$ are all equivalent (using topological orders if needed):
  \[ \begin{array}{ccccccc}& &\text{\hyl{mMP}{moLMP}}&\iff& \text{\hyl{mMP}{mrFP}}  \\
    &\iff& \text{\hyl{SEP}{SEP}} &\iff& \text{\hyl{lusSEP}{lusSEP}} \\
    &\iff& \text{\hyl{mMP}{mdLMP}}  &\iff& \text{\hyl{mMP}{mdGMP}} &\iff\;& \text{\hyl{mMP}{mgdGMP}}.
\end{array}\]
\begin{proof}
  See \ref{oLMP-rFP}, \ref{SEP-equiv}, \ref{mDAG-mdLMP-SEP}, \ref{lsSEP-mdLMP}, \ref{dGMP-gdGMP}, \ref{mdGMP-dGMP}, \ref{oLMP-dLMP-dGMP}.
\end{proof}
\end{Cor}

The next result is the mentioned generalization of \cite{Spirtes94} Thm. 5 from directed graphs under restrictive assumptions to the case of arbitrary \HEDG{}es under a solvability assumption on the structural equations (\hyl{csSEP}{csSEP}).

\begin{Thm}
\label{csSEP-smgdGMP}
For a \HEDG{} $G=(V,E,H)$ and a probability distribution $\Pr_V$ on $\Xcal_V$ we have the following implications for $(G,\Pr_V)$:
\[ \xymatrix{ 
      \text{\hyl{csSEP}{csSEP}}\ar@{=>}[r]& \text{\hyl{smgdGMP}{smgdGMP}.} 
} \]
  In words: For all \HEDG{}es the component-wisely solvable structural equations property (\hyl{csSEP}{csSEP}) implies the strong marginal general directed global Markov property (\hyl{smgdGMP}{smgdGMP}).
\begin{proof}
  Only consider the functions $\tilde{g}_v$, $v \in V$, given by \hyl{csSEP}{csSEP}. Every $\tilde{g}_v$ depends on the variables $X_w$ and error terms $E_F$ going into $\Sc^G(v)$. 
So the dependency structure of these functions is given by the DAG $G^\acag:=(V^\aug, (E^\aug)^\acy)$.
So the probability distribution $\Pr_{V^\aug}:=\Pr^{\lp(X_v)_{v \in V}, (E_F)_{F \in \tilde{H}}\rp}$ together with the DAG $G^\acag$ satisfies the
structural equations property (\hyl{SEP}{SEP}) with those $\tilde{g}_v$, $v \in V$, as functions.
  So we get the directed global Markov property (\hyl{dGMP}{dGMP}) for $(G^\acag, \Pr_{V^\aug})$ by the DAG case \ref{dag-markov-prop} or with \ref{mDAG-all}, which was proven in this paper. This shows the claim.
\end{proof}
\end{Thm}

The next theorem fixes the proof of \cite{PearlDechter96} that roughly stated that structural equations on directed graphs with discrete distributions with finite domain imply the directed global Markov property (\hyl{dGMP}{dGMP}). 
To avoid the counterexample from \cite{Neal00} we make the assumption of ancestrally unique solvability (\hyl{ausSEP}{ausSEP}).
Furthermore, we generalize the claim to work also for arbitrary \HEDG{}es with discrete distributions (not necessarily of finite domain)
and at the same time get as a stronger result the marginal ancestral factorization property (\hyl{mMP}{maFP}).
Note that the latter even implies the marginal directed global Markov property (\hyl{mMP}{mdGMP}), the ordinary ancestral factorization property (\hyl{aFP}{aFP}) and the ordinary directed global Markov property (\hyl{dGMP}{dGMP}) by \ref{aFP-dGMP}, \ref{mdGMP-dGMP}, \ref{dGMP-auPMP}.

\begin{Thm}
\label{ausSEP-wmaFP}
  Let $G=(V,E,H)$ be a \HEDG{} and $\Pr_V$ a discrete (N.B.) probability distribution on $\Xcal_V$.
Then we have the implications:
\[ \xymatrix{ 
      \text{\hyl{ausSEP}{ausSEP}}\ar@{=>}[r]& \text{\hyl{mMP}{maFP}} %
} \]
In words: For all \HEDG{}es with discrete probability distributions the ancestrally uniquely solvable structural equations property
 implies the marginal ancestral factorization property.%
\begin{proof}
We put $\nu_F:=\Pr^{E_F}$ as the measure on $\Ecal_F$ for every $F \in \tilde{H}$. Let $\mu_v$ be a counting measure ($=$ sum of Dirac measures)
on $\Xcal_v$ for every $v \in V$ such that $\Pr_V$ has a density/mass function $p_V$ w.r.t.\ the product measure $\mu_V:=\otimes_{v \in V} \mu_v$.
Let $A' \ins V^\aug$ be an ancestral sub-\HEDG{} and $A = A' \cap V$. %
The function:
\[ q_A: \Xcal_{A} \x \Ecal_{A} \to \R_{\ge 0},\qquad \text{ with }\qquad q_A(x,e):= \I_{x}(g_A(e)), \]
is then a density of $\Pr_{A'}$ w.r.t.\ the product measure $\mu_A \otimes \nu_A$.
Indeed for measurable $B \ins \Xcal_A$ and $D \ins \Ecal_A$ we have:
\begin{eqnarray*}
&& \Pr(X_A \in B, E_A \in D) \\
&=& \int_{D} \Pr(X_A \in B|E_A=e) \, d \Pr^{E_A}(e) \\
&\stackrel{\text{discr.}}{=}& \sum_{x \in B} \int_{D} \Pr(X_A =x|E_A=e) \, d \Pr^{E_A}(e) \\
&\stackrel{X_A=g_A(E_A)}{=}&  \sum_{x \in B} \int_{D} \Pr(g_A(E_A) =x|E_A=e) \, d \Pr^{E_A}(e) \\
&\stackrel{}{=}&  \sum_{x \in B} \int_{D} \I_{x}(g_A(e)) \, d \Pr^{E_A}(e) \\
&\stackrel{}{=}&  \int_{B} \int_{D} q_A(x,e) \, d \nu_A(e) \,d\mu_A(x).
 \end{eqnarray*}
  Furthermore, for $x_A=(x_v)_{v \in A} \in \Xcal_A$ by the unique solvability (\hyl{ausSEP}{ausSEP}) we have the equivalence for $\Pr^{E_A}$-almost-all $e_A=(e_v)_{v \in A}$:
\[ \forall v \in A:\;g_v(e_A)=x_v \quad \iff \quad \forall v \in A:\; x_v=f_v(x_{\Pa^G(v)},e_v).\]
From this we get the factorization for the density:
\[ q_A(x_A,e_A) = \I_{x_A}(g_A(e_A)) = \prod_{v \in A} \I_{x_v=f_v(x_{\Pa^G(v)},e_v)}. \]
  Note that the last factors are functions in the variables from $\{v\} \cup \Pa^{A'}(v)$ and the latter set is a complete subgraph of $(A')^\moral$ by definition. This shows \hyl{mMP}{maFP}.
\end{proof}
\end{Thm}

\begin{Rem}
\label{cusSEP-mFP}
  Note that in the proof of \ref{ausSEP-wmaFP} we do not get the marginal factorization property (\hyl{mFP}{mFP}). 
For this we would need to show in addition that for every strongly connected component $S \in \Scal(G)$ and almost every value $x_{\Pa^G(S)\sm S}, e_S$ the equation:
\[ \int_{\Xcal_S} \prod_{v \in S} \I_{x_v=f_v(x_{\Pa^G(v)},e_v)} \, d\mu_S(x_S) =1  \]
holds. The integral counts the solutions on each strongly connected component $S$ given the values $x_{\Pa^G(S)\sm S}$, $e_S$.
  Since we only enforced unique solvability  for ancestral subsets $A$ (\hyl{ausSEP}{ausSEP}) we do not get the above equation for $S$ i.g..
  By assuming unique solutions on every strongly connected component (\hyl{cusSEP}{cusSEP}), the equation follows and we also get the marginal factorization property (\hyl{mFP}{mFP}).
\end{Rem}

In the rest of this subsection we will investigate the question under which assumption on the structural equations we will get the ancestral factorization property (\hyl{aFP}{aFP}) and thus the directed global Markov property (\hyl{dGMP}{dGMP}) for not necessarily discrete distributions $\Pr_V$ and not necessarily linear equations. Since the claim is quite strong the needed assumptions are lengthy and quite restrictive. Nonetheless, these assumptions are not empty in general as they are generalizations of the linear model, which we will get as an example.

\begin{Def}
\label{def-SEPwared}
Let $G=(V,E,H)$ be a \HEDG{} and $\Pr_V$ a probability distribution $\Pr_V$ on 
$\Xcal_V=\prod_{v \in V} \Xcal_v$. Then we also have the following property relating $G$ to $\Pr_V$:
\begin{enumerate}
  \item[8.] The \emph{structural equations property with ancestrally reconstructible errors and density (\hyt{SEPwared}{SEPwared})}: We have:
     \begin{enumerate}
         \item a probability space $(\Omega, \fa, \Pr)$,
				 \item standard Borel spaces $\Ecal_v$ for every $v \in V$ and random variables $E_v: (\Omega, \fa, \Pr) \to \Ecal_v$,				     
	       \item  standard Borel spaces $\Ecal_F$ and random variables $E_F: (\Omega, \fa, \Pr) \to \Ecal_F$ for $F \in \tilde{H}$ such that:
				       \begin{enumerate}
	                \item $(E_F)_{F \in \tilde{H}}$  are jointly $\Pr$-independent,
	                \item there are measurable functions $\pi_v: \prod_{v \in F \in \tilde{H}} \Ecal_ F \to \Ecal_v$ for $v \in V$ such that:
						            \[E_v=\pi_v((E_F)_{v \in F \in \tilde{H}}), \qquad \Pr\text{-a.s.,}\]
							\end{enumerate}
				\item a %
				measure $\nu_v$ on $\Ecal_v$ for every $v \in V$, such that $\Pr^{(E_v)_{v \in V}}$ 
				     has a density $p_{\Ecal,V}$ w.r.t.\ the product measure $\nu_V=\otimes_{v \in V} \nu_v$,
				\item a %
				measure $\mu_v$ on $\Xcal_v$ for every $v \in V$,			
         \item measurable functions $f_v: \Xcal_{\Pa^G(v)} \x \Ecal_v \to \Xcal_v$ for $v \in V$,
	       \item random variables $X_v: (\Omega, \fa, \Pr) \to \Xcal_v$ for $v \in V$ that satisfy:
	          \begin{enumerate}
          	   \item $X_v=f_v(X_{\Pa^G(v)}, E_v)$ $\Pr$-a.s., 
	             \item the joint distribution of $(X_v)_{v \in V}$ under $\Pr$ is given by:
	           		\[ \Pr^{(X_v)_{v \in V}} = \Pr_V.  \]
	           \end{enumerate}
	       \item  measurable functions $h_v: \Xcal_{D(v) \cup \Pa^G(D(v)) } \to \Ecal_ v$ for $v \in V$,  where $D(v):=\Dist^{\Anc^G(v)}(v)$ 
				   is the district of $v$ in $\Anc^G(v)$, such that:
					\begin{enumerate}
          	   \item $  h_v(X_{D(v) \cup \Pa^G(D(v))})=E_v \qquad \Pr\text{-a.s.}, $
	             \item for all ancestral sub-\HEDG{}es $A \ins G$ the induced map $h_A: \prod_{v \in A} \Xcal_v \to \prod_{v \in A} \Ecal_v$ with 
							   components $(h_v)_{v \in A}$ is a bijection 
								with inverse function $g_A:=h_A^{-1}: \prod_{v \in A} \Ecal_v \to \prod_{v \in A} \Xcal_v$
								(which is then measurable by \cite{Kec95} 15.A),
							such that the image measure $\nu_A^{g_A}=g_{A,*}\lp \otimes_{v \in A} \nu_v \rp$ 
							has a density $|h_A'|:=\frac{d\nu_A^{g_A}}{d \mu_A}$ w.r.t.\ the product measure $\mu_A=\otimes_{v \in A} \mu_v$ and factorizes according to the complete subgraphs of $A^\moral$, i.e., for $\mu_A$-almost all $x_A \in \Xcal_A$:
                \[|h_A'|(x_A) = \prod_{C \in \mathcal{C}_A} k_{A,C}(x_C)\]
                for a finite set $\mathcal{C}_A$ of complete subgraphs of $A^\moral$ and integrable functions $k_{A,C} : \Xcal_C \to \R_{\ge 0}$. 
	           \end{enumerate}	
   \end{enumerate}
\end{enumerate}

\end{Def}

\begin{Thm}
\label{SEPwared-aFP}
Let $G=(V,E,H)$ be a \HEDG{} and $\Pr_V$ a probability distribution on $\Xcal_V$.
Then we have the implication:
\[ \xymatrix{ 
      \text{\hyl{SEPwared}{SEPwared}}\ar@{=>}[r]& \text{\hyl{aFP}{aFP}} \; \&\; \text{\hyl{asSEP}{asSEP}}. 
} \]
\begin{proof}
For $A \ins G$ ancestral we have $h_A(X_A)=E_A$ $\Pr$-a.s. with $E_A=(E_v)_{v \in A}$. It follows that $X_A=g_A(E_A)$ $\Pr$-a.s. and:
\[ \Pr^{X_A} = g_{A,*} \Pr^{E_A}.  \] %
The latter has a density w.r.t.\ to $\nu_A^{g_A}=g_{A,*} \nu_A$ given by:
\[ \tilde{p}_A(x_A)=g_{A,*}(p_{\Ecal,A})(x_A)= p_{\Ecal,A}(h_A(x_A)).  \]
Since $(E_F)_{F \in \tilde{H}}$ is jointly $\Pr$-independent and $E_v=\pi_v((E_F)_{v \in F} )$ we also have that
\[ E_{D_1}:=(E_v)_{v \in D_1}, \dots, E_{D_r}:=(E_v)_{v \in D_r}   \] 
are jointly $\Pr$-independent, where $D_1,\dots, D_r \ins A$ are the districts of $A$.
Indeed, $E_{D_i}$ is a function of $(E_F)_{\substack{F \in \tilde{H}\\ F \cap D_i \neq \emptyset}}$ and for indices $i \neq j$ the sets $\{ F \in \tilde{H} | F \cap D_i \neq \emptyset\}$ and $\{ F \in \tilde{H} | F \cap D_j \neq \emptyset\}$ are disjoint (by definition of districts).\\
It then follows that the density $p_{\Ecal,A}(e_A)$ factors as $p_{\Ecal,D_1}(e_{D_1})\cdots p_{\Ecal,D_r}(e_{D_r})$. %
So for $\tilde{p}_A$ this means:
\[ \tilde{p}_A(x_A)= p_{\Ecal,D_1}(h_{D_1}(x_{D_1\cup \Pa^G(D_1)}))\cdots p_{\Ecal,D_r}(h_{D_r}(x_{D_r\cup \Pa^G(D_r)})).  \]
Note that the components of $h_{D_i}$ by definition really only depend on $\Xcal_{D_i\cup \Pa^G(D_i)}$ because for $v \in D_i$ we have $\Anc^G(v) \ins A$ and thus $D(v) \ins D_i$. 
  Since for every district $D \ins A$ the set of nodes $D \cup \Pa^G(D)$ build a complete subgraph in $A^\moral$ the density $\tilde{p}_A$ factors according to $A^\moral$. \\
Furthermore, by assumption we have that $\nu_A^{g_A}$ has the density $|h'_A|$ w.r.t.\ $\mu_A$ that factors according to $A^\moral$.
It follows that $\Pr^{X_A}$ has a density w.r.t.\ $\mu_A$ given by:
\[p_A(x_A) = \frac{d \Pr^{X_A} }{d \mu_A} (x_A)= \frac{d \Pr^{X_A} }{d \nu_A^{g_A}} (x_A) \cdot \frac{d \nu_A^{g_A}}{d \mu_A} (x_A) = \tilde{p}_A(x_A) \cdot |h'_A|(x_A), \]
which then factors according to the complete subgraphs of $A^\moral$.
  This shows \hyl{aFP}{aFP}. The \hyl{asSEP}{asSEP} follows directly with $g_v$ given by the projection of $g_A$ onto the $v$-component for $A=\Anc^G(v)$.
\end{proof}
\end{Thm}

\begin{Rem}
  For the implication ``\hyl{SEPwared}{SEPwared} $\implies$ \hyl{aFP}{aFP}'' in \ref{SEPwared-aFP} the random variables $E_F$ for $F \in \tilde{H}$ were only used to encode the independence structure between the $E_v$ for $v \in V$. So instead of the existence of the $E_F$ in \hyl{SEPwared}{SEPwared} we could only ask for:
\begin{enumerate}
\item[(c')] For all  subsets $X,Y \ins V$ we have:
 \[  X \Indep_{(V,\emptyset,H_2)}^d Y \qquad \implies \qquad (E_v)_{v \in X} \Indep_\Pr (E_v)_{v \in Y}.   \]
 Or even weaker:
\item[(c'')] For every ancestral sub-\HEDG{} $A \ins G$ with its districts $D_1,\dots,D_r$ we have that $E_{D_1},\dots,E_{D_r}$ are jointly $\Pr$-independent, where we put $E_{D_i}:=(E_v)_{v \in D_i}$.
\end{enumerate}
  Then the implication ``\hyl{SEPwared}{SEPwared} $\implies$ \hyl{aFP}{aFP}'' for $(G,\Pr_V)$ holds as well. Note that these weaker versions only need bidirected edges (not the generality of hyperedges).
\end{Rem}

As an example and corollary of \ref{SEPwared-aFP} we get the linear model with density back (see \cite{Spirtes93}, \cite{Spirtes95}, \cite{Koster96}, \cite{foygel2012}, a.o.). Note that we do not need to assume Gaussianity and that we work with dependencies and not just correlations as other authors do. 

\begin{Eg}[Solvable linear model with density]
\label{linear-SEP}
Let $G=(V,E,H)$ be a \HEDG{}. Assume the following:
\begin{enumerate}
 \item There is an underlying probability space $(\Omega, \fa, \Pr)$.
 \item For every $v \in V$ we have a $1$-dimensional real-valued random variable $X_v$, which together satisfy the linear equations: 
        \[ X_v = \sum_{w \in \Pa^G(v)} b_{v,w} \cdot X_w + E_v,\] %
	 with real-valued error variables $E_v$, $v \in V$, such that:
 \item the distribution of the tuple $(E_v)_{v \in V}$ has a joint density w.r.t.\ the Lebesgue measure $\lambda$, and: %
	\item the (in)dependence structure (not just correlations!) is assumed to be given by:
\[  X \Indep_{(V,\emptyset,H_2)}^d Y \qquad \implies \qquad (E_v)_{v \in X} \Indep_\Pr (E_v)_{v \in Y}.   \]		
 \item[] For instance, this holds if for every $F \in \tilde{H}$ there is a $\#F$-dimensional random variable $E_F=(E_{F,v})_{v \in F}$ 
		such that $(E_F)_{F \in \tilde{H}}$ is jointly $\Pr$-independent and
 \[ E_v= \pi_v((E_F)_{v \in F \in \tilde{H}}):= \sum_{\substack{v \in F \in \tilde{H}}} c_{F,v} \cdot E_{F,v}\]
 for some real numbers $c_{F,v}$.
 \item The following matrix is invertible (i.e.\ assuming solvability):
				 \[ I_V-B_V := \lp \I_{v}(w)- b_{v,w} \cdot \I_{\Pa^G(v)}(w) \rp_{v,w \in V}.    \]
\end{enumerate}
Then for all ancestral $A \ins G$ the matrix 
\[ I_A-B_A := \lp \I_{v}(w)- b_{v,w} \cdot \I_{\Pa^G(v)}(w) \rp_{v,w \in A}    \]
is also invertible.
Furthermore, we can put:
 \begin{enumerate}
  \item $h_v(x_{\{v\} \cup \Pa^G(v)}):=x_v - \sum_{w \in \Pa^G(v)} b_{v,w} \cdot x_w$ leading to
	 \[  h_v(X_{\{v\} \cup \Pa^G(v)}) = E_v, \]
	\item $h_A(x_A)=(I_A-B_A) \cdot x_A$ and
  \item $g_A(e_A):=(I_A-B_A)^{-1} \cdot e_A$, which implies
	\item $|h'_A| (x_A) = \frac{d\lambda^{g_A}_A}{d\lambda_A}(x_A)= |\frac{d h_A }{d x_A}|=|\det(I_A-B_A)| =$ constant $ \neq 0$ (by the standard transformation rule for the Lebesgue measure).
 \end{enumerate}
  It follows that $\Pr^{(X_v)_{v \in V}}$ satisfies \hyl{SEPwared}{SEPwared} and thus the ancestral factorization property (\hyl{aFP}{aFP}) and the directed global Markov property (\hyl{dGMP}{dGMP}).\\
	Also note that by directly working with $E_v$ for $v \in V$ we only need to use bidirected edges (and not general hyperedges).
	
\end{Eg}

\section{Graphical Models}

In this section we will shortly and explicitely indicate how to use the Markov properties for directed graphs with hyperedges (\HEDG{}es) to create model classes.

\subsection{Probabilistic Models for \HEDG{}es}

In this subsection we shortly define probabilistic graphical models attached to some directed graph with hyperedges (\HEDG{}) w.r.t.\ some Markov property.

\begin{Def}
\label{def-prob-model}
Let $V$ be finite set of nodes, $\Xcal_v$ for $v \in V$ standard Borel spaces and $\Xcal_V:=\prod_{v\in V} \Xcal_v$ the product space.
Furthermore, let $G=(V,E,H)$ be a directed graph with hyperedges (\HEDG{}) with set of nodes $V$.
We then define:
$$\Pcal(V) := \{ \Pr_V \,|\, \Pr_V \text{ a probability distribution on } \Xcal_V \}$$
 to be the set of all probability distributions on $\Xcal_V$.\\
Furthermore, let \emph{MP} denote a fixed Markov property from the last section (see overview in figure \ref{fig:overview}).
Then we define the model (class) of $G$ w.r.t.\ \emph{MP} as:
\[ \Mcal_{\mathrm{MP}}(G):=\{\Pr_V \in \Pcal(V)\,|\, (G,\Pr_V) \text{ satisfies MP}\}  \]
the set of all probability distributions on $\Xcal_V$ that satisfies the Markov property \emph{MP} w.r.t.\ $G$.
\end{Def}

\begin{Rem}
\label{rem-nmp}
Let the notations be like in \ref{def-prob-model}. 
\begin{enumerate}
\item Every implication from the last section (see figure \ref{fig:overview}) leads to an inclusion of models.
E.g.\ ``\hyl{mMP}{mdGMP} $\implies$ \hyl{dGMP}{dGMP}'' (see \ref{mdGMP-dGMP}) induces: 
$$\Mcal_{\textrm{mdGMP}}(G) \ins \Mcal_{\textrm{dGMP}}(G),$$
the inclusion of the marginal directed global Markov model into the (ordinary) directed global Markov model.
\item If $G$ is an mDAG (i.e.\ $G^\acy$ of any \HEDG{} $G$) then also an intermediate model called ``\emph{nested Markov model}'' $\Mcal_{\textrm{nMP}}(G)$ was introduced and analysed in \cite{Eva15}, \cite{SERR14}, \cite{RERS17}:
$$\Mcal_{\textrm{mdGMP}}(G) \ins \Mcal_{\textrm{nMP}}(G) \ins \Mcal_{\textrm{dGMP}}(G).$$
\end{enumerate}
\end{Rem}

\begin{Eg}
\label{prob-model-eg}
Again, let the notations be like in \ref{def-prob-model}. As a result of the previous section we then have the following inclusions of models:
\begin{enumerate}
\item For every \HEDG{} $G$ we have (see \ref{aFP-dGMP}, \ref{dGMP-auPMP}, \ref{dGMP-gdGMP}):
\[ \Mcal_{\textrm{aFP}}(G) \ins \Mcal_{\textrm{auGMP}}(G)=\Mcal_{\textrm{dGMP}}(G) \ins \Mcal_{\textrm{gdGMP}}(G), \]
highlighting the most important ``ordinary'' Markov models (for general \HEDG{}es) that behave well under marginalizations.
\item Note that these properties do only depend on the induced directed mixed graph $G_2$ (using directed and bidirected edges only). So, for instance, we have:
\[ \Mcal_{\textrm{dGMP}}(G)= \Mcal_{\textrm{dGMP}}(G_2) \ins \Mcal_{\textrm{gdGMP}}(G_2) = \Mcal_{\textrm{gdGMP}}(G). \]
\item For every \HEDG{} $G$ we also have the inclusions of the marginal versions (see \ref{mFP-maFP}, \ref{mdGMP-dGMP}):
\[ \Mcal_{\textrm{mFP}}(G) \ins\Mcal_{\textrm{maFP}}(G) \ins \Mcal_{\textrm{mdGMP}}(G) \ins \Mcal_{\textrm{mgdGMP}}(G), \]
which then i.g.\ will depend on the hyperedges (and not only on the induced bidirected edges as shown in example 
\ref{eg-dGMP-mdGMP}, figure \ref{fig:eg-dGMP-mdGMP}.).
\item For every \HEDG{} $G$ we also have (see \ref{s-sep-d-sep-eg}, \ref{dGMP-gdGMP}):
\[ \Mcal_{\textrm{gdGMP}}(G) = \Mcal_{\textrm{gdGMP}}(G^\acy) = \Mcal_{\textrm{dGMP}}(G^\acy),\] %
showing the equivalence between the general directed global Markov model of a \HEDG{} $G$ and the ``ordinary'' directed global Markov model of its acyclification $G^\acy$.
\item Note that for any \HEDG{} $G$ by \ref{rem-nmp} we have the inclusions:
\[ \Mcal_{\textrm{mdGMP}}(G^\acy) %
\stackrel{\ref{rem-nmp}}{\ins} \Mcal_{\textrm{nMP}}(G^\acy) \stackrel{\ref{rem-nmp}}{\ins} \Mcal_{\textrm{dGMP}}(G^\acy), \]
indicating the relation between the marginal directed global Markov model, the nested Markov model of \cite{SERR14}, \cite{RERS17} and the ``ordinary'' directed global Markov model.
\item For every \HEDG{} $G$ we furthermore have (see \ref{SEP-equiv}, \ref{csSEP-smgdGMP}, \ref{smgdGMP-mdGMP-acy}, \ref{smgdGMP-mgdGMP}):
\[\begin{array}{ccccc}
 \Mcal_{\textrm{lsSEP}}(G) &\ins& \Mcal_{\textrm{csSEP}}(G) &\ins& \Mcal_{\textrm{smgdGMP}}(G)\\
 &\ins&\Mcal_{\textrm{mdGMP}}(G^\acy) &\cap & \Mcal_{\textrm{mgdGMP}}(G)\\
 &\ins&\Mcal_{\textrm{gdGMP}}(G),
\end{array}\]
focusing on the fact that the (loop-wisely solvable) structural equations model is included in the (strong marginal) general directed global Markov model encoding conditional independence relations using $\sigma$-separation.
\item If $G$ is a DAG and $<$ a fixed topological order on $G$ then we have the equalities (see \ref{dag-markov-prop}):
\[ \Mcal_{\textrm{dGMP}}(G) = \Mcal_{\textrm{rFP}}(G) = \Mcal_{\textrm{dLMP}}(G) = \Mcal_{\textrm{oLMP}}^<(G) = \Mcal_{\textrm{SEP}}(G).   \]
\item If $G$ is an mDAG we then have the equalities (see \ref{mDAG-all}):
\[ \Mcal_{\textrm{mdGMP}}(G) = \Mcal_{\textrm{mgdGMP}}(G) = \Mcal_{\textrm{lusSEP}}(G) =  \Mcal_{\textrm{SEP}}(G).   \]

\end{enumerate}
\end{Eg}

\subsection{Causal Models for \HEDG{}es}

In this subsection we shortly want to indicate how to use a more rigid version of the loop-wisely solvable structural equations property (\hyl{lsSEP}{lsSEP})
to model causal relations that allow for interactions, feedback loops, latent confounding and non-linear functions. The advantages for using this approach are that it is stable under marginalization and, surprisingly, 
also stable under interventions. In other words, all marginalizations of all interventions also satisfy the \hyl{lsSEP}{lsSEP} and thus the general directed global Markov property (\hyl{gdGMP}{gdGMP}) based on $\sigma$-separation. It reduces to the usual models in the acyclic cases (mDAGs, ADMGs, DAGs). 
It follows the idea that every (cyclic) subsystem is already determined by its joint causal parents outside of the subsystem. Furthermore, it abstains from modelling self-loops, which cannot even be detected by interventions.
More extensive work on causality and structural causal models can be found in \cite{Pearl09}, \cite{SGS00}, \cite{PJS17}, \cite{BPSM16} a.o..

\begin{Def}[Loops of a \HEDG{}]
Let $G=(V,E,H)$ be a \HEDG{}.
\begin{enumerate}
  \item Remember that a \emph{strongly connected induced subgraph} of $G$ is a set of nodes $S \ins V$ that is strongly connected in the induced subgraph structure of $G$. As an abbreviation (in this sub-\HEDG{} structure) we will call $S$ a \emph{loop of $G$}.
\item Let $\Lcal(G):=\{ S \ins G \;|\, S \text{ a loop of } G  \}$ be the \emph{loop set of $G$}. 
\end{enumerate}
\end{Def}

\begin{Rem}
Note that the loop set $\Lcal(G)$ contains all single element loops $\{v\} \in \Lcal(G)$, $v \in V$, as the smallest loops and all strongly connected components $\Sc^G(v) \in \Lcal(G)$, $v \in V$, as the largest loops, but also all non-trivial intermediate loops $S$ with $\{v\} \subsetneq S \subsetneq \Sc^G(v)$ inside the strongly connected components (if existent). If $G$ is an mDAG (e.g.\ an ADMG or DAG) then $\Lcal(G)$ only consists of the single element loops: $\Lcal(G)=\{\{v\} \,|\, v \in V\}$.
\end{Rem}

\begin{Def}[Compatible system of structural equations]
\label{compatible-str-eq}
Let $G=(V,E,H)$ be a \HEDG{} with its loop set $\Lcal(G)$.
For every $v \in V$, $F \in \tilde{H}$ resp., fix a standard Borel space $\Xcal_v$, $\Ecal_F$ resp.. 
For a subset $S \ins V$ we put: 
\[\Xcal_S:=\prod_{v \in S} \Xcal_v \qquad \text{and} \qquad \Ecal_S := \prod_{\substack{F \in \tilde{H}\\ F \cap S \neq \emptyset}} \Ecal_F.\]
We now consider a family of measurable functions $(g_S)_{S\in \Lcal(G)}$ indexed by $\Lcal(G)$ as follows:
\[ g_S: \; \Xcal_{\Pa^G(S)\sm S} \x \Ecal_S \to \Xcal_S.\]
We say that the system $(g_S)_{S\in\Lcal(G)}$ is \emph{compatible} if for all $S, S' \in \Lcal(G)$ with $S' \ins S$ and all 
$x_{\Pa^G(S)\cup S} \in \Xcal_{\Pa^G(S) \cup S}$ and $e_S \in \Ecal_S$ we have:
\[ x_S = g_S(x_{\Pa^G(S)\sm S},e_S) \; \implies \; x_{S'} = g_{S'}(x_{\Pa^G(S')\sm S'},e_{S'}),\]
where $x_S$, $e_{S'}$ etc.\ are the corresponding components of $x_{\Pa^G(S)\cup S}$ and $e_S$ resp..
Note that with $S' \ins S$ we also have $\Pa^G(S') \cup S' \ins \Pa^G(S) \cup S$ and
 $\{ F \in \tilde{H}\, |\, F \cap S' \neq \emptyset\} \ins \{ F \in \tilde{H}\, |\, F \cap S \neq \emptyset\}$.\\
In this case we will call $(g_S)_{S \in \Lcal(G)}$ a \emph{compatible system of structural equations of $G$} or just a \emph{compatible system of equations of $G$}.
\end{Def}

\begin{Def}[Modular structural causal model (mSCM)]
\label{mSCM-def}
A \emph{modular structural causal model (mSCM)} by definition consists of:
\begin{enumerate}
\item A \HEDG{} $G=(V,E,H)$,
\item a standard Borel space $\Xcal_v$ for every $v \in V$,
\item a standard Borel space $\Ecal_F$ for every $F \in \tilde{H}$,
\item a probability measure $\Pr_{\Ecal_F}$ on $\Ecal_F$ for every $F \in \tilde{H}$,
\item a compatible system of structural equations $(g_S)_{S \in \Lcal(G)}$:
 \[ g_S: \; \Xcal_{\Pa^G(S)\sm S} \x \Ecal_S \to \Xcal_S,\]
with their components $g_S=(g_{S,v})_{v \in S}$.
\end{enumerate}
We make the following abbreviations:
\[\begin{array}{ccll}
 \Xcal&:=&\Xcal_V:=\prod_{v \in V} \Xcal_v, \\
 \Ecal&:=&\Ecal_V:=\prod_{F \in \tilde{H}} \Ecal_F,\\ 
 \Pr&:=& \otimes_{F \in \tilde{H}} \Pr_{\Ecal_F} & \text{ the product measure on } \Ecal,  \\
 E_S &:& \Ecal \to \Ecal_S & \text{ the natural projection},\\
 g &=&(g_S)_{S \in \Lcal(G)}.
\end{array}\]
The mSCM is then given by the tuple $M=(G,\Xcal,\Ecal,\Pr,g)$.
\end{Def}

\begin{Rem}
Let $M=(G, \Xcal, \Ecal, \Pr, g)$ be a mSCM and $S_1,\dots,S_r \in \Lcal(G)$ the strongly connected components of $G$ ordered according to a topological order of the DAG of strongly connected components of $G$.
Then inductively we can define the random variables:
\[ \begin{array}{ccll}
X_v &:=& g_{S_1,v}(E_{S_1}) & \text{ for all } v \in S_1,\\
X_v &:=& g_{S_i,v}(X_{\Pa^G(S_i)\sm S_i},E_{S_i}) & \text{ for all } v \in S_i,
\end{array}\]
for all $i>1$ etc.. By the compatibility of the system $g$ we then also have for every $S \in \Lcal(G)$:
\[ X_S = g_S(X_{\Pa^G(S)\sm S},E_S).   \]
Furthermore, if we put $\Pr_V:=\Pr^{(X_v)_{v \in V}}$ to be the induced distribution on $\Xcal=\Xcal_V$ then clearly $(G,\Pr_V)$ satisfies the loop-wisely solvable structural equations property (\hyl{lsSEP}{lsSEP}) from \ref{def-sep} and thus the (strong marginal) general directed global Markov property (\hyl{smgdGMP}{smgdGMP}) by \ref{csSEP-smgdGMP}. In short: 
$$\Pr_V \in \Mcal_{\mathrm{lsSEP}}(G) \ins \Mcal_{\mathrm{smgdGMP}}(G) \ins \Mcal_{\mathrm{mdGMP}}(G^\acy) \cap \Mcal_{\mathrm{mgdGMP}}(G).$$
\end{Rem}

\begin{Def}[Marginalization of a mSCM]
\label{mSCM-marg}
Let $M=(G, \Xcal, \Ecal, \Pr, g)$ be a mSCM with \HEDG{} $G=(V,E,H)$ and $W \ins V$ a subset.
Fix a map $\varphi: \tilde{H} \to \tilde{H}^\marg(W)$ with $\varphi(F) \supseteq F^\marg$ for all $F \in \tilde{H}$ (with $F^\marg$ as in \ref{HEDG-marg-def}). 
Then we can define a \emph{marginalization} 
\[M^{\marg(W)} =(G^{\marg(W)}, \Xcal^{\marg(W)}, \Ecal^{\marg(W)}, \Pr^{\marg(W)}, g^{\marg(W)})\]
of $M$ as the mSCM given as follows: 
\begin{enumerate}
\item $G^{\marg(W)}=(W,E^{\marg(W)},H^{\marg(W)})$ the marginalized \HEDG{} from \ref{HEDG-marg-def},
\item $\Xcal_v^{\marg(W)}:= \Xcal_v$ for $v \in W$,
\item $\Xcal^{\marg(W)}=\Xcal_W= \prod_{v \in W} \Xcal_v$,
\item $\Ecal_{F'}^{\marg(W)} := \prod_{\substack{F \in \tilde{H}\\\varphi(F)=F'}} \Ecal_F$ for every $F' \in \tilde{H}^{\marg(W)}$,
\item $\Ecal^{\marg(W)} = \prod_{F' \in \tilde{H}^{\marg(W)}} \Ecal_{F'} ^{\marg(W)} = 
\Ecal$,
\item $\Pr_{\Ecal_{F'}} ^{\marg(W)} = \otimes_{\substack{F \in \tilde{H}\\\varphi(F)=F'}} \Pr_{\Ecal_F}$,
\item $\Pr^{\marg(W)}= \Pr$ on $\Ecal^{\marg(W)}=\Ecal$.
\item $E_{F'}^{\marg(W)} =(E_F)_{\substack{F \in \tilde{H}\\\varphi(F)=F'}}$,
\item For $C \in \Lcal(G^{\marg(W)})$ we construct $g^{\marg(W)}_{C}$ from the components of $g_S$ corresponding to $v \in C$, where $S$ is the strongly connected component of $C$ in $C \cup W$ in the induced sub-\HEDG{} structure of $G$, inductively by marginalizing out one element $u \in V \sm W$ at a time. The latter can be done by substituting $x_u$ with $g_{\{u\}}(x_{\Pa^G(u)\sm \{u\}},e_u)$ in any function if it occurs (like in the proof of \ref{lsSEP-marg}). Since substituting for a node $u$ commutes with substiting for a node $u'$ the final procedure for the whole set $V \sm W$ does not depend on the marginalization order of the nodes.
Also the compatibility condition follows directly from this construction.
\end{enumerate}
\end{Def}

\begin{Def}[Stochastic Interventions on a mSCM]
\label{mSCM-intv}
Let $M=(G, \Xcal, \Ecal, \Pr, g)$ be a mSCM with \HEDG{} $G=(V,E,H)$ and $I \ins V$ a set of nodes (the intervention target). 
For every $v \in I$ fix a probability measure $\Pr'_v$ on $\Xcal_v$ (e.g.\ $\delta_{x_v}$ for an $x_v \in \Xcal_v$).
Then we can define a \emph{post-intervention mSCM}
\[M^{\doi(I)} =(G^{\doi(I)}, \Xcal^{\doi(I)}, \Ecal^{\doi(I)},\Pr^{\doi(I)},g^{\doi(I)})\] 
of $M$ as follows:
\begin{enumerate}
\item Put $G^{\doi(I)}:=(V^{\doi(I)},E^{\doi(I)},H^{\doi(I)})$ with:
   \begin{enumerate}
	   \item $V^{\doi(I)}:=V$,
		 \item $E^{\doi(I)}:=E \sm \{ v \to w\, | \,w \in I\}$,
		 \item $H^{\doi(I)}:=\{F \sm I \,|\, F \in H\} \cup \{ \{v\}\,|\,v \in I\}$.
	\end{enumerate}
\item $\Xcal^{\doi(I)}_v := \Xcal_v$ for every $v \in V$, and thus $\Xcal^{\doi(I)} = \Xcal$.
\item Fix a map $\varphi: \{ F \in \tilde{H}\,|\, F \sm I \neq \emptyset \} \to \tilde{H}^{\doi(I)} \sm \{ \{v\}\,|\,v \in I\}$ with $\varphi(F) \supseteq F \sm I$, and put:
\item $\Ecal_{F'}^{\doi(I)}:= \prod_{\varphi(F)=F'} \Ecal_F$ for $F' \in \tilde{H}^{\doi(I)} \sm \{ \{v\}\,|\,v \in I\}$
together with the product measure $\Pr^{\doi(I)}_{\Ecal_{F'}^{\doi(I)}} := \otimes_{\varphi(F)=F'} \Pr_{\Ecal_F}$, and
\item $\Ecal^{\doi(I)}_{\{v\}}:=\Xcal_v$ for $v \in I$ together with the measure $\Pr^{\doi(I)}_{\Ecal_{\{v\}}^{\doi(I)}}:=\Pr'_v$.
\item For $v \in I$ put $g^{\doi(I)}_{\{v\}}:=\id_{\Xcal_v}$ and
\item for every other $S \in \Lcal(G^{\doi(I)})$ put $g^{\doi(I)}_S:=g_S$. Note that if $S$ is a loop in  $G^{\doi(I)}$ then $S$ is also a loop in $G$.
\end{enumerate}
Note that the construction is dependent on the choice of $\varphi$ and the measures $\Pr_v'$ for $v \in I$.
\begin{proof}
For the very last point let $S \in \Lcal(G^{\doi(I)})$ be a strongly connected induced sub-\HEDG{} of $G^{\doi(I)}$. Let $S'$ be the induced sub-\HEDG{} of the nodes of $S$ in $G$. If $S'$ had an edge $v \to w$ that was not in $S$ then $w \in I$, which was excluded. So $S'=S$ is a loop in $G$.
\end{proof}
\end{Def}

\begin{Rem}
Let $M=(G, \Xcal, \Ecal, \Pr, g)$ be a mSCM. Then all marginalizations of all possible interventions are defined. Every such mSCM has then well-defined random variables $X_v$ whose induced probability distribution satisfies the loop-wisely solvable structural equations property (\hyl{lsSEP}{lsSEP}) and thus the (strong marginal) general directed global Markov property (\hyl{smgdGMP}{smgdGMP}) by \ref{csSEP-smgdGMP}.
\end{Rem}

\begin{Rem}[Reduced modular structural causal model (rmSCM)]
If we want to have a causal model using only bidirected edges (instead of more general hyperedges) one could make the following reductions in the definition of mSCM \ref{mSCM-def}:
\begin{enumerate}
 \item Take $G=(V,E,H_2)$ with only bidirected edges (e.g.\ the induced directed mixed graph of a \HEDG{} $(V,E,H)$).
 \item Directly assume the existence of $\Ecal_v$  for every $v \in V$ and put $\Ecal_S:=\prod_{v \in S} \Ecal_v$ 
    (instead of $\prod_{\substack{F \in \tilde{H}\\ F \cap S \neq \emptyset}} \Ecal_F$).
 \item With this replacement take the corresponding definitions for compatible functions, etc..
 \item Assume a probability measure $\Pr$ on $\Ecal_V=\prod_{v \in V} \Ecal_v$ with independence structure of the canonical projections $E_v: \Ecal_V \to \Ecal_v$ for $v \in V$ given by the rule:
   \[  X \Indep_{(V,\emptyset,H_2)}^d Y \qquad \implies \qquad (E_v)_{v \in X} \Indep_\Pr (E_v)_{v \in Y}.   \]
 \item Marginalization and interventions can then be defined as the induced directed mixed graph versions of the mSCM version.	
 \item Furthermore, if one even wanted to have the directed global Markov property (\hyl{dGMP}{dGMP}) in such a rmSCM one could enforce the conditions of \hyl{SEPwared}{SEPwared}, i.e., roughly speaking, that we assume densities and that we only allow for functions such that
for every ancestral $A \ins G$ the inverse to $g_A$ exists on a reduction of the error spaces and that its Jacobian determinant factorizes according to $A^\moral$.
\end{enumerate}

\end{Rem}

\section{Discussion}

We introduced \emph{directed graphs with hyperedges} (\HEDG{}es) as a combination and generalization of marginal directed acyclic graphs (mDAGs) and directed mixed graphs (DMGs) allowing for cycles.
We connected \HEDG{}es with probability distributions in different ways by introducing several \emph{Markov properties}. 
We formally justified the use of hyperedges as a summary of the whole latent world for the most important Markov properties.
We analysed the logical relations between all found Markov properties. An overview is given in figure \ref{fig:overview}.
We generalized central results from 
\cite{Pearl86c}, \cite{Lau90}, \cite{Verma93}, \cite{Richardson03}, \cite{Eva15}, \cite{Spirtes94}, \cite{PearlDechter96}, \cite{Koster96} in a unified way by either allowing for cycles, hyperedges or both. 
The relations between the most important Markov properties for \HEDG{}es that are stable under marginalizations are given as follows:
\[   \hyl{aFP}{aFP} \;\implies\; \hyl{dGMP}{dGMP} \;\implies\; \hyl{gdGMP}{gdGMP} \;\Longleftarrow\; \hyl{lsSEP}{lsSEP}. \]
For marginal versions and reverse implications see figure \ref{fig:overview}.
Based on the different types of Markov properties we introduced different types of \emph{probabilistic graphical Markov models} for \HEDG{}es, which might have both cycles and hyperedges.

Finally, we used our insights into structural equations (\hyl{lsSEP}{lsSEP}) for \HEDG{}es to define the well-behaved class of \emph{modular structural causal models} (mSCM) that allow for all at once: non-linear functional relations, interactions, feedback loops and latent confounding. We showed that all interventions, all marginalizations and all combinations of these exist and are well-defined. Each of the models will then encode conditional independence relations based on our introduced $\sigma$-separation criterion (\hyl{gdGMP}{gdGMP}). 
Since structural causal models with cycles i.g.\ do not encode conditional independence relations based on d/m/m*-separation (see \ref{main-example}) we suggest using $\sigma$-separation for causal modelling and causal discovery algorithms instead. Note that $\sigma$-separation reduces to the usual d/m/m*-separation in the acyclic case (i.e.\ for DAGs, ADMGs, mDAGs, etc.)\ and is thus a better generalization of d/m/m*-separation to the \HEDG{} case in the presence of cycles and non-linear functions (and latent confounders).

Future work might address the following topics:
\begin{enumerate}
\item Extending \HEDG{}es to graphical structures with additional (hyper)edge types, e.g.\ \emph{undirected} edges to represent \emph{selection biases}, etc..
\item Merging \HEDG{}es with \emph{acyclic graphs} (see \cite{Lau16}).
\item Adjusting Markov properties of \HEDG{}es for \emph{deterministic relations}. 
\item \emph{Parametrisation} of a \HEDG{} w.r.t.\ a fixed Markov property and \emph{fitting} such models to data.
\item Proving \emph{completeness} results for some of the found Markov properties for \HEDG{}es.
\item Improving conditional independence \emph{constraints based causal discovery} algorithms (e.g.\ \cite{HEJ14}, \cite{HHEJ13}, a.o.)\ 
by using \emph{$\sigma$-separation} in \HEDG{}es instead of d/m/m*-separation.
\item Using \HEDG{}es for \emph{score based causal discovery} algorithms.
\item Extending the \emph{nested Markov properties} from mDAGs to \HEDG{}es and analysing their corresponding models (see \cite{Eva15}, \cite{SERR14}, \cite{RERS17}).
\end{enumerate}

\section{Acknowledgments}

This work was supported by the European Research Council (ERC) under the European Union's Horizon 2020 research and innovation programme (grant agreement 639466).

\bibliographystyle{amsalpha} 

\providecommand{\bysame}{\leavevmode\hbox to3em{\hrulefill}\thinspace}
\providecommand{\MR}{\relax\ifhmode\unskip\space\fi MR }
\providecommand{\MRhref}[2]{%
  \href{http://www.ams.org/mathscinet-getitem?mr=#1}{#2}
}
\providecommand{\href}[2]{#2}

\end{document}